\documentclass[12pt, leqno]{amsart}

\usepackage{graphicx}
\usepackage{amssymb,amsmath,amsthm,amscd,accents}
\usepackage{stackengine,scalerel}
\usepackage{mathrsfs}
\usepackage{lscape}
\usepackage{enumerate}
\usepackage{upgreek}
\usepackage{stmaryrd}
\usepackage[usenames,dvipsnames]{color}
\usepackage[colorlinks=true, pdfstartview=FitV,
 linkcolor=blue,citecolor=blue,urlcolor=blue]{hyperref}
\usepackage{tikz}
\usetikzlibrary{calc,tikzmark,decorations.pathmorphing}
\tikzset{dynkdot/.style={circle,draw,scale=.38}}
\usetikzlibrary{arrows.meta,arrows}
\usepackage{bm}
\usepackage{verbatim}
\usepackage{xspace}

\usepackage[all]{xy}

\allowdisplaybreaks[3]

\newcommand{\nc}{\newcommand}
\numberwithin{equation}{section}

\nc{\eq}{\begin{myequation}} \nc{\eneq}{\end{myequation}}
\nc{\eqn}{\begin{myequationn}} \nc{\eneqn}{\end{myequationn}}

\nc{\ms}{\mspace}

\nc{\qR}[1]{\ttq_{\mspace{-2mu}\raisebox{-.8ex}{${\scriptstyle{#1}}$}}}

\theoremstyle{plain}
\newtheorem{lemma}{Lemma}[section]
\newtheorem{proposition}[lemma]{Proposition}
\newtheorem{theorem}[lemma]{Theorem}

\newtheorem{corollary}[lemma]{Corollary}
\newtheorem{conjecture}[lemma]{Conjecture}

\newtheorem*{thmx}{Main Theorem}

\theoremstyle{definition}
\newtheorem{remark}[lemma]{Remark}
\newtheorem{example}[lemma]{Example}

\newtheorem{definition}[lemma]{Definition}

\nc{\Def}{\begin{definition}}
\nc{\edf}{\end{definition}}
\nc{\mTh}{\begin{mtheorem}}
\nc{\enmth}{\end{mtheorem}}
\nc{\Th}{\begin{theorem}}
\nc{\enth}{\end{theorem}}
\nc{\Prop}{\begin{proposition}}
\nc{\enprop}{\end{proposition}}
\nc{\Lemma}{\begin{lemma}}
\nc{\enlemma}{\end{lemma}}
\nc{\Cor}{\begin{corollary}}
\nc{\encor}{\end{corollary}}
\nc{\Rem}{\begin{remark}}
\nc{\enrem}{\end{remark}}
\nc{\Conj}{\begin{conjecture}}
\nc{\enconj}{\end{conjecture}}

\nc{\ledot}{\mathrel{\le\ms{-11mu}\raisebox{.2ex}{$\cdot$}}}
\nc{\predot}{\mathrel{\preceq\ms{-9mu}\raisebox{.35ex}{$\centerdot$}}}

\newcommand{\hxi}{{\widehat{\xi}}}

\renewcommand{\le}{\leqslant}
\renewcommand{\ge}{\geqslant}
\renewcommand{\preceq}{\preccurlyeq}

\newcommand{\seteq}{\mathbin{:=}}

\newcommand{\soplus}{\mathop{\mbox{\normalsize$\bigoplus$}}\limits}

\newcommand{\ev}{{\operatorname{ev}}}

\newcommand{\gmod}{\text{-}\mathrm{gmod}}

\newcommand{\ex}{\mathrm{ex}}
\newcommand{\fr}{\mathrm{fr}}

\newcommand{\g}{\mathfrak{g}}

\newcommand{\n}{\mathfrak{n}}

\newcommand{\C}{\mathbb{C}}
\newcommand{\Q}{\mathbb{Q}}
\newcommand{\Z}{\mathbb{Z}\ms{1mu}}
\newcommand{\N}{\mathbb{N}\ms{1mu}}

\newcommand{\al}{{\ms{1mu}\alpha}}

\newcommand{\la}{\lambda}
\newcommand{\be}{{\ms{1mu}\beta}}

\newcommand{\ve}{\boldsymbol{\varepsilon}}
\newcommand{\La}{\Lambda}

\newcommand{\wt}{{\rm wt}}

\newcommand{\Fold}{{\rm Fold}}
\newcommand{\bx}{{\rm box}}
\newcommand{\re}{{\rm re}}
\newcommand{\loc}{{\rm loc}}

\newcommand{\ch}{{\rm ch}}
\newcommand{\rev}{{\rm rev}}
\newcommand{\up}{{\rm up}}
\newcommand{\full}{{\rm full}}

\newcommand{\Dynkin}{\triangle}  

\newcommand{\im}{\imath}
\newcommand{\jm}{\jmath}

\newcommand{\ii}{ \textbf{\textit{i}}}
\newcommand{\jj}{ \textbf{\textit{j}}}

\newcommand*\circled[1]{ \fontsize{6}{6}\selectfont \tikz[baseline=(char.base)]{
  \node[shape=circle,draw,inner sep=0.4pt] (char) {#1};} \fontsize{12}{12}\selectfont }

\newcommand*\bcircled[1]{\fontsize{6}{6}\selectfont \tikz[baseline=(char.base)]{
    \node[shape=circle, fill= black, draw=black, text=white,  inner sep=0.4pt] (char) {#1};} \fontsize{12}{12}\selectfont}

\newcommand*\mcircled[1]{\fontsize{6}{6}\selectfont \tikz[baseline=(char.base)]{
    \node[shape=circle, fill= magenta, draw=black, text=white, inner sep=0.4pt] (char) {#1};} \fontsize{12}{12}\selectfont}

\newcommand*\rcircled[1]{\fontsize{6}{6}\selectfont \tikz[baseline=(char.base)]{
    \node[shape=circle, fill= blue, draw=black, text=white, inner sep=0.4pt] (char) {#1};} \fontsize{12}{12}\selectfont}

\newcommand{\Hom}{\operatorname{Hom}}

\newcommand{\Fib}{\operatorname{Fib}}
\newcommand{\eff}{\operatorname{eff}}

\newcommand{\tG}{\widetilde{G}}

\newcommand{\tb}{\widetilde{b}}

\newcommand{\tx}{\widetilde{x}}
\newcommand{\ty}{\widetilde{y}}

\newcommand{\tii}{{\widetilde{\ii}}}
\newcommand{\htii}{{\widehat{\tii}}}

\newcommand{\hB}{\widehat{B}}

\newcommand{\hI}{\widehat{I}}
\newcommand{\hJ}{\widehat{J}}

\newcommand{\hM}{\widehat{M}}
\newcommand{\hN}{\widehat{N}}

\newcommand{\hP}{\widehat{P}}
\newcommand{\hQ}{\widehat{Q}}

\newcommand{\hZ}{\widehat{Z}}

\newcommand{\ha}{\widehat{a}}
\newcommand{\hb}{\widehat{b}}
\newcommand{\hc}{\widehat{c}}

\newcommand{\hhm}{\widehat{m}}
\newcommand{\hg}{{\widehat{\g}}}

\newcommand{\hk}{\widehat{k}}
\newcommand{\hl}{\widehat{l}}

\newcommand{\hs}{\widehat{s}}

\newcommand{\hr}{\widehat{r}}

\newcommand{\hw}{\widehat{w}}
\newcommand{\hx}{\widehat{x}}
\newcommand{\hy}{\widehat{y}}
\newcommand{\hz}{\widehat{z}}

\newcommand{\hal}{\widehat{\al}}
\newcommand{\hbe}{\widehat{\be}}

\newcommand{\heta}{\widehat{\eta}}
\newcommand{\htau}{\widehat{\tau}}

\newcommand{\hii}{{\widehat{\ii}}}
\newcommand{\hjj}{{\widehat{\jj}}}

\newcommand{\hve}{{\widehat{\ve}}}

\newcommand{\oT}{\overline{T}}

\newcommand{\ocalL}{\overline{\calL}}
\newcommand{\ocalR}{\overline{\calR}}

\newcommand{\frakB}{\mathfrak{B}}
\newcommand{\frakC}{\mathfrak{C}}
\newcommand{\frakD}{\mathfrak{D}}

\newcommand{\frakK}{\mathfrak{K}}

\newcommand{\frakS}{\mathfrak{S}}
\newcommand{\frakT}{\mathfrak{T}}

\newcommand{\fraka}{\mathfrak{a}}

\newcommand{\frakc}{\mathfrak{c}}
\newcommand{\frakd}{\mathfrak{d}}

\newcommand{\frakf}{\mathfrak{f}}

\newcommand{\frakn}{\mathfrak{n}}

\newcommand{\fraks}{\mathfrak{s}}

\newcommand{\fraku}{\mathfrak{u}}

\newcommand{\tfrakC}{\widetilde{\frakC}}
\newcommand{\tfrakc}{\widetilde{\frakc}}

\newcommand{\hfrakT}{{\widehat{\frakT}}}
\newcommand{\hfrakC}{\widehat{\frakC}}
\newcommand{\hfrakc}{\widehat{\frakc}}
\newcommand{\hfrakf}{\widehat{\frakf}}
\newcommand{\hfrakn}{\widehat{\frakn}}
\newcommand{\hfrakd}{\widehat{\frakd}}

\newcommand{\htfrakC}{\widehat{\tfrakC}}
 
\newcommand{\sfA}{\mathsf{A}}
\newcommand{\sfB}{\mathsf{B}}
\newcommand{\sfC}{\mathsf{C}}
\newcommand{\sfD}{\mathsf{D}}

\newcommand{\sfL}{\mathsf{L}}

\newcommand{\sfP}{\mathsf{P}}
\newcommand{\sfQ}{\mathsf{Q}}

\newcommand{\sfT}{\mathsf{T}}

\newcommand{\sfW}{\mathsf{W}}

\newcommand{\sff}{\mathsf{f}}

\newcommand{\sfh}{\mathsf{h}}

\newcommand{\tsfA}{\widetilde{\sfA}}
\newcommand{\tsfB}{\widetilde{\sfB}}

\newcommand{\hsfC}{\widehat{\sfC}}

\newcommand{\hsff}{\widehat{\sff}}

\newcommand{\hsfQ}{\widehat{\sfQ}}

\newcommand{\hsfT}{\widehat{\sfT}}

\newcommand{\hsfW}{\widehat{\sfW}}

\newcommand{\bbA}{\mathbb{A}}

\newcommand{\bbD}{\mathbb{D}}

\newcommand{\bbF}{\mathbb{F}}

\newcommand{\bbZ}{\mathbb{Z}}

\newcommand{\bmu}{{\boldsymbol{\mu}}}

\newcommand{\bfB}{\mathbf{B}}

\newcommand{\bfG}{\mathbf{G}}

\newcommand{\bfL}{\mathbf{L}}

\newcommand{\bfV}{\mathbf{V}}

\newcommand{\bfX}{\mathbf{X}}
\newcommand{\bfY}{\mathbf{Y}}

\newcommand{\hbfB}{\widehat{\bfB}}

\newcommand{\hbfn}{\widehat{\bfn}}
\newcommand{\hbfg}{\widehat{\bfg}}

\newcommand{\tbfG}{\widetilde{\bfG}}

\newcommand{\bfa}{\mathbf{a}}
\newcommand{\bfb}{\mathbf{b}}

\newcommand{\bfg}{\mathbf{g}}

\newcommand{\bfi}{\mathbf{i}}

\newcommand{\bfn}{\mathbf{n}}

\newcommand{\bfp}{\mathbf{p}}
\newcommand{\bfq}{\mathbf{q}}

\newcommand{\bfu}{\mathbf{u}}
\newcommand{\bfv}{\mathbf{v}}

\newcommand{\bfx}{\mathbf{x}}
\newcommand{\bfy}{\mathbf{y}}

\newcommand{\calA}{\mathcal{A}}
\newcommand{\calB}{\mathcal{B}}

\newcommand{\calD}{\mathcal{D}}
\newcommand{\calE}{\mathcal{E}}
\newcommand{\calF}{\mathcal{F}}
\newcommand{\calG}{\mathcal{G}}

\newcommand{\calK}{\mathcal{K}}
\newcommand{\calL}{\mathcal{L}}

\newcommand{\calN}{\mathcal{N}}
\newcommand{\calO}{\mathcal{O}}
\newcommand{\calP}{\mathcal{P}}
\newcommand{\calQ}{\mathcal{Q}}
\newcommand{\calR}{\mathcal{R}}
\newcommand{\calS}{\mathcal{S}}
\newcommand{\calT}{\mathcal{T}}
\newcommand{\calU}{\mathcal{U}}

\newcommand{\calY}{\mathcal{Y}}
\newcommand{\calZ}{\mathcal{Z}}

\newcommand{\hcalB}{\widehat{\calB}}

\newcommand{\hcalE}{\widehat{\calE}}
\newcommand{\hcalF}{{\widehat{\calF}}}
\newcommand{\hcalG}{{\widehat{\calG}}}

\newcommand{\tcalA}{\widetilde{\calA}}

\newcommand{\scrA}{\mathscr{A}}

\newcommand{\scrC}{\mathscr{C}}

\newcommand{\ttB}{\mathtt{B}}

\newcommand{\ttb}{\mathtt{b}}

\newcommand{\ttg}{\mathtt{g}}

\newcommand{\ttq}{\mathtt{q}}

\nc{\bg}{\sigma}

\newcommand{\httB}{\widehat{\ttB}}
\newcommand{\httb}{\widehat{\ttb}}

\newcommand{\hfrakB}{\widehat{\frakB}}

\newlength{\mylength}
\newcommand*{\para}{\rlap{\rotatebox{-30}{\rule[.05ex]{.4pt}{.77em}}}\kern.04em\rlap{\kern.36em\raisebox{0.649519052835em}{\rule{.6em}{.4pt}}}\rule{.6em}{.4pt}\kern-.04em\rotatebox{-30}{\rule[.05ex]{.4pt}{.77em}}}

\newcommand{\hDynkin}{{\widehat{\Dynkin}}}

\newcommand{\tDynkin}{{\widetilde{\Dynkin}}}

\newcommand{\diag}{\mathrm{diag}}

\newcommand{\wl}{\sfP}
\newcommand{\rl}{\sfQ}
\newcommand{\weyl}{\sfW}
\newcommand{\lan}{\langle}
\newcommand{\ran}{\rangle}

\newcommand{\qtq}[1][{and}]{\quad\text{{#1}}\quad}
\newcommand{\qtqh}[1][{and hence}]{\quad\text{{#1}}\quad}
\newcommand{\qtoq}[1][{or}]{\quad\text{{#1}}\quad}

\newcommand{\ee}{\end{enumerate}}
\newcommand{\bitem}{\begin{itemize}}
\newcommand{\eitem}{\end{itemize}}
\newcommand{\ben}{\begin{enumerate}[{\rm (1)}]}
\newcommand{\bnum}{\begin{enumerate}[{\rm (i)}]}
\newcommand{\bnump}{\begin{enumerate}[{\rm (i)$'$}]}
\newcommand{\bna}{\begin{enumerate}[{\rm (a)}]}
\newcommand{\bnA}{\begin{enumerate}[{\rm (A)}]}
\newcommand{\bc}{\begin{cases}}
\newcommand{\ec}{\end{cases}}
\newcommand{\ba}{\begin{array}}
\newcommand{\ea}{\end{array}}

\newenvironment{myequation}
{\relax\setlength{\arraycolsep}{1pt}\begin{eqnarray}}
{\end{eqnarray}}
\newenvironment{myequationn}
{\relax\setlength{\arraycolsep}{1pt}\begin{eqnarray*}}
{\end{eqnarray*}}

\nc{\eqs}[1]{\underset{\raisebox{.4ex}[.7ex][0ex]{$\scriptstyle{#1}$}}{=}}

\nc{\seq}[1]{ \boldsymbol{(} {#1} \boldsymbol{)}   }
\nc{\Bseq}[1]{ \boldsymbol{\left(} {#1} \boldsymbol{\right)}   }

\newcommand{\Dcan}{{\bbD_\can}}

\newcommand{\Ang}[1]{  \bigl\lan #1 \bigr\ran  }
\newcommand{\ang}[1]{\langle#1\rangle}

\newcommand{\longepito}[1][]{\xymatrix@C=4ex{{}\ar@{->>}[r]^{#1}&{}}}

\newcommand{\oprod}{\prod^{\xrightarrow{}}}

\newcommand{\uii}{{\underline{\boldsymbol{\im}}}}

\nc{\ake}[1][2ex]{\rule[-.5ex]{0ex}{#1}}

\newcommand{\ujj}{{\underline{\boldsymbol{\jm}}}}
\newcommand{\TT}{ \textbf{\textit{T}}}

\newcommand{\Red}{\mathrm{Red}}
\newcommand{\Seq}{\mathrm{Seq}}

\newcommand{\tPsi}{{\widetilde{\Psi}}}

\nc{\col}{\colon}
\nc{\ord}{\mathrm{ord}}
\nc{\catQ}{\scrC_\calQ}
\nc{\catD}{\scrC_\bbD}
\nc{\vs}{\vspace*}
\nc{\D}{\mathscr{D}}
\nc{\Proof}{\begin{proof}}
  \nc{\QED}{\end{proof}}

\nc{\tLa}{\widetilde{\Lambda}}
\nc{\hLa}{\widehat{\Lambda}}
\nc{\hSigma}{\widehat{\Sigma}}
\nc{\ro}{{\rm(}}
\nc{\rf}{{\rm)}\xspace}
\nc{\Aut}{\mathrm{Aut}}
\nc{\can}{\mathrm{can}}
\nc{\Dc}{{\bbD_{\can}}}
\nc{\Cgo}{\scrC_{\mathfrak{g}}^{0}}
\nc{\bb}{\mathtt{b}}
\nc{\catC}{\scrC}
\nc{\braid}{\ttB}
\nc{\Ass}[1][\bbD]{\calE_{#1}}
\nc{\Ci}{C^\uii}
\nc{\Cj}{C^\ujj}
\nc{\Di}[1][{\uii}]{{\bbD,{#1}}}
\nc{\Vi}[1][\uii]{V^{#1}}
\nc{\Pii}[1][\uii]{P^{#1}}
\nc{\Vdi}[1][{\Di}]{V^{#1}}
\nc{\Pdi}[1][{\Di}]{P^{#1}} 
\nc{\PBix}[1][{[l,r]}]{\Z_{\ge0}^{\oplus{#1}}}
\nc{\mcf}[1][{$[a,b]$}]{maximal commuting family of $\bfi$-boxes in {#1}\xspace}
\nc{\nn}{\nonumber}
\nc{\bwr}{\mbox{\large$\wr$}}
\newcommand{\tch}{\widetilde{\ch}}
\nc{\tchDcan}{\ms{2mu}\tch_{\mspace{.1mu}\raisebox{-.4ex}{${\scriptstyle{\Dcan}}$}}}
\nc{\monoto}[1][]{\xymatrix@C=2ex{\ar@{>->}[r]^-{{#1}}&}\ms{-8mu}}

\title[Folding of cluster algebras]{Folding of cluster algebras}

\author[W.-S. Jung]{Woo-Seok Jung}
\thanks{The research of W.-S. Jung was supported by the National Research Foundation of Korea (NRF) grants funded by the Korea government (Ministry of Education and MSIT) (RS-2024-00451844 and RS-2023-00273425).}
\address[W.-S. Jung]{Department of Mathematics, University of Seoul, Seoul 02504, South Korea}
\email[W.-S. Jung]{jungws@uos.ac.kr}

\author[S.-j. Oh]{Se-jin Oh}
\thanks{ The research of S.-j.\ Oh was supported by the National Research Foundation of
	Korea (NRF) Grant funded by the Korea government(MSIT) (NRF-2022R1A2C1004045).}
\address[S.-j. Oh]{ Department of Mathematics, Sungkyunkwan University, Suwon, South Korea}
\email[S.-j. Oh]{sejin092@gmail.com}

\date{\today}

\begin{document}

\begin{abstract}
We develop a general framework for folding cluster algebras from
simply-laced Cartan data equipped with admissible finite-group actions.
For arbitrary finite sequences, we establish folding of root data,
compatible pairs, $\bfi$-box seeds, and their individual cluster variables,
without reducedness, finite-type, or acyclicity assumptions.
In the reduced case, the folding is realized through quantum unipotent
coordinate algebras and is compatible with dual PBW generators.
It sends cluster monomials associated with canonical lifted
$\bfi$-box seeds to elements of the classical dual canonical basis
 for the quotient Cartan datum.
For arbitrary positive braids in finite Cartan type, the folding is
realized through quantum virtual Grothendieck rings and, after
specialization at $q=1$, is induced by a single ambient
Grothendieck-ring homomorphism.
Finally, for arbitrary symmetrizable Cartan data and arbitrary
positive braids, we prove that  the Laurent homomorphism identifying corresponding seed variables restricts
to a surjective algebra homomorphism between the corresponding
cluster algebras.
We also extend the Kimura--Oya quantum twist automorphisms
to arbitrary positive braids after localization at the frozen
variables, and prove that their commutative specializations
commute with the folding homomorphisms.
\end{abstract}

\maketitle
\tableofcontents

\section*{Introduction}

Global and canonical bases and their dual bases are fundamental
objects in the representation theory of quantum groups. Their
positivity and integrality properties are accompanied by a more subtle
multiplicative problem: one would like to recognize when products of
distinguished basis elements are again distinguished. Cluster algebras and their quantum analogues were introduced in part
to provide a combinatorial framework for this problem
\cite{FZ02,BZ05}. In coordinate rings arising from Lie theory and
their quantizations, cluster variables and cluster monomials are
expected to be closely related to upper global or dual canonical basis
elements, while exchange relations make part of their multiplicative
structure explicit.

\smallskip

A fundamental model for this interaction is the quantum unipotent
coordinate algebra
\[
\sfA_q(\n(w)),
\]
where $w$ is a Weyl group element. A reduced sequence $\ii$ of $w$
determines dual PBW root vectors, and different reduced
sequences are related by braid moves. The algebra
$\sfA_q(\n(w))$ carries quantum cluster algebra structures
\cite{GLS11,GLS13,GY17,GY21}; its cluster variables include quantum
minors, and its cluster monomials are related to the dual canonical/upper global 
basis \cite{Kimura12,KKKO18}. In symmetric type, modules over quiver
Hecke algebras provide a \emph{monoidal categorification}: Under monoidal categorification, cluster variables correspond to \emph{real
prime} simple modules; the modules attached to a cluster form a \emph{maximal
commuting family}, and mutations are encoded by \emph{short exact sequences}
\cite{KKKO18}. Thus quantum unipotent coordinate algebras form a
setting in which quantum-group-theoretic, cluster-theoretic, and
categorical structures meet.

\smallskip

Recent work extends this picture from Weyl group elements to positive
braids. Bosonic extensions of quantum groups have braid group
actions, and to a symmetrizable Cartan datum $\g$ and a positive braid
$\ttb\in\ttB_\g^+$ one associates a positive-braid subalgebra
\[
\scrA_q(\ttb)
\]
of the corresponding bosonic extension
\cite{HL16,FHOO,JLO2,OP24,KKOP24,KKOP24B,KKOP25M}. 
When $\ttb$ is represented by a reduced expression of
$w\in\sfW_\g$, this construction recovers $\sfA_q(\n(w))$. Thus
positive-braid algebras extend quantum unipotent coordinate algebras
to arbitrary positive braids, whose associated root
sequences may contain repeated or negative roots. In finite Cartan type, they admit
realizations inside quantum (virtual) Grothendieck rings; in simply-laced
type $ADE$, the corresponding realizations are related to quantum
Grothendieck rings of finite-dimensional representations of quantum
loop algebras. Positive braids of finite type also occur naturally in the geometry of
braid varieties and their cluster structures \cite{CGGLSS25}.

\smallskip

To make this representation-theoretic connection more explicit, let
$\scrC_\g$ be the category of finite-dimensional representations of
$U_q(L\g)$.\footnote{Here $t$ denotes the deformation parameter of the quantum
Grothendieck ring; in the main body, this parameter is denoted by
$q$. In simply-laced type, we also identify the usual Grothendieck
ring $\calK(\g)$ with the corresponding virtual Grothendieck ring
$\frakK(\g)$.} The $q$-character theory of Frenkel and Reshetikhin
\cite{FR98,FR99} realizes its Grothendieck ring as a subring of an
ambient Laurent polynomial ring:
\begin{align} \label{eq: in Laurent}
\calK_\Z(\g)\hookrightarrow \Z[y_{i,p}^{\pm1}].
\end{align}
In simply-laced type, Nakajima's quiver varieties give geometric
constructions of $q$-characters and their quantum deformations
\cite{Nak01,Nak04,VV03}. The corresponding quantum Grothendieck rings
$\calK_t(\g)$ have canonical or global
bases given by the classes of simple modules, while distinguished
monoidal subcategories of $\scrC_\g$ provide monoidal
categorifications of cluster algebras
\cite{HL10,HL16,KKOP24A,KKOP25M}. In particular, for the seeds associated
with adapted sequences, the initial cluster variables are given by $q$-characters 
of Kirillov--Reshetikhin modules \cite{HL16,FHOO,JLO1}.

Even when $\g$ is non-simply-laced, however, these cluster algebras $\calK_\Z(\g)$
are of skew-symmetric type and are naturally related to simply-laced
bosonic extensions \cite{HL16,FHOO,B21}. This
simply-laced nature of the usual Grothendieck ring led to the
introduction of the \emph{quantum virtual Grothendieck ring}
$\frakK_q(\g)$ as a genuinely non-simply-laced counterpart
\cite{KO23,JLO1}. Its quantum cluster algebra structure is
skew-symmetrizable and reflects the Cartan datum of $\g$ itself
\cite{JLO1}, as does its braid-group action \cite{JLO2}.

The Laurent realizations~\eqref{eq: in Laurent}  make the corresponding
folding particularly natural.
 Suppose that the Cartan matrix of $\g$ is the quotient Cartan matrix
associated with an admissible finite-group action on the Dynkin
diagram of a simply-laced algebra $\hg$, and write $I_i$ for the
orbit corresponding to vertex $i$ of $\g$
(see Subsection~\ref{subsec: admissible folding}).
At $q=1$, identifying the variables in each Dynkin-diagram orbit defines
\[
\uppi:\C[y_{\im,p}^{\pm1}]
\twoheadrightarrow
\C[y_{i,p}^{\pm1}],
\qquad
\xymatrix@C=4em{
y_{\im,p}\ar@{|->}[r]^{\ |I_i|:1\ }&y_{i,p}
}.
\]
A nontrivial feature of this elementary identification is that it
restricts to a folding homomorphism
\[
\uppi:\calK(\hg)\twoheadrightarrow\frakK(\g).
\]

At the level of heart subrings \cite{HL15,JLO2}, this map is \emph{somewhat} natural
in view of the identifications with the unipotent coordinate rings of
the simply-laced and folded types:
\[
\xymatrix@C=3em@R=0.7em{
N_{\hg} && \C[N_{\hg}] \; \ar@{->>}[dd] \ar@{^{(}->}[r]
& \calK(\hg) \ar@{->>}[dd]^{\uppi} \ar@{^{(}->}[r]
& \C[y_{\im,p}^{\pm1}] \ar@{->>}[dd]^{\uppi} \\
& \ar@{}[r]|{\leadsto \qquad \quad } & \\
N_{\g} \ar@{^{(}->}[uu]
&& \C[N_{\g}] \; \ar@{^{(}->}[r]
& \frakK(\g) \ar@{^{(}->}[r]
& \C[y_{i,p}^{\pm1}]
}
\]
\emph{What is much less apparent is that this elementary folding is
highly compatible with the corresponding cluster structures.}

For a non-simply-laced Cartan datum, analogous geometric or monoidal
structures are not always available at the same level of generality.
This motivates the use of an admissible Dynkin folding 
\[
(\hg,\Gamma)\longrightarrow\g,
\]
where $\hg$ is simply-laced and $\Gamma$ acts admissibly on its Dynkin
diagram. Our broader goal is to understand how the geometric,
representation-theoretic, and categorical structures in simply-laced
types are related to those in non-simply-laced types through folding. 
 We construct surjective homomorphisms from cluster algebras associated with $\hg$ to cluster algebras associated with $\g$.{}
The direction of this motivation may be summarized as
\begin{align} \label{eq: philosophy}
\begin{array}{c}
\text{well-understood simply-laced upstairs structures}\\
\downarrow\ \text{folding}\\
\text{generally less-understood downstairs structures}.
\end{array}
\end{align}
The approach described in \eqref{eq: philosophy} has been
explored in various settings under additional assumptions
\cite{Dupont08,Demonet11,FST12Unfoldings,HL18Unfolding,CHL22,NS14}.
In contrast, our approach applies to arbitrary positive braids
over symmetrizable Cartan data without these restrictions;
see Remark~\ref{rem: comparison with earlier folding constructions}.

\smallskip

At the level of Cartan matrices, admissible folding is elementary. If
$I_i$ is the orbit corresponding to $i\in I$, then (see \S~\ref{subsec: admissible folding})
\[
c_{ij}=\sum_{a\in I_i}\hc_{a,b}
\qquad(b\in I_j),
\qquad
d_i|I_i|=m.
\]
The corresponding problem for \emph{quantum} cluster seeds is more delicate:
one must fold not only the exchange matrices $B$, but also the
compatible skew-symmetric matrices $\Lambda$, together with mutations
and individual cluster variables. A further difficulty is that identifying
the variables in each orbit does not in general define a homomorphism
between the corresponding quantum tori.

One of our main observations is that the Cartan-matrix quotient rule
also governs the root-theoretic construction of the corresponding
quantum cluster seeds.
For every finite sequence $\ii$ and every symmetrizable
Cartan datum $\g$, we
give an \emph{explicit} root-theoretic formula for the
skew-symmetric matrix $\Lambda^\ii$ and prove that
\[
(B^\ii,\Lambda^\ii)
\quad\text{is a compatible pair} \quad
\text{(Proposition~\ref{prop: compatibility arbitrary sequence}).}
\]
For the blockwise lift $\hii$ of $\ii$, the associated root sequence
folds to that of $\ii$, and the bilinear and skew-symmetric pairings
upstairs recover their downstairs counterparts by explicit orbit sums.
Schematically,
\[
\widehat C\longrightarrow C
\quad\Longrightarrow\quad
\bigl(\text{upstairs roots and pairings}\bigr)
\longrightarrow
\bigl(\text{downstairs roots and pairings}\bigr),
\]
and hence
\[
\xymatrix@C=5em{
(B^{\hii},\Lambda^{\hii})
\ar[r]^{\text{folding}}
&
(B^\ii,\Lambda^\ii).
}
\]
By Lemmas~\ref{lem: root folding}, \ref{lem: root pairing folding},
and~\ref{lem: lambda pairing folding}, this holds for arbitrary
sequences $\ii$, without assuming reducedness or finite Cartan type.

More importantly, this construction extends uniformly from individual
positions to $\bfi$-boxes. These boxes were introduced in
\cite{KKOP24B,KK24} as a combinatorial tool for the monoidal
categorification of cluster algebras associated with arbitrary
sequences. A maximal commuting family $\calF$ of $\bfi$-boxes
determines an exchange matrix $B^\calF$. Using the box vectors \eqref{eq: big epsilon box} and the
skew-symmetric pairing $\lambda^\ii$ \eqref{eq: lambda pairing}, we define $\Lambda^\calF$ and
prove that
\[
(B^\calF,\Lambda^\calF)
\quad\text{is a compatible pair} \qquad\text{(Theorem~\ref{thm: ibox compatibility})}
\]
for every maximal commuting family $\calF$. Thus the compatible-pair
construction for $(B^\ii,\Lambda^\ii)$ extends naturally to all
$\bfi$-box seeds.

We then compare this construction under folding.
Let $\frakC$ be an admissible chain of $\bfi$-boxes realizing
$\calF$ (see Subsection~\ref{subsec: ibox compatible pair}). 
The canonical lift
$\hcalF$ of $\calF$ is organized into  fibers over the boxes of $\calF$, and the entire compatible pair folds along these fibers:
\[
\xymatrix@C=5em{
(B^{\hcalF},\Lambda^{\hcalF})
\ar[r]^{\text{folding}}
&
(B^\calF,\Lambda^\calF).
}
\]
More precisely, for a fixed representative
$\hfrakd_t\in\Fib^{\frakC}(\frakd)$,
\[
b^\calF_{\frakc,\frakd}
=
\sum_{\hfrakc\in\Fib^{\frakC}(\frakc)}
b^{\hcalF}_{\hfrakc,\hfrakd_t},
\qquad
\Lambda^\calF_{\frakc,\frakd}
=
d_{i_\frakc}
\sum_{\hfrakc\in\Fib^{\frakC}(\frakc)}
\Lambda^{\hcalF}_{\hfrakc,\hfrakd_t} \quad \text{(Theorem~\ref{thm: ibox pair folding})}.
\]
Thus the quotient Cartan rule governs not only the root data and the
compatible pair associated with the sequence itself, but also the
exchange $B$ and commutation matrices $\La$ of every $\bfi$-box seed.

\smallskip

 The results above compare the compatible pairs
$(B^\hcalF,\Lambda^\hcalF)$ and $(B^\calF,\Lambda^\calF)$
at a fixed $\bfi$-box seed. Sections~\ref{sec: braid}
and~\ref{sec: box} establish the same matrix identities for
the seeds obtained after braid moves and box moves.
Using the mutation--permutation realizations of these moves,
we compare the downstairs seeds before and after each move with
the upstairs seeds before and after the corresponding sequence of moves.
The Laurent homomorphism sending each upstairs variable in a fiber
before the move to the corresponding downstairs variable also sends
every upstairs variable in a fiber after the move to the corresponding
downstairs cluster variable.{} For a box move $\frakC\longmapsto\frakC'$, the individual-variable statement{}
takes the uniform form
\[
\phi^\frakC(\hx'_{\hfrakc'})=x'_{\frakc'}
\qquad
\bigl(
\hfrakc'\in\Fib^{\frakC'}(\frakc')
\bigr) \quad \text{(Theorem~\ref{thm: ibox endpoint folding})}.
\]
The corresponding
$g$-vector statements are Proposition~\ref{prop: g vector folding}
and Corollary~\ref{cor: ibox g-vector folding}.
{}

The two kinds of moves serve complementary purposes. Braid moves
 compare the seeds{} attached to different sequences, whereas box
moves compare maximal commuting families for a fixed sequence. {}

\smallskip

We combine this machinery in Section~\ref{sec: global folding}. For a
positive braid $\ttb$, a sequence $\ii\in\Seq(\ttb)$, and a maximal
commuting family $\calF$ of $\bfi$-boxes, let $\hcalF$ be its
canonical lift.
 The fibers over the boxes of $\calF$ define  a surjective Laurent polynomial homomorphism
\[
\phi^\calF:\calL(\Sigma^\hcalF)\longrightarrow
\calL(\Sigma^\calF),
\qquad
\hx_{\hfrakc}\longmapsto x_\frakc
\quad\bigl(\hfrakc\in\Fib^\frakC(\frakc)\bigr).
\]
Theorem~\ref{thm: global cluster monomial folding} shows that this one map folds all cluster variables
and cluster monomials arising from the box seeds connected by the braid
and box transformations under consideration. The exponents of a
cluster monomial are added along each fiber, and the corresponding
$g$-vectors fold by  the same summation rule.
Explicitly,
\[
\phi^\calF\left(
\prod_{\hfrakd\in\hcalG}
(\hx_{\hfrakd}^{\hcalG})^{a_{\hfrakd}}
\right)
=
\prod_{\frakd\in\calG}
(x_{\frakd}^{\calG})^{
\sum_{\hfrakd\in\Fib^{\frakD}(\frakd)}a_{\hfrakd}} \qquad\text{(Theorem~\ref{thm: global cluster monomial folding})}.
\]

The Laurent phenomenon, together with the individual-variable folding
for $\bfi$-box seeds, naturally suggests that this fixed Laurent
homomorphism should restrict to the corresponding cluster algebras.
Our main result confirms this expectation in a stronger form: in the
commutative setting, the restriction is a surjective algebra
homomorphism.

\begin{thmx}
Fix an admissible Dynkin folding
$(\hg,\Gamma)\to\mathfrak g$.
Let $\ttb\in \ttB_{\g}^{+}$ be an arbitrary positive braid. 
For every maximal commuting family $\calF$ of
$\bfi$-boxes and an admissible chain realizing $\calF$, the fixed
Laurent homomorphism $\phi^\calF$ restricts to a surjective
$\C$-algebra homomorphism
\[
\phi^\calF_\ttb:
\calA(\httb)\twoheadrightarrow\calA(\ttb)
\]
\end{thmx}
\noindent
(Theorem~\ref{thm: cluster algebra folding}, proved in full
generality in Section~\ref{sec: arbitrary positive braid folding}).

\smallskip

\noindent
$\boxed{1}$ We first consider reduced sequences, where the cluster
algebras are realized by quantum unipotent coordinate algebras
$\sfA_q(\n(w))$ and have established structures in quantum group theory and related categories.

\smallskip

\noindent
$\boxed{2}$ We then allow arbitrary positive braids while assuming
that the Cartan datum is of finite type. In this setting, the cluster algebras are realized by the
positive-braid algebras $\scrA_q(\ttb)$ inside quantum virtual
Grothendieck rings and carry interpretations in terms of bosonic
extensions and representation theory.

\smallskip

\noindent
$\boxed{3}$ Finally, we remove both reducedness and the finite Cartan-type assumption.
For an arbitrary symmetrizable Cartan datum and an arbitrary positive
braid, no additional quantum-group realization or categorification is used: 
the cluster algebra is considered on its own. Section~\ref{sec: arbitrary positive braid folding} shows that
the folding machinery developed above is sufficient to prove the
theorem in this full generality.

\medskip
\noindent
\medskip
\noindent
$\boxed{1}^*$\emph{Reduced sequences and quantum unipotent coordinate
algebras.}
Let $w\in\sfW_\g$ and  $\ii=\seq{i_1,\ldots,i_r}\in\Red(w)$.  Its blockwise lift $\hii$ is{}
again reduced and represents an upstairs Weyl group element $\hw$
(Corollary~\ref{cor: reduced lift}).
The associated quantum cluster algebras are realized by
\[
\sfA_q(\n(\hw)) \qquad\text{and}\qquad \sfA_q(\n(w)).
\]
Thus, in the reduced case, the folding homomorphism can be described
in terms of dual PBW generators and dual canonical or upper global
basis elements{}.

The cluster variables associated with singleton $\bfi$-boxes
correspond to the dual PBW generators. Both upstairs and downstairs, these
generators are categorified by cuspidal modules over the corresponding
quiver Hecke algebras (Lemma~\ref{lem: singleton PBW generators}). After
specialization at $q=1$, the fixed Laurent homomorphism
constructed above restricts to a surjective $\C$-algebra homomorphism
\[
\phi^\calF_w: \calA(\hw)\twoheadrightarrow\calA(w).
\]
This restriction is displayed by the commutative diagram
\[
\xymatrix@C=5em@R=1.5em{
\calA(\hw) \ar@{^{(}->}[r] \ar@{->>}[d]_{\phi^\calF_w} & \calL(\Sigma^{\hcalF}) \ar@{->>}[d]^{\phi^\calF}\\
\calA(w) \ar@{^{(}->}[r] & \calL(\Sigma^\calF).
}
\]
Under the polynomial presentations by dual PBW generators, the
restriction is determined by
\[
\begin{aligned}
\calA(\hw)&=\C[\hsff_{k,s}\mid 1\le k\le r,\ 1\le s\le e_{i_k}],\\
\calA(w) &= \C[\sff_1^\ii,\ldots,\sff_r^\ii], \qquad \phi^\calF_w(\hsff_{k,s})=\sff_k^\ii
\end{aligned}
\]
(Proposition~\ref{prop: folding PBW generators} and
Theorem~\ref{thm: reduced cluster algebra folding}).
Thus every individual dual PBW generator in an upstairs fiber 
folds to the corresponding downstairs
generator. Moreover, the folding homomorphism sends the cluster monomials
associated with  canonical lifted {}$\bfi$-box seeds to elements of the classical dual
canonical basis $\bfB^{\rm cl}(w)$
(Corollary~\ref{cor: classical dual canonical basis folding}).

\medskip
\noindent
$\boxed{2}^*$\emph{Arbitrary positive braids in finite Cartan type and
quantum $($virtual$)$ Grothendieck rings.}
We next assume that $\g$ is of finite type, while the positive braid
$\ttb\in\ttB_\g^+$ is arbitrary. The positive-braid algebras
\[
\scrA_q(\httb)
\qquad\text{and}\qquad
\scrA_q(\ttb)
\]
are defined through the braid-group actions on the corresponding
bosonic extensions. 
Put $\bbA=\Z[q^{\pm1/2}]$. 
Their $\bbA$-forms are identified with the quantum
cluster algebras
\[
\scrA_{\bbA}(\httb)\simeq\calA_q(\httb),
\qquad 
\scrA_{\bbA}(\ttb) \simeq\calA_q(\ttb) \quad
\text{(Corollaries~\ref{cor: upstairs cluster algebra bosonic}
and~\ref{cor: positive braid cluster algebra})}.
\]
 The algebra $\scrA_{\bbA}(\httb)$ admits a monoidal categorification
by finite-dimensional representations of the quantum loop algebra
of $\hg$, while $\scrA_{\bbA}(\ttb)$ is realized in the quantum
virtual Grothendieck ring of $\g$.

After specialization at $q=1$, the ambient Grothendieck-ring folding
homomorphism
\[
\uppi:\calK(\hg)\twoheadrightarrow\frakK(\g)
\]
restricts to the positive-braid cluster algebras. More precisely, we
obtain the commutative diagram
\[
\xymatrix@C=4em@R=1.5em{
\calA(\httb)
\ar@{^{(}->}[r]
\ar@{->>}[d]_{\uppi_{\ttb}}
&
\calK(\hg)
\ar@{^{(}->}[r]
\ar@{->>}[d]^{\uppi}
&
\C[y_{\im,p}^{\pm1}]
\ar@{->>}[d]^{\uppi}
\\
\calA(\ttb)
\ar@{^{(}->}[r]
&
\frakK(\g)
\ar@{^{(}->}[r]
&
\C[y_{i,p}^{\pm1}].
}
\]
Thus the folding homomorphism for every finite-type positive braid is
the restriction of the single global homomorphism $\uppi$.

Moreover, for every sequence $\ii\in\Seq(\ttb)$ and every maximal
commuting family $\calF$ of $\bfi$-boxes, this restriction agrees with
 $\phi^\calF_\ttb$:
\[
\uppi_\ttb=\phi^\calF_\ttb \qquad\text{(Theorem~\ref{thm: surjectivity arbitrary sequences})}.
\]
Hence the folding homomorphisms arising from all choices of sequences
and maximal commuting families are induced by the same ambient map
$\uppi$.

 The same realization shows that cluster monomials associated with
canonical lifted $\bfi$-box seeds map to elements of the classical global basis of $\frakK(\g)${}
(Corollary~\ref{cor: classical global basis folding}).  The homomorphism $\uppi${} is also compatible with the
substitution isomorphisms  on the ambient rational function fields associated with Dynkin quivers: the upstairs substitution formulas specialize under
$\uppi$ to their downstairs counterparts
(Corollary~\ref{cor: substitution folding}). Thus, in finite type, the
cluster-theoretic folding is realized by a global homomorphism already
present in the quantum virtual Grothendieck-ring framework.

\medskip
\noindent
$\boxed{3}^*$\emph{Arbitrary positive braids over symmetrizable Cartan
data.}
Finally, let both the symmetrizable Cartan datum $\g$ and the positive
braid $\ttb$ be arbitrary. {}

Proposition~\ref{lem: Qin box interval identification} identifies
the cluster variables attached to $\bfi$-boxes with the interval
variables of Qin.
Qin's standard-monomial basis then gives the polynomial presentation
in Proposition~\ref{prop: singleton polynomial presentation}.
Write $\ii=\seq{i_1,\ldots,i_l}$ and $e_i=|I_i|$.
For $1\le k\le l$ and $1\le r\le e_{i_k}$, let $\kappa(k,r)$
be the position in $\hii$ of the $r$th letter replacing $i_k$.
Applying this proposition to $\ii$ and $\hii$, we obtain
\[
\calA(\ttb)=\C[x_{[1]},\ldots,x_{[l]}],
\qquad
\calA(\httb)=
\C[\hx_{[\kappa(k,r)]}\mid
1\le k\le l,\ 1\le r\le e_{i_k}].
\]
The folding homomorphism sends each upstairs singleton variable to
the corresponding downstairs singleton variable, as shown in
\eqref{eq: singleton generator folding}:
\[
\phi^\calF\bigl(\hx_{[\kappa(k,r)]}\bigr)=x_{[k]}.
\]
Consequently, the fixed Laurent homomorphism restricts to the entire
cluster algebra, completing the commutative diagram
\[
\xymatrix@C=5em@R=1.5em{
\calA(\httb)
\ar@{^{(}->}[r]
\ar@{->>}[d]_{\phi^\calF_\ttb}
&
\calL(\Sigma^{\hcalF})
\ar@{->>}[d]^{\phi^\calF}
\\
\calA(\ttb)
\ar@{^{(}->}[r]
&
\calL(\Sigma^\calF).
}
\]
This proves Theorem~\ref{thm: cluster algebra folding} without
reducedness or finite-type assumptions.

\smallskip

As a further application of the preceding framework, let
$\tcalA_q(\ttb)$ denote the localization of $\calA_q(\ttb)$ at its
frozen variables. In Section~\ref{sec: arbitrary positive word DT twists},
we explicitly construct automorphisms $\tau_{\ttb,q}$ of
$\tcalA_q(\ttb)$ for arbitrary positive braids $\ttb$
over symmetrizable Cartan data, using the quantum cluster
variables attached to $\bfi$-boxes and the skew-symmetric
pairing on box vectors.
These automorphisms recover the quantum twist automorphisms
of Kimura--Oya \cite[Theorem~6.1]{KimuraOya21}
when $\ttb$ is the positive braid lift of a Weyl group element.
The automorphism $\tau_{\ttb,q}$ is a twist automorphism of
Donaldson--Thomas type in the sense of Kimura--Qin--Wei
\cite{KimuraQinWei} and permutes the common triangular basis of
$\tcalA_q(\ttb)$ in the sense of
\cite[Definitions~6.2.1--6.2.2]{Qin20}
(Corollary~\ref{prop: positive braid twist common triangular basis}).
Let $\tcalA(\ttb)$ denote the localization of $\calA(\ttb)$
at its frozen variables, and let $\tau_\ttb$ denote the automorphism
of $\tcalA(\ttb)$ obtained by specializing $\tau_{\ttb,q}$ at $q^{1/2}=1$.
The folding homomorphism $\phi^\calF_\ttb$ extends to a surjective homomorphism
$\phi^{\calF,\loc}_\ttb:
\tcalA(\httb)\twoheadrightarrow\tcalA(\ttb)$ satisfying
\[
\phi^{\calF,\loc}_\ttb\circ\tau_{\httb}
=\tau_\ttb\circ\phi^{\calF,\loc}_\ttb
\qquad
\text{(Theorem~\ref{thm: arbitrary positive word DT folding compatibility}).}
\]

The folding homomorphisms constructed here also suggest a natural
direction for further study. It would be interesting to construct
their categorical realizations, possibly using module categories over
quiver Hecke algebras, so that the induced maps on Grothendieck rings
recover these surjective homomorphisms. Such realizations could provide
a mechanism for relating the well-established simply-laced structures
to their counterparts for the quotient Cartan datum. More broadly, the
folding machinery developed here provides a starting point for studying
  global and canonical{} bases, and their categorical realizations, under
folding.

\smallskip

We conclude with the organization of the paper.
Sections~\ref{sec: Cartan data}--\ref{sec: compatible pair folding}
develop the folding of Cartan and root data and the compatible-pair
machinery for $\bfi$-box seeds. Sections~\ref{sec: braid} and
\ref{sec: box} establish its compatibility with braid moves, box
moves, and individual cluster variables.
 Section~\ref{sec: global folding} combines these results to prove
the folding theorem for cluster monomials.

Sections~\ref{sec: reduced} and \ref{sec: finite case} give,
respectively, the quantum-unipotent and finite-type positive-braid
interpretations. Section~\ref{sec: arbitrary positive braid folding}
proves the general commutative folding theorem. 
Section~\ref{sec: arbitrary positive word DT twists} constructs
quantum twist automorphisms for arbitrary positive braids, proves
that they preserve the common triangular bases, and establishes
compatibility with folding after commutative specialization.
Finally,
Appendix~\ref{app: affine braid folding} contains the additional
  calculations for braid moves with orbit sizes $(3,3)$ and $(4,2)$. 

\section*{Acknowledgements}

The authors thank Masaki Kashiwara, Euiyong Park, and Myungho Kim
for helpful discussions. 
They also thank George Lusztig and Nathan Reading for their helpful comments.
They used SageMath~\cite{SageMath}
for computations in the examples.

\subsection*{Convention}
Throughout this paper, we use the following conventions.
\bitem
\item For a statement $P$, we set $\delta(P)$ to be $1$ or $0$ depending on whether $P$ is true or not. In particular, $\delta_{i,j}=\delta(i =j)$
is the Kronecker delta.
\item For integers $a\le b$, we write $[a,b]=\{a,a+1,\ldots,b\}$.
We understand $[a,b] = \emptyset$ if $b<a$ and set $[a]\seteq [a,a]$.
\eitem

\section{Cartan data and admissible Dynkin foldings} \label{sec: Cartan data}
In this section, we fix the notation for Cartan data and admissible Dynkin foldings that will be used throughout this paper.

\subsection{Symmetrizable Cartan data}

Let $I$ be a finite index set. A \emph{Cartan datum} $(\sfC,\wl,\Pi,\wl^\vee,\Pi^\vee)$ consists of \bna
\item a symmetrizable Cartan matrix
$\sfC=(c_{ij})_{i,j\in I}$; i.e., $\sfD\sfC$ is symmetric for a diagonal matrix $\sfD={\rm diag}(d_i \in \Z_{\ge1} \ | \ i \in I)$,
\item  a free abelian group $\wl$, called the \emph{weight lattice},
\item $\Pi=\{ \al_i \ | \ i \in I \}$, called the set of \emph{simple roots},
\item  $\Pi^\vee=\{ h_i \ | \ i \in I \} \subset \wl^\vee \seteq \Hom(\wl,\Z)$ , called the set of \emph{simple coroots},
\item a $\Q$-valued symmetric bilinear form $( \cdot,\cdot)$ on $\wl$,
\ee
satisfying the standard property (see \cite[\S 1.1]{KKKO18} for instance). Here we take $\sfD$ such that
\begin{align} \label{eq: d_i and pairing}
d_i=(\al_i,\al_i)/2,
\qquad
(\al_i,\al_j)=d_i c_{ij}
\qquad (i,j\in I),
\end{align}
and normalize it by $\gcd_{i\in I}(d_i)=1$.

For $i \in I$, we choose $\varpi_i \in \wl$ such that $\lan h_i,\varpi_j \ran =\delta_{ij}$ for any $j \in I$ and call it the \emph{$i$-th fundamental weight.} We denote by
$\wl^+ \seteq \{\lambda\in \wl\mid \langle h_i,\lambda\rangle\geq 0 \ \text{ for all } i\in I \}$ the set of \emph{dominant integral weights}.
The free abelian group $\rl \seteq \bigoplus_{i \in I} \Z \al_i$ is called the \emph{root lattice} and we set $\rl^+ = \sum_{i \in I}\Z_{\ge 0} \al_i \subset \rl$ and $\rl^- = \sum_{i \in I}\Z_{\le 0} \al_i \subset \rl$. We denote by $\Phi$ the set of roots and by $\Phi^\pm$ the set of \emph{positive} roots (resp. \emph{negative} roots). We also denote by $\Phi_{{\rm re}}$ the set of \emph{real} roots.
We set $\Phi_{\re}^{\pm}:=\Phi_{\re}\cap\Phi^\pm$.

Let $\g$ be the Kac-Moody algebra associated with the Cartan datum $(\sfC,\wl,\Pi,\wl^\vee,\Pi^\vee)$.
We say that $\g$ is \emph{simply-laced} if $c_{ij}<0$ implies $c_{ij}=c_{ji}=-1$ for all $i\ne j$.
We denote by $\weyl$
the \emph{Weyl group} of $\g$. It is generated by the simple reflections $s_i\in\Aut(\wl)$ $(i \in I)$ defined by $s_i(\la) = \la - \lan h_i,\la \ran \al_i$ for
$\la \in \wl$.
In particular,
\begin{align} \label{eq: s_i action on simple}
s_i(\al_j)=\al_j-c_{ij}\al_i
\qquad (i,j\in I).
\end{align}

We denote by $\ttB_\g = \Ang{r_i^{\pm 1}}_{i \in I}$ the braid group associated with $\g$; i.e, it is the group generated by $\{ r_i^{\pm 1} \}_{i \in I}$
subject to the following relations:
\begin{align} \label{eq: braid relation}
{\rm (i)} \ r_{i}r_{i}^{-1} = r_{i}^{-1}r_{i}=1 \qtq
{\rm (ii)} \ \underbrace{r_{i}r_{j}\cdots}_{m_{i,j}\text{-times}} =\underbrace{r_{j}r_{i}\cdots}_{m_{i,j}\text{-times}}
\text{ for } i \ne j \in I,
\end{align}
where
\begin{align*}
m_{i,j} \seteq \bc
c_{i,j}c_{j,i} +2 & \text{ if } c_{i,j}c_{j,i} \le 2, \\
6  & \text{ if } c_{i,j}c_{j,i} =3, \\
\infty & \text{ otherwise.}
\ec
\end{align*}
We call the relations in~\eqref{eq: braid relation}~{\rm (ii)} the \emph{braid relations} of $\ttB_\g$.
The Weyl group $W_\g$ is the quotient of $\ttB_\g$ obtained by imposing
$r_i^2=1$ for all $i\in I$; we denote the image of $r_i$ in $W_\g$ by $s_i$.

We denote by $\ttB^+_\g$ the \emph{positive braid monoid} generated by
$\{r_i\}_{i\in I}$.  There is a natural surjective monoid homomorphism
\[
{\rm pr}_\sfW:\ttB^+_\g\longrightarrow \sfW_\g, \qquad r_i\longmapsto s_i.
\]

For a sequence $\ii=(i_1,\ldots,i_r)\in I^r$, we write
\[
r_\ii:=r_{i_1}\cdots r_{i_r}\in\ttB^+_\g, \qquad s_\ii:=s_{i_1}\cdots s_{i_r}\in W_\g.
\]
For $w\in\weyl$ and $\ttb\in\ttB^+_\g$, set
\begin{align*}
\Red(w) &\seteq \bigl\{\ii=\seq{i_1,\ldots,i_r}\in I^r \mid s_\ii=w,\ r=\ell(w)\bigr\} \quad \text{ and} \allowdisplaybreaks\\
\Seq(\ttb) &\seteq \bigl\{\ii=\seq{i_1,\ldots,i_r}\in I^r \mid r_\ii=\ttb\bigr\}.
\end{align*}
Thus $\Red(w)$ is the set of reduced sequences of $w$, while
$\Seq(\ttb)$ is the set of all finite sequences realizing $\ttb$.
Thus, an arbitrary sequence $\ii$ may always be regarded as an element
$r_\ii$ of the positive braid monoid $\ttB^+_\g$.
If $\ii$ is reduced, then $r_\ii$ depends only on $s_\ii$.
Thus, in this case, we shall identify $s_\ii\in W_\g$ with its
positive braid lift $r_\ii\in\ttB^+_\g$ whenever no confusion can arise.
We usually drop the subscript $\g$ from the notation when there is no
danger of confusion.

\subsection{Admissible Dynkin foldings and quotient Cartan data}
\label{subsec: admissible folding}

Let $\g$ be a symmetrizable Kac-Moody algebra. We denote by $\Dynkin$ the \emph{Dynkin diagram} of $\g$.

Let $\Gamma$ be a finite group acting on $\Dynkin$ by diagram automorphisms, i.e.,
$$\text{ $c_{g(i),g(j)}=c_{ij}$ for all $g\in \Gamma$ and $i,j\in I$.}$$
We call the action of $\Gamma$ \emph{admissible} if
\[
c_{i,g(i)}=0
\qquad
\text{for all $i\in I$ and $g\in \Gamma$ with $g(i)\ne i$}.
\]

From now on, let $\hg$ be a simply-laced Kac--Moody algebra
with index set $\hI$ and Cartan matrix
$\hsfC=(\hc_{a,b})_{a,b\in\hI}$.
We normalize the bilinear form~\eqref{eq: d_i and pairing} of $\hg$ so that $(\al_i,\al_i)_\hg=2$ for all $i \in \hI$.
Assume that a finite group $\Gamma$ acts admissibly on the Dynkin diagram
of $\hg$.

We identify the orbit set $\hI/\Gamma$ with an index set $I$ and
denote by
\[
\pi:\hI\longrightarrow I
\]
the canonical projection.  For $i\in I$, we set
\begin{align} \label{eq: I_i and e_i}
I_i:=\pi^{-1}(i)
\qtq
e_i:=|I_i|.
\end{align}

\begin{definition}\label{def: admissible quotient Cartan datum}
The \emph{quotient Cartan matrix} is the matrix \begin{align} \label{eq: quotient Cartan matrix}
\sfC=(c_{ij})_{i,j\in I}
\quad
\text{ defined by } \ \
c_{ij} \seteq \sum_{a\in I_i}\hc_{a,b},
\qquad b\in I_j.
\end{align}
\end{definition}

\begin{lemma}
\label{lem: quotient Cartan}
The matrix $\sfC=(c_{ij})_{i,j\in I}$ in
\eqref{eq: quotient Cartan matrix} is a symmetrizable Cartan matrix.
More precisely, there exist positive integers $d_i$ $(i\in I)$
with
$\gcd_{i\in I}\{d_i \}=1$
such that
\[
d_i c_{ij}=d_j c_{ji} \qquad (i,j\in I).
\]
Moreover, there exists a positive integer $m$ such that
\begin{align}
\label{eq: d_i e_i constant}
d_i e_i=m
\qquad (i\in I).
\end{align}
\end{lemma}

\begin{proof}
The definition of $c_{ij}$ is independent of the choice of
$b\in I_j$, since $\Gamma$ acts transitively on each orbit $I_j$ and preserves $\hsfC$.

By admissibility, for $b\in I_i$,
$c_{ii} = \sum_{a\in I_i}\hc_{a,b}  = \hc_{b,b} =
2$. For $i\ne j$, we have $c_{ij}\in\Z_{\le0}$, and clearly
$c_{ij}=0$ if and only if $c_{ji}=0$.  Hence $\sfC$ is a generalized Cartan matrix.

Since $\hsfC$ is symmetric, counting the edges between $I_i$ and $I_j$ in two ways gives
\begin{align}
\label{eq: orbit-size-cartan}
e_j(-c_{ij})=e_i(-c_{ji})
\iff
e_jc_{ij}=e_ic_{ji}.
\end{align}
Thus $\sfC$ is symmetrizable.

Since $(L/e_i)_{i\in I}$, where
$L=\operatorname{lcm}\{e_i\mid i\in I\}$, is an integral symmetrizer,
set $g=\gcd_{i\in I}(L/e_i)$, $d_i=(L/e_i)/g$, and $m=L/g$.
Then $\gcd_{i\in I}d_i=1$ and
$d_i e_i=m$ $(i\in I)$.
\end{proof}

\begin{remark}\label{rem: every symmetrizable datum is a quotient}
Lusztig proved that every Cartan datum arises by folding a finite
graph along the orbits of an admissible automorphism
\cite[Proposition~14.1.2]{LusztigBook}.
For any symmetrizable Cartan matrix $\sfC$, applying Salisbury and
Scrimshaw's modification of this construction
\cite[Proposition~3.2]{SS17} to $\sfC^{{\rm tr}}$ gives a
simply-laced Dynkin diagram with an admissible automorphism whose
quotient Cartan matrix is $\sfC$ under~\eqref{eq: quotient Cartan matrix}.
\end{remark}

In what follows, we mainly consider the pair $(\hg,\g)$ arising from
an admissible $\Gamma$-action on $\hg$, where $\g$ denotes the Kac--Moody
algebra associated with the quotient Cartan matrix $\sfC$ in
\eqref{eq: quotient Cartan matrix}.

\subsection{Examples} \label{subsec: Examples}

In this subsection, we collect several examples of pairs $(\hg,\g)$ arising from admissible Dynkin foldings, which will be used throughout
the paper.  We describe them mainly by their Dynkin diagrams, together with the $\Gamma$-orbits and the corresponding quotient diagrams.  Besides
the standard finite and affine foldings, we include examples illustrating non-cyclic group actions, fibers of larger cardinality, and indefinite
Cartan data.

We first consider the classical folding pair
$(A_{2n-1},C_n)$.
Let $\sigma$ be the involutive diagram automorphism of $A_{2n-1}$ given by
\[
\sigma(i)=2n-i \qquad (1\le i\le 2n-1).
\]
Then the corresponding admissible folding is
\begin{align}
A_{2n-1} \ \ \xymatrix@R=0.5ex@C=5ex{ *{\circled{2}}<3pt> \ar@{-}[r]_<{ 1 \ \  } \ar@/^1.5pc/@{<.>}[rrrr]|{\sigma}
& *{\circled{2}}<3pt> \ar@{.}[r]_<{ 2 \ \  } \ar@/^1pc/@{<.>}[rr]|{\sigma}
& *{\circled{2}}<3pt> \ar@{}[r]_<{ n  \ \ } \ar@{.}[r]
&*{\circled{2}}<3pt> \ar@{-}[r]_>{ \ \ 2n-1} &*{\circled{2}}<3pt> \ar@{}[l]^>{   2n-2} }
\overset{/\ang{\sigma}}{\Longrightarrow}
\quad
C_{n} \ \  \xymatrix@R=0.5ex@C=3.5ex{ *{\circled{2}}<3pt> \ar@{-}[r]_<{ 1 \ \  } & *{\circled{2}}<3pt> \ar@{-}[r]_<{ 2 \ \  } & *{\circled{2}}<3pt> \ar@{.}[r]_<{ 3  \ \  }
&*{\circled{2}}<3pt> \ar@{<=}[r]_>{ \ \ \raisebox{-0.16cm}{\scriptsize $n$}} &*{\bcircled{4}}<3pt> \ar@{}[l]^>{   n-1} }.
\end{align}
Here $ \xymatrix@R=0.5ex@C=3.5ex{    *{\circled{t}}<3pt>}_k$ means that $(\al_k,\al_k)=t$ and dotted arrows indicate the action of the
diagram automorphism $\sigma$. In this case, the $\Gamma=\ang{\sigma}$-orbits are
$I_i=\{i,2n-i\}$ $(1\le i\le n-1)$,
and  $I_n=\{n\}$,
so that
\[
(e_1,\ldots,e_n)=(2,\ldots,2,1).
\]
The quotient Cartan matrix is of type $C_n$, with normalized symmetrizer
$(d_1,\ldots,d_n)=(1,\ldots,1,2)$.
In particular, $d_i e_i=2$ $(1\le i\le n)$.

The other standard finite-type admissible foldings are depicted below:
\begin{align*}
 D_{n+1}\ \
\raisebox{1em}{
\xymatrix@R=2ex@C=4ex{
&&& *{\circled{2}}<3pt>    \ar@{-}[d]
\ar@/^.6pc/@{<.>}[dr]|{\sigma} \ar@{}[r]^<{n \ \ } & \\
*{\circled{2}}<3pt>  \ar@{-}[r]_<{1\ \ }
& *{\circled{2}}<3pt> \ar@{-}[r]_<{2\ \ }
& *{\circled{2}}<3pt> \ar@{.}[r]
& *{\circled{2}}<3pt> \ar@{-}[r]_>{\ \ n+1}
& *{\circled{2}}<3pt> \ar@{}[l]^>{n-1}
}}
\quad
&\overset{/\ang{\sigma}}{\Longrightarrow}
\quad
B_n\ \
\xymatrix@R=0.5ex@C=3.5ex{
*{\bcircled{4}}<3pt>\ar@{-}[r]_<{1\ \ }
& *{\bcircled{4}}<3pt> \ar@{-}[r]_<{2\ \ }
& *{\bcircled{4}}<3pt> \ar@{.}[r]_<{3\ \ }
& *{\bcircled{4}}<3pt>  \ar@{=>}[r]_>{\ \ \raisebox{-0.16cm}{\scriptsize $n$}}
& *{\circled{2}}<3pt> \ar@{}[l]^>{n-1}
},
\allowdisplaybreaks\\
 E_6\ \
\raisebox{1em}{ \xymatrix@R=2ex@C=3.5ex{
&& *{\circled{2}}<3pt> \ar@{-}[d]^<{2\ \ }\\
*{\circled{2}}<3pt> \ar@{-}[r]_<{1\ \ }
      \ar@/_1.5pc/@{<.>}[rrrr]|{\sigma}
& *{\circled{2}}<3pt> \ar@{-}[r]_<{3\ \ }
      \ar@/_1pc/@{<.>}[rr]|{\sigma}
& *{\circled{2}}<3pt> \ar@{-}[r]_<{4\ \ }
& *{\circled{2}}<3pt> \ar@{-}[r]_>{\ \ 6}
& *{\circled{2}}<3pt> \ar@{}[l]^>{5}
}}
\quad
& \overset{/\ang{\sigma}}{\Longrightarrow}
\quad
F_4\ \
\xymatrix@R=0.5ex@C=3.5ex{
*{\bcircled{4}}<3pt> \ar@{-}[r]_<{1\ \ }
& *{\bcircled{4}}<3pt> \ar@{=>}[r]_<{2\ \ }
& *{\circled{2}}<3pt> \ar@{-}[r]_>{\ \ 4}
& *{\circled{2}}<3pt>
      \ar@{}[l]^>{3}
},
\allowdisplaybreaks\\
 D_4\ \
\raisebox{1em}{
\xymatrix@R=2ex@C=4.5ex{
& *{\circled{2}}<3pt> \ar@{-}[d]_<{3}
      \ar@/^.6pc/@{<.}[dr]|{\sigma}
      \ar@/_.6pc/@{.>}[dl]|{\sigma}
\\
*{\circled{2}}<3pt>
      \ar@{-}[r]_<{1\ \ }
      \ar@/_1.5pc/@{.>}[rr]|{\sigma}
&
*{\circled{2}}<3pt>
      \ar@{-}[r]_>{\ \ 4}
&
*{\circled{2}}<3pt>
      \ar@{}[l]^>{2}
}}
\quad
&\overset{/\ang{\sigma}}{\Longrightarrow}
\quad
G_2\ \
\xymatrix@R=0.5ex@C=3.5ex{
*{\circled{2}}<3pt>
      \ar@{-}[r]_<{1\ \ }
&
*{\rcircled{6}}<3pt>
      \ar@{=>}[l]^<{\ \ 2}
}.
\end{align*}
For later reference, the orbit fibers and their cardinalities in the
above three foldings are summarized as follows:
\[
\begin{array}{c|c|c}
(\hg,\g) & \text{orbit fibers } I_i & (e_i)_{i\in I}
\\ \hline
(D_{n+1},B_n)
&
\begin{array}{l}
I_i=\{i\}\quad(1\le i\le n-1),\\
I_n=\{n,n+1\}
\end{array}
&
(1,\ldots,1,2)
\\[3ex]
(E_6,F_4)
&
\begin{array}{ll}
I_1=\{2\}, & I_2=\{4\},\\
I_3=\{3,5\}, & I_4=\{1,6\}
\end{array}
&
(1,1,2,2)
\\[3ex]
(D_4,G_2)
&
\begin{array}{l}
I_1=\{1,3,4\},\\
I_2=\{2\}
\end{array}
&
(3,1)
\end{array}
\]

Several affine admissible foldings which will be useful later are
depicted below:
\begin{align*}
D_4^{(1)}\ \
\raisebox{1.5em}{
\xymatrix@R=1ex@C=4.2ex{
&
*{\circled{2}}<3pt>
   \ar@{-}[d]_>{1}
   \ar@/_0.5pc/@{.>}[dl]|{\sigma}
&
\\
*{\circled{2}}<3pt>
   \ar@{-}[r]_<{3\ \ }
   \ar@/_0.5pc/@{.>}[rr]|{\sigma}
&
*{\circled{2}}<3pt>
   \ar@{-}[r]_>{\ \ 2}
   \ar@{-}[d]
&
*{\circled{2}}<3pt>
   \ar@/_0.5pc/@{.>}[ul]_>{4}|{\sigma}
\\
&
*{\circled{2}}<3pt>
   \ar@{}^>{0}
&
}}
\quad
&\overset{/\ang{\sigma}}{\Longrightarrow}
\quad
G_2^{(1)}\ \
\xymatrix@R=0.5ex@C=3.8ex{
*{\circled{2}}<3pt>
   \ar@{<=}[r]_<{2\ \ }
&
*{\rcircled{6}}<3pt>
   \ar@{-}[r]_<{1\ \ } \ar@{-}[l]
&
*{\rcircled{6}}<3pt>
   \ar@{}[l]^>{\qquad \qquad 0}
},
\\
E_6^{(1)}\ \
\raisebox{2.1em}{
\xymatrix@R=1ex@C=4.2ex{
&& *{\circled{2}}<3pt>
\ar@{}[r]_<{ \ 0 }
   \ar@{-}[d] \ar@/^0.5pc/@{.>}[ddrr]|{\sigma} & \\
&& *{\circled{2}}<3pt> \ar@{-}[d]
   \ar@{-}[d] \ar@/^0.5pc/@{.>}[dr]|{\sigma}
   \ar@{}[r]_<{ \ 2 } & \\
*{\circled{2}}<3pt>\ar@{-}[r]_<{  1 \ }
\ar@/^0.5pc/@{.>}[uurr]|{\sigma}
& *{\circled{2}}<3pt> \ar@{-}[r]_<{ 3 \ }
\ar@/^0.5pc/@{.>}[ur]|{\sigma}
& *{\circled{2}}<3pt> \ar@{-}[r]_<{ 4 \ } & *{\circled{2}}<3pt> \ar@{-}[r]_<{ 5 \ }
\ar@/^1pc/@{.>}[ll]|{\sigma}
& *{\circled{2}}<3pt> \ar@{}[l]^<{ \ 6  }
\ar@/^2pc/@{.>}[llll]|{\sigma}
}}
\quad
&\overset{/\ang{\sigma}}{\Longrightarrow}
\quad
D_4^{(3)}\ \
\xymatrix@R=0.5ex@C=3.8ex{
*{\circled{2}}<3pt>
   \ar@{-}[r]_<{0\ \ }
&
*{\circled{2}}<3pt>
   \ar@{<=}[r]_<{1\ \ }
&
*{\rcircled{6}}<3pt>
   \ar@{}[l]^>{\qquad \qquad 2}
   \ar@{-}[l]
},
\\
D_{2n+2}^{(1)}\ \
\raisebox{1.5em}{
\xymatrix@R=1ex@C=3.8ex{
*{\circled{2}}<3pt>
   \ar@{-}[dr]
   \ar@{}[r]^<{0\ \ }
&&&&&&
*{\circled{2}}<3pt>
   \ar@{-}[dl]
   \ar@/_0.3pc/@{<.}[llllll]|{\sigma}
   \ar@/_2.3pc/@{.>}[ddllllll]|{\sigma}
      \ar@{}[l]_>{ \qquad \qquad 2n+1}
\\
&
*{\circled{2}}<3pt>
   \ar@{-}[r]_<{2\ \ }
   \ar@/^1pc/@{<.>}[rrrr]|{\sigma}
&
*{\circled{2}}<3pt>
   \ar@{-}[r]_<{3\ \ }
   \ar@/_1pc/@{<.>}[rr]|{\sigma}
&
\cdots
   \ar@{-}[r]_<{\qquad \qquad 2n-1 }
&
*{\circled{2}}<3pt>
   \ar@{-}[r]_>{\ \ 2n}
&
*{\circled{2}}<3pt>
   \ar@{-}[ur]
   \ar@{-}[dr]
\\
*{\circled{2}}<3pt>
   \ar@{-}[ur]
   \ar@{}[r]_>{1 \qquad \qquad }
&&&&&&
*{\circled{2}}<3pt>
   \ar@{-}[ul]
   \ar@{}[l]^>{ \qquad \qquad 2n+2}
   \ar@/^0.3pc/@{<.}[llllll]|{\sigma}
   \ar@/^2.3pc/@{.>}[uullllll]|{\sigma}
}}
\quad
&\underset{n\ge1}{\overset{/(C_4\simeq\ang{\sigma})}{\Longrightarrow}}
\
A_{2n}^{(2)}\
\xymatrix@R=0.5ex@C=3.4ex{
*{\circled{2}}<3pt>
   \ar@{<=}[r]_<{0\ \ }
&
*{\bcircled{4}}<3pt>
   \ar@{-}[r]_<{1\ \ }
&
\cdots
   \ar@{-}[r]
&
*{\bcircled{4}}<3pt>
   \ar@{<=}[r]_>{\ \ n}
&
*{\mcircled{8}}<3pt>
   \ar@{}[l]^>{n-1}
}
\\
D_n^{(1)}\ \
\raisebox{1.5em}{
\xymatrix@R=1ex@C=3.8ex{
*{\circled{2}}<3pt>
   \ar@{-}[dr]
   \ar@{<.>}[dd]|{\sigma_1}
   \ar@{}[r]^<{0\ \ }
&&&&&&
*{\circled{2}}<3pt>
   \ar@{-}[dl]
   \ar@{<.>}[dd]|{\sigma_2}
      \ar@{}[l]_>{ \qquad \qquad n-1}
\\
&
*{\circled{2}}<3pt>
   \ar@{-}[r]_<{2\ \ }
&
*{\circled{2}}<3pt>
   \ar@{-}[r]_<{3\ \ }
&
\cdots
   \ar@{-}[r]_<{\qquad \qquad n-3 }
&
*{\circled{2}}<3pt>
   \ar@{-}[r]_>{\ \ n-2}
&
*{\circled{2}}<3pt>
   \ar@{-}[ur]
   \ar@{-}[dr]
\\
*{\circled{2}}<3pt>
   \ar@{-}[ur]
   \ar@{}[r]_>{1 \qquad \qquad }
&&&&&&
*{\circled{2}}<3pt>
   \ar@{-}[ul]
   \ar@{}[l]^>{ \qquad \qquad n}
}}
\quad
&\underset{n\ge4}{\overset{/\Gamma}{\Longrightarrow}}
\quad
D_{n-1}^{(2)}\ \
\xymatrix@R=0.5ex@C=3.4ex{
*{\circled{2}}<3pt>
   \ar@{<=}[r]_<{0\ \ }
&
*{\bcircled{4}}<3pt>
   \ar@{-}[r]_<{1\ \ }
&
\cdots
   \ar@{-}[r]
&
*{\bcircled{4}}<3pt>
   \ar@{=>}[r]_>{\ \  n-2  }
&
*{\circled{2}}<3pt>
   \ar@{}[l]^>{  n-3  }
}.
\end{align*}
Unlike the preceding three examples, the last folding is induced by a
non-cyclic group action.  More precisely,
\[
\Gamma=\ang{\sigma_1,\sigma_2}\simeq C_2\times C_2,
\]
and the orbit fibers are
\[
I_0=\{0,1\},\qquad
I_i=\{i+1\}\quad(1\le i\le n-3),\qquad
I_{n-2}=\{n-1,n\}.
\]
Hence
\[
(e_0,\ldots,e_{n-2})=(2,1,\ldots,1,2),\qquad
(d_0,\ldots,d_{n-2})=(1,2,\ldots,2,1),
\]
so that $m=2$ in \eqref{eq: d_i e_i constant}.

Notice that in the folding
$E_6^{(1)}\to D_4^{(3)}$, the two adjacent quotient vertices $0$ and $1$ of  $D_4^{(3)}$ 
have fibers of cardinality three.
Some further standard admissible affine foldings are
\[
(A_3^{(1)},A_1^{(1)}), \ \
(A_{2n-1}^{(1)},C_{n}^{(1)})_{n\ge2}, \ \
(D_{2n}^{(1)},A_{2n-1}^{(2)})_{n\ge2}, \ \
(D_{n+1}^{(1)},B_n^{(1)})_{n\ge3}, \ \
(E_6^{(1)},F_4^{(1)}),
\]
where
$$\text{$\Gamma\simeq C_2$ in all five cases}.$$

Finally, admissible folding is not restricted to finite or affine
Cartan data.  The following is a minimal simply-laced indefinite
example with a nontrivial admissible action:
\begin{align} \label{eq: indefinite folding}
\raisebox{1.2em}{
\xymatrix@R=1.5ex@C=4.5ex{
*{\circled{2}}<3pt>
   \ar@{-}[d]
   \ar@{}[r]^>{0\qquad \qquad \ \ }
   \ar@{-}[dr]
   \ar@{<.>}[r]^{\sigma}
&
*{\circled{2}}<3pt>
   \ar@{-}[dl]
   \ar@{-}[d]
   \ar@{}[l]_>{\qquad \qquad \ \  1}
\\
*{\circled{2}}<3pt>
   \ar@{-}[r]_>{2\qquad \qquad \ \ }
&
*{\circled{2}}<3pt>
   \ar@{-}[l]^>{ \qquad \qquad \ \ 3}
}}
\quad
&\overset{/\ang{\sigma}}{\Longrightarrow}
\quad
\raisebox{1.2em}{
\xymatrix@R=2.2ex@C=4.5ex{
&
*{\circled{2}}<3pt>
   \ar@{}[r]^>{0\qquad\qquad \ \ }
&
\\
*{\bcircled{4}}<3pt>
   \ar@{=>}[ur]
 \ar@{-}[rr]_>{1\qquad\qquad\quad\qquad\qquad \ \ }
&&
*{\bcircled{4}}<3pt>
   \ar@{=>}[ul]
   \ar@{-}[l]^>{\qquad\qquad \ \ 2}
}}.
\end{align}

\begin{remark}
Neither the simply-laced Cartan datum $\hg$ nor the quotient Cartan
datum $\g$ is determined by the other alone; the quotient depends also
on the choice of an admissible group action $\Gamma$.
Conversely, a given Cartan datum $\g$ may admit different admissible
simply-laced realizations.
The folding results below apply to any such realization and give the
same statements for the quotient datum $\g$.

For example, $C_2^{(1)}$ admits admissible simply-laced realizations
of types $A_7^{(1)}$ and $A_3^{(1)}$, with $m=4$ and $2$,
respectively.
\begin{align*}
\raisebox{2.5em}{
\xymatrix@R=1ex@C=3.5ex{
&& *{\circled{2}}<3pt>
\ar@{-}[dl]
\ar@{-}[dr]
\ar@{}[r]^<{ 2 \ \ }
&
\\
&
*{\circled{2}}<3pt>
\ar@{}[r]^<{1 \ \ }
&
&
*{\circled{2}}<3pt>
\ar@{}[l]_<{ \ \ 3}
\\
*{\circled{2}}<3pt>
\ar@{}[r]_<{0 \ \ }
\ar@{-}[dr]
\ar@{-}[ur]
&&&&
*{\circled{2}}<3pt>
\ar@{-}[dl]
\ar@{-}[ul]
\ar@{}[l]^<{ \ \ 4}
\\
&
*{\circled{2}}<3pt>
\ar@{}[r]_<{7 \ \ }
&
&
*{\circled{2}}<3pt>
\ar@{}[l]^<{ \ \ 5}
\\
&&
*{\circled{2}}<3pt>
\ar@{-}[ul]
\ar@{-}[ur]
\ar@{}[r]_<{6 \ \ }
&
}}
\ \
&\overset{/C_2 \times C_2}{\Longrightarrow}
\ \
\xymatrix@R=0.5ex@C=3.5ex{
*{\bcircled{4}}<3pt>
      \ar@{=>}[r]_<{0\ \ }
&
*{\circled{2}}<3pt>
       \ar@{}[r]_<{1\ \ }
&
*{\bcircled{4}}<3pt>
      \ar@{=>}[l]^<{ \ \ 2 }
}
\ \
\overset{/C_2}{\Longleftarrow}
\ \
\raisebox{1.5em}{
\xymatrix@R=1ex@C=3.5ex{
& *{\circled{2}}<3pt>
\ar@{-}[dl] \ar@{-}[dr]
\ar@{}[r]^<{ 1 \ }
\ar@{<.>}[dd]
&
\\
*{\circled{2}}<3pt>
\ar@{}[r]_<{0 \ \ }
&
&
*{\circled{2}}<3pt>
\ar@{}[l]^<{ \ \ 2}
\\
&*{\circled{2}}<3pt>
\ar@{-}[ul] \ar@{-}[ur]
\ar@{}[r]_<{3 \ \ }
&
}}.
\end{align*}
Here $C_2\times C_2\simeq\langle\tau,\sigma\rangle$, where
$\sigma=(1\,7)(2\,6)(3\,5)$ and $\tau=(0\,4)(1\,3)(5\,7)$.
For the left realization, the three fibers over $0,1,2$ are
$\{0,4\}$, $\{1,3,5,7\}$, and $\{2,6\}$, respectively.

Throughout this paper, in finite and affine types, we use the realizations fixed in the examples above.
\end{remark}

\section{Folding of root data} \label{sec: folding of roots}

Throughout this section, we retain the notation and conventions of
Section~\ref{sec: Cartan data}.
We establish the folding properties of root sequences and their associated pairings.

\subsection{Weyl groups and root lattices under folding}

Let $\hsfW$ and $\sfW$ be  the Weyl groups of $\hg$ and $\g$, respectively, and let
$\httB^+$ and $\ttB^+$ be the corresponding positive braid monoids.
For each $i\in I$, set
\[
\hs_i:=\prod_{a\in I_i}s_a\in\hsfW \qtq \hr_i:=\prod_{a\in I_i}r_a\in\httB^+.
\]
By admissibility, the factors in the above product commute.

We use the same symbol $\pi$ for the induced $\Z$-linear map
\[
\pi:\hsfQ\longrightarrow \sfQ, \qquad \pi(\hal_a)=\alpha_i \quad (a\in I_i).
\]
We also define the $\Q$-linear averaging map
\begin{align} \label{eq: averaging}
\jm:\sfQ_{\Q}\longrightarrow\hsfQ_{\Q}, \qquad \jmath(\alpha_i) = \frac{1}{e_i}\sum_{a\in I_i}\hal_a.
\end{align}

\begin{lemma}\label{lem: Weyl intertwining}
For every $i\in I$, we have
\[
\pi\circ\hs_i=s_i\circ\pi
\qtq
\jm\circ s_i=\hs_i\circ\jm.
\]
\end{lemma}

\begin{proof}
Let $b\in I_j$. By admissibility, we have
\( \hs_i(\hal_b) = \hal_b-\sum_{a\in I_i}\hc_{a,b}\hal_a \).
Hence
\[
\pi\bigl(\hs_i(\hal_b)\bigr)
= \alpha_j- \left(\sum_{a\in I_i}\hc_{a,b}\right)\alpha_i
= \alpha_j-c_{ij}\alpha_i
= s_i(\alpha_j),
\]
which proves $\pi\circ\hs_i=s_i\circ\pi$.

On the other hand,
\begin{align*}
\hs_i\bigl(\jm(\alpha_j)\bigr) &= \frac1{e_j}\sum_{b\in I_j} \left( \hal_b-\sum_{a\in I_i}\hc_{a,b}\hal_a \right)
=\jm(\alpha_j) -\frac{c_{ji}}{e_j}\sum_{a\in I_i}\hal_a \allowdisplaybreaks\\
&\overset{*}{=} \jm(\alpha_j) -\frac{c_{ij}}{e_i}\sum_{a\in I_i}\hal_a
= \jm(\alpha_j-c_{ij}\alpha_i) = \jm\bigl(s_i(\alpha_j)\bigr),
\end{align*}
where $\overset{*}{=}$ holds by $e_jc_{ij}=e_ic_{ji}$. \end{proof}

Recall from Lemma~\ref{lem: quotient Cartan} that
\[
d_i e_i=m\qquad(i\in I)
\]
for a positive integer $m$ independent of $i$.

\begin{lemma}\label{lem: form scaling}
For $\alpha,\beta\in\sfQ$, we have
\[
(\alpha,\beta)_\g = m\bigl(\jm(\alpha),\jm(\beta)\bigr)_\hg.
\]
\end{lemma}

\begin{proof}
For $i,j\in I$, we have
\[
m\bigl(\jm(\alpha_i),\jm(\alpha_j)\bigr)_\hg
= \frac{m}{e_i e_j} \sum_{a\in I_i}\sum_{b\in I_j} (\hal_a,\hal_b)_\hg
= \frac{m}{e_i}c_{ij} = d_i c_{ij} = (\alpha_i,\alpha_j)_\g. \qedhere
\]
\end{proof}

\begin{remark}
Throughout this paper, we fix the finite and affine realizations listed below, with the corresponding values of \(m\): 
\[
\begin{aligned}
m=2 \quad &\text{for }
(A_{2n-1},C_n),\ (D_{n+1},B_n),\ (E_6,F_4), \ (A_{2n-1}^{(1)},C_n^{(1)}), \   (A_{3}^{(1)},A_1^{(1)}),  \allowdisplaybreaks\\
&  \hspace{1.5em}  (D_n^{(1)},D_{n-1}^{(2)}),\ (D_{2n}^{(1)},A_{2n-1}^{(2)}),\   (D_{n+1}^{(1)},B_n^{(1)}),\  (E_6^{(1)},F_4^{(1)}),\allowdisplaybreaks\\
m=3 \quad &\text{for } (D_4,G_2),\ (D_4^{(1)},G_2^{(1)}),\ (E_6^{(1)},D_4^{(3)}),\allowdisplaybreaks\\
m=4 \quad &\text{for } (D_{2n+2}^{(1)},A_{2n}^{(2)}).
\end{aligned}
\]
In particular, $m$ need not be the order of $\Gamma$:
for $(D_n^{(1)},D_{n-1}^{(2)})$, we have
$\Gamma\simeq C_2\times C_2$ but $m=2$.
\end{remark}

\subsection{Arbitrary sequences and their lifts} \label{subsec: seq lift}

For each $i\in I_\g$, fix an ordering
\[
I_i=\seq{\im_1,\ldots,\im_{e_i}} \subset I_\hg.
\]
We use $\im$ for an index of $\hg$ when we wish to distinguish it
from an index $i$ of $\g$.

Let $\ii=\seq{i_1,\ldots,i_l}\in I^l_\g$
be an arbitrary finite sequence. We do not assume that $\ii$ is reduced.
For each $1\le k\le l$, write $I_{i_k}=\seq{\im_{k,1},\ldots,\im_{k,e_{i_k}}}$.
We define the blockwise lift of $\ii$ by
\begin{align} \label{eq: hii}
\hii \seteq I_{i_1}\star\cdots\star I_{i_l} = \seq{\im_{1,1},\ldots,\im_{1,e_{i_1}},\ldots, \im_{l,1},\ldots,\im_{l,e_{i_l}}}
\in I_\hg^{N_l},
\end{align}
where $N_s \seteq \sum_{k=1}^{s}e_{i_k}$ for $1 \le s \le l$. Here $\star$ denotes the concatenation operator.

Set $N_0\seteq0$ and, for $1\le k\le l$ and $1\le r\le e_{i_k}$, define
\begin{align} \label{eq: kappa k,r}
\kappa(k,r)\seteq N_{k-1}+r.
\end{align}
Thus $\kappa(k,r)$ is the position of $\im_{k,r}$ in $\hii$.

We associate with $\ii$ the root sequence
\[
\be_k^\ii \seteq s_{i_1}\cdots s_{i_{k-1}}(\alpha_{i_k}) \in \Phi_{\re,\g} \qquad(1\le k\le l).
\]
For its lift $\hii$, we set
\[
\hbe_{k,r} \seteq \be_{\kappa(k,r)}^{\hii} \in\Phi_{\re,\hg} \qquad (1\le k\le l,\ 1\le r\le e_{i_k}).
\]

By admissibility, we have
\begin{align} \label{eq: lifted block roots}
\hbe_{k,r}
= \hs_{i_1}\cdots\hs_{i_{k-1}} s_{\im_{k,1}}\cdots s_{\im_{k,r-1}} (\hal_{\im_{k,r}})
= \hs_{i_1}\cdots\hs_{i_{k-1}} (\hal_{\im_{k,r}}).
\end{align}

\begin{lemma}\label{lem: root folding}
For any $1\le k\le l$ and $1\le r\le e_{i_k}$, we have
\[
\pi(\hbe_{k,r})=\be_k^\ii \qtq \jm(\be_k^\ii) = \frac{1}{e_{i_k}} \sum_{t=1}^{e_{i_k}}\hbe_{k,t}.
\]
\end{lemma}

\begin{proof}
By Lemma~\ref{lem: Weyl intertwining} and \eqref{eq: lifted block roots}, we have
\[
\pi(\hbe_{k,r}) = s_{i_1}\cdots s_{i_{k-1}} \bigl(\pi(\hal_{\im_{k,r}})\bigr) = s_{i_1}\cdots s_{i_{k-1}}(\alpha_{i_k})
= \be_k^\ii.
\]
Similarly, we have
\[
\jm(\be_k^\ii) = \hs_{i_1}\cdots\hs_{i_{k-1}} \bigl(\jm(\alpha_{i_k})\bigr)
= \dfrac{1}{e_{i_k}} \sum_{t=1}^{e_{i_k}} \hs_{i_1}\cdots\hs_{i_{k-1}} (\hal_{\im_{k,t}})
= \dfrac{1}{e_{i_k}} \sum_{t=1}^{e_{i_k}}\hbe_{k,t}.
\qedhere
\]
\end{proof}

\begin{corollary}\label{cor: reduced lift}
The sequence $\ii$ is reduced if and only if $\hii$ is reduced.
\end{corollary}

\begin{proof}
Recall that $\ii$ is reduced if and only if $\be_k^\ii\in\Phi_{\re,\g}^+$ $(1\le k\le l)$.
By Lemma~\ref{lem: root folding}, we have $\pi(\hbe_{k,r})=\be_k^\ii$.
Since $\pi(\Phi_{\re,\hg}^+)\subset\sfQ^+$ and $\pi(\Phi_{\re,\hg}^-)\subset\sfQ^-$, we have
\[
\be_k^\ii\in\Phi_{\re,\g}^+ \ \ \text{ if and only if} \ \ \hbe_{k,r}\in\Phi_{\re,\hg}^+ \qquad(1\le r\le e_{i_k}).
\]
Thus all roots in the root sequence of $\ii$ are positive if and only if
all roots in the root sequence of $\hii$ are positive.
\end{proof}

By Corollary~\ref{cor: reduced lift} and
\cite[Theorem~28]{Castella08}, the assignments
$s_i\longmapsto\hs_i$ and $r_i\longmapsto\hr_i$
induce injective homomorphisms
\[
\iota:\sfW\hookrightarrow\hsfW
\qtq
\iota:\ttB^+\hookrightarrow\httB^+.
\]
Moreover,
\[
\iota(\sfW)=\hsfW^\Gamma
\qtq
\iota(\ttB^+)=(\httB^+)^\Gamma.
\]
In particular, if $\ii=\seq{i_1,\ldots,i_l}$ is reduced, then
$\hii$ is reduced and
\[
\ell_{\hsfW}\bigl(\hs_{i_1}\cdots\hs_{i_l}\bigr) = N_l = \sum_{k=1}^{l}e_{i_k}.
\]

\subsection{Folding of root pairings and $\lambda$-pairings} \label{subsec: lambda pairing}

The $\Gamma$-action on $I_\hg$ preserves each fiber $I_{i_k}$ and permutes its entries.
It induces the same permutation on the
 roots $\hbe_{k,r}$; namely,
\[
g(\im_{k,r})=\im_{k,r'} \qtq g(\hbe_{k,r})=\hbe_{k,r'}.
\]

\begin{lemma}\label{lem: root pairing folding}
For any $1\le k,p\le l$ and $1\le t\le e_{i_p}$, we have
\[
(\be_k^\ii,\be_p^\ii)_\g = \frac{m}{e_{i_k}} \sum_{r=1}^{e_{i_k}} (\hbe_{k,r},\hbe_{p,t})_\hg.
\]
In particular, the right-hand side is independent of $t$.
\end{lemma}

\begin{proof}
By Lemmas~\ref{lem: form scaling} and~\ref{lem: root folding},
\begin{align*}
(\be_k^\ii,\be_p^\ii)_\g
= m\bigl(\jm(\be_k^\ii),\jm(\be_p^\ii)\bigr)_\hg
= \frac{m}{e_{i_k}e_{i_p}} \sum_{r=1}^{e_{i_k}}\sum_{s=1}^{e_{i_p}} (\hbe_{k,r},\hbe_{p,s})_\hg.
\end{align*}
For $1\le t,t'\le e_{i_p}$, choose $g\in \Gamma$ such that
$g(\im_{p,t})=\im_{p,t'}$. Since $g$ permutes each fiber and preserves the bilinear form,
\begin{align*}
\sum_{r=1}^{e_{i_k}} (\hbe_{k,r},\hbe_{p,t})_\hg = \sum_{r=1}^{e_{i_k}} \bigl(g(\hbe_{k,r}),g(\hbe_{p,t})\bigr)_\hg
= \sum_{r=1}^{e_{i_k}} (\hbe_{k,r},\hbe_{p,t'})_\hg.
\end{align*}
Hence
\[
(\be_k^\ii,\be_p^\ii)_\g = \frac{m}{e_{i_k}e_{i_p}} \sum_{s=1}^{e_{i_p}}\sum_{r=1}^{e_{i_k}} (\hbe_{k,r},\hbe_{p,s})_\hg
= \frac{m}{e_{i_k}} \sum_{r=1}^{e_{i_k}} (\hbe_{k,r},\hbe_{p,t})_\hg.
\qedhere
\]
\end{proof}

Let
\[
L^\ii \seteq \soplus_{k=1}^{l}\Z\ve_k,
\]
where $\ve_k$ is labeled by the root $\be_k^\ii$.
With the sign opposite to that of \cite[(5.4)]{KKOP23F}, we define the skew-symmetric bilinear form
\[
\lambda^\ii:L^\ii\times L^\ii\longrightarrow\Z
\]
by
\begin{align}\label{eq: lambda pairing}
\lambda^\ii(\ve_k,\ve_p) \seteq (-1)^{\delta(k<p)} \delta(k\neq p) (\be_k^\ii,\be_p^\ii)_\g.
\end{align}

Likewise, let
\[
L^\hii \seteq \soplus_{k=1}^{l}\soplus_{r=1}^{e_{i_k}} \Z\hve_{k,r},
\]
where $\hve_{k,r}$ is labeled by $\hbe_{k,r}$.
For $r\neq s$, we have
\[
(\hbe_{k,r},\hbe_{k,s})_\hg = (\hal_{\im_{k,r}},\hal_{\im_{k,s}})_\hg  =0
\]
by admissibility. Hence we define
\[
\lambda^\hii:L^\hii\times L^\hii\longrightarrow\Z
\]
by
\begin{align}\label{eq: lifted lambda pairing}
\lambda^\hii(\hve_{k,r},\hve_{p,t}) \seteq (-1)^{\delta(k<p)} \delta\bigl((k,r)\neq(p,t)\bigr) (\hbe_{k,r},\hbe_{p,t})_\hg.
\end{align}

\begin{lemma}\label{lem: lambda pairing folding}
For any $1\le k,p\le l$ and $1\le t\le e_{i_p}$, we have
\[
\lambda^\ii(\ve_k,\ve_p)
= \frac{m}{e_{i_k}} \sum_{r=1}^{e_{i_k}} \lambda^\hii(\hve_{k,r},\hve_{p,t})
= d_{i_k} \sum_{r=1}^{e_{i_k}} \lambda^\hii(\hve_{k,r},\hve_{p,t}).
\]
In particular, the right-hand side is independent of $t$.
\end{lemma}

\begin{proof}
Since
$(\hbe_{k,r},\hbe_{k,t})_\hg=0$ for $r\neq t$,
we have
\begin{align*}
\frac{m}{e_{i_k}} \sum_{r=1}^{e_{i_k}} \lambda^\hii(\hve_{k,r},\hve_{p,t})
&= (-1)^{\delta(k<p)}\delta(k\neq p) \frac{m}{e_{i_k}} \sum_{r=1}^{e_{i_k}} (\hbe_{k,r},\hbe_{p,t})_\hg\\
& \overset{*}{=} (-1)^{\delta(k<p)}\delta(k\neq p) (\be_k^\ii,\be_p^\ii)_\g = \lambda^\ii(\ve_k,\ve_p),
\end{align*}
where $\overset{*}{=}$ follows from Lemma~\ref{lem: root pairing folding}.
The remaining equality follows from $m/e_{i_k}=d_{i_k}$.
\end{proof}

\section{Folding of compatible pairs} \label{sec: compatible pair folding}

We study compatible pairs associated with arbitrary finite sequences
and $\bfi$-box families, and establish their folding properties.
We slightly generalize the relevant constructions from
\cite{GLS06,GLS08,GLS11,FHOO2,KKOP24A,KK24} to the present setting.

\subsection{Quantum cluster algebras and pointed elements}

We briefly recall the notation for quantum cluster algebras that will be
used throughout this paper. We refer to \cite{FZ02,BZ05} for the general theory
and follow the conventions of \cite{KKOP24A,FHOO2}.

Let $J$ be a finite index set and let $J_{\ex}\subset J$ be the set of
exchangeable indices. Let
$B=(b_{ij})_{i\in J,\,j\in J_{\ex}}$
be an integer matrix and let
$\Lambda=(\lambda_{ij})_{i,j\in J}$
be a skew-symmetric integer matrix. We say that $(B,\Lambda)$ is a
compatible pair if
\begin{equation}\label{eq: compatible pair}
\sum_{p\in J} b_{pk}\lambda_{pj} =2d_k\delta_{k,j} \qquad(k\in J_{\ex},\ j\in J)
\end{equation}
for some positive integers $d_k$.

Associated with $\Lambda$ is the quantum torus $\calT(\Lambda)$
over $\Z[q^{\pm1/2}]$, generated by $X_i^{\pm1}$ $(i\in J)$ with
\[
X_iX_j=q^{\lambda_{ij}}X_jX_i.
\]
Fix a total order on $J$.
For $\bfa=(a_i)_{i\in J}\in\Z^J$, the
\emph{based quantum-torus monomial} $X^{\bfa}$ is defined by
\[
X^{\bfa}\seteq q^{c(\bfa)}
\prod_{i\in J}X_i^{a_i},
\qquad
c(\bfa)\seteq
-\frac{1}{2}\sum_{\substack{i,j\in J\\i<j}}\lambda_{ij}a_i a_j,
\]
where the product is taken in increasing order.
This definition is independent of the chosen total order. We write
\[
\Sigma_q=(\bfX,B,\Lambda), \qquad \bfX=(X_i\mid i\in J),
\]
for the corresponding quantum seed, whose based quantum torus is
$\calT(\Lambda)$, and denote its quantum cluster algebra by
$\calA_q(\Sigma_q)$.
We refer to \cite[Section~2.1]{KKOP23F} for the quantum cluster algebra
conventions used here. By specializing $q^{1/2}$ to $1$, a quantum cluster algebra specializes
to the corresponding (commutative) cluster algebra (see for instance~\cite[\S A.4]{FHOO}). We refer to these
two settings as the \emph{quantum setting} and the \emph{commutative
setting}, respectively.

In the commutative setting, a seed is a pair $\Sigma=(\bfx,B)$,
where $\bfx=(x_i\mid i\in J)$ is called the cluster of $\Sigma$.
Throughout this paper, all commutative cluster algebras and classical
Grothendieck rings are considered over $\C$.\footnote{Although cluster
algebras and Grothendieck rings have natural integral forms over $\Z$,
we work at the classical level over $\C$ for their later realization
as coordinate rings $\C[N(w)]$.} Whenever we move from the quantum
setting to the commutative cluster-algebra setting by specializing at
$q^{1/2}=1$, we implicitly extend scalars from $\Z$ to $\C$. We denote by
\[
\calL(\Sigma)\seteq\C[x_i^{\pm1}\mid i\in J]
\]
the Laurent polynomial algebra associated with $\Sigma$.

For $k\in J_{\ex}$, the mutation
\[
\mu_k(B,\Lambda)=(B',\Lambda')
\]
is defined by
\[
b'_{ij} =
\begin{cases}
-b_{ij}, & i=k \text{ or } j=k,\\
b_{ij}+[b_{ik}]_+[b_{kj}]_+ -[-b_{ik}]_+[-b_{kj}]_+, & \text{otherwise},
\end{cases}
\]
where $[a]_+\seteq\max(a,0)$, and
\[
\lambda'_{ij}=\lambda_{ij}\qquad(i,j\neq k),
\qquad\lambda'_{kk}=0,
\]
\[
\lambda'_{ik} = -\lambda_{ik} +\sum_{p\in J}[-b_{pk}]_+\lambda_{ip}, \qquad \lambda'_{ki}=-\lambda'_{ik}\qquad(i\neq k).
\]
Then $(B',\Lambda')$ is again a compatible pair with  the same positive integers $d_k$.

At the commutative level, for a cluster
$\bfx=(x_i\mid i\in J)$, set
\[
M_k^\pm(\bfx) \seteq \prod_{i\in J}x_i^{[\pm b_{ik}]_+}.
\]
The mutation at $k\in J_{\ex}$ is given by
\[
x'_i=x_i\qquad(i\neq k) \qtq x_kx'_k=M_k^+(\bfx)+M_k^-(\bfx).
\]

The cluster algebra $\calA(\Sigma)$ is the $\C$-subalgebra of
$\calL(\Sigma)$ generated by all cluster variables obtained from
$\Sigma$ by iterated mutations. By the Laurent phenomenon~\cite{FZ02}
and its quantum analogue~\cite{BZ05}, for every seed $\Sigma'$ mutation-equivalent
to $\Sigma$ and every quantum seed
$\Sigma_q'=(\bfX',B',\Lambda')$ mutation-equivalent to $\Sigma_q$, we have
\[
\calA(\Sigma)\subseteq\calL(\Sigma'),   \qquad   \calA_q(\Sigma_q)\subseteq\calT(\Lambda').
\]

We recall pointed elements and $g$-vectors (\cite{FZ07,Qin17}). An element
$z\in\calT(\Lambda)$ is said to be pointed at
$\bfg\in\Z^J$ if it has an expansion of the form
\begin{equation}\label{eq: pointed element}
z = q^a X^{\bfg} + \sum_{\bfn\in\N^{J_{\ex}}\setminus\{0\}} p_{\bfn}X^{\bfg+B\bfn}.
\end{equation}
The vector $\bfg$ is called the \emph{degree}, or the
\emph{$g$-vector}, of $z$ with respect to $\Sigma_q$,
and we write
\[
\deg_{\Sigma_q}(z)=\bfg.
\]
The same terminology will be used in the commutative setting by
specializing $q^{1/2}$ to $1$. Cluster monomials are pointed and hence
have well-defined $g$-vectors; see \cite[\S 5]{GHKK18}.

For later use, we recall the mutation rule for $g$-vectors.
Let $\Sigma_q'=\mu_k(\Sigma_q)$, and let $z$ be an element
compatibly pointed at $\Sigma_q$ and $\Sigma_q'$. Write
\[
\deg_{\Sigma_q}(z)=\bfg=(g_j)_{j\in J},
\qquad
\deg_{\Sigma_q'}(z)=\bfg'=(g'_j)_{j\in J}.
\] 
Then
\begin{equation}\label{eq: g-vector mutation}
g_j=
\begin{cases}
-g'_k, & j=k,\\[1mm]
g'_j+[b'_{jk}]_+g'_k,  & j\neq k,\ g'_k\geq0,\\[1mm]
g'_j+[-b'_{jk}]_+g'_k,  & j\neq k,\ g'_k\leq0,
\end{cases}
\end{equation}
where $B'=(b'_{ij})$ is the exchange matrix of the mutated seed
(see~\cite{FZ07,GHKK18,FHOO2} for more details).

\subsection{Compatible pairs for arbitrary sequences}

Let $\ii=\seq{i_1,\ldots,i_l}\in I^l$ be an arbitrary finite sequence.
For $1\le k\le l$ and $j\in I$, set
\begin{equation}\label{eq: k pm j}
\begin{aligned}
k_+ &\seteq \min\bigl(\{p\in[1,l]\mid k<p,\ i_p=i_k\}\cup\{\infty\}\bigr),\\
k_- &\seteq \max\bigl(\{p\in[1,l]\mid p<k,\ i_p=i_k\}\cup\{-\infty\}\bigr),\\
k(j)^+ &\seteq \min\bigl(\{p\in[1,l]\mid k\le p,\ i_p=j\}\cup\{\infty\}\bigr),\\
k(j)^- &\seteq \max\bigl(\{p\in[1,l]\mid p\le k,\ i_p=j\}\cup\{-\infty\}\bigr).
\end{aligned}
\end{equation}

Set
\[
J^\ii\seteq[1,l] \qtq J^\ii_{\ex}\seteq\{k\in[1,l]\mid k_+<\infty\}.
\]
Following the Geiss--Leclerc--Schr\"oer construction~\cite{GLS08}, we define the exchange matrix
\begin{equation}\label{eq: GLS exchange matrix}
B^\ii=(b^\ii_{p,k})_{p\in J^\ii,\,k\in J^\ii_{\ex}} \quad \text{ by } \ \
b^\ii_{p,k}\seteq
\begin{cases}
{\rm (i\text{-}1)} \ \ 1, & k=p_+,\\
{\rm (i\text{-}2)} \ \ -1, & p=k_+,\\
{\rm (ii\text{-}1)} \ \ c_{i_p,i_k}, & p<k<p_+<k_+,\\
{\rm (ii\text{-}2)} \ \ -c_{i_p,i_k}, & k<p<k_+<p_+,\\
{\rm (iii)} \ \  0, & \text{otherwise}.
\end{cases}
\end{equation}

For $1\le k\le l$, set
\[
k_{\min}\seteq \min\{p\in[1,k]\mid i_p=i_k\}.
\]
Thus $[k_{\min},k]$ is the $\bfi$-box determined by the position $k$;
see Subsection~\ref{subsec: ibox compatible pair} below.

Extending the construction for reduced sequences in
\cite{KKOP24A,KKOP23F}, we define the skew-symmetric matrix
\begin{equation}\label{eq: GLS Lambda matrix}
\Lambda^\ii=(\Lambda^\ii_{p,k})_{p,k\in J^\ii} \quad
\text{ by } \ \
\Lambda^\ii_{p,k} \seteq \sum_{\substack{p_{\min}\le a\le p\\ i_a=i_p}} \sum_{\substack{k_{\min}\le b\le k\\ i_b=i_k}}
\lambda^\ii(\ve_a,\ve_b).
\end{equation}
We call $(B^\ii,\Lambda^\ii)$ the \emph{GLS pair} associated with
$\ii$.
We denote by $\Sigma^\ii$ the commutative seed associated with the GLS pair
$(B^\ii,\Lambda^\ii)$, retaining $\Lambda^\ii$ as part of the seed data, and write
\[
  \bfx^\ii=(x_1,\ldots,x_l)
\]
for its cluster.
 We call $\Sigma^\ii$ the \emph{GLS seed} associated with $\ii$.

\begin{proposition}\label{prop: compatibility arbitrary sequence}
For any finite sequence $\ii=\seq{i_1,\ldots,i_l}$,
the GLS pair $(B^\ii,\Lambda^\ii)$ is compatible. More precisely,
\[
\sum_{a\in J^\ii} b^\ii_{a,k}\Lambda^\ii_{a,p}
= 2d_{i_k}\delta_{k,p} \qquad (k\in J^\ii_{\ex},\ p\in J^\ii).
\]
\end{proposition}

\begin{proof}
For $a\in J^\ii$, set
\begin{equation}\label{eq: big epsilon}
\calE_a \seteq \sum_{\substack{a_{\min}\le r\le a\\ i_r=i_a}}\ve_r \quad \text{ so that } \quad
\Lambda^\ii_{a,p} = \lambda^\ii(\calE_a,\calE_p).
\end{equation}
We use the convention $\calE_{-\infty}=0$.

Fix $k\in J^\ii_{\ex}$ and put $i\seteq i_k=i_{k_+}$.
In the $k$-th column of $B^\ii$, the positions of color $i$
contribute
\[
\calE_{k_-}-\calE_{k_+}=-\ve_k-\ve_{k_+}.
\]
For $j\neq i$, let $a_1<\cdots<a_s$ be the positions of color $j$
lying strictly between $k$ and $k_+$.
By the conditions~{\rm (ii)} in \eqref{eq: GLS exchange matrix}, the
contribution of the $j$-colored positions to
$\sum_{a\in J^\ii}b^\ii_{a,k}\calE_a$ is
\[
-c_{j,i}(\ve_{a_1}+\cdots+\ve_{a_s}).
\]

Summing over all colors $j\neq i$, we obtain
\begin{equation}\label{eq: B column telescoping}
-\xi_k \seteq \sum_{a\in J^\ii}b^\ii_{a,k}\calE_a = -\left( \ve_k+\ve_{k_+} +\sum_{k<a<k_+}c_{i_a,i}\ve_a \right)
= -\sum_{k\le a\le k_+}u_a\ve_a,
\end{equation}
where $u_k=u_{k_+}=1$ and  $u_a=c_{i_a,i}$ $(k<a<k_+)$.

Let $w_a=s_{i_1}\cdots s_{i_a}$, with $w_0=1$. Since $w_a(\alpha_i)-w_{a-1}(\alpha_i) =-c_{i_a,i}\be_a^\ii$,
summing over $k\le a<k_+$ gives
\[
\be_{k_+}^\ii-\be_k^\ii = -\sum_{a=k}^{k_+-1}c_{i_a,i}\be_a^\ii = -2\be_k^\ii -\sum_{k<a<k_+}c_{i_a,i}\be_a^\ii,
\]
which is equivalent to
\begin{equation}\label{eq: root telescoping compatibility}
\be_k^\ii+\be_{k_+}^\ii
+\sum_{k<a<k_+}c_{i_a,i}\be_a^\ii
=0.
\end{equation}
In terms of the coefficients $u_a$, this is
\begin{equation}\label{eq: root telescoping u}
\sum_{k\le a\le k_+}u_a\be_a^\ii=0.
\end{equation}

We now compute $\lambda^\ii(\xi_k,\ve_q)$. Recall that
\[
\lambda^\ii(\ve_a,\ve_q) = \delta(a \ne q) (-1)^{\delta(a<q)} (\be_a^\ii,\be_q^\ii).
\]

Thus~\eqref{eq: root telescoping u} implies that
\[
\lambda^\ii(\xi_k,\ve_q) = 0 \qquad \text{ for $q \not\in [k,k_+]$}.
\]
Next, suppose that $k<q<k_+$ and  $j\seteq i_q$. Summing
$w_a(\alpha_i)-w_{a-1}(\alpha_i) = -c_{i_a,i}\be_a^\ii$
over $k\le a<q$, we obtain
\[
w_{q-1}(\alpha_i)-w_{k-1}(\alpha_i) = -2\be_k^\ii -\sum_{k<a<q}c_{i_a,i}\be_a^\ii.
\]
Since $w_{k-1}(\alpha_i)=\be_k^\ii$, it follows that
\begin{equation}\label{eq: partial root telescoping}
\sum_{k\le a<q}u_a\be_a^\ii = \be_k^\ii +\sum_{k<a<q}c_{i_a,i}\be_a^\ii  = -w_{q-1}(\alpha_i).
\end{equation}
On the other hand, since $u_q=c_{j,i}$,
\eqref{eq: root telescoping u} and
\eqref{eq: partial root telescoping} give
\[
\sum_{q<a\le k_+}u_a\be_a^\ii
= -\sum_{k\le a<q}u_a\be_a^\ii -u_q\be_q^\ii
= w_{q-1}(\alpha_i)-c_{j,i}\be_q^\ii.
\]
Therefore, by the definition of $\lambda^\ii$,
\[
\lambda^\ii(\xi_k,\ve_q) =
-\left(\sum_{k\le a<q}u_a\be_a^\ii, \be_q^\ii \right) + \left( \sum_{q<a\le k_+}u_a\be_a^\ii, \be_q^\ii \right)
= \bigl( 2w_{q-1}(\alpha_i)-c_{j,i}\be_q^\ii, \be_q^\ii \bigr).
\]
Since $\be_q^\ii=w_{q-1}(\alpha_j)$ and the bilinear form is $\sfW$-invariant, we have
\[
\bigl(w_{q-1}(\alpha_i),\be_q^\ii\bigr) = \bigl(w_{q-1}(\alpha_i),w_{q-1}(\alpha_j)\bigr) = (\alpha_i,\alpha_j)
\qtq
(\be_q^\ii,\be_q^\ii) = (\alpha_j,\alpha_j).
\]
Consequently,
\[
\lambda^\ii(\xi_k,\ve_q) = 2(\alpha_i,\alpha_j) -c_{j,i}(\alpha_j,\alpha_j)  = 2d_ic_{ij}-2d_jc_{ji} =0,
\]
where the last equality follows from $d_ic_{ij}=d_jc_{ji}$.

It remains to consider $q=k$ and $q=k_+$. For $q=k$, all terms of
$\xi_k$ other than $\ve_k$ lie to the right of $k$, and hence
\[
\lambda^\ii(\xi_k,\ve_k) = \left( \sum_{k<a\le k_+}u_a\be_a^\ii, \be_k^\ii \right) = (-\be_k^\ii,\be_k^\ii) = -(\alpha_i,\alpha_i) = -2d_i.
\]
Here we used \eqref{eq: root telescoping u} in the second equality.

For $q=k_+$, all terms other than $\ve_{k_+}$ lie to the left of
$k_+$, so
\[
\lambda^\ii(\xi_k,\ve_{k_+}) = -\left( \sum_{k\le a<k_+}u_a\be_a^\ii, \be_{k_+}^\ii \right) = (\be_{k_+}^\ii,\be_{k_+}^\ii) = (\alpha_i,\alpha_i) = 2d_i.
\]
Thus, for every $q\in J^\ii$,
\begin{equation}\label{eq: lambda xi basis}
\lambda^\ii(\xi_k,\ve_q)
=
\begin{cases}
-2d_i,&q=k,\\
2d_i,&q=k_+,\\
0,&\text{otherwise}.
\end{cases}
\end{equation}

We now return to the compatibility relation. By
\eqref{eq: B column telescoping} and
\eqref{eq: big epsilon}, we have
\[
\begin{aligned}
\sum_{a\in J^\ii}b^\ii_{a,k}\Lambda^\ii_{a,p} = \sum_{a\in J^\ii} b^\ii_{a,k}\lambda^\ii(\calE_a,\calE_p)
= \lambda^\ii \left( \sum_{a\in J^\ii}b^\ii_{a,k}\calE_a, \calE_p \right)
= -\lambda^\ii(\xi_k,\calE_p).
\end{aligned}
\]

Since $\calE_p$ is the sum of the basis vectors corresponding to
the occurrences of the color $i_p$ up to $p$, and $k_+$ is the next
occurrence of the color $i_k$ after $k$, \eqref{eq: lambda xi basis}
gives
\[
\lambda^\ii(\xi_k,\calE_p) = -2d_i\delta_{k,p}.
\]
Therefore,
\[
\sum_{a\in J^\ii}b^\ii_{a,k}\Lambda^\ii_{a,p} = -\lambda^\ii(\xi_k,\calE_p) = 2d_i\delta_{k,p} = 2d_{i_k}\delta_{k,p}. \qedhere
\]
\end{proof}

\begin{remark}
For reduced sequences, compatible pairs of this type originate from
the Geiss--Leclerc--Schr\"oer construction
\cite{GLS06,GLS08,GLS11} for the symmetric setting.
Their counterparts in the symmetrizable setting appear in the
quantum cluster structures of Goodearl--Yakimov \cite{GY17,GY21}.
The formulation used here for arbitrary sequences in finite type is
essentially that of \cite[Proposition~1.2, Corollary~1.3]{FHOO2}.
The proof above uses neither reducedness nor finiteness of the Weyl
group. Thus, the same argument applies to the present symmetrizable
Kac--Moody setting.
\end{remark}

\subsection{$\bfi$-box families and compatible pairs associated with them}
\label{subsec: ibox compatible pair}

We recall the $\bfi$-box combinatorics of \cite{KK24}; see also
\cite{KKOP24A}. Let $\ii=\seq{i_1,\ldots,i_l}\in I^l$ be an arbitrary finite sequence.
An \emph{$\bfi$-box} is an interval
\[
\frakc=[a,b]\subset[1,l] \ \ \text{ such that $i_a=i_b$}.
\]
We write $i_{\frakc}\seteq i_a=i_b$ for its \emph{color}. When $a=b=k$, we simply write $[k]$ for the singleton
$\ii$-box.

Associated with an $\bfi$-box $\frakc=[a,b]$, we define
\begin{equation}\label{eq: big epsilon box}
\calE_{\frakc} \seteq \sum_{\substack{a\le r\le b\\ i_r=i_{\frakc}}}\ve_r.
\end{equation}
We call $\calE_{\frakc}$ the \emph{box vector} associated with
$\frakc$. We set $\calE_\emptyset \seteq 0$.
In particular, for the $\bfi$-box $\frakc_k=[k_{\min},k]$, we have $\calE_{\frakc_k}=\calE_k$ by \eqref{eq: big epsilon}.

Two $\bfi$-boxes $\frakc=[a,b]$ and $\frakd=[c,d]$
are said to \emph{commute} if
\[
a_-<c\le d<b_+ \qtoq  c_-<a\le b<d_+.
\]
A collection of pairwise commuting $\bfi$-boxes is called a
\emph{commuting family}. A commuting family of $\bfi$-boxes in an interval $[a,b]$ is called
\emph{maximal} if it is maximal among the commuting families of
$\bfi$-boxes in $[a,b]$.

For an interval $[a,b]\subset[1,l]$, set
\[
[a,b\}\seteq[a,b(i_a)^-] \qtq \{a,b]\seteq[a(i_b)^+,b].
\]

A sequence
\[
\frakC=\seq{\frakc_1,\ldots,\frakc_r}, \qquad \frakc_k=[x_k,y_k],
\]
of $\bfi$-boxes is called an \emph{admissible chain} if, for each
$1\le k\le r$, the union
\[
\tfrakc_k \seteq \bigcup_{1\le s\le k}[x_s,y_s]
\]
is an interval of cardinality $k$, and
\[
\frakc_k=[\tx_k,\ty_k\} \quad\text{or}\quad \frakc_k=\{\tx_k,\ty_k],
\]
where $\tfrakc_k=[\tx_k,\ty_k]$. We call $\tfrakc_k$ the \emph{envelope} of $\frakc_k$, and $\tfrakc_r$ the \emph{extent} of $\frakC$.

The underlying family
$$\calF(\frakC) \seteq \{\frakc_1,\ldots,\frakc_r\} \text{ of $\frakC$ }$$
is a maximal commuting family of $\bfi$-boxes in the extent of
$\frakC$. Conversely, every maximal commuting family of
$\bfi$-boxes in an interval $[a,b]$ arises in this way from an
admissible chain of extent $[a,b]$ (\cite[Proposition~5.3]{KKOP24A} and \cite[Corollary~2.10]{KK24}).
In particular, every such maximal commuting family has cardinality
$b-a+1$.

Let $\calF$ be a maximal commuting family of $\bfi$-boxes in
$[a,b]$, and choose an admissible chain
$\frakC=\seq{\frakc_1,\ldots,\frakc_r}$
of extent $[a,b]$ such that $\calF=\calF(\frakC)$. Note that such an admissible chain $\frakC$ need not be unique.
Since $|\tfrakc_k|=|\tfrakc_{k-1}|+1$, there exists a unique
\[
z_k\in\tfrakc_k\setminus\tfrakc_{k-1} \qtq z_k\in\{x_k,y_k\},
\]
where $\frakc_k=[x_k,y_k]$ and $\tfrakc_0=\emptyset$.

By \cite[Proposition~2.12]{KK24}, the endpoint $z_k$ depends only on
$\calF$ and $\frakc_k$, and not on the choice of an admissible chain
realizing $\calF$. We call $z_k$ the \emph{effective end} of
$\frakc_k$ in $\calF$, and write
\[
\eff_{\calF}(\frakc_k)\seteq z_k.
\]
Thus
\begin{align} \label{eq: eff_F}
\eff_{\calF}:\calF\longrightarrow[a,b] \text{ is a bijection.}
\end{align}

We now associate an exchange matrix with a maximal commuting family
$\calF$ of $\bfi$-boxes. We first define the full exchange matrix
\[
B^\calF_{\full} = \bigl(b^\calF_{\frakc,\frakd}\bigr)_{\frakc,\frakd\in\calF},
\]
with skew-symmetrizer $d_{\frakc}\seteq d_{i_{\frakc}} \qquad(\frakc\in\calF)$.

Let $\frakc=[x,y]$ and  $\frakd=[x',y']$ be distinct $\bfi$-boxes in $\calF$. The \emph{negative} entries of
$B^\calF_{\full}$ are defined as follows:
\begin{equation}\label{eq: full ibox exchange matrix}
b^\calF_{\frakc,\frakd}
\seteq
\begin{cases}
-1, & \ {\rm (i\text{-}1)} \ x=x' \text{ and } y'=y_-, \text{ or } {\rm (i\text{-}2)} \ y=y' \text{ and } x'=x_-,  \\
c_{i_{\frakc},i_{\frakd}}, &
\begin{array}{l}
{\rm (ii)} \ c_{i_{\frakc},i_{\frakd}}<0
\text{ and one of {\rm (a)--(d)} holds:}
\\[1mm]
\ \ {\rm (ii\text{-}a)}\ [x,y_+]\in\calF,\  \eff_\calF(\frakc)=x,\  x'_-<x<x',\ y'<y_+<y'_+, \\
\ \ {\rm (ii\text{-}b)}\ [x,y_+]\in\calF,\  \eff_\calF(\frakd)=y',\ x'_-<x,\  y<y'<y_+<y'_+, \\
\ \ {\rm (ii\text{-}c)}\ [x'_-,y']\in\calF,\  \eff_\calF(\frakd)=y',\ x_-<x'_-<x,\ y<y'<y_+, \\
\ \ {\rm (ii\text{-}d)}\ [x'_-,y']\in\calF,\  \eff_\calF(\frakc)=x,\ x_-<x'_-<x<x',\ y'<y_+.
\end{array}
\end{cases}
\end{equation}
We call the conditions {\rm (i\text{-}1)} and {\rm (i\text{-}2)}
the \emph{horizontal relations}, and the conditions
{\rm (ii\text{-}a)}--{\rm (ii\text{-}d)} the \emph{vertical relations}.
The remaining entries are uniquely determined by
\[
d_{\frakc}b^\calF_{\frakc,\frakd} = -d_{\frakd}b^\calF_{\frakd,\frakc}, \qquad b^\calF_{\frakc,\frakc}=0.
\]
Our convention is opposite to that of \cite[Section~3.2]{KK24}; namely, if $B_{\rm KK}(\calF)$ denotes the full exchange matrix
defined in \cite[Section~3.2]{KK24},
\[
B^\calF_{\full}=-B_{\rm KK}(\calF).
\]

An $\bfi$-box $\frakc=[x,y]\in\calF$ is called \emph{frozen} if
\[
x_-=-\infty \qtq y_+=\infty,
\]
and \emph{exchangeable} otherwise. Set
\[
\calF_{\ex}\seteq \{\frakc\in\calF\mid \frakc\text{ is exchangeable}\}.
\]
The exchange matrix associated with $\calF$ is the restriction
\[
B^\calF \seteq \left.B^\calF_{\full}\right|_{\calF\times\calF_{\ex}}
= \bigl(b^\calF_{\frakc,\frakd}\bigr)_{\frakc\in\calF,\,\frakd\in\calF_{\ex}}.
\]

For a maximal commuting family $\calF$, we define the skew-symmetric matrix
\begin{equation}\label{eq: ibox Lambda matrix}
\Lambda^\calF = \bigl(\Lambda^\calF_{\frakc,\frakd}\bigr)_{\frakc,\frakd\in\calF}
\ \  \text{ by } \ \   \Lambda^\calF_{\frakc,\frakd} \seteq \lambda^\ii(\calE_{\frakc},\calE_{\frakd}).
\end{equation}

For $1\le k\le l$, set
\[
\frakc_k\seteq[k_{\min},k] \qtq \calF^\ii\seteq\{\frakc_k\mid 1\le k\le l\}.
\]
Then $\calF^\ii$ is a maximal commuting family of $\bfi$-boxes, and
\[
B^{\calF^\ii}=B^\ii \qtq \Lambda^{\calF^\ii}=\Lambda^\ii.
\]
Indeed, the equality for the $\Lambda$-matrices follows immediately
from \eqref{eq: big epsilon} and \eqref{eq: ibox Lambda matrix}.
For the exchange matrices, the horizontal relations in {\rm (i)} of
\eqref{eq: full ibox exchange matrix} give precisely the entries in
{\rm (i)} of \eqref{eq: GLS exchange matrix}, while the vertical
relations in {\rm (ii)} of \eqref{eq: full ibox exchange matrix}
reduce exactly to the conditions in {\rm (ii)} of
\eqref{eq: GLS exchange matrix}.

We recall the $\{\calL,\calR\}$-sequence associated with an admissible
chain. Let $\frakC=\seq{\frakc_1,\ldots,\frakc_r}$
be an admissible chain of extent $[a,b]$, and write
$\frakc_k=[x_k,y_k]$ and $\tfrakc_k=[\tx_k,\ty_k]$.
Note that $\frakC$ is uniquely determined by its sequence of envelopes
$\seq{\tfrakc_1,\ldots,\tfrakc_r}$. Indeed, for $1<k\le r$, we have
\begin{equation}\label{eq: box from envelopes}
\frakc_k=[x_k,y_k] =
\begin{cases}
[\tx_k,\ty_{k-1}\}, & \tx_k=\tx_{k-1}-1,\\
\{\tx_{k-1},\ty_k], & \ty_k=\ty_{k-1}+1.
\end{cases}
\end{equation}
In the first case, we set $T_{k-1}\seteq \calL$, while in the second case we set
$T_{k-1}\seteq  \calR$. Thus, to $\frakC$ we associate an $\{\calL,\calR\}$-sequence
\[
\frakT(\frakC)=\seq{T_1,\ldots,T_{r-1}}\in\{\calL,\calR\}^{r-1}.
\]

Conversely, an $\{\calL,\calR\}$-sequence
\[
\frakT=\seq{T_1,\ldots,T_{r-1}}\in\{\calL,\calR\}^{r-1}
\]
together with the extent $[a,b]$ uniquely determines an admissible
chain. More precisely, its initial envelope is the singleton
\[
\tfrakc_1=[x_1],
\qquad
x_1=a+\#\{k\mid T_k=\calL\},
\]
and the subsequent envelopes are determined recursively by
\[
\tfrakc_{k+1} =
\begin{cases}
[\tx_k-1,\ty_k], & T_k=\calL,\\[1mm]
[\tx_k,\ty_k+1], & T_k=\calR.
\end{cases}
\qquad(1\le k<r).
\]
The corresponding $\bfi$-boxes are then recovered from
\eqref{eq: box from envelopes}.
Hence admissible chains of extent $[a,b]$ are in bijection with
$\{\calL,\calR\}$-sequences of length $r-1=b-a$.

For $1\le k<r$, we say that $\frakc_k$ is \emph{movable} if
$k=1$ or $T_{k-1}\neq T_k \quad (k\ge2)$.
For a movable $\frakc_k$ in $\frakC$, the \emph{box move} of
$\frakC$ at $k$ is the admissible chain $\frakB_k(\frakC)$
determined by $x'_1$ and
$\frakT'=\seq{T'_1,\ldots,T'_{r-1}}$ as follows:
\[
\begin{array}{ll}
{\rm (i)}&
\begin{cases}
x_1'=x_1\pm1, & \text{if $k=1$ and $T_1=\calR$ (resp.\ $\calL$)},\\
x_1'=x_1, & \text{if $k>1$,}
\end{cases}
\\
{\rm (ii)} & T'_s=T_s \quad\text{for }s\notin\{k-1,k\},  \qquad T'_s\neq T_s \quad\text{for }s\in\{k-1,k\}.
\end{array}
\]
Thus $\frakB_k(\frakC)$ is obtained from $\frakC$ by shifting
$\tfrakc_k$ by one step to the right or to the left inside
$\tfrakc_{k+1}$.

\begin{proposition} [{\cite[Proposition 5.6, Proposition 5.7]{KKOP24A}}] \label{prop:boxmove}
Let $\frakC=\seq{\frakc_k}_{1 \le k \le l}$ be an admissible chain of
$\bfi$-boxes and let $k_0$ be the index of 
a movable $\bfi$-box $(1\le k_0<l)$. Set
$\frakB_{k_0}(\frakC) = (\frakc_k')_{1 \le k \le l}$.
\bnum
\item \label{it: per boxmove}
Assume that $i_{\tx_{k_0+1}} \ne i_{\ty_{k_0+1}}$, i.e., $\tfrakc_{k_0+1}$ is
not an $\bfi$-box. Then we have
$$\frakc_k' = \frakc_{\fraks_{k_0}(k)} \quad \text{for}\quad 1\le k \le l, $$
where $\fraks_{k_0} \in \frakS_l$ is the transposition of $k_0$
and $k_0+1$.
\item \label{it: mutation box-move} Assume that $i_{\tx_{k_0+1}}=i_{\ty_{k_0+1}}$, i.e., $\tfrakc_{k_0+1}$ is an $\bfi$-box.
Set $\frakc_{k_0+1}=[x,y]$ with $x=\tx_{k_0+1}$ and $y=\ty_{k_0+1}$.

\noindent
$\bullet$ If $k_0>1$, then
\bna
\item $\frakc_{k_0}=[x_+,y]$ and $\frakc'_{k_0}=[x,y_-]$ if $T_{k_0-1}=\calR$,
\item $\frakc_{k_0}=[x,y_-]$ and $\frakc'_{k_0}=[x_+,y]$ if $T_{k_0-1}=\calL$.
\ee
\noindent
$\bullet$ If $k_0=1$, then $y=x+1$ and
\[
\begin{array}{lll}
T_1=\calL:& \frakc_1=[y],& \frakc'_1=[x],\\
T_1=\calR:& \frakc_1=[x],& \frakc'_1=[y].
\end{array}
\]
\ee
\end{proposition}
A box move as in~\eqref{it: per boxmove} is called a \emph{permutation-type box move},
while a box move as in~\eqref{it: mutation box-move} is called a \emph{mutation-type box move}.

By abuse of notation, we also write $\frakB_k\calF$ for the maximal
commuting family obtained from $\calF$ by the box move $\frakB_k$.

\begin{proposition}\label{prop: ibox pair box move}
Let $\calF'=\frakB_k\calF$.
\bna
\item If $\frakB_k$ is of permutation type, then
\[
B^{\calF'} = P B^\calF (P_{\ex})^T \qtq \Lambda^{\calF'} = P\Lambda^\calF P^T,
\]
where $P$ is the induced permutation matrix and $P_{\ex}$ is its
restriction to the exchangeable indices.
\item If $\frakB_k$ is of mutation type at $\frakc\in\calF$, then
\[
B^{\calF'}=\mu_{\frakc}(B^\calF) \qtq \Lambda^{\calF'}=\mu_{\frakc}(\Lambda^\calF).
\]
\ee
\end{proposition}

\begin{proof}
The permutation case follows immediately from the definitions.

Assume that $\frakB_k$ is of mutation type at
$\frakc\in\calF$, and let $\frakc'$ be the new $\bfi$-box.
The equality
\[
B^{\calF'}=\mu_{\frakc}(B^\calF)
\]
is the local box-move formula of \cite{KK24}, with our opposite sign
convention. We explain below how the relevant entries of the
$\frakc$-column follow from \eqref{eq: full ibox exchange matrix}. It remains to prove that
\[
\Lambda^{\calF'}=\mu_{\frakc}(\Lambda^\calF).
\]

We consider the two orientations separately.

\smallskip

\noindent
{\rm (i)} Suppose that
\[
\frakc=[x_+,y] \qtq \frakc'=[x,y_-], \qquad i\seteq i_x=i_y.
\]
For $k>1$, one has $T_{k-1}=\calR$; for $k=1$, one has
$\frakc=[y]$. Thus
\[
\eff_{\calF}(\frakc)=y.
\]

\noindent
$\bullet$ For each $j\neq i$, set
\[
\fraka_j=[x(j)^+,y(j)^-]
\]
whenever this $\bfi$-box is nonempty. By
Proposition~\ref{prop:boxmove}~\eqref{it: mutation box-move},
$\calF$ contains the $\bfi$-box $[x,y]$.

For the ordered pair $(\fraka_j,\frakc)$, condition~{\rm (ii-c)} of
\eqref{eq: full ibox exchange matrix} implies that
\[
b^\calF_{\fraka_j,\frakc}=c_{j,i}\leq0.
\]

\noindent
$\bullet$ For boxes of color $i$, no horizontal relation contributes
to the negative part of the $\frakc$-column in this orientation.

Thus, apart from the boxes $\fraka_j$ above, all boxes contribute
zero to the negative part of the $\frakc$-column. Consequently,
\begin{equation}\label{eq: first box move negative column}
[-b^\calF_{\fraka_j,\frakc}]_+=-c_{j,i}
\qquad(j\neq i),
\end{equation}
where empty boxes $\fraka_j$ are omitted.

We next compare the corresponding box vectors. By definition,
\[
\calE_{\frakc}
=
\sum_{\substack{x_+\le p\le y\\i_p=i}}\ve_p
\qtq
\calE_{\frakc'}
=
\sum_{\substack{x\le p\le y_-\\i_p=i}}\ve_p.
\]
Thus their sum contains $\ve_x$ and $\ve_y$ once and every
$i$-colored position strictly between $x$ and $y$ twice. Since
$c_{i,i}=2$, we have
\[
\calE_{\frakc'}+\calE_{\frakc}
=
\ve_x+\ve_y
+
\sum_{\substack{x<p<y\\i_p=i}}
c_{i_p,i}\ve_p.
\]
For $j\neq i$,  the positions of color $j$ in $\fraka_j$ are precisely
those lying strictly between $x$ and $y$.{} Hence
\[
c_{j,i}\calE_{\fraka_j}
=
\sum_{\substack{x<p<y\\i_p=j}}
c_{i_p,i}\ve_p.
\]
Summing over all $j\neq i$ gives
\begin{equation}\label{eq: first box move vector identity}
\calE_{\frakc'}+\calE_{\frakc}
-
\sum_{j\neq i}(-c_{j,i})\calE_{\fraka_j}
=
\ve_x+\ve_y+\sum_{x<p<y}c_{i_p,i}\ve_p.
\end{equation}

Since $i_x=i_y=i$, summing $w_p(\alpha_i)-w_{p-1}(\alpha_i) = -c_{i_p,i}\be_p^\ii$
over $x\le p<y$ gives
\[
\be_x^\ii+\be_y^\ii +\sum_{x<p<y}c_{i_p,i}\be_p^\ii =0.
\]
 Summing \eqref{eq: lambda xi basis} over the positions $p$ with
$x\le p<y$ and $i_p=i$ gives{}
\begin{equation}\label{eq: local lambda box move}
\lambda^\ii
\left(
\ve_x+\ve_y+\sum_{x<p<y}c_{i_p,i}\ve_p,
\ve_q
\right)
=
2d_i(\delta_{q,y}-\delta_{q,x}).
\end{equation}

Now let $\frakd\in\calF\cap\calF'$. By
\eqref{eq: ibox Lambda matrix},
\eqref{eq: first box move negative column}, and
\eqref{eq: first box move vector identity}, we have
\[
\Lambda^{\calF'}_{\frakc',\frakd} +\Lambda^\calF_{\frakc,\frakd} -\sum_{\fraka\in\calF} [-b^\calF_{\fraka,\frakc}]_+ \Lambda^\calF_{\fraka,\frakd}
= \lambda^\ii \left( \ve_x+\ve_y+\sum_{x<p<y}c_{i_p,i}\ve_p, \calE_{\frakd} \right).
\]
Using \eqref{eq: local lambda box move} and the definition of
$\calE_{\frakd}$, the right-hand side is
\[
2d_i
\sum_{\substack{q\in\frakd\\i_q=i_{\frakd}}}
(\delta_{q,y}-\delta_{q,x}).
\]
If $i_{\frakd}\neq i$, both delta terms vanish. Suppose that
$i_{\frakd}=i$. Since $\frakd\in\calF\cap\calF'$, the box $\frakd$ commutes with
both $\frakc=[x_+,y]$ and $\frakc'=[x,y_-]$.
The commuting condition gives
\[
x\in\frakd
\quad\Longleftrightarrow\quad
y\in\frakd.
\]
Thus the two delta contributions cancel, and therefore
\[
\Lambda^{\calF'}_{\frakc',\frakd}
=
-\Lambda^\calF_{\frakc,\frakd}
+
\sum_{\fraka\in\calF}
[-b^\calF_{\fraka,\frakc}]_+
\Lambda^\calF_{\fraka,\frakd}.
\]

\smallskip

\noindent
{\rm (ii)}
Suppose that
\[
\frakc=[x,y_-] \qtq \frakc'=[x_+,y].
\]
For $k>1$, one has $T_{k-1}=\calL$; for $k=1$, one has
$\frakc=[x]$. Thus
\[
\eff_{\calF}(\frakc)=x.
\]

By Proposition~\ref{prop:boxmove}, the box $[x,y]$ belongs to
$\calF$. It has the same left endpoint as $\frakc=[x,y_-].$
Hence condition~{\rm (i-1)} of \eqref{eq: full ibox exchange matrix} gives
\[
b^\calF_{[x,y],\frakc}=-1.
\]

If $[x_+,y_-]$ is nonempty, it also belongs to $\calF$. This box
has the same right endpoint as $\frakc$. Hence condition~{\rm (i-2)} gives
\[
b^\calF_{[x_+,y_-],\frakc}=-1.
\]

These are negative entries in the $\frakc$-column.
Therefore,
\begin{align} \label{eq: 1 and 1}
[-b^\calF_{[x,y],\frakc}]_+=1
\quad\text{and}\quad
[-b^\calF_{[x_+,y_-],\frakc}]_+=\delta([x_+,y_-] \ne \emptyset).
\end{align}
Note that the effective-end conditions in the vertical cases
{\rm (ii-a)}--{\rm (ii-d)} are not satisfied for any additional
negative entry in this orientation. Thus the only nonzero coefficients
$[-b^\calF_{\fraka,\frakc}]_+$ are given in~\eqref{eq: 1 and 1}.

The corresponding box vectors satisfy
\begin{equation}\label{eq: second box move vector identity}
\calE_{\frakc'}+\calE_{\frakc} = \calE_{[x_+,y_-]}+\calE_{[x,y]}.
\end{equation}
Indeed, both sides count $\ve_x$ and $\ve_y$
once and every $i$-colored position strictly between them twice.

Let $\frakd\in\calF\cap\calF'$. By
\eqref{eq: ibox Lambda matrix},
\eqref{eq: 1 and 1}, and
\eqref{eq: second box move vector identity}, we obtain
\[
\Lambda^{\calF'}_{\frakc',\frakd} +\Lambda^\calF_{\frakc,\frakd}
= \Lambda^\calF_{[x_+,y_-],\frakd} +\Lambda^\calF_{[x,y],\frakd}
= \sum_{\fraka\in\calF} [-b^\calF_{\fraka,\frakc}]_+ \Lambda^\calF_{\fraka,\frakd}.
\]
Hence
\[
\Lambda^{\calF'}_{\frakc',\frakd}= -\Lambda^\calF_{\frakc,\frakd} + \sum_{\fraka\in\calF} [-b^\calF_{\fraka,\frakc}]_+ \Lambda^\calF_{\fraka,\frakd}.
\]

Thus, in both orientations,
\[
\Lambda^{\calF'}=\mu_{\frakc}(\Lambda^\calF).
\qedhere
\]
\end{proof}

\begin{theorem}\label{thm: ibox compatibility}
Let $\ii=\seq{i_1,\ldots,i_l}\in I^l$ be an arbitrary finite sequence,
and let $\calF$ be a maximal commuting family of $\bfi$-boxes in
$[1,l]$. Then $(B^\calF,\Lambda^\calF)$ is a compatible pair. More
precisely,
\[
\sum_{\frakc\in\calF} b^\calF_{\frakc,\frakd}\Lambda^\calF_{\frakc,\fraku}
= 2d_{i_{\frakd}}\delta_{\frakd,\fraku}
\qquad
(\frakd\in\calF_{\ex},\ \fraku\in\calF).
\]
\end{theorem}

\begin{proof}
Let $\calF^\ii$ be the maximal commuting family corresponding to the
constant $\calR$-sequence. Since
\[
B^{\calF^\ii}=B^\ii \qtq \Lambda^{\calF^\ii}=\Lambda^\ii,
\]
Proposition~\ref{prop: compatibility arbitrary sequence} gives
\[
\sum_{\frakc\in\calF^\ii} b^{\calF^\ii}_{\frakc,\frakd} \Lambda^{\calF^\ii}_{\frakc,\fraku}
= 2d_{i_{\frakd}}\delta_{\frakd,\fraku}
\qquad
(\frakd\in\calF^\ii_{\ex},\ \fraku\in\calF^\ii).
\]

Every maximal commuting family $\calF$ is connected to $\calF^\ii$
by a sequence of box moves. By Proposition~\ref{prop: ibox pair box move},
at each step the pair $(B,\Lambda)$ is transformed by either a
simultaneous permutation or a simultaneous matrix mutation.
Since compatibility is preserved under both operations, we obtain
\[
\sum_{\frakc\in\calF} b^\calF_{\frakc,\frakd} \Lambda^\calF_{\frakc,\fraku}
= 2d_{i_{\frakd}}\delta_{\frakd,\fraku} \qquad (\frakd\in\calF_{\ex},\ \fraku\in\calF).
\qedhere
\]
\end{proof}

For reduced sequences, the corresponding compatibility appears in
the quiver Hecke setting of \cite{KK24}; Theorem~\ref{thm: ibox compatibility} applies to an
arbitrary finite sequence.

\smallskip

For a maximal commuting family $\calG$ of $\bfi$-boxes, let
$\Sigma^\calG=(\bfx^\calG,B^\calG)$ be the associated commutative
seed, where $\bfx^\calG=(x^\calG_\frakc\mid\frakc\in\calG)$.

\begin{proposition}
\label{lem: box classical T system}
Let $\ii=\seq{i_1,\ldots,i_l}\in I^l$ be an arbitrary finite sequence,
let $\calG$ be a maximal commuting family of $\bfi$-boxes in $[1,l]$,
and let $\calG'=\frakB_k\calG$ be obtained by a mutation-type box move
replacing $[x_+,y]$ with $[x,y_-]$.
Set $i\seteq i_x=i_y$ and $x^\calG_\emptyset\seteq1$.
The corresponding mutation from $\Sigma^\calG$ to $\Sigma^{\calG'}$
has exchange matrix $B^{\calG'}$ and exchange relation
\begin{equation}
x^{\calG}_{[x_+,y]}x^{\calG'}_{[x,y_-]}
=x^\calG_{[x,y]}x^\calG_{[x_+,y_-]}+
\prod_{\substack{j\in I\setminus\{i\},\ c_{j,i}<0\\ x(j)^+\le y(j)^-}}
\bigl(x^\calG_{[x(j)^+,y(j)^-]}\bigr)^{-c_{j,i}}.
\label{eq: box Qin classical T system}
\end{equation}
\end{proposition}

\begin{proof}
By Proposition~\ref{prop:boxmove}\eqref{it: mutation box-move}
and~\eqref{eq: box from envelopes}, the box $[x,y]$ belongs to
$\calG\cap\calG'$. For $k>1$, the boxes preceding the mutated box
form a maximal commuting family in $[x+1,y-1]$. Every nonempty box
$[x_+,y_-]$ or $[x(j)^+,y(j)^-]$ commutes with this family and hence
belongs to $\calG\cap\calG'$ by maximality.

By~\eqref{eq: first box move negative column}, the negative entries
of the mutation column are $c_{j,i}$ at $[x(j)^+,y(j)^-]$ for the
indicated colors $j$. By~\eqref{eq: 1 and 1} and
Proposition~\ref{prop: ibox pair box move}, the positive entries
are $1$ at $[x,y]$ and, when nonempty, $[x_+,y_-]$.
All other entries vanish, so the cluster mutation formula gives
\eqref{eq: box Qin classical T system}.
For $k=1$, all interior boxes are empty, and the same argument
applies with the empty-box convention.
\end{proof}

\begin{remark}
For finite symmetric type, the exchange relation in
the proposition is the specialization at $q^{1/2}=1$ of the categorical $T$-system in
\cite{KKOP24A,KKOP25M,KK24}. 
The proposition gives the same classical exchange relation for arbitrary finite sequences and symmetrizable Cartan data through the
exchange matrices of $\bfi$-box seeds. 
\end{remark}

\subsection{Folding of compatible pairs}

We compare the compatible pairs associated with an arbitrary
finite sequence $\ii$ and its lift $\hii$. We first establish the
folding formulas for the exchange and $\Lambda$-matrices, and then
extend them to maximal commuting $\bfi$-box families.
Throughout this subsection, we refer to the data associated with
$\ii$ and $\hii$ as the \emph{downstairs} and \emph{upstairs} data,
respectively.

Retain the notation of Subsections~\ref{subsec: seq lift}--\ref{subsec: lambda pairing}, and identify the position
$\kappa(k,r)$ of $\hii$ with $(k,r)$.

\begin{proposition}\label{prop: GLS pair folding}
For $1\le k,p\le l$ and $1\le t\le e_{i_p}$, we have
\[
b^\ii_{k,p} = \sum_{r=1}^{e_{i_k}} b^\hii_{(k,r),(p,t)}  \quad (p\in J^\ii_{\ex}),
\qtq
\Lambda^\ii_{k,p} = d_{i_k}\sum_{r=1}^{e_{i_k}} \Lambda^\hii_{(k,r),(p,t)}.
\]
In particular, both right-hand sides are independent of the choice of
$t$.
\end{proposition}

\begin{proof}
Fix $t$. We first prove the assertion for the exchange matrices.
If $k=p_+$, then
\[
b^\hii_{(k,r),(p,t)} =
\begin{cases}
-1 ,&r=t,\\
0,&r\neq t,
\end{cases}
\qtqh
\sum_{r=1}^{e_{i_k}} b^\hii_{(k,r),(p,t)} =-1=b^\ii_{k,p}.
\]
Similarly, if $p=k_+$, the fiber sum is $1=b^\ii_{k,p}$.

For the vertical relations, the blockwise lift preserves the strict
interlacing conditions, and
\[
\sum_{r=1}^{e_{i_k}} \hc_{\im_{k,r},\im_{p,t}} = c_{i_k,i_p}.
\]
Hence the corresponding fiber sum is respectively
$c_{i_k,i_p}$ or $-c_{i_k,i_p}$, as in
\eqref{eq: GLS exchange matrix}. In all remaining cases both sides
vanish.

For the $\Lambda$-matrices, set
\[
\hcalE_{k,r} \seteq \sum_{\substack{1\le a\le k\\ i_a=i_k}}\hve_{a,r}.
\]
By \eqref{eq: big epsilon} and
Lemma~\ref{lem: lambda pairing folding},
\[
\Lambda^\ii_{k,p} = \lambda^\ii(\calE_k,\calE_p)
= d_{i_k}\sum_{r=1}^{e_{i_k}} \lambda^\hii(\hcalE_{k,r},\hcalE_{p,t})
= d_{i_k}\sum_{r=1}^{e_{i_k}} \Lambda^\hii_{(k,r),(p,t)}.
\qedhere
\]
\end{proof}

\begin{remark}
The exchange-matrix part of Proposition~\ref{prop: GLS pair folding} is reminiscent of the
classical folding of exchange matrices; see, for example,
\cite[\S4.4]{FWZ2}, where global foldability is also considered.
To our knowledge, the simultaneous folding of the exchange and
compatible $\Lambda$-matrices for GLS seeds associated with sequences,
as above, has not appeared before.
\end{remark}

\begin{lemma}\label{lem: row fiber weak sign}
For $1\le k\le l$, $p\in J^\ii_{\ex}$, and
$1\le t\le e_{i_p}$, the integers
\[
b^\hii_{\kappa(k,r),\kappa(p,t)}   \qquad   (1\le r\le e_{i_k})
\]
have the same weak sign. Consequently,
\[
\sum_{r=1}^{e_{i_k}} \bigl[b^\hii_{\kappa(k,r),\kappa(p,t)}\bigr]_+
= \bigl[b^\ii_{k,p}\bigr]_+,
\qtq
\sum_{r=1}^{e_{i_k}}  \bigl[-b^\hii_{\kappa(k,r),\kappa(p,t)}\bigr]_+
= \bigl[-b^\ii_{k,p}\bigr]_+.
\]
\end{lemma}

\begin{proof}
The weak-sign assertion follows directly from the case-by-case
description in the proof of Proposition~\ref{prop: GLS pair folding}.
Together with
\[
  b^\ii_{k,p}   =   \sum_{r=1}^{e_{i_k}}   b^\hii_{\kappa(k,r),\kappa(p,t)},
\]
it gives the two identities involving positive parts.
\end{proof}

For the exchange matrices, the weak-sign property in Lemma~\ref{lem: row fiber weak sign}
corresponds to the sign condition in $\Gamma$-admissibility; see
\cite[Definition~4.4.1 and Lemma~4.4.4]{FWZ2}.

\smallskip

We next describe the lifting of maximal commuting families.
Assume first that $l\ge2$.  
Let $\frakC=\seq{\frakc_1,\ldots,\frakc_l}$ be an admissible chain of extent $[1,l]$, with
$\calF=\calF(\frakC)$ and associated  $\{\calL,\calR\}$-sequence{} $\frakT=\seq{T_1,\ldots,T_{l-1}}$.
For $1\le k\le l$, let $p_k$ be the unique position such that
$\{p_k\}=\tfrakc_k\setminus\tfrakc_{k-1}$ $(\tfrakc_0=\emptyset)$.
Thus $p_k=\eff_\calF(\frakc_k)$. We call
\[
\bfp(\frakC)\seteq\seq{p_1,\ldots,p_l} \ \ \text{ the \emph{effective-end sequence} of $\frakC$}.
\]

For $T\in\{\calL,\calR\}$, let $\oT$ denote the opposite symbol; i.e., $\ocalL=\calR$ and $\ocalR=\calL$.
We define the \emph{canonical lift} $\hfrakT$ of $\frakT$ by
\begin{align} \label{eq: hfrakT}
\hfrakT \seteq \seq{ \oT_1^{\,e_{i_{p_1}}-1},  T_1^{\,e_{i_{p_2}}}, T_2^{\,e_{i_{p_3}}},  \cdots T_{l-1}^{\,e_{i_{p_l}}} }.
\end{align}
Thus the initial effective end $p_1$ is replaced by an internal block
of length $e_{i_{p_1}}-1$ with the opposite orientation, while, for
$2\le k\le l$, the step adjoining $p_k$ is replaced by a constant
block of length $e_{i_{p_k}}$ with the same orientation.

Since $\bfp(\frakC)$ is a permutation of $[1,l]$, the sequence
$\hfrakT$ has length
\[
\sum_{k=1}^l e_{i_{p_k}}-1=N_l-1.
\]

Let
$\hfrakC \seteq\frakC(\hii,\hfrakT)$
be the corresponding admissible chain of $\bfi$-boxes, and set
$\hcalF\seteq\calF(\hfrakC)$.

For $l=1$, set
\[
\hfrakT\seteq\seq{\calR^{\,e_{i_1}-1}}.
\]
The canonical lift of $\frakC=\seq{[1]}$ is then
\[
\hfrakC=\seq{[\kappa(1,1)],\ldots,[\kappa(1,e_{i_1})]},
\qquad
\hcalF\seteq\calF(\hfrakC).
\]
For the empty sequence, the chain, the family, and their lifts are empty.

\begin{example}\label{ex: lifted LR A3 C2}
Consider the folding $A_3\to C_2$ with $I_1=(1,3)$, $I_2=(2)$,
and let
\[
\ii=(1,2,1,2),
\qquad
\hii=(1,3,2,1,3,2).
\]

First take $\frakT=(\calR,\calR,\calR)$.
Then
\[
\bfp(\frakC)=(1,2,3,4).
\]
Since $p_1=1$ and $e_{i_{p_1}}=e_1=2$,  the block replacing $p_1${}
contributes the internal symbol $\ocalR=\calL$.
The remaining steps are lifted by \eqref{eq: hfrakT}, and hence \[
\hfrakT = \seq{\calL,\calR,\calR,\calR,\calR}.
\]
Notice that $\seq{\calR,\calR,\calR,\calR,\calR}$
 determines the same lifted family, with the first two singleton boxes
in the opposite order{}, but
 the rule in \eqref{eq: hfrakT} selects $\hfrakT$ above.

Now take $\frakT=\seq{\calL,\calR,\calR}$.
Then
\[
\bfp(\frakC)=\seq{2,1,3,4}.
\]
Here $p_1=2$ and $e_{i_{p_1}}=e_2=1$, so  the first block in
\eqref{eq: hfrakT} is empty{}. Since
\[
\seq{p_2,p_3,p_4}=\seq{1,3,4},
\]
\eqref{eq: hfrakT} gives \[
\hfrakT = \seq{\calL,\calL,\calR,\calR,\calR}.
\]
\end{example}

By construction,  each effective end $p$ of $\frakC$  is replaced by the
 positions
\[
\{\kappa(p,r)\mid 1\le r\le e_{i_p}\}
\]
in the effective-end sequence of $\hfrakC$.
We therefore define the \emph{effective-end folding map}
\[
\pi_{\eff}:
\{\text{effective ends of }\hcalF\}
\longrightarrow [1,l]
\ \ \text{ by } \ \
\pi_{\eff}(\kappa(p,r))\seteq p.
\]
Thus
\begin{equation}\label{eq: effective end fiber}
\bigl|\pi_{\eff}^{-1}(p)\bigr|=e_{i_p}
\qquad (1\le p\le l).
\end{equation}

Since the effective-end maps
\[
\eff_\calF:\calF\longrightarrow[1,l]
\qtq
\eff_{\hcalF}:\hcalF\longrightarrow[1,N_l]
\]
are bijections, we define
\begin{equation}\label{eq: box folding map}
\pi_{\bx}^{\frakC}:
\hcalF\longrightarrow\calF
\ \ \text{ by } \ \
\pi_{\bx}^{\frakC}
\seteq
\eff_\calF^{-1}\circ\pi_{\eff}\circ\eff_{\hcalF}.
\end{equation}

\begin{example}\label{ex: box folding A5 C3}
Consider the folding $A_5\to C_3$ with
\[
I_1=(1,5),\qquad I_2=(2,4),\qquad I_3=(3).
\]
Let $\ii=\seq{1,2,3,2,1,3,2,3}$,
so that
\[
\hii=\seq{1,5,2,4,3,2,4,1,5,3,2,4,3}.
\]

Consider the maximal commuting families
\[
\calF=
\{[1,5],[2,4],[3],[4], [5],[3,6],[2,7],[3,8] \}
\]
and
\[
\hcalF= \{ [1,8],[2,9],[3,6],[4,7],[5],[6],[7], [8],[9],[5,10],[3,11],[4,12],[5,13] \},
\]
ordered by their effective ends. Thus
\[
\begin{array}{c|cccccccc}
\frakc
&[1,5]&[2,4]&[3]&[4]&[5]&[3,6]&[2,7]&[3,8]\\
\hline
\eff_\calF(\frakc)
&1&2&3&4&5&6&7&8
\end{array}
\]
and
\[
\begin{array}{c|ccccccccccccc}
\hfrakc
&[1,8]&[2,9]&[3,6]&[4,7]&[5]&[6]&[7]
&[8]&[9]&[5,10]&[3,11]&[4,12]&[5,13]\\
\hline
\eff_{\hcalF}(\hfrakc)
&1&2&3&4&5&6&7&8&9&10&11&12&13 .
\end{array}
\]

The effective-end folding is
\[
\begin{array}{cccccccc}
\{1,2\}&\{3,4\}&\{5\}&\{6,7\}&
\{8,9\}&\{10\}&\{11,12\}&\{13\}\\
\downarrow&\downarrow&\downarrow&\downarrow&
\downarrow&\downarrow&\downarrow&\downarrow\\
1&2&3&4&5&6&7&8 .
\end{array}
\]
Hence the box folding map
\[
\pi_{\bx}^{\frakC}
=
\eff_\calF^{-1}\circ\pi_{\eff}\circ\eff_{\hcalF}
\]
is given by
\[
\begin{aligned}
\{[1,8],[2,9]\}&\longmapsto[1,5],&
\{[3,6],[4,7]\}&\longmapsto[2,4],\\
\{[5]\}&\longmapsto[3],&
\{[6],[7]\}&\longmapsto[4],\\
\{[8],[9]\}&\longmapsto[5],&
\{[5,10]\}&\longmapsto[3,6],\\
\{[3,11],[4,12]\}&\longmapsto[2,7],&
\{[5,13]\}&\longmapsto[3,8].
\end{aligned}
\]
\end{example}

For $\frakc\in\calF$, we define the \emph{box fiber} over $\frakc$ by
\[
\Fib^{\frakC}(\frakc) \seteq (\pi_{\bx}^{\frakC})^{-1}(\frakc)
\subset\hcalF.
\]
We say that the $\bfi$-boxes in $\Fib^{\frakC}(\frakc)$ fold to $\frakc$. By \eqref{eq: effective end fiber},
\[
\bigl|\Fib^{\frakC}(\frakc)\bigr|=e_{i_{\frakc}}.
\]

\begin{lemma}\label{lem: big epsilon box fiber}
Let $\frakc=[a,b]\in\calF$, and write
$\Fib^{\frakC}(\frakc) = \seq{\hfrakc_1,\ldots,\hfrakc_{e_{i_{\frakc}}}}$
in  increasing order of their effective ends{}. Then
\[
\calE_{\hfrakc_r}
=
\sum_{\substack{a\le p\le b\\ i_p=i_{\frakc}}}
\hve_{p,r}
\qquad
(1\le r\le e_{i_{\frakc}}).
\]
In particular,
\[
\hfrakc_r=[\kappa(a,r),\kappa(b,r)]
\qquad(1\le r\le e_{i_{\frakc}}).
\] 
\end{lemma}

\begin{proof}
The lemma follows directly from the fiberwise construction of
$\hfrakC$.
By \eqref{eq: big epsilon box}, the endpoints of a box are the
positions of the first and last nonzero coordinates of its box vector.
\end{proof}

\begin{lemma}\label{lem: G action box fibers}
The $\Gamma$-action on the effective ends of $\hcalF$ induces an action on
$\hcalF$ by
\[
\ttg \cdot \hfrakc  \seteq \eff_{\hcalF}^{-1} \bigl(\ttg \cdot \eff_{\hcalF}(\hfrakc) \bigr) \qquad (\ttg \in \Gamma).
\]
For every $\frakc\in\calF$ and $\hfrakc\in\Fib^{\frakC}(\frakc)$,
we have
\[
\Gamma \cdot \hfrakc=\Fib^{\frakC}(\frakc).
\]
Moreover, for $\hfrakc,\hfrakd\in\hcalF$ and $\ttg\in \Gamma$,
\[
b^{\hcalF}_{\ttg \cdot\hfrakc ,\ttg \cdot\hfrakd } = b^{\hcalF}_{\hfrakc,\hfrakd} \qtq
\Lambda^{\hcalF}_{\ttg \cdot\hfrakc,\ttg \cdot\hfrakd } = \Lambda^{\hcalF}_{\hfrakc,\hfrakd}.
\]
\end{lemma}

\begin{proof}
By construction,
\[
\eff_{\hcalF}\bigl(\Fib^{\frakC}(\frakc)\bigr) = \pi_{\eff}^{-1}\bigl(\eff_\calF(\frakc)\bigr),
\]
and the latter is a single $\Gamma$-orbit. Hence
\[
\Gamma\cdot \hfrakc=\Fib^{\frakC}(\frakc).
\]
The remaining assertions follow from the $\Gamma$-invariance of the
upstairs Cartan data and of $\lambda^\hii$.
\end{proof}

\begin{theorem}\label{thm: ibox pair folding}
Let $\frakc,\frakd\in\calF$,  and write
$\Fib^{\frakC}(\frakc) =\seq{\hfrakc_1,\ldots,\hfrakc_{e_{i_\frakc}}}$. For any $\hfrakd_t\in\Fib^{\frakC}(\frakd)$, we have
\[
b^\calF_{\frakc,\frakd} = \sum_{\hfrakc\in\Fib^{\frakC}(\frakc)}b^{\hcalF}_{\hfrakc,\hfrakd_t},
\qtq
\Lambda^\calF_{\frakc,\frakd} = d_{i_{\frakc}} \sum_{\hfrakc\in\Fib^{\frakC}(\frakc)} \Lambda^{\hcalF}_{\hfrakc,\hfrakd_t}.
\]
In particular, both right-hand sides are independent of the choice of $\hfrakd_t$.
\end{theorem}

\begin{proof}
We first prove the assertion for $B$.  The construction of $\hfrakC$ preserves the effective-end and endpoint conditions
in \eqref{eq: full ibox exchange matrix}.

Suppose first that $b^\calF_{\frakc,\frakd}<0$. If this entry comes
from a horizontal relation, then, for the fixed
$\hfrakd_t\in\Fib^{\frakC}(\frakd)$, exactly one
$\hfrakc\in\Fib^{\frakC}(\frakc)$ satisfies the corresponding
upstairs horizontal relation. Its contribution is $-1$, while all
other members of the fiber contribute $0$.

If $b^\calF_{\frakc,\frakd}$ comes from a vertical relation, then the
same effective-end and endpoint conditions hold fiberwise upstairs.
Hence
\[
\sum_{\hfrakc\in\Fib^{\frakC}(\frakc)} b^{\hcalF}_{\hfrakc,\hfrakd_t}
= \sum_{\hfrakc\in\Fib^{\frakC}(\frakc)} \hc_{i_{\hfrakc},i_{\hfrakd_t}}
= c_{i_{\frakc},i_{\frakd}} = b^\calF_{\frakc,\frakd},
\]
where the second equality is the definition of the quotient Cartan
matrix. If $b^\calF_{\frakc,\frakd}=0$, none of the corresponding
upstairs horizontal or vertical relations occurs, and every summand
is zero.

The case $b^\calF_{\frakc,\frakd}>0$ follows from the negative case
by skew-symmetrizability, together with $d_i e_i=m$ and the
$\Gamma$-invariance of the upstairs exchange matrix. Therefore,
\[
b^\calF_{\frakc,\frakd} = \sum_{\hfrakc\in\Fib^{\frakC}(\frakc)} b^{\hcalF}_{\hfrakc,\hfrakd_t}.
\]

Write
$\Fib^{\frakC}(\frakc) = \seq{\hfrakc_1,\ldots,\hfrakc_{e_{i_\frakc}}}$
in the fixed fiber order. By Lemma~\ref{lem: lambda pairing folding} and Lemma~\ref{lem: big epsilon box fiber},
\begin{align*}
\Lambda^\calF_{\frakc,\frakd} = \lambda^\ii(\calE_\frakc,\calE_\frakd)
= d_{i_\frakc} \sum_{r=1}^{e_{i_\frakc}}  \lambda^\hii(\calE_{\hfrakc_r},\calE_{\hfrakd_t}) =
d_{i_\frakc} \sum_{r=1}^{e_{i_\frakc}} \Lambda^{\hcalF}_{\hfrakc_r,\hfrakd_t}.
\end{align*}
The independence of $\hfrakd_t$ follows from the $\Gamma$-invariance above.
\end{proof}

\subsection{Folding of $g$-vectors}

We introduce the folding map on exponent lattices that will be used
later for the comparison of $g$-vectors.

Define
\begin{equation}\label{eq: g vector folding map}
\Fold^\ii:\mathbb Z^{J^\hii}\longrightarrow\mathbb Z^{J^\ii}
\ \ \text{ by } \ \
\bigl(\Fold^\ii(\hbfg)\bigr)_k \seteq \sum_{r=1}^{e_{i_k}}\hbfg_{(k,r)} \qquad (k\in J^\ii).
\end{equation}
We denote by
\[
\Fold^\ii_{\ex}: \mathbb Z^{J^\hii_{\ex}}\longrightarrow \mathbb Z^{J^\ii_{\ex}}
\]
its restriction to the exchangeable coordinates.

\begin{proposition}\label{prop: Fold B intertwining}
We have
\[
\Fold^\ii\circ B^\hii = B^\ii\circ\Fold^\ii_{\ex}.
\]
\end{proposition}

\begin{proof}
For the standard basis vector $\ve_{(p,t)}$,
Proposition~\ref{prop: GLS pair folding} gives
\[
\bigl(\Fold^\ii(B^\hii\ve_{(p,t)})\bigr)_k
= \sum_{r=1}^{e_{i_k}} b^\hii_{(k,r),(p,t)}
= b^\ii_{k,p}.
\]
On the other hand,
\[
B^\ii\Fold^\ii_{\ex}(\ve_{(p,t)}) = B^\ii\ve_p.
\]
Hence the assertion follows.
\end{proof}

\section{Braid moves and folding} \label{sec: braid}

 In this section, we lift finite braid moves to sequences of $2$-
and $3$-moves and realize these sequences by mutations and
permutations. We prove the folding identities for the exchange
and commutation matrices, cluster variables, and $g$-vectors
before and after each finite braid move and its lifted sequence.{}

\subsection{A motivating example: the $4$-move $A_3\to C_2$} \label{subsec: motivating ex braid}

We begin with the $4$-move
\[
\ii=1212\longrightarrow \ii'=2121 \text{ of type $C_2$.}
\]
With our fixed ordering of the fibers,
\[
\hii=132132,\qquad \hii'=213213.
\]
 The positions in $\hii$ corresponding to the successive positions of $\ii$ are
\[
\seq{1,2},\ \seq3,\ \seq{4,5},\ \seq6,
\]
 while the positions in $\hii'$ corresponding to those in $\ii'$ are
\[
\seq1,\ \seq{2,3},\ \seq4,\ \seq{5,6}.
\]

Let $\bfx=(x_1,x_2,x_3,x_4)$ and $\bfy=(y_1,\ldots,y_6)$ be  the initial downstairs and upstairs clusters, respectively.
Define
\[
\phi:\bbZ[y_1^{\pm1},\ldots,y_6^{\pm1}]
\longrightarrow
\bbZ[x_1^{\pm1},\ldots,x_4^{\pm1}]
\]
by
\[
\phi(y_1)=\phi(y_2)=x_1,\qquad \phi(y_3)=x_2,\qquad \phi(y_4)=\phi(y_5)=x_3,\qquad \phi(y_6)=x_4.
\]

\smallskip
\noindent\emph{The downstairs calculation.}
The initial compatible pair is
\[
B^{(0)} = \left(
\begin{array}{cc}
0&-2\\
1&0\\ \hline
-1&2\\
0&-1
\end{array}
\right),
\qquad
\Lambda^{(0)}
=
\left(
\begin{array}{cc|cc}
0&-2&0&0\\
2&0&0&0\\ \hline
0&0&0&0\\
0&0&0&0
\end{array}
\right).
\]
The exchangeable positions are $1,2$.
Let $P_p$ denote the permutation matrix exchanging positions $p$ and $p+1$.
The $4$-move is realized by
\[
\mu_1,\ \mu_2,\ \mu_1,\ P_3,\ P_1
\]
from left to right (see Subsection~\ref{subsec: mu+per} below).

Starting with $\bfx^{(0)}=(x_1,x_2,x_3,x_4)$, the first mutation gives
\[
x_1a=x_2+x_3, \qquad a=\frac{x_2+x_3}{x_1} \qtq
B^{(1)} =
\begin{pmatrix}
0&2\\
-1&0\\
1&0\\
0&-1
\end{pmatrix}.
\]
The second mutation gives
\[
x_2b=a^2+x_4, \qquad b=\frac{x_1^2x_4+(x_2+x_3)^2}{x_1^2x_2}
\quad \text{with} \quad
B^{(2)}
=
\begin{pmatrix}
0&-2\\
1&0\\
1&0\\
-1&1
\end{pmatrix}.
\]
Finally,
\[
ac=bx_3+x_4, \qquad c=\frac{x_1^2x_4+x_2x_3+x_3^2}{x_1x_2}
\qtq
B^{(3)} =
\begin{pmatrix}
0&2\\
-1&0\\
-1&0\\
1&-1
\end{pmatrix}.
\]
Applying $P_3$ and then $P_1$, we obtain
\[
\bfx'=(b,c,x_4,x_3),
\]
that is,
\[
x'_1=\frac{x_1^2x_4+(x_2+x_3)^2}{x_1^2x_2}, \qquad x'_2=\frac{x_1^2x_4+x_2x_3+x_3^2}{x_1x_2},
\qquad x'_3=x_4,\quad x'_4=x_3.
\]
The final compatible pair is
\[
B^{\ii'}
=
\left(
\begin{array}{cc}
0&-1\\
2&0\\ \hline
-1&1\\
0&-1
\end{array}
\right),
\qquad
\Lambda^{\ii'}
=
\left(
\begin{array}{cc|cc}
0&-2&0&0\\
2&0&0&0\\ \hline
0&0&0&0\\
0&0&0&0
\end{array}
\right).
\]
\smallskip
\noindent\emph{The upstairs calculation.}
The initial compatible pair associated with
$\hii=132132$ is
\[
B^{\hii,(0)}
=
\left(
\begin{array}{ccc}
0&0&-1\\
0&0&-1\\ \hline
1&1&0\\ \hline
-1&0&1\\
0&-1&1\\ \hline
0&0&-1
\end{array}
\right),
\qquad
\Lambda^{\hii,(0)}
=
\left(
\begin{array}{cc|c|cc|c}
0&0&-1&1&-1&0\\
0&0&-1&-1&1&0\\ \hline
1&1&0&0&0&0\\ \hline
-1&1&0&0&0&0\\
1&-1&0&0&0&0\\ \hline
0&0&0&0&0&0
\end{array}
\right).
\]
The exchangeable positions are $1,2,3$.  The following sequence of $2$- and $3$-moves connects $132132$ to $213213$:
\[
\underline{13}2132\to3\underline{121}32\to321\underline{232}
\to 32\underline{13}23 \to \underline{323}123
\to23\underline{212}3\to2\underline{31}213\to213213,
\]
with corresponding transformations
\[
P_1,\quad \mu_2P_3,\quad \mu_4P_5,\quad
P_3,\quad \mu_1P_2,\quad \mu_3P_4,\quad P_2
\]
from left to right (see Subsection~\ref{subsec: mu+per} below).

Starting with $\bfy^{(0)}=\seq{y_1,\ldots,y_6}$, the first move gives
\[
\bfy^{(1)}=\seq{y_2,y_1,y_3,y_4,y_5,y_6}.
\]
Next,
\[
y_1u=y_3+y_4, \qquad u=\frac{y_3+y_4}{y_1},
\]
and after $P_3$,
\[
\bfy^{(3)}=\seq{y_2,u,y_4,y_3,y_5,y_6},
\qquad
B^{\hii,(3)}
=
\begin{pmatrix}
0&0&-1\\
0&0&1\\
0&1&0\\
1&-1&0\\
-1&0&1\\
0&0&-1
\end{pmatrix}.
\]

The next mutation gives
\[
y_3v=uy_5+y_2y_6,
\qquad
v=\frac{y_1y_2y_6+y_3y_5+y_4y_5}{y_1y_3},
\]
and after $P_5$ and then $P_3$,
\[
\bfy^{(6)}=\seq{y_2,u,v,y_4,y_6,y_5},
\qquad
B^{\hii,(6)}
=
\begin{pmatrix}
0&-1&1\\
1&0&-1\\
-1&1&0\\
0&1&0\\
0&-1&1\\
0&0&-1
\end{pmatrix}.
\]

Next,
\[
y_2w=u+v,
\qquad
w= \frac{y_1y_2y_6+y_3^2+y_3y_4+y_3y_5+y_4y_5}{y_1y_2y_3},
\]
and after $P_2$,
\[
\bfy^{(8)}=\seq{w,v,u,y_4,y_6,y_5},
\qquad
B^{\hii,(8)}
=
\begin{pmatrix}
0&-1&1\\
1&0&0\\
-1&0&0\\
0&0&1\\
0&1&-1\\
0&-1&0
\end{pmatrix}.
\]

Finally,
\[
uz=wy_4+y_6,
\qquad
z=
\frac{y_1y_2y_6+y_3y_4+y_4y_5}{y_2y_3}.
\]
After $P_4$ and $P_2$, we obtain
\[
\bfy'=\seq{w,z,v,y_6,y_4,y_5},
\ \
B^{\hii'} =
\left(
\begin{array}{ccc}
0&-1&-1\\ \hline
1&0&0\\
1&0&0\\ \hline
-1&1&1\\ \hline
0&-1&0\\
0&0&-1
\end{array}
\right)
\qtq
\Lambda^{\hii'}
=
\left(
\begin{array}{c|cc|c|cc}
0&-1&-1&0&0&0\\ \hline
1&0&0&0&1&-1\\
1&0&0&0&-1&1\\ \hline
0&0&0&0&0&0\\ \hline
0&-1&1&0&0&0\\
0&1&-1&0&0&0
\end{array}
\right).
\]
Thus
\[
\begin{aligned}
y'_1&= \frac{y_1y_2y_6+y_3^2+y_3y_4+y_3y_5+y_4y_5}{y_1y_2y_3},\\
y'_2&= \frac{y_1y_2y_6+y_3y_4+y_4y_5}{y_2y_3},
\qquad
y'_3= \frac{y_1y_2y_6+y_3y_5+y_4y_5}{y_1y_3},\\
y'_4&=y_6,\qquad y'_5=y_4,\qquad y'_6=y_5.
\end{aligned}
\]

\smallskip
\noindent\emph{Comparison under folding.}
We now apply the initial folding homomorphism $\phi$ to  the final upstairs cluster $\bfy'$. The positions in $\hii'$ corresponding to the successive positions of $\ii'$ are
\[
\seq1,\qquad \seq{2,3},\qquad \seq4,\qquad \seq{5,6}.
\]

\begin{proposition}\label{prop: motivating C2 variable folding}
With the notation above,
\[
\phi(y'_1)=x'_1,\qquad \phi(y'_2)=\phi(y'_3)=x'_2,\qquad \phi(y'_4)=x'_3,\qquad  \phi(y'_5)=\phi(y'_6)=x'_4.
\]
 Thus $\phi$ sends every variable in $\bfy'$ to the corresponding variable in $\bfx'$.
\end{proposition}

\begin{proof}
Using
\[
\phi(y_1)=\phi(y_2)=x_1,\qquad \phi(y_3)=x_2,\qquad  \phi(y_4)=\phi(y_5)=x_3,\qquad \phi(y_6)=x_4,
\]
we obtain
\begin{align*}
\phi(y'_1) = \frac{x_1^2x_4+x_2^2+x_2x_3+x_2x_3+x_3^2}{x_1^2x_2}
= \frac{x_1^2x_4+(x_2+x_3)^2}{x_1^2x_2} =x'_1.
\end{align*}
Notice that the coefficient $2$ in $x'_1$ comes from the two distinct
upstairs monomials $y_3y_4$ and $y_3y_5$, both of which fold to
$x_2x_3$.
Similarly,
\[
\phi(y'_2) = \frac{x_1^2x_4+x_2x_3+x_3^2}{x_1x_2} = \phi(y'_3) = x'_2.
\]
For the remaining variables,
\[
\phi(y'_4)=\phi(y_6)=x_4=x'_3, \qquad \phi(y'_5)=\phi(y_4)=x_3=x'_4, \qquad \phi(y'_6)=\phi(y_5)=x_3=x'_4.
\]
Hence the assertion follows.
\end{proof}

\subsection{Finite braid moves and their unfoldings}
\label{subsec: braid move}
We describe how a finite braid move lifts to the upstairs sequence. Let $i,j\in I_\g$ be distinct vertices such that
\[
c_{i,j}c_{j,i}<4,
\]
and write
\[
I_i=\seq{\im_{i,1},\ldots,\im_{i,e_i}}, \qquad I_j=\seq{\im_{j,1},\ldots,\im_{j,e_j}}
\]
for the corresponding $\Gamma$-orbits in $I_\hg$.

Since the $\Gamma$-action is admissible, there are no edges between two
vertices in the same orbit. Since $\hg$ is simply-laced, the Cartan
folding formula~\eqref{eq: quotient Cartan matrix} gives
\[
-c_{i,j} = \#\{\im\in I_i\mid \im\sim\jm\},
\qquad
-c_{j,i} = \#\{\jm\in I_j\mid \jm\sim\im\},
\]
for $\im\in I_i$ and $\jm\in I_j$. Hence
\[
c_{i,j}c_{j,i}\in\{0,1,2,3\}, \qquad m_{i,j}\in\{2,3,4,6\}.
\]

\begin{lemma}\label{lem: connected finite braid configurations}
Suppose that $0<c_{ij}c_{ji}<4$. Then every connected component of
the subdiagram on $I_i\sqcup I_j$ consists of one vertex on one side
joined to respectively one, two, or three vertices on the other side,
according as
\[
c_{ij}c_{ji}=1,\ 2,\ \text{or }3.
\]
Moreover, all connected components are isomorphic.
\end{lemma}

\begin{proof}
Every vertex of $I_j$ has $-c_{ij}$ neighbors in $I_i$, and every
vertex of $I_i$ has $-c_{ji}$ neighbors in $I_j$. Since
\[
c_{ij}c_{ji}\in\{1,2,3\},
\]
one of these two numbers is $1$. Hence each connected component has
one vertex on one side and respectively one, two, or three vertices
on the other side. The last assertion follows from the transitivity
of the $\Gamma$-action on $I_i$.
\end{proof}

 The subdiagram on $I_i\sqcup I_j$ is therefore as follows:
\begin{eqnarray} &&
\parbox{85ex}{
\bna
\item If $m_{i,j}=2$, there are no edges between $I_i$ and $I_j$.
\item If $m_{i,j}=3$, every vertex in either orbit has exactly one
neighbor in the other orbit. Hence $e_i=e_j$, and the edges give a
bijection between $I_i$ and $I_j$.
\item If $m_{i,j}=4$, up to interchanging $i$ and $j$, we have
$(c_{i,j},c_{j,i})=(-1,-2)$. Thus every vertex of $I_i$ has two
neighbors in $I_j$, while every vertex of $I_j$ has one neighbor in
$I_i$.

\item If $m_{i,j}=6$, up to interchanging $i$ and $j$, we have
$(c_{i,j},c_{j,i})=(-1,-3)$. Thus every vertex of $I_i$ has three
neighbors in $I_j$, while every vertex of $I_j$ has one neighbor in
$I_i$.
\ee
}\label{eq: braid lift classification}
\end{eqnarray}

We shall use the following terminology repeatedly. Let $\ii=\seq{i_1,\ldots,i_l}$ and let $k$ be a position.
\bnum
\item A \emph{$2$-move at $k$} replaces $(i_k,i_{k+1})=(i,j)$ by $(j,i)$, when $c_{i,j}=0$.
\item A \emph{$3$-move at $k$} replaces $(i_k,i_{k+1},i_{k+2})=(i,j,i)$ by $(j,i,j)$, when $c_{i,j}c_{j,i}=1$.
\item A \emph{$4$-move at $k$} replaces $(i,j,i,j)$ by $(j,i,j,i)$, when $c_{i,j}c_{j,i}=2$.
\item A \emph{$6$-move at $k$} replaces $(i,j,i,j,i,j)$ by $(j,i,j,i,j,i)$, when $c_{i,j}c_{j,i}=3$.
\ee
We denote an $m$-move starting at position $k$ by $(m;k)$.

\begin{example}\label{ex: E6 D43 rank two orbit}
Consider the folding
\[
E_6^{(1)}\ \
\raisebox{2.1em}{
\xymatrix@R=1ex@C=4.2ex{
&& *{\circled{2}}<3pt>
\ar@{}[r]_<{ \ 0 }
   \ar@{-}[d] \ar@/^0.5pc/@{.>}[ddrr]|{\sigma} & \\
&& *{\circled{2}}<3pt> \ar@{-}[d]
   \ar@{-}[d] \ar@/^0.5pc/@{.>}[dr]|{\sigma}
   \ar@{}[r]_<{ \ 2 } & \\
*{\circled{2}}<3pt>\ar@{-}[r]_<{  1 \ }
\ar@/^0.5pc/@{.>}[uurr]|{\sigma}
& *{\circled{2}}<3pt> \ar@{-}[r]_<{ 3 \ }
\ar@/^0.5pc/@{.>}[ur]|{\sigma}
& *{\circled{2}}<3pt> \ar@{-}[r]_<{ 6 \ } & *{\circled{2}}<3pt> \ar@{-}[r]_<{ 5 \ }
\ar@/^1pc/@{.>}[ll]|{\sigma}
& *{\circled{2}}<3pt> \ar@{}[l]^<{ \ 4  }
\ar@/^2pc/@{.>}[llll]|{\sigma}
}}
\quad
\overset{/\ang{\sigma}}{\Longrightarrow}
\quad
D_4^{(3)}\ \
\xymatrix@R=0.5ex@C=3.8ex{
*{\circled{2}}<3pt>  \ar@{-}[r]_<{0\ \ }
&
*{\circled{2}}<3pt> \ar@{<=}[r]_<{1\ \ }
&
*{\rcircled{6}}<3pt> \ar@{}[l]^>{\qquad \qquad 2}
},
\]
by the order-three diagram automorphism. With the numbering above and our fixed ordering of the fibers, we have
\[
I_0=\seq{1,0,4},\qquad I_1=\seq{3,2,5},\qquad I_2=\seq6.
\]
Thus $e_0=e_1=3$ and $e_2=1$.

The edges between $I_0$ and $I_1$ are
\[
1\sim3,\quad 0\sim2,\quad 4\sim5, \ \ \text{ while } \ \ 3\sim6,\quad 2\sim6,\quad 5\sim6.
\]
There are no other edges between distinct fibers. Hence, by
\eqref{eq: quotient Cartan matrix},
\[
\sfC=
\begin{pmatrix}
2&-1&0\\
-1&2&-3\\
0&-1&2
\end{pmatrix}.
\]
Indeed,
$c_{0,1}=c_{1,0}=-1$,  $c_{1,2}=-3$, and $c_{2,1}=-1$.

The pair $0,1$ therefore gives a $3$-move, with
$e_0=e_1=3$.
 The subdiagram on $I_0\sqcup I_1$ consists of the three disjoint edges
\[
1\sim3,\qquad 0\sim2,\qquad 4\sim5.
\]
In particular, the two sides of the braid relation lift to
\[
010 \longmapsto (1,0,4)(3,2,5)(1,0,4), \qtq 101 \longmapsto (3,2,5)(1,0,4)(3,2,5).
\]

On the other hand, the pair $1,2$ gives a $6$-move. Here
\[
(c_{2,1},c_{1,2})=(-1,-3),
\]
and the unique vertex $6\in I_2$ is adjacent to all three vertices
$3,2,5\in I_1$. Thus this single affine folding exhibits both the
$3$-move with $e_0=e_1=3$ and the $6$-move with
$(e_1,e_2)=(3,1)$.
\end{example}

Let $\ii=\seq{i_1,\ldots,i_l}$, and suppose that a finite braid move
of length $m=m_{i,j}$ starts at position $k$.
We call the subword $\seq{i_k,\ldots,i_{k+m-1}}$ the
\emph{braid block} of this move and its blockwise lift in $\hii$
the \emph{lifted braid block}.
The subword $I_{i_k}$ starts at position $\kappa(k,1)$ in $\hii$.

\begin{proposition}\label{prop: canonical finite braid lift}
Let $\ii'$ be obtained from $\ii$ by a finite braid move at position
$k$. Then $\hii'$ is obtained from $\hii$ by a finite sequence of
$2$- and $3$-moves  depending only on the subdiagram on $I_i\sqcup I_j$ described in
\eqref{eq: braid lift classification}. The reverse braid move is
obtained by reversing the same sequence.
\end{proposition}

\begin{proof}
Since the lift is blockwise, it is enough to consider the subword
involving $I_i$ and $I_j$. In what follows,
$\overset{(2)}{\longrightarrow}$ and
$\overset{(3)}{\longrightarrow}$ denote sequences of $2$-moves and
$3$-moves, respectively, while
$\overset{(2,3)}{\longrightarrow}$ denotes a sequence involving both.

\smallskip
\noindent
{\rm (i) $m_{i,j}=2$.}
Write
\[
I_i=\seq{a_1,\ldots,a_p}, \qquad I_j=\seq{b_1,\ldots,b_q}.
\]
Since $a_r\not\sim b_s$ for all $r,s$,
\[
\seq{a_1\cdots a_p}\seq{b_1\cdots b_q} \overset{(2)}{\longrightarrow} \seq{b_1\cdots b_q}\seq{a_1\cdots a_p}.
\]

\smallskip
\noindent
{\rm (ii) $m_{i,j}=3$.}
Put $r=e_i=e_j$. Using $2$-moves inside the fibers, we write
$I_i=\seq{a_1,\ldots,a_r}$ and
$I_j=\seq{b_1,\ldots,b_r}$
so that
\[
a_s\sim b_s,
\qquad
a_s\not\sim b_t\quad(s\neq t).
\]
Then
\begin{align*}
\seq{a_1\cdots a_r}\seq{b_1\cdots b_r}\seq{a_1\cdots a_r}
&\overset{(2)}{\longrightarrow}
\seq{a_1b_1a_1}\cdots\seq{a_rb_ra_r}
\overset{(3)}{\longrightarrow}
\seq{b_1a_1b_1}\cdots\seq{b_ra_rb_r} \allowdisplaybreaks \\
&\overset{(2)}{\longrightarrow}
\seq{b_1\cdots b_r}\seq{a_1\cdots a_r}\seq{b_1\cdots b_r}.
\end{align*}

\smallskip
\noindent
{\rm (iii) $m_{i,j}=4$.}
Up to interchanging $i$ and $j$,
\[
(c_{i,j},c_{j,i})=(-1,-2).
\]
Writing $r=e_i$, we have $e_j=2r$. Using $2$-moves inside the fibers,
we write
$I_i=\seq{c_1,\ldots,c_r}$,
and
$I_j=\seq{a_1,b_1,\ldots,a_r,b_r}$,
where
\[
a_s\sim c_s\sim b_s, \qquad a_s\not\sim b_s,
\]
and vertices with distinct subscripts are mutually nonadjacent. Hence
\[
I_iI_jI_iI_j
\overset{(2)}{\longrightarrow}
\prod_{s=1}^r c_sa_sb_sc_sa_sb_s.
\]
For each $s$, suppressing the subscript, we have
\begin{align*}
cabcab
\overset{(2)}{\longrightarrow}cbacab
\overset{(3)}{\longrightarrow}cbcacb
\overset{(3)}{\longrightarrow}bcbacb
\overset{(2)}{\longrightarrow}bcabcb
\overset{(3)}{\longrightarrow}bcacbc
\overset{(3)}{\longrightarrow}bacabc
\overset{(2)}{\longrightarrow}abcabc .
\end{align*}
Therefore
\[
I_iI_jI_iI_j \overset{(2,3)}{\longrightarrow} \prod_{s=1}^r a_sb_sc_sa_sb_sc_s \overset{(2)}{\longrightarrow} I_jI_iI_jI_i.
\]

\smallskip
\noindent
{\rm (iv) $m_{i,j}=6$.}
Up to interchanging $i$ and $j$,
\[
(c_{i,j},c_{j,i})=(-1,-3).
\]
Writing $r=e_i$, we have $e_j=3r$. Using $2$-moves inside the fibers,
we write
$I_i=\seq{d_1,\ldots,d_r}$,
and
$I_j=\seq{a_1,b_1,c_1,\ldots,a_r,b_r,c_r}$,
where
\[
a_s\sim d_s,\qquad b_s\sim d_s,\qquad c_s\sim d_s,
\]
the vertices $a_s,b_s,c_s$ are pairwise nonadjacent, and vertices with
distinct subscripts are mutually nonadjacent. Thus
\[
I_iI_jI_iI_jI_iI_j \overset{(2)}{\longrightarrow} \prod_{s=1}^r d_sa_sb_sc_sd_sa_sb_sc_sd_sa_sb_sc_s.
\]

For one component, suppressing the subscript, the fixed sequence
\begin{align*}
&(2;5),(2;6),(2;5),(3;3),(3;7),(3;5),(2;4),(2;7),\allowdisplaybreaks\\
&(3;2),(3;4),(2;6),(3;8),(3;10),(2;9),\allowdisplaybreaks\\
&(3;7),(3;5),(2;4),(2;3),(3;9),(2;7),\allowdisplaybreaks\\
&(3;5),(3;1),(3;3),(2;8),(2;5),
(3;6),(3;8),(3;4),(2;6),(2;7),(2;6)
\end{align*}
gives
\[
abcdabcdabcd \overset{(2,3)}{\longrightarrow} dabcdabcdabc.
\]
Hence, by reversing this sequence,
\[
dabcdabcdabc \overset{(2,3)}{\longrightarrow} abcdabcdabcd.
\]
Applying this to each component gives
\begin{align*}
I_iI_jI_iI_jI_iI_j &\overset{(2,3)}{\longrightarrow}
\prod_{s=1}^r
a_sb_sc_sd_sa_sb_sc_sd_sa_sb_sc_sd_s \overset{(2)}{\longrightarrow}
I_jI_iI_jI_iI_jI_i.
\end{align*}

Thus, in every case,  the two lifted words  are connected
by a finite sequence of $2$- and $3$-moves. Since both moves are
involutive, the reverse braid move is obtained by reversing the same
sequence.
\end{proof}

\subsection{Mutation--permutation realizations and 
folding
} 
\label{subsec: mu+per}

We recall the mutation--permutation realizations of finite braid moves.
The $2$- and $3$-moves were treated in \cite{FHOO2}, and the $4$- and
$6$-moves in \cite{LO24}  (see also \cite{SW21}).

As in Subsection~\ref{subsec: motivating ex braid}, mutations and permutations below are applied
from left to right.

\begin{proposition} \cite{FHOO2,LO24}\label{prop: mutation permutation braid moves}
Let $\ii'$ be obtained from $\ii$ by an $(m;k)$-move.
Then the compatible pair $(B^{\ii'},\Lambda^{\ii'})$ is obtained from
$(B^\ii,\Lambda^\ii)$ by the following mutation--permutation sequence$:$
\[
\begin{array}{c|c}
m & \text{mutation--permutation sequence} \\ \hline
2 & P_k \\[1mm]
3 & \mu_k,\ P_{k+1} \\[1mm]
4 & \mu_k,\ \mu_{k+1},\ \mu_k,\ P_k,\ P_{k+2} \\[1mm]
6 & \mu_k,\ \mu_{k+1},\ \mu_{k+2},\ \mu_k,\ \mu_{k+3},\
\mu_{k+1},\ \mu_k,\ \mu_{k+2},\ \mu_{k+3},\ \mu_k,\
P_k,\ P_{k+2},\ P_{k+4}.
\end{array}
\]
\end{proposition}

For an $(m;k)$-move, denote by $\sfT$ the corresponding
mutation--permutation sequence in Proposition~\ref{prop: mutation permutation braid moves}.
 For the sequence of $2$- and $3$-moves{}
\[
\hii=\hjj^{(0)}\longrightarrow\hjj^{(1)} \longrightarrow\cdots\longrightarrow\hjj^{(s)}=\hii',
\]
let $\hsfT_a$ be the mutation--permutation sequence
corresponding to
\[
\hjj^{(a-1)}\longrightarrow\hjj^{(a)}   \qquad(1\le a\le s).
\]
We call the sequence from $\hii$ to $\hii'$ a \emph{lifted path}
of the $(m;k)$-move and $\seq{\hsfT_1,\ldots,\hsfT_s}$ its
mutation--permutation realization.

\begin{proposition}\label{prop: endpoint braid folding}
Let $\ii'$ be obtained from $\ii$ by an $(m;k)$-move, and let
\[
\hii=\hjj^{(0)}   \longrightarrow\hjj^{(1)}   \longrightarrow\cdots\longrightarrow  \hjj^{(s)}=\hii'
\]
 be the sequence of $2$- and $3$-moves in
Proposition~\ref{prop: canonical finite braid lift}.{}
Then
\[
\xymatrix@C=60pt@R=35pt{
(B^\hii,\Lambda^\hii)   \ar[r]^{\hsfT_1,\ldots,\hsfT_s}   \ar[d]_{\mathrm{Fold}} & (B^{\hii'},\Lambda^{\hii'})   \ar[d]^{\mathrm{Fold}}\\
(B^\ii,\Lambda^\ii)   \ar[r]_{ \sfT} & (B^{\ii'},\Lambda^{\ii'})
}
\]
where the vertical arrows are given by
Proposition~\ref{prop: GLS pair folding}.
\end{proposition}

\begin{proof}
By Proposition~\ref{prop: mutation permutation braid moves},
\[
(B^\ii,\Lambda^\ii) \xrightarrow{\ (m;k)\text{-move}\ } (B^{\ii'},\Lambda^{\ii'}).
\]
Applying the same proposition successively to
\[
\hjj^{(a-1)}\longrightarrow\hjj^{(a)} \quad(1\le a\le s)
\ \text{ gives } \
(B^\hii,\Lambda^\hii) \xrightarrow{\ \hsfT_1,\ldots,\hsfT_s\ } (B^{\hii'},\Lambda^{\hii'}).
\]
 The two vertical arrows are given by{}
Proposition~\ref{prop: GLS pair folding}, applied to
$(\hii,\ii)$ and $(\hii',\ii')$, respectively.
\end{proof}

\subsection{Folding of cluster variables and $g$-vectors}

We work in the commutative setting throughout this subsection.
Let $\Sigma^\ii$ and $\Sigma^\hii$ be  the initial seeds associated with 
$\ii$ and $\hii$, respectively. Define the algebra homomorphism between the seed Laurent polynomial rings
\[
\phi^\ii:\calL(\Sigma^\hii)\longrightarrow\calL(\Sigma^\ii),  \qquad
\hx_{(k,r)}\longmapsto x_k  \quad  (1\le k\le l,\ 1\le r\le e_{i_k}).
\]
We show that $\phi^\ii$ identifies  the cluster variables of $\Sigma^{\hii'}$ and $\Sigma^{\ii'}$
obtained by the mutation--permutation sequences of
Proposition~\ref{prop: endpoint braid folding}.

\begin{remark}\label{rem: larger fibers braid}
 For $m=3,4,6$, the connected components of the subdiagram on $I_i\sqcup I_j$
are classified in Lemma~\ref{lem: connected finite braid configurations}.
The case $m=2$ involves only permutations.
The additional affine calculations are recorded in
Appendix~\ref{app: affine braid folding}, and
Proposition~\ref{prop: reduction to connected configurations}
reduces \eqref{eq: folding of cluster var} to the calculation for each connected component.
The indefinite example in \eqref{eq: indefinite folding} introduces
no additional case of a finite braid move.
\end{remark}

\noindent
{}
 By the Laurent phenomenon, the resulting cluster variables belong
to $\calL(\Sigma^\hii)$ and $\calL(\Sigma^\ii)$, respectively.  We shall prove
\[
\xymatrix@C=65pt@R=38pt{
\bigl(\hx_{\kappa(p,r)}\bigr)_{p,r} \ar[r]^{\hsfT_1,\ldots,\hsfT_s}  \ar[d]_{\phi^\ii}
& \bigl(\hx'_{\kappa'(p,r)}\bigr)_{p,r}  \ar[d]^{\phi^\ii}
\\
\bigl(x_p\bigr)_p \ar[r]_{\sfT} & \bigl(x'_p\bigr)_p .
}
\]
Notice that the right vertical arrow is  the \emph{same} Laurent homomorphism $\phi^\ii$, rather than
the Laurent homomorphism $\phi^{\ii'}$ defined from $\Sigma^{\hii'}$. Equivalently, we shall prove
\begin{align} \label{eq: folding of cluster var}
\phi^\ii\bigl(\hx'_{\kappa'(p,r)}\bigr)=x'_p \qquad  (1\le p\le l,\ 1\le r\le e_{i'_p}),
\end{align}
 where $\kappa'$ is defined for $\ii'$ by
\eqref{eq: kappa k,r}.

\medskip
\noindent
$\bullet$ \ \textbf{The $2$-move.}
 By Proposition~\ref{prop: canonical finite braid lift}, the lifted path consists only of $2$-moves, each of which acts on seeds by a permutation.

\begin{proposition}\label{prop: cluster variable folding 2 move}
Suppose that $\ii'$ is obtained from $\ii$ by a $2$-move.
Then \eqref{eq: folding of cluster var} holds.
\end{proposition}

\smallskip
\noindent
$\bullet$\ \textbf{The $3$-move.}
We first consider the case $e_i=e_j=2$.  This case occurs in the folding
$A_{2n-1}\to C_n$ away from the multiple edge.
Write $I_i=(a,A)$ and $I_j=(b,B)$, where $a\sim b$, $A\sim B$, $a\nsim B$ and $A\nsim b$.
Thus the braid move
\[
iji\longrightarrow jij \quad \text{lifts to} \quad  aAbBaA\longrightarrow bBaAbB.
\]
By Proposition~\ref{prop: canonical finite braid lift}, the lifted
$3$-move decomposes into $e_i=e_j$ simply-laced $3$-moves, together
with $2$-moves rearranging the fibers.  Thus, in the present case
$e_i=e_j=2$, exactly two $3$-moves occur.
For the calculation below, we choose the following  sequence of $2$- and $3$-moves:{}
\begin{equation} \label{eq: seq of lifted 3-move 2,2}
\begin{aligned}
a\underline{Ab}BaA
&\longrightarrow abA\underline{Ba}A
 \longrightarrow ab\underline{Aa}BA
 \longrightarrow \underline{aba}ABA\\
&\longrightarrow bab\underline{ABA}
 \longrightarrow ba\underline{bB}AB
 \longrightarrow b\underline{aB}bAB
 \longrightarrow bBa\underline{bA}B
 \longrightarrow bBaAbB.
\end{aligned}
\end{equation}
Here the underlined subwords indicate the successive braid moves.
More precisely, the successive moves are
\begin{align*}
 &(2;\kappa(k,1)+1),\quad (2;\kappa(k,1)+3),\quad (2;\kappa(k,1)+2),\quad (3;\kappa(k,1)),\allowdisplaybreaks\\
 &(3;\kappa(k,1)+3),\quad (2;\kappa(k,1)+2),\quad (2;\kappa(k,1)+1),\quad (2;\kappa(k,1)+3).
\end{align*}
The two $3$-moves occur at positions $\kappa(k,1)$ and $\kappa(k,1)+3$; all the
remaining moves are $2$-moves.

We compare the exchange relations associated with these mutations.
For the three positions of the downstairs braid block, write
\[
x_1 \seteq x_k,\qquad x_2 \seteq x_{k+1},\qquad x_3 \seteq x_{k+2}.
\]
By Proposition~\ref{prop: mutation permutation braid moves}, the
downstairs $3$-move is realized by $\seq{\mu_k,\ P_{k+1}}$ in order.
 Let $D^+$ and $D^-$ denote the contributions of the rows indexed by
$r\notin\{k,k+1,k+2\}$ to the positive and negative exchange monomials
at the mutation $\mu_k$.  Then
\[
x_1x_1'=D^+x_2+D^-x_3,
\]
and hence, after $P_{k+1}$,
\[
x'_1=\frac{D^+x_2+D^-x_3}{x_1},   \qquad   x'_2=x_3,   \qquad   x'_3=x_2.
\]

  For the lifted braid block $aAbBaA$, write 
\[
\hx_1,\hx_2,\hx_3,\hx_4,\hx_5,\hx_6
\]
for the corresponding cluster variables, in this order. Thus
\[
\phi^\ii(\hx_1)=\phi^\ii(\hx_2)=x_1,\qquad \phi^\ii(\hx_3)=\phi^\ii(\hx_4)=x_2,\qquad \phi^\ii(\hx_5)=\phi^\ii(\hx_6)=x_3.
\]

For $a=1,2$, let $U_a^+$ and $U_a^-$ denote the contributions of the   rows with indices outside $[\kappa(k,1),\kappa(k+2,2)]$ 
to the positive and negative exchange monomials  at the $a$th upstairs mutation, computed after all preceding
mutations and permutations in \eqref{eq: seq of lifted 3-move 2,2}.{} Thus $U_1^\pm$ are computed after the three preceding $2$-moves,
whereas $U_2^\pm$ are computed after all preceding moves, including
the first $3$-move.

\begin{lemma}\label{lem: coefficient folding three move}
With the notation above, we have
\[
\phi^\ii(U_1^\pm)=\phi^\ii(U_2^\pm)=D^\pm.
\]
\end{lemma}

\begin{proof}
Set
\[
k_a \seteq \kappa(k,a)\qquad(a=1,2).
\]
Let $\widehat B^{(a)}$ denote the exchange matrix immediately before
the $a$-th upstairs mutation, and let $\hk_a$ be the
corresponding mutation position.

By \eqref{eq: seq of lifted 3-move 2,2}, the two mutation vertices are
the positions of $k_1$ and $k_2$ after the preceding permutations,
respectively. Since $k_1$
and $k_2$ belong to the same fiber,
\[
b^\hii_{k_1,k_2}=0.
\]
Therefore the first mutation does not change the $k_2$-column:
\[
(\mu_{k_1}B^\hii)_{r,k_2}
= b^\hii_{r,k_2} + [b^\hii_{r,k_1}]_+[b^\hii_{k_1,k_2}]_+ - [-b^\hii_{r,k_1}]_+[-b^\hii_{k_1,k_2}]_+
= b^\hii_{r,k_2}.
\]
 All moves other than the two $3$-moves in
\eqref{eq: seq of lifted 3-move 2,2} are $2$-moves and only permute
rows and columns.  Hence, up to these permutations, the mutation
column of $\widehat B^{(a)}$ is the initial $k_a$-column of $B^\hii$,
for $a=1,2$.

Now fix $r\notin\{k,k+1,k+2\}$, and put $r_s \seteq \kappa(r,s)$  for $(1\le s\le e_{i_r})$.
By Proposition~\ref{prop: GLS pair folding} and
Lemma~\ref{lem: row fiber weak sign},
\[
\sum_{s=1}^{e_{i_r}}   \bigl[\pm b^\hii_{r_s,k_a}\bigr]_+   =   \bigl[\pm b^\ii_{r,k}\bigr]_+   \qquad(a=1,2).
\]
Following the corresponding rows and columns through
\eqref{eq: seq of lifted 3-move 2,2}, we obtain
\[
\sum_{s=1}^{e_{i_r}} \bigl[\pm\hb^{(a)}_{\hr_s,\hk_a} \bigr]_+
=\bigl[\pm b^\ii_{r,k}\bigr]_+,   \qquad (a=1,2),
\]
where $\hr_s$ is the position of $r_s$ immediately
before the $a$-th mutation.

Since no position   outside $[\kappa(k,1),\kappa(k+2,2)]$  is mutated, the corresponding
cluster variables are only permuted, and hence
\[
\phi^\ii(\hx_{\hr_s})=x_r.
\]
It follows that  the contribution of the rows indexed by $\hr_s$ $(1\le s\le e_{i_r})$ to
$U_a^\pm$ satisfies
\[
\phi^\ii\left(\prod_{s=1}^{e_{i_r}}\hx_{\hr_s}^{[\pm\hb^{(a)}_{\hr_s,\hk_a}]_+} \right)
= x_r^{[\pm b^\ii_{r,k}]_+}.
\]
Multiplying over all
$r\notin\{k,k+1,k+2\}$ gives
\[
\phi^\ii(U_a^\pm)=D^\pm  \qquad(a=1,2). \qedhere
\]
\end{proof}

\begin{proposition}\label{prop: cluster variable folding three move 22}
Suppose that $\ii'$ is obtained from $\ii$ by a $3$-move
\[
iji\longrightarrow jij \quad \text{with $e_i=e_j=2$.}
\]
Then \eqref{eq: folding of cluster var} holds.
\end{proposition}

\begin{proof}
 With the local notation introduced above, the final downstairs variables are
\[
x'_1=\frac{D^+x_2+D^-x_3}{x_1}, \qquad x'_2=x_3, \qquad x'_3=x_2.
\]

We now follow the two mutations in
\eqref{eq: seq of lifted 3-move 2,2}.
For the first $3$-move, the exchange relation is
\[
\hx_1\hx'_1=U_1^+\hx_3+U_1^-\hx_5,
\]
whereas for the second $3$-move it is
\[
\hx_2\hx'_2=U_2^+\hx_4+U_2^-\hx_6.
\]
All other steps in
\eqref{eq: seq of lifted 3-move 2,2}
are $2$-moves and hence only permute the cluster variables. Therefore
\[
\hx'_3=\hx_5,\qquad   \hx'_4=\hx_6,\qquad   \hx'_5=\hx_3,\qquad   \hx'_6=\hx_4.
\]
Thus
\[
\hx'_1= \frac{U_1^+\hx_3+U_1^-\hx_5}{\hx_1},
\qquad
\hx'_2= \frac{U_2^+\hx_4+U_2^-\hx_6}{\hx_2}.
\]

By Lemma~\ref{lem: coefficient folding three move},
\[
\phi^\ii(U_1^\pm) =\phi^\ii(U_2^\pm) = D^\pm.
\]
Using
\[
\phi^\ii(\hx_1)=\phi^\ii(\hx_2)=x_1,\qquad
\phi^\ii(\hx_3)=\phi^\ii(\hx_4)=x_2,\qquad
\phi^\ii(\hx_5)=\phi^\ii(\hx_6)=x_3,
\]
we obtain
\[
\begin{aligned}
\phi^\ii(\hx'_1) =   \frac{D^+x_2+D^-x_3}{x_1} =x'_1,\quad \phi^\ii(\hx'_2) =  \frac{D^+x_2+D^-x_3}{x_1} =x'_1.
\end{aligned}
\]
For the remaining variables,
\[
\phi^\ii(\hx'_3)=\phi^\ii(\hx'_4)=x_3=x'_2,
\qtq
\phi^\ii(\hx'_5)=\phi^\ii(\hx'_6)=x_2=x'_3.
\]
 Therefore
\eqref{eq: folding of cluster var} holds.
\end{proof}

We next consider the case $e_i=e_j=1$.
Write $I_i=(a)$ and $I_j=(b)$, with $a\sim b$.  Then
\[
  iji\longrightarrow jij \quad\text{lifts simply to}\quad  aba\longrightarrow bab.
\]
Thus the lifted path consists of a single $3$-move.
 This case occurs, for example, in the
folding $D_{n+1}\to B_n$.

\begin{proposition}\label{prop: cluster variable folding three move 11}
Suppose that $\ii'$ is obtained from $\ii$ by a $3$-move with
$e_i=e_j=1$. Then \eqref{eq: folding of cluster var} holds.
\end{proposition}

\begin{proof}
This follows immediately from
Proposition~\ref{prop: GLS pair folding} and
Lemma~\ref{lem: row fiber weak sign}.
\end{proof}

The additional affine case $e_i=e_j=3$ is recorded in
Appendix~\ref{app: affine braid folding}, where the lifted path and
the identities for the final cluster variables yielding \eqref{eq: folding of cluster var}
are given.

\medskip
\noindent
$\bullet$\ \textbf{The $4$-move.}
We first consider the case $e_i=1$ and $e_j=2$.
 This case occurs in the folding $A_{2n-1}\longrightarrow C_n$ at the multiple edge; in this case,
$i=n$ and $j=n-1$. Write
$I_i=(c)$ and $I_j=(a,b)$,
where $a\sim c\sim b$ and $a\nsim b$.
Thus the braid move
\[
ijij\longrightarrow jiji \quad\text{lifts to}\quad   cabcab\longrightarrow abcabc.
\]

By Proposition~\ref{prop: canonical finite braid lift}, the lifted
$4$-move is obtained by four simply-laced $3$-moves together with
three $2$-moves.  For the calculation below, we choose the following
 sequence of $2$- and $3$-moves:{}
\begin{equation}\label{eq: seq of lifted 4-move 1,2}
\begin{aligned}
c\underline{ab}cab &\longrightarrow cb\underline{aca}b \longrightarrow \underline{cbc}acb \longrightarrow bc\underline{ba}cb
\longrightarrow bca\underline{bcb} \longrightarrow b\underline{cac}bc  \longrightarrow \underline{ba}cabc \longrightarrow abcabc.
\end{aligned}
\end{equation}
Here the underlined subwords indicate the successive braid moves.
More precisely, the successive moves are
\begin{align*}
 &(2;\kappa(k,1)+1),\quad (3;\kappa(k,1)+2),\quad
 (3;\kappa(k,1)),\quad (2;\kappa(k,1)+2),\allowdisplaybreaks\\
 &(3;\kappa(k,1)+3),\quad (3;\kappa(k,1)+1),\quad
 (2;\kappa(k,1)).
\end{align*}
Thus the four $3$-moves occur at positions
\[
  \kappa(k,1)+2,\qquad \kappa(k,1),\qquad \kappa(k,1)+3,\qquad \kappa(k,1)+1,
\]
and the remaining three moves are $2$-moves.

\begin{lemma}\label{lem: GLS four move outside rows 12}
For $r\notin\{k,k+1,k+2,k+3\}$, the entries in columns $k,k+1$ satisfy
\begin{equation}\label{eq: four move outside row patterns 12}
\bigl(b^\ii_{r,k},b^\ii_{r,k+1}\bigr)=
\begin{cases}
(-2,1),&r_+=k+1,\\
(1,0),&r_+=k,\\
(0,0),&\text{otherwise}.
\end{cases}
\end{equation}
In particular, for the mutation sequence
\[
\seq{\mu_k,\mu_{k+1},\mu_k},
\]
the relevant initial data are
\[
\rho_{\rm down}(r) \seteq  \bigl(b^\ii_{r,k},b^\ii_{r,k+1}\bigr)   \in  \{(-2,1),(1,0),(0,0)\}.
\]
We call the tuple of entries of a row in a specified ordered list
of columns its \emph{row pattern}. The three patterns above define
the classes $\calP,\calQ,\calZ$ of row indices, respectively.
The entries of  the rows of $B^\hii$ indexed by $\kappa(r,s)$ $(1\le s\le e_{i_r})$ in columns $\kappa(k,1),\kappa(k+1,1),\kappa(k+1,2)$, in this order, are
\begin{equation}\label{eq: four move active row data 12}
\begin{array}{c|ccc}
 & \calP & \calQ & \calZ\\ \hline
\rho_{\rm down}  & (-2,1) & (1,0) & (0,0)\\[1mm]
\rho_{\rm up}  & (-1,1,0),\ (-1,0,1)  & (1,0,0)  & (0,0,0).
\end{array}
\end{equation}
Here $\calP$ occurs only when $i_r=j$ and hence $e_j=2$, whereas $\calQ$ occurs only when $i_r=i$ and hence $e_i=1$.
The class $\calZ$ consists of the remaining cases with
$\rho_{\rm down}(r)=(0,0)$;  every row indexed by $\kappa(r,s)$ $(1\le s\le e_{i_r})$ has entries $(0,0,0)$ in the three columns specified above. \end{lemma}

\begin{proof}
By the defining formula for the GLS exchange matrix, a nonzero entry
in columns $k,k+1$ can occur only from ~{\rm (i)}~\eqref{eq: GLS exchange matrix} or from~{\rm (ii)}~\eqref{eq: GLS exchange matrix}. Since
\[
k_+=k+2, \qquad (k+1)_+=k+3,
\]
a direct case-by-case check, according as $r<k$ or $r>k+3$, gives
\eqref{eq: four move outside row patterns 12}.

For the mutation sequence $\seq{\mu_k,\mu_{k+1},\mu_k}$,
the initial entries in columns $k,k+1$ are \[
(-2,1),\qquad (1,0),\qquad (0,0),
\]
which give the classes $\calP,\calQ,\calZ$, respectively.

  For the subword $cabcab$ of $\hii$,  the simply-laced GLS formula gives,
for a row of class $\calP$, the two triples in columns $\kappa(k,1),\kappa(k+1,1),\kappa(k+1,2)$ \[
(-1,1,0),  \qquad  (-1,0,1),
\]
corresponding to  the rows indexed by $\kappa(r,1)$ and $\kappa(r,2)$.  For $r\in\calQ$, the row indexed by $\kappa(r,1)$ has the triple \[
(1,0,0).
\]
 For $r\in\calZ$, every row indexed by $\kappa(r,s)$ $(1\le s\le e_{i_r})$ has the triple \[
(0,0,0). \qedhere
\]
\end{proof}

\begin{example}\label{ex: folding and fiber}
For the folding $A_5\longrightarrow C_3$, consider the reduced sequence
\[
\ii=\seq{2,1,3,2,3,2}
\]
of type $C_3$, whose $4$-braid block $\seq{3,2,3,2}$ starts at $k=3$.
Its exchange matrix is
\[
B^\ii=
\begin{pmatrix}
0&-2&1\\
1&0&0\\
1&0&-1\\
-1&2&0\\
0&-1&1\\
0&0&-1
\end{pmatrix}.
\]
The columns of $B^\ii$ are indexed by $1,3,4$, in this order.
Thus, for the outside row $r=1$,
\[
\rho_{\rm down}(1)=(-2,1),
\]
so that $r=1$ belongs to the class $\calP$.

Its lift is
\[
\hii=\seq{2,4,1,5,3,2,4,3,2,4},
\]
with exchange matrix
\[
B^\hii=
\left(
\begin{array}{cc|c|cc}
 0& 0&-1& 1& 0\\
 0& 0&-1& 0& 1\\ \hline
 1& 0& 0& 0& 0\\
 0& 1& 0& 0& 0\\ \hline
 1& 1& 0&-1&-1\\ \hline
-1& 0& 1& 0& 0\\
 0&-1& 1& 0& 0\\ \hline
 0& 0&-1& 1& 1\\ \hline
 0& 0& 0&-1& 0\\
 0& 0& 0& 0&-1
\end{array}
\right).
\]
The columns of $B^\hii$ are indexed by $1,2,5,6,7$, in this order.
 The rows corresponding to $r=1$ are indexed by $1,2$, 
and we consider columns $5,6,7$.
The entries in these columns are
\[
(-1,1,0),\qquad (-1,0,1),
\]
as in \eqref{eq: four move active row data 12}. Moreover,
\[
b^\hii_{1,5}+b^\hii_{2,5}=-2,   \qquad   b^\hii_{1,6}+b^\hii_{2,6}   =   b^\hii_{1,7}+b^\hii_{2,7}=1,
\]
which recovers
\[
\rho_{\rm down}(1)=(-2,1)
\]
under folding.
\end{example}

We introduce notation for the outside variables in the $4$-move. Set
\[
\calO \seteq \{r\notin\{k,k+1,k+2,k+3\}\},
\]
and
\[
\begin{aligned}
\calO_{\calP} &\seteq \{r\in\calO\mid \rho_{\rm down}(r)=(-2,1)\}, \quad
\calO_{\calQ} &\seteq \{r\in\calO\mid \rho_{\rm down}(r)=(1,0)\}.
\end{aligned}
\]
Set
\[
P \seteq \prod_{r\in\calO_{\calP}}x_r, \qquad  Q \seteq \prod_{r\in\calO_{\calQ}}x_r,
\]
where an empty product is understood to be $1$.

By Lemma~\ref{lem: GLS four move outside rows 12},  for $r\in\calO_{\calP}$, let $\hr_A,\hr_B$ be the indices
$\kappa(r,1),\kappa(r,2)$, ordered so that the corresponding triples in
  \eqref{eq: four move active row data 12}  are
\[
  (-1,1,0),  \qquad   (-1,0,1),
\]
respectively. Set
\[
\hP_A \seteq \prod_{r\in\calO_{\calP}}\hx_{\hr_A},
\qquad
\hP_B \seteq \prod_{r\in\calO_{\calP}}\hx_{\hr_B}.
\]
 Similarly, writing $\hr=\kappa(r,1)$ for $r\in\calO_{\calQ}$, set
\[
\hQ \seteq \prod_{r\in\calO_{\calQ}}\hx_{\hr}.
\]
Then
\[
\phi^\ii(\hP_A) = \phi^\ii(\hP_B) = P, \qquad \phi^\ii(\hQ)=Q.
\]

\begin{lemma}\label{lem: tropical four move 12}
In the setting of Lemma~\ref{lem: GLS four move outside rows 12}, let
$D_s^\pm$ $(s=1,2,3)$ be the factors from  rows indexed by $r\in\calO$ in the
successive downstairs exchange relations. Let $U_a^\pm$ $(a=1,\ldots,4)$
be those for the successive upstairs mutations, determined immediately
before each mutation after all preceding mutations and permutations
 in the mutation--permutation sequence realizing \eqref{eq: seq of lifted 4-move 1,2}. Then
\[
(D_1^+,D_1^-)=(Q,P^2),\qquad
(D_2^+,D_2^-)=(1,P),\qquad
(D_3^+,D_3^-)=(1,Q),
\]
whereas
\begin{align*}
(U_1^+,U_1^-)&=(\hP_A,1),&
(U_2^+,U_2^-)&=(\hQ,\hP_B),\allowdisplaybreaks\\
(U_3^+,U_3^-)&=(1,1),&
(U_4^+,U_4^-)&=(\hQ,\hP_A).
\end{align*}
\end{lemma}

\begin{proof}
For  $r\in\calO$, set
\[
(\alpha,\beta)\seteq(b^\ii_{r,k},b^\ii_{r,k+1}),
\qquad
B^{(0)}=B^\ii,\quad B^{(1)}=\mu_kB^{(0)},\quad B^{(2)}=\mu_{k+1}B^{(1)},
\]
and
\[
\eta_1=b^{(0)}_{r,k}, \qquad \eta_2=b^{(1)}_{r,k+1}, \qquad \eta_3=b^{(2)}_{r,k}.
\]
Thus $\eta_a$ is the entry in row $r$ of the exchange column used at the $a$-th mutation.

Since  the final permutations  occur only after the three mutations,
they play no role here. By~\eqref{eq: GLS exchange matrix},
\[
b^{(0)}_{k,k+1}=-1,\qquad b^{(0)}_{k+1,k}=2,
\]
and hence $\eta_1=\alpha$.

For the second mutation happening at $k+1$,
\[
\begin{aligned}
\eta_2 =b^{(1)}_{r,k+1} =\beta  +[\alpha]_+[b^{(0)}_{k,k+1}]_+  -[-\alpha]_+[-b^{(0)}_{k,k+1}]_+ =\beta-[-\alpha]_+.
\end{aligned}
\]
Note that we have
$b^{(1)}_{r,k}=-\alpha$ and $b^{(1)}_{k+1,k} =-b^{(0)}_{k+1,k} =-2$.
Therefore, for the third mutation happening at $k$,
\[
\eta_3 =b^{(2)}_{r,k} =-\alpha +[\eta_2]_+[b^{(1)}_{k+1,k}]_+ -[-\eta_2]_+[-b^{(1)}_{k+1,k}]_+ =-\alpha-2[-\eta_2]_+.
\]

Thus Lemma~\ref{lem: GLS four move outside rows 12} gives
\begin{align}\label{eq: PQZ matrix}
\begin{array}{c|ccc}
 &\eta_1&\eta_2&\eta_3\\ \hline
\calP&-2&-1&0\\
\calQ&1&0&-1\\
\calZ&0&0&0
\end{array}
\end{align}
and hence, reading the three mutation columns,
\[
(D_1^+,D_1^-)=(Q,P^2),\qquad
(D_2^+,D_2^-)=(1,P),\qquad
(D_3^+,D_3^-)=(1,Q).
\]

 We now consider the mutation--permutation sequence realizing
\eqref{eq: seq of lifted 4-move 1,2}.
 Number the six positions of the subword $cabcab$
by $1,\ldots,6$ from left to right. With respect to these local
positions, the four mutations occur at
\[
(m_1,m_2,m_3,m_4)=\seq{3,1,4,2}.
\]
 For $r\in\calO$ and $1\le s\le e_{i_r}$, put
$r_s=\kappa(r,s)$ and let
\[
(\alpha,\beta,\gamma)
\]
denote its entries in local columns $1,2,3$ of the initial
  exchange matrix, as in \eqref{eq: four move active row data 12}.
Let $\hB^{[a]}$ be the exchange matrix immediately before the
$a$-th mutation, after all preceding mutations and permutations, and
let $r_s^{[a]}$ be the current position of $r_s$. Set
\[
\heta_a(r_s)\seteq
\hB^{[a]}_{r_s^{[a]},m_a}.
\]

The words immediately before the four mutations 
in
\eqref{eq: seq of lifted 4-move 1,2} are
\begin{align*}
W_1&=\seq{c,b,\underline{\boxed{a},c,a},b},&
W_2&=\seq{\underline{\boxed{c},b,c},a,c,b},\allowdisplaybreaks\\
W_3&=\seq{b,c,a,\underline{\boxed{b},c,b}},&
W_4&=\seq{b,\underline{\boxed{c},a,c},b,c}.
\end{align*}
Thus the mutation colors are $\seq{a,c,b,c}$, respectively.

For the word 
$W_1$,~\eqref{eq: GLS exchange matrix} gives
\[
\hB^{[1]}_{3,1}=1, \qquad \hB^{[1]}_{3,2}=0.
\]
Since the initial permutation gives
\[
(\alpha,\beta,\gamma)
\longmapsto
(\alpha,\gamma,\boxed{\beta}),
\]
and the mutation happens at the third position,
\[
\heta_1(r_s)=\beta.
\]
Mutation at $m_1=3$ gives
\begin{align}
\hB'_{r_s,1} &= \alpha+[\beta]_+[1]_+ -[-\beta]_+[-1]_+ = \alpha+[\beta]_+ \quad (\because \hB^{[1]}_{3,1}=1), \nonumber\\
\hB'_{r_s,2} &= \gamma+[\beta]_+[0]_+ -[-\beta]_+[0]_+ = \gamma, \qquad (\because \hB^{[1]}_{3,2}=0) \label{eq: gamma in 2nd} \\
\hB'_{r_s,3} &=-\beta.  \label{eq: -beta in 3rd}
\end{align}

The following permutation does not affect positions $1,2,3$.
Since the second mutation happens at the first position,
\[
\heta_2(r_s) = \alpha+[\beta]_+ = \alpha+[\heta_1(r_s)]_+.
\]

At $W_2$,~\eqref{eq: GLS exchange matrix} gives
\[
\hB^{[2]}_{1,2}=-1, \qquad \hB^{[2]}_{1,3}=1.
\]
The entry $\gamma$ is still in position $2$ by
\eqref{eq: gamma in 2nd}, so after the second mutation it becomes
\[
\gamma+[\heta_2(r_s)]_+[-1]_+ -[-\heta_2(r_s)]_+[1]_+
= \gamma-[-\heta_2(r_s)]_+.
\]
After the following two permutations, the position of this entry changes as
\[
2\longmapsto3\longmapsto4.
\]
Since the third mutation happens at the fourth position,
\[
\heta_3(r_s) = \gamma-[-\heta_2(r_s)]_+.
\]

At $W_2$, the entry originating from $\beta$ has already changed to
$-\beta$ by \eqref{eq: -beta in 3rd}, and
$\hB^{[2]}_{1,3}=1$.
Thus, after the second mutation, it becomes
\[
-\beta \longmapsto -\beta+[\heta_2(r_s)]_+[1]_+ -[-\heta_2(r_s)]_+[-1]_+ = -\beta+[\heta_2(r_s)]_+ \ \ (\because \hB^{[2]}_{1,3}=1).
\]
The subsequent permutation moves this entry from position $3$
to position $2$, while the next one leaves position $2$ fixed.

At $W_3$,~\eqref{eq: GLS exchange matrix} gives
\[
\hB^{[3]}_{4,2}=1.
\]
Thus the third mutation gives
\begin{align*}
-\beta+[\heta_2(r_s)]_+ &\longmapsto -\beta+[\heta_2(r_s)]_+ +[\heta_3(r_s)]_+[1]_+ -[-\heta_3(r_s)]_+[-1]_+\allowdisplaybreaks\\
&= -\beta+[\heta_2(r_s)]_++[\heta_3(r_s)]_+.
\end{align*}
The following permutation leaves position $2$ fixed. Since the last
mutation is at $m_4=2$,
\[
\heta_4(r_s) = -\beta+[\heta_2(r_s)]_++[\heta_3(r_s)]_+.
\]

Consequently,
\begin{align*}
\heta_1(r_s)&=\beta,
&\qquad
\heta_2(r_s)&=\alpha+[\heta_1(r_s)]_+, \allowdisplaybreaks\\
\heta_3(r_s)&=\gamma-[-\heta_2(r_s)]_+,
&
\heta_4(r_s)&=-\beta+[\heta_2(r_s)]_++[\heta_3(r_s)]_+.
\end{align*}

 For $r\in\calP$, the rows indexed by $\kappa(r,1)$ and $\kappa(r,2)$ have initial triples \[
(-1,1,0), \qquad (-1,0,1),
\]
and hence
\[
\seq{\heta_1,\heta_2,\heta_3,\heta_4} = \seq{1,0,0,-1}, \qquad \seq{0,-1,0,0}.
\]
For $\calQ$,
\[
(\alpha,\beta,\gamma)=(1,0,0), \qquad \seq{\heta_1,\heta_2,\heta_3,\heta_4} = \seq{0,1,0,1}.
\]
 For every row indexed by $\kappa(r,s)$ with $r\in\calZ$ and $1\le s\le e_{i_r}$,
\[
(\alpha,\beta,\gamma)=(0,0,0), \qquad  \heta_a=0 \quad(1\le a\le4).
\]
Therefore
\[
\begin{array}{c|cccc}
 &\heta_1&\heta_2&\heta_3&\heta_4\\ \hline
\calP_A&1&0&0&-1\\
\calP_B&0&-1&0&0\\
\calQ&0&1&0&1\\
\calZ&0&0&0&0
\end{array}
\]
and hence
\begin{align*}
(U_1^+,U_1^-)&=(\hP_A,1),&
(U_2^+,U_2^-)&=(\hQ,\hP_B),\allowdisplaybreaks\\
(U_3^+,U_3^-)&=(1,1),&
(U_4^+,U_4^-)&=(\hQ,\hP_A).
\end{align*}
Hence the assertion follows.
\end{proof}

\begin{proposition}\label{prop: cluster folding four move 12}
Suppose that $\ii'$ is obtained from $\ii$ by a $4$-move
\[
ijij\longrightarrow jiji \quad\text{with $e_i=1$ and $e_j=2$.}
\]
Then \eqref{eq: folding of cluster var} holds.
\end{proposition}

\begin{proof}
Write
\[
x_1=x_k,\qquad x_2=x_{k+1},\qquad
x_3=x_{k+2},\qquad x_4=x_{k+3}.
\]
By Lemma~\ref{lem: tropical four move 12},
\[
(D_1^+,D_1^-)=(Q,P^2),\quad
(D_2^+,D_2^-)=(1,P),\quad
(D_3^+,D_3^-)=(1,Q).
\]
Hence the three downstairs exchange relations are
\[
x_1x_1^{(1)}=Qx_2^2+P^2x_3,\quad
x_2x_2^{(1)}=x_1^{(1)}+Px_4,\quad
x_1^{(1)}x_1^{(2)}=(x_2^{(1)})^2x_3+Qx_4^2.
\]
Set
\[
T \seteq Qx_2^2+P^2x_3,\qquad
W \seteq T+Px_1x_4.
\]
Then
\[
x_1^{(1)}=\frac{T}{x_1},\qquad
x_2^{(1)}=\frac{W}{x_1x_2},\qquad
x_1^{(2)}
=\frac{x_3W^2+Qx_1^2x_2^2x_4^2}{x_1x_2^2T}.
\]
After the final permutations,
\[
x_k'=\frac{W}{x_1x_2},\quad
x_{k+1}'=\frac{x_3W^2+Qx_1^2x_2^2x_4^2}{x_1x_2^2T},
\quad x_{k+2}'=x_4,\quad x_{k+3}'=x_3.
\]

  For the lifted braid block $\seq{c,a,b,c,a,b}$, write 
$\seq{\hx_1,\ldots,\hx_6}$ for the corresponding cluster variables.
Then
\[
\phi^\ii(\hx_1)=x_1,\quad
\phi^\ii(\hx_2)=\phi^\ii(\hx_3)=x_2,\quad
\phi^\ii(\hx_4)=x_3,\quad
\phi^\ii(\hx_5)=\phi^\ii(\hx_6)=x_4.
\]
The lifted path in~\eqref{eq: seq of lifted 4-move 1,2},
\begin{align*}
\boldsymbol{(}&(2;\kappa(k,1)+1),(3;\kappa(k,1)+2),
(3;\kappa(k,1)),(2;\kappa(k,1)+2), \allowdisplaybreaks \\
&(3;\kappa(k,1)+3),(3;\kappa(k,1)+1),
(2;\kappa(k,1))\boldsymbol{)},
\end{align*}
 has a mutation--permutation realization with four mutations at the current local positions $3,1,4,2$.
Using Lemma~\ref{lem: tropical four move 12}, set
\[
R \seteq U_1^+\hx_4+U_1^-\hx_1\hx_5,\quad
C \seteq U_2^+\hx_2\hx_3+U_2^-R,\quad
E \seteq U_3^+\hx_4C+U_3^-\hx_1\hx_6R.
\]
Eliminating the intermediate mutated variables  in this realization gives
\[
\hx_1'=\frac{U_4^-E+U_4^+\hx_1\hx_2\hx_3\hx_5}
{\hx_1\hx_3R},\quad
\hx_2'=\frac{C}{\hx_1\hx_2},\quad
\hx_3'=\frac{E}{\hx_1\hx_2\hx_3},
\]
and
\[
\seq{\hx_4',\hx_5',\hx_6'}=\seq{\hx_5,\hx_6,\hx_4}.
\]

Again by Lemma~\ref{lem: tropical four move 12} and
\[
\phi^\ii(\hP_A)=\phi^\ii(\hP_B)=P,
\qquad
\phi^\ii(\hQ)=Q,
\]
we have
\[
\phi^\ii(U_1^+,U_1^-)=(P,1),\quad
\phi^\ii(U_2^+,U_2^-)=(Q,P),\quad
\phi^\ii(U_3^+,U_3^-)=(1,1),\quad
\phi^\ii(U_4^+,U_4^-)=(Q,P).
\]
Consequently
\[
\phi^\ii(R)=Px_3+x_1x_4,\qquad
\phi^\ii(C)=Qx_2^2+P\phi^\ii(R)=W,\qquad
\phi^\ii(E)=x_3W+x_1x_4\phi^\ii(R).
\]
Moreover,
\[
P\phi^\ii(E)+Qx_1x_2^2x_4=W\phi^\ii(R),\qquad
T\phi^\ii(E)=x_3W^2+Qx_1^2x_2^2x_4^2,
\]
where the first follows from
$W=Qx_2^2+P\phi^\ii(R)$ and the second from
$T=Qx_2^2+P^2x_3$ and $W=T+Px_1x_4$.

Therefore
\[
\phi^\ii(\hx_1') =\frac{P\phi^\ii(E)+Qx_1x_2^2x_4}{x_1x_2\phi^\ii(R)} =\frac{W}{x_1x_2}
=x_k' =\phi^\ii(\hx_2'),
\]
and
\[
\phi^\ii(\hx_3')
=\frac{\phi^\ii(E)}{x_1x_2^2}
=\frac{x_3W^2+Qx_1^2x_2^2x_4^2}{x_1x_2^2T}
=x_{k+1}'.
\]
Finally,
\[
\phi^\ii(\hx_4')=\phi^\ii(\hx_5')=x_4=x_{k+2}',
\qquad
\phi^\ii(\hx_6')=x_3=x_{k+3}'.
\]
Thus  the final variables at positions $1,\ldots,6$ satisfy
\[
\phi^\ii(\hx_1')=\phi^\ii(\hx_2')=x_k',\quad
\phi^\ii(\hx_3')=x_{k+1}',\quad
\phi^\ii(\hx_4')=\phi^\ii(\hx_5')=x_{k+2}',\quad
\phi^\ii(\hx_6')=x_{k+3}'.
\]
Hence \eqref{eq: folding of cluster var} holds.
\end{proof}

We next consider the case $e_i=2$ and $e_j=1$.
 This case occurs in the folding
$D_{n+1}\longrightarrow B_n$ at the multiple edge; in this case,
$i=n$ and $j=n-1$. Write $I_i=(a,b)$ and $I_j=(c)$,
where $a\sim c\sim b$ and $a\nsim b$. Thus the braid move
\[
ijij\longrightarrow jiji
\quad\text{lifts to}\quad
abcabc\longrightarrow cabcab.
\]

For the calculation below, we choose the following  sequence of $2$- and $3$-moves:{}
\begin{equation}\label{eq: seq of lifted 4-move 2,1}
\begin{aligned}
\underline{ab}cabc
&\longrightarrow b\underline{aca}bc
 \longrightarrow bca\underline{cbc}
 \longrightarrow bc\underline{ab}cb\\
&\longrightarrow \underline{bcb}acb
 \longrightarrow cb\underline{cac}b
 \longrightarrow c\underline{ba}cab
\longrightarrow cabcab.
\end{aligned}
\end{equation}
More precisely, the successive moves are
\begin{align*}
&(2;\kappa(k,1)),\quad (3;\kappa(k,1)+1),\quad
(3;\kappa(k,1)+3),\quad (2;\kappa(k,1)+2),\allowdisplaybreaks\\
&(3;\kappa(k,1)),\quad (3;\kappa(k,1)+2),\quad
(2;\kappa(k,1)+1).
\end{align*}
Thus the four $3$-moves occur at positions
\[
\kappa(k,1)+1,\qquad \kappa(k,1)+3,\qquad \kappa(k,1),\qquad \kappa(k,1)+2.
\]

The following calculation is entirely parallel to the case $e_i=1$ and
$e_j=2$ above, so we record only the data and omit the proof.

\begin{lemma}\label{lem: GLS four move outside rows 21}
For $r\notin\{k,k+1,k+2,k+3\}$, 
the entries in columns $k,k+1$ satisfy
\begin{equation}\label{eq: four move outside row patterns 21}
\bigl(b^\ii_{r,k},b^\ii_{r,k+1}\bigr)=
\begin{cases}
(-1,1),&r_+=k+1,\\
(1,0),&r_+=k,\\
(0,0),&\text{otherwise}.
\end{cases}
\end{equation}
In particular,
\[
\rho_{\rm down}(r)\seteq
\bigl(b^\ii_{r,k},b^\ii_{r,k+1}\bigr)
\in\{(-1,1),(1,0),(0,0)\},
\]
defining the classes $\calN,\calP,\calZ$, respectively.
The entries of  the rows of $B^\hii$ indexed by $\kappa(r,s)$ $(1\le s\le e_{i_r})$ in columns $\kappa(k,1),\kappa(k,2),\kappa(k+1,1)$, in this order, are
\begin{equation}\label{eq: four move active row data 21}
\begin{array}{c|ccc}
 &\calN&\calP&\calZ\\ \hline
\rho_{\rm down}&(-1,1)&(1,0)&(0,0)\\[1mm]
\rho_{\rm up}&(-1,-1,1)&(1,0,0),\ (0,1,0)&(0,0,0).
\end{array}
\end{equation}
Here $\calN$ occurs only when $i_r=j$, and hence  $e_j=1$, whereas $\calP$ occurs only when $i_r=i$,
and hence $e_i=2$.
\end{lemma}

We introduce notation for the outside variables in the $4$-move. Set
\[
\calO\seteq\{r\notin\{k,k+1,k+2,k+3\}\},
\]
and
\begin{align*}
\calO_{\calN}&\seteq \{r\in\calO\mid\rho_{\rm down}(r)=(-1,1)\},\quad
\calO_{\calP}&\seteq \{r\in\calO\mid\rho_{\rm down}(r)=(1,0)\}.
\end{align*}
Set
\[
N\seteq\prod_{r\in\calO_{\calN}}x_r, \qquad P\seteq\prod_{r\in\calO_{\calP}}x_r,
\]
where an empty product is understood to be $1$.

By Lemma~\ref{lem: GLS four move outside rows 21},  for $r\in\calO_{\calN}$, put $\hr=\kappa(r,1)$ and set
\[
\hN\seteq\prod_{r\in\calO_{\calN}}\hx_{\hr}.
\]
 For $r\in\calO_{\calP}$, let $\hr_A,\hr_B$ be the indices
$\kappa(r,1),\kappa(r,2)$, ordered so that the corresponding triples in \eqref{eq: four move active row data 21}
are
\[
(1,0,0),\qquad (0,1,0),
\]
respectively, and set
\[
\hP_A\seteq\prod_{r\in\calO_{\calP}}\hx_{\hr_A}, \qquad \hP_B\seteq\prod_{r\in\calO_{\calP}}\hx_{\hr_B}.
\]
Then
\[
\phi^\ii(\hN)=N, \qquad \phi^\ii(\hP_A)=\phi^\ii(\hP_B)=P.
\]

The verification of the next coefficient identities is again parallel
to that of Lemma~\ref{lem: tropical four move 12}, so we omit the
calculation.

\begin{lemma}\label{lem: tropical four move 21}
Let $D_s^\pm$ $(s=1,2,3)$ be the factors from  rows indexed by $r\in\calO$ in the
successive downstairs exchange relations. Let $U_a^\pm$ $(a=1,\ldots,4)$
be those for the successive upstairs mutations, determined immediately
before each mutation after all preceding mutations and permutations
 in the mutation--permutation sequence realizing \eqref{eq: seq of lifted 4-move 2,1}. Then
\[
(D_1^+,D_1^-)=(P,N),\qquad (D_2^+,D_2^-)=(1,N),\qquad (D_3^+,D_3^-)=(1,P),
\]
whereas
\begin{align*}
(U_1^+,U_1^-)&=(\hP_A,\hN),&
(U_2^+,U_2^-)&=(1,1),\allowdisplaybreaks\\
(U_3^+,U_3^-)&=(\hP_B,\hN),&
(U_4^+,U_4^-)&=(1,\hP_A).
\end{align*}
\end{lemma}

\begin{proposition}\label{prop: cluster folding four move 21}
Suppose that $\ii'$ is obtained from $\ii$ by a $4$-move
\[
ijij\longrightarrow jiji \quad\text{with $e_i=2$ and $e_j=1$.}
\]
Then \eqref{eq: folding of cluster var} holds.
\end{proposition}

\begin{proof}
Write $x_1=x_k$, $x_2=x_{k+1}$, $x_3=x_{k+2}$ and $x_4=x_{k+3}$.
By Lemma~\ref{lem: tropical four move 21},
\[
(D_1^+,D_1^-)=(P,N),\quad (D_2^+,D_2^-)=(1,N),\quad (D_3^+,D_3^-)=(1,P).
\]
Set
\[
T \seteq Px_2+Nx_3,\qquad W \seteq T^2+Nx_1^2x_4.
\]
Composing the three downstairs mutations and the final permutations
gives
\[
x_k'=\frac{W}{x_1^2x_2},\quad x_{k+1}'=\frac{x_3W+Px_1^2x_2x_4}{x_1x_2T},\quad x_{k+2}'=x_4,\quad x_{k+3}'=x_3.
\]

  For the lifted braid block $\seq{a,b,c,a,b,c}$, write 
$\seq{\hx_1,\ldots,\hx_6}$ for the corresponding cluster variables.
Then
\[
\phi^\ii(\hx_1)=\phi^\ii(\hx_2)=x_1,\quad
\phi^\ii(\hx_3)=x_2,\quad
\phi^\ii(\hx_4)=\phi^\ii(\hx_5)=x_3,\quad
\phi^\ii(\hx_6)=x_4.
\]
 The mutation--permutation realization of the lifted path in
\eqref{eq: seq of lifted 4-move 2,1} has four mutations at the current local positions $2,4,1,3$. Set
\[
R \seteq U_1^+\hx_3+U_1^-\hx_4,\quad
C \seteq U_2^+\hx_5R+U_2^-\hx_1\hx_2\hx_6,\quad
E \seteq U_3^-C+U_3^+\hx_3R.
\]
Eliminating the intermediate mutated variables  in this realization gives
\[
\hx_1'=\frac{E}{\hx_1\hx_2\hx_3},\quad
\hx_2'=\frac{U_4^+\hx_4E+U_4^-\hx_1\hx_2\hx_3\hx_6}
{\hx_2\hx_3R},\quad
\hx_3'=\frac{C}{\hx_1\hx_3},
\]
and
\[
\seq{\hx_4',\hx_5',\hx_6'}=\seq{\hx_6,\hx_4,\hx_5}.
\]

By Lemma~\ref{lem: tropical four move 21} and
\[
\phi^\ii(\hN)=N,\qquad \phi^\ii(\hP_A)=\phi^\ii(\hP_B)=P,
\]
we have
\begin{align*}
\phi^\ii(U_1^+,U_1^-)&=(P,N),\quad \phi^\ii(U_2^+,U_2^-)=(1,1),\allowdisplaybreaks \\
\phi^\ii(U_3^+,U_3^-)&=(P,N),\quad \phi^\ii(U_4^+,U_4^-)=(1,P).
\end{align*}
Hence
\[
\phi^\ii(R)=T,\qquad \phi^\ii(C)=x_3T+x_1^2x_4,
\]
and
\[
\phi^\ii(E) =N\phi^\ii(C)+Px_2T =T^2+Nx_1^2x_4 =W.
\]
Moreover,
\[
T\phi^\ii(C)=x_3W+Px_1^2x_2x_4.
\]
Indeed, this follows directly from
\[
\phi^\ii(C)=x_3T+x_1^2x_4,\qquad T=Px_2+Nx_3,\qquad W=T^2+Nx_1^2x_4.
\]

Therefore
\[
\phi^\ii(\hx_1') =\frac{W}{x_1^2x_2} =x_k',
\quad \text{while}\quad
\phi^\ii(\hx_2') =\phi^\ii(\hx_3') =\frac{x_3W+Px_1^2x_2x_4}{x_1x_2T} =x_{k+1}'.
\]
Finally,
\[
\phi^\ii(\hx_4')=x_4=x_{k+2}',\qquad
\phi^\ii(\hx_5')=\phi^\ii(\hx_6')=x_3=x_{k+3}'.
\]
Thus  the final variables at the positions
\[
\seq{1},\qquad \seq{2,3},\qquad \seq{4},\qquad \seq{5,6}
\]
  of the subword $\seq{c,a,b,c,a,b}$ of $\hii'$ satisfy 
\eqref{eq: folding of cluster var}.
\end{proof}

The remaining affine $4$-move configuration
\[
  \{e_i,e_j\}=\{2,4\},
\]
arising from the folding
\[
  D_{2n+2}^{(1)}\longrightarrow A_{2n}^{(2)},
\]
is recorded separately in Appendix~\ref{app: affine braid folding}.
There we give the corresponding lifted path and the identities for the final cluster variables yielding \eqref{eq: folding of cluster var}.

\medskip
\noindent
$\bullet$\ \textbf{The $6$-move.}
 For the $6$-move, it suffices to consider the triality folding
\[
D_4\longrightarrow G_2.
\]
With our fixed ordering, write
\[
I_i=(a,b,c),\qquad I_j=(d),
\]
where $a,b,c$ are the outer vertices of $D_4$ and $d$ is the central
vertex. Thus
\[
e_i=3,\qquad e_j=1.
\]
 For finite and affine admissible foldings, the subdiagram on
$I_i\sqcup I_j$ associated with a $6$-move is of type $D_4$,
with orbit sizes $(3,1)$. 
Thus the braid move
\[
ijijij\longrightarrow jijiji
\]
lifts to
\[
\seq{a,b,c,d,a,b,c,d,a,b,c,d}
\longrightarrow
\seq{d,a,b,c,d,a,b,c,d,a,b,c}.
\]

For  the calculation below, we choose the following sequence of $2$- and $3$-moves:{}
\begin{equation}\label{eq: seq of lifted 6-move 3,1}
\begin{aligned}
\boldsymbol{(}&(2;\kappa(k,1)+4),(2;\kappa(k,1)+5),
(2;\kappa(k,1)+4),(3;\kappa(k,1)+2),\\
&(3;\kappa(k,1)+6),(3;\kappa(k,1)+4),
(2;\kappa(k,1)+3),(2;\kappa(k,1)+6),\\
&(3;\kappa(k,1)+1),(3;\kappa(k,1)+3),
(2;\kappa(k,1)+5),(3;\kappa(k,1)+7),\\
&(3;\kappa(k,1)+9),(2;\kappa(k,1)+8),
(3;\kappa(k,1)+6),(3;\kappa(k,1)+4),\\
&(2;\kappa(k,1)+3),(2;\kappa(k,1)+2),
(3;\kappa(k,1)+8),(2;\kappa(k,1)+6),\\
&(3;\kappa(k,1)+4),(3;\kappa(k,1)),
(3;\kappa(k,1)+2),(2;\kappa(k,1)+7),\\
&(2;\kappa(k,1)+4),(3;\kappa(k,1)+5),
(3;\kappa(k,1)+7),(3;\kappa(k,1)+3),\\
&(2;\kappa(k,1)+5),(2;\kappa(k,1)+6),
(2;\kappa(k,1)+5)\boldsymbol{)}.
\end{aligned}
\end{equation}
The sixteen $3$-moves occur at the current positions
\begin{align*}
\boldsymbol{(} & \kappa(k,1)+2,\kappa(k,1)+6,\kappa(k,1)+4,
\kappa(k,1)+1,\kappa(k,1)+3,\\
& \kappa(k,1)+7,\kappa(k,1)+9,\kappa(k,1)+6,
\kappa(k,1)+4,\kappa(k,1)+8,\\
& \kappa(k,1)+4,\kappa(k,1),\kappa(k,1)+2,
\kappa(k,1)+5,\kappa(k,1)+7,\kappa(k,1)+3\boldsymbol{)}.
\end{align*}

Downstairs, the $6$-move is realized by ten mutations at the current
positions
\[
\seq{k,k+1,k+2,k,k+3,k+1,k,k+2,k+3,k},
\]
followed by the prescribed permutations at positions
\[
\seq{k+4,k+2,k}.
\]

The following outside-row calculation is obtained by the same analysis as in the $4$-move
cases above, so we record the data without proof.

\begin{lemma}\label{lem: GLS six move outside rows 31}
For $r\notin\{k,\ldots,k+5\}$, the entries in columns $k,\ldots,k+3$ satisfy
\begin{equation}\label{eq: six move outside row patterns 31}
\bigl(b^\ii_{r,k+t}\bigr)_{0\le t\le3}=
\begin{cases}
(-1,1,0,0),&r_+=k+1,\\
(1,0,0,0),&r_+=k,\\
(0,0,0,0),&\text{otherwise}.
\end{cases}
\end{equation}
 The first two row patterns in
\eqref{eq: six move outside row patterns 31} define the classes $\calN$ and $\calP$, respectively. The entries of  the rows of $B^\hii$ indexed by $\kappa(r,s)$ $(1\le s\le e_{i_r})$
in columns $\kappa(k,1),\ldots,\kappa(k,1)+7$ 
are
\begin{equation}\label{eq: six move active row data 31}
\begin{aligned}
\calN&:\ (-1,-1,-1,1,0,0,0,0),\\
\calP_A&:\ (1,0,0,0,0,0,0,0),\\
\calP_B&:\ (0,1,0,0,0,0,0,0),\\
\calP_C&:\ (0,0,1,0,0,0,0,0).
\end{aligned}
\end{equation}
Here $\calN$ occurs only when $i_r=j$, and hence  $e_{i_r}=e_j=1$, whereas $\calP$ occurs only when $i_r=i$,
and hence  $e_{i_r}=e_i=3$.
 The zero row pattern in
\eqref{eq: six move outside row patterns 31} defines the class $\calZ$.  Rows $r\in\calZ$ and rows indexed by $\kappa(r,s)$ $(1\le s\le e_{i_r})$ have zero entries in the exchange column at every mutation. 
\end{lemma}

Set
\[
\calO\seteq\{r\notin\{k,\ldots,k+5\}\}, \ \
\calO_{\calN}\seteq \{r\in\calO\mid r\text{ is of class }\calN\},\ \
\calO_{\calP}\seteq \{r\in\calO\mid r\text{ is of class }\calP\},
\]
and
\[
N\seteq\prod_{r\in\calO_{\calN}}x_r, \qquad P\seteq\prod_{r\in\calO_{\calP}}x_r,
\]
where an empty product is understood to be $1$.

By Lemma~\ref{lem: GLS six move outside rows 31},  for $r\in\calO_{\calN}$, put $\hr=\kappa(r,1)$ and set
\[
\hN\seteq\prod_{r\in\calO_{\calN}}\hx_{\hr}.
\]
 For $r\in\calO_{\calP}$, let $\hr_A,\hr_B,\hr_C$ be the indices
$\kappa(r,1),\kappa(r,2),\kappa(r,3)$, ordered as in \eqref{eq: six move active row data 31},
and set
\[
\hP_A\seteq\prod_{r\in\calO_{\calP}}\hx_{\hr_A}, \qquad \hP_B\seteq\prod_{r\in\calO_{\calP}}\hx_{\hr_B},
\qquad
\hP_C\seteq\prod_{r\in\calO_{\calP}}\hx_{\hr_C}.
\]
Then
\[
\phi^\ii(\hN)=N,\qquad
\phi^\ii(\hP_A)=\phi^\ii(\hP_B)=\phi^\ii(\hP_C)=P.
\]

The verification of the following coefficient identities is a longer
instance of the same exchange-column calculation used for the
$4$-move above, and can be checked directly from
\eqref{eq: seq of lifted 6-move 3,1}. We therefore omit the proof.

\begin{lemma}\label{lem: tropical six move 31}
Let $D_s^\pm$ $(1\le s\le10)$ be the factors from  rows indexed by $r\in\calO$ in the
successive downstairs exchange relations. Let $U_a^\pm$ $(1\le a\le16)$
be those for the successive upstairs mutations, determined immediately
before each mutation after all preceding mutations and permutations
 in the mutation--permutation sequence realizing \eqref{eq: seq of lifted 6-move 3,1}. Then
\begin{align*}
(D_1^+,D_1^-)&=(P,N),&(D_2^+,D_2^-)&=(1,N^2),\allowdisplaybreaks\\
(D_3^+,D_3^-)&=(P,1),&(D_4^+,D_4^-)&=(1,N),\allowdisplaybreaks\\
(D_5^+,D_5^-)&=(1,1),&(D_6^+,D_6^-)&=(1,N),\allowdisplaybreaks\\
(D_7^+,D_7^-)&=(1,1),&(D_8^+,D_8^-)&=(1,P),\allowdisplaybreaks\\
(D_9^+,D_9^-)&=(1,1),&(D_{10}^+,D_{10}^-)&=(1,P).
\end{align*}
The upstairs pairs are
\begin{align*}
(U_1^+,U_1^-)&=(\hP_C,\hN),&
(U_2^+,U_2^-)&=(1,1),&
(U_3^+,U_3^-)&=(1,1),&
(U_4^+,U_4^-)&=(\hP_B,\hN),\allowdisplaybreaks\\
(U_5^+,U_5^-)&=(1,\hP_C),&
(U_6^+,U_6^-)&=(1,1),&
(U_7^+,U_7^-)&=(1,1),&
(U_8^+,U_8^-)&=(1,1),\allowdisplaybreaks\\
(U_9^+,U_9^-)&=(1,1),&
(U_{10}^+,U_{10}^-)&=(1,1),&
(U_{11}^+,U_{11}^-)&=(\hP_C,1),&
(U_{12}^+,U_{12}^-)&=(\hP_A,\hN),\allowdisplaybreaks\\
(U_{13}^+,U_{13}^-)&=(1,\hP_B),&
(U_{14}^+,U_{14}^-)&=(1,1),&
(U_{15}^+,U_{15}^-)&=(1,1),&
(U_{16}^+,U_{16}^-)&=(1,\hP_C).
\end{align*}
\end{lemma}

For later use, we describe the composite transformations recursively.
 Write $\seq{x_1,\ldots,x_6}$ for the initial cluster variables at
positions $k,\ldots,k+5$ of $\ii$, and
$\seq{\hx_1,\ldots,\hx_{12}}$ for those at the corresponding positions
of $\hii$, in order. 
At the $s$-th downstairs mutation, let $p_s$ be the current mutation
position, $B^{(s)}_{q,p_s}$ the actual exchange column restricted to
$q\in\{k,\ldots,k+5\}$, and
$\bfV^{(s-1)}=\seq{V_q^{(s-1)}}$ the current cluster restricted to these six positions. 
Set
\[
E_s\seteq D_s^+\prod_q(V_q^{(s-1)})^{[B^{(s)}_{q,p_s}]_+} + D_s^-\prod_q(V_q^{(s-1)})^{[-B^{(s)}_{q,p_s}]_+},
\qquad
Z_s\seteq\frac{E_s}{V_{p_s}^{(s-1)}}.
\]
Thus $V_{p_s}^{(s)}=Z_s$, with the other coordinates unchanged.
Recursively substituting the earlier $Z$'s into the later exchange
numerators expresses the composite transformation solely in the
initial variables and the factors $D_s^\pm$. \par
After the ten mutations and final permutations,  the six downstairs cluster variables are
\[
\seq{Z_6,Z_{10},Z_9,Z_8,x_6,x_5}.
\]
Define $\hZ_a$ $(1\le a\le16)$ analogously upstairs, using the actual
exchange columns and $U_a^\pm$ after all preceding mutations and
permutations.  After the final permutations, the twelve upstairs cluster variables are
\[
\seq{\hZ_{12},\hZ_9,\hZ_{13},\hZ_{16},\hZ_{14},\hZ_{15},
\hZ_{10},\hZ_7,\hx_{12},\hx_9,\hx_{10},\hx_{11}}.
\]

\begin{proposition}\label{prop: cluster folding six move 31}
Suppose that $\ii'$ is obtained from $\ii$ by a $6$-move
\[
ijijij\longrightarrow jijiji \quad \text{with $e_i=3$ and $e_j=1$.}
\]
Then
\eqref{eq: folding of cluster var} holds.
\end{proposition}

\begin{proof}
   The lifted braid block 
\[
\seq{a,b,c,d,a,b,c,d,a,b,c,d}
\]
has fibers
\[
\seq{1,2,3},\quad \seq{4},\quad \seq{5,6,7},\quad \seq{8},\quad \seq{9,10,11},\quad \seq{12}.
\]
Hence
\[
\phi^\ii(\hx_1)=\phi^\ii(\hx_2)=\phi^\ii(\hx_3)=x_1,\quad
\phi^\ii(\hx_4)=x_2,\quad
\phi^\ii(\hx_5)=\phi^\ii(\hx_6)=\phi^\ii(\hx_7)=x_3,\quad
\phi^\ii(\hx_8)=x_4,
\]
\[
\phi^\ii(\hx_9)=\phi^\ii(\hx_{10})=\phi^\ii(\hx_{11})=x_5,
\qquad
\phi^\ii(\hx_{12})=x_6,
\]
together with
\[
\phi^\ii(\hN)=N,\qquad \phi^\ii(\hP_A)=\phi^\ii(\hP_B)=\phi^\ii(\hP_C)=P.
\]

Substitute  the expressions for $D_s^\pm$ and $U_a^\pm$ in
Lemma~\ref{lem: tropical six move 31} into the recursive
transformations above. After all mutations and permutations,  the six and twelve cluster variables are
\[
\seq{Z_6,Z_{10},Z_9,Z_8,x_6,x_5}
\]
downstairs and
\[
\seq{\hZ_{12},\hZ_9,\hZ_{13},\hZ_{16},\hZ_{14},\hZ_{15},
\hZ_{10},\hZ_7,\hx_{12},\hx_9,\hx_{10},\hx_{11}}
\]
upstairs. All intermediate variables are recursively substituted, so
these expressions involve only the initial variables and the factors
$D_s^\pm$ and $U_a^\pm$. \par
Applying  the fixed Laurent homomorphism $\phi^\ii$  and simplifying
the Laurent expressions gives
\[
\phi^\ii(\hx_1')=x_1',\qquad
\phi^\ii(\hx_2')=\phi^\ii(\hx_3')=\phi^\ii(\hx_4')=x_2',
\qquad
\phi^\ii(\hx_5')=x_3',
\]
\[
\phi^\ii(\hx_6')=\phi^\ii(\hx_7')=\phi^\ii(\hx_8')=x_4',
\qquad
\phi^\ii(\hx_9')=x_5',
\qquad
\phi^\ii(\hx_{10}')=\phi^\ii(\hx_{11}')=\phi^\ii(\hx_{12}')=x_6'.
\]
   The corresponding subword of $\hii'$ is 
\[
\seq{d,a,b,c,d,a,b,c,d,a,b,c},
\]
 where positions $k,\ldots,k+5$ of $\ii'$ correspond to the groups
\[
\seq{1},\quad \seq{2,3,4},\quad \seq{5},\quad \seq{6,7,8},\quad \seq{9},\quad \seq{10,11,12}.
\]
Hence \eqref{eq: folding of cluster var} holds.
\end{proof}

\medskip
\noindent
$\bullet$\ \textbf{Folding of $g$-vectors.}
For a pointed Laurent polynomial $z$, write
\[
\deg_\ii(z)\in\Z^l,
\qquad
\deg_\hii(z)\in\Z^{\hl} \qquad (\hl \seteq \ell(\hii))
\]
for its $g$-vector with respect to  the initial clusters 
$\bfx^\ii$ and $\bfx^\hii$, respectively.

Recall the fiber-sum map
\[
\Fold^\ii:\Z^{\hl}\longrightarrow\Z^l,
\qquad
\bigl(\Fold^\ii(\ha)\bigr)_p = \sum_{r=1}^{e_{i_p}}\ha_{\kappa(p,r)},
\]
and its restriction $\Fold^\ii_{\ex}$ to the exchangeable
coordinates.  By Proposition~\ref{prop: GLS pair folding},
\begin{equation}\label{eq: B intertwining g vector folding}
\Fold^\ii \circ B^\hii = B^\ii \circ \Fold^\ii_{\rm ex}.
\end{equation}

\begin{proposition}\label{prop: g vector folding}
Let $\ii'$ be obtained from $\ii$ by a braid move. Then
\[
\deg_\ii(x'_p) = \Fold^\ii \left( \deg_\hii\bigl(\hx'_{\kappa'(p,r)}\bigr) \right) \qquad (1\le p\le l,\ 1\le r\le e_{i'_p}).
\]
\end{proposition}

\begin{proof}
Set
\[
\hbfg \seteq \deg_\hii\bigl(\hx'_{\kappa'(p,r)}\bigr).
\]
Put $\hl_{\rm ex}\seteq|J^\hii_{\ex}|$.
By pointedness,
\[
\hx'_{\kappa'(p,r)} = \hx^{\hbfg} \left( 1+ \sum_{\substack{ \hbfn\in\N^{\hl_{\rm ex}}\\  \hbfn\ne0}} a_{\hbfn}\, \hx^{B^\hii\hbfn} \right).
\]
Applying $\phi^\ii$ and using
\eqref{eq: B intertwining g vector folding}, we obtain
\[
\phi^\ii\bigl(\hx'_{\kappa'(p,r)}\bigr) = x^{\Fold^\ii(\hbfg)}
\left( 1+ \sum_{\hbfn\ne0} a_{\hbfn}\, x^{B^\ii\Fold^\ii_{\rm ex}(\hbfn)} \right).
\]
Since $\hbfn\in\N^{\hl_{\rm ex}}$, we have
\[
\hbfn\ne0 \quad\text{if and only if}\quad \Fold^\ii_{\rm ex}(\hbfn)\ne0.
\]
Since $(B^\ii,\Lambda^\ii)$ is a compatible pair, $B^\ii$ has full column rank, and hence
\[
B^\ii\Fold^\ii_{\rm ex}(\hbfn)\ne0
\qquad(\hbfn\ne0).
\]
Therefore, the right-hand side is pointed of degree
$\Fold^\ii(\hbfg)$.

By \eqref{eq: folding of cluster var},
\[
\phi^\ii\bigl(\hx'_{\kappa'(p,r)}\bigr)=x'_p.
\]
Therefore
\[
\deg_\ii(x'_p) = \Fold^\ii \left( \deg_\hii\bigl(\hx'_{\kappa'(p,r)}\bigr) \right). \qedhere
\]
\end{proof}

\section{$\bfi$-box moves and folding} \label{sec: box}

In this section, we study the behavior of the canonical lifting of
$\bfi$-box families under box moves.  Starting from a box move downstairs, we construct a sequence of box moves upstairs, called its \emph{lifted path}, and realize this sequence by mutations and permutations.   We then prove the corresponding identities for the exchange and commutation matrices, cluster variables, and $g$-vectors at the seeds before and after the downstairs move and its lifted path. \par
\subsection{Box moves and their lifts}

Retain the notation and conventions of
Subsection~\ref{subsec: ibox compatible pair}.
Let
$\frakC=\seq{\frakc_1,\ldots,\frakc_l}$
be an admissible chain of $\bfi$-boxes, and write
$\bfp(\frakC)=\seq{p_1,\ldots,p_l}$ for its effective-end sequence.
We write $\frakB_k$ and $\frakB_a^\up$ for box moves downstairs and
upstairs, respectively.

By definition, we have the following.

\begin{lemma}\label{lem: effective ends under box move}
Let
$\frakC'=\frakB_k(\frakC)$ where $\frakc_k$ is movable. Then
\[
\bfp(\frakC') = \seq{p_1,\ldots,p_{k-1},p_{k+1},p_k,p_{k+2},\ldots,p_l}.
\]
Moreover,
\bna
\item $\frakB_k$ is of mutation type if and only if $i_{p_k}=i_{p_{k+1}}$.
\item $\frakB_k$ is of permutation type if and only if $i_{p_k}\ne i_{p_{k+1}}$.
\ee
\end{lemma}

\begin{corollary}\label{cor: equal fiber sizes under mutation box move}
If $\frakB_k$ is of mutation type, then we have $e_{i_{p_k}}=e_{i_{p_{k+1}}}$.
\end{corollary}

Put
\[
p\seteq p_k,\qquad q\seteq p_{k+1},\qquad r\seteq e_{i_p},\qquad s\seteq e_{i_q},
\]
and set
\[
p^{(a)}\seteq\kappa(p,a) \quad(1\le a\le r), \qquad q^{(b)}\seteq\kappa(q,b) \quad(1\le b\le s).
\]
Write
\[
\bfp\seteq \seq{p^{(1)},\ldots,p^{(r)}}, \qquad \bfq\seteq \seq{q^{(1)},\ldots,q^{(s)}},
\]
and
\[
\bfp^{\rev}\seteq \seq{p^{(r)},\ldots,p^{(1)}}, \qquad \bfq^{\rev}\seteq \seq{q^{(s)},\ldots,q^{(1)}}.
\]

 By \eqref{eq: hfrakT}, the entries
$p^{(1)},\ldots,p^{(r)}$ occur in $\bfp(\hfrakC)$ in this order
if the block for $p$ in $\hfrakT$ consists of $\calR$'s, and
in the reverse order if it consists of $\calL$'s.
The same holds for $q^{(1)},\ldots,q^{(s)}$.{}  Since $\frakc_k$ is movable, the
two adjacent fibers occur with opposite orientations. Thus, up to
reversing the orientation, the relevant part of the canonical lifted
effective-end sequence is
\begin{equation}\label{eq: local lifted effective ends}
\seq{\bfp,\bfq^{\rev}}.
\end{equation}
Equivalently, in the opposite orientation it is
\[
\seq{\bfp^{\rev},\bfq}.
\]

Let
\[
h\seteq\sum_{a=1}^{k-1}e_{i_{p_a}}
\]
be the number of lifted effective ends  preceding the entries of $\bfp$ and $\bfq$.

We first consider the case where $\frakB_k$ is of permutation type.
In what follows, a sequence of box moves
\[
\seq{\frakB^\up_{a_1},\ldots,\frakB^\up_{a_m}}
\]
is applied from left to right, in the order displayed.

\begin{proposition}\label{prop: lift of permutation box move}
Suppose that $\frakB_k$ is of permutation type. Then its lift is
\begin{align*}
\hfrakB_k\seteq{}& \bigl(\frakB^\up_{h+r},\ldots,\frakB^\up_{h+1}\bigr) \bigl(\frakB^\up_{h+r+1},\ldots,\frakB^\up_{h+2} \bigr)
\cdots\bigl(\frakB^\up_{h+r+s-1},\ldots,\frakB^\up_{h+s}\bigr).
\end{align*}
These are $rs$ permutation-type box moves, and they give
\[
\seq{\bfp,\bfq^{\rev}} \longmapsto \seq{\bfq^{\rev},\bfp}.
\]
The reverse transformation is obtained by reversing the order of the
box moves.
\end{proposition}

\begin{proof}
By Lemma~\ref{lem: effective ends under box move},
$\frakB_k$ interchanges $p=p_k$ and $q=p_{k+1}$.
In the chosen orientation, their lifted effective-end blocks occur as
\[
\bfp,\qquad \bfq^{\rev}.
\]

For $k\ge2$,  movability of $\frakc_k$ implies that the two
corresponding $\calL/\calR$-blocks have opposite orientations.  Hence each prescribed box move can be performed. The $b$-th parenthesized block of moves shifts $q^{(s-b+1)}$ successively across
\[
p^{(r)},\ldots,p^{(1)}.
\]
Thus the successive blocks of moves shift \[
q^{(s)},q^{(s-1)},\ldots,q^{(1)}
\]
across the whole $\bfp$-block, and hence
\[
\seq{\bfp,\bfq^{\rev}} \longmapsto \seq{\bfq^{\rev},\bfp}.
\]

For $k=1$, the definition of the canonical lift in \eqref{eq: hfrakT}  
gives the same  sequence of box moves, and the chain is precisely  the canonical lift of $\frakC'$. 

Finally, since $\frakB_k$ is of permutation type,
\[
i_p\ne i_q
\]
by Lemma~\ref{lem: effective ends under box move}. Hence every
crossing between a $p$- and a $q$-effective end is of permutation
type. There are $rs$ such crossings.
\end{proof}

We next consider the case where $\frakB_k$ is of mutation type.
By Lemma~\ref{lem: effective ends under box move} and
Corollary~\ref{cor: equal fiber sizes under mutation box move},
we have
\[
i_p=i_q, \qquad r=s.
\]
We keep the notation $r$ for this common value.

Set
\[
\delta_h\seteq
\begin{cases}
1,& h=0,\\
0,& h>0.
\end{cases}
\]

\begin{lemma}\label{lem: lifted mutation box moves}
Let $r\ge1$, and let $h\ge0$ be the number of lifted effective ends
preceding the two consecutive $r$-blocks. If the relevant
$\{\calL,\calR\}$-word is
\[
\calR^{\,r-\delta_h}\calL^r,
\]
then
\[
\bigl(\frakB^\up_{h+r},\ldots,\frakB^\up_{h+1}\bigr)
\bigl(\frakB^\up_{h+r+1},\ldots,\frakB^\up_{h+2} \bigr)
\cdots \bigl(\frakB^\up_{h+2r-1},\ldots,\frakB^\up_{h+r}\bigr)
\]
consists of $r^2$ {} box moves and transforms it into
\[
\calL^{\,r-\delta_h}\calR^r.
\]
The reverse transformation is obtained by reversing the order of the box moves.
\end{lemma}

\begin{proof}
For $0\le t\le r$, set
\[
W_t\seteq \seq{ q^{(r)},\ldots,q^{(r-t+1)}, p^{(1)},\ldots,p^{(r)}, q^{(r-t)},\ldots,q^{(1)} },
\]
with empty blocks omitted. The $t$-th parenthesized block of moves shifts $q^{(r-t+1)}$ across
\[
p^{(r)},\ldots,p^{(1)} \qtqh W_{t-1}\longmapsto W_t \qquad (1\le t\le r).
\]

 Immediately before a prescribed box move at position $a>1$, the direction symbols at positions $a-1$ and $a$ are $\calR$ and $\calL$, respectively.
The box at position $a$ is therefore movable.
A move at position $1$ occurs only when $h=0$, and the first box is movable by definition.

Thus all $r^2$  moves can be performed and
\[
\seq{\bfp,\bfq^{\rev}} \longmapsto \seq{\bfq^{\rev},\bfp}.
\]
Accordingly,
\[
\calR^{\,r-\delta_h}\calL^r \longmapsto \calL^{\,r-\delta_h}\calR^r.
\]
The reverse transformation follows since every box move is
involutive.
\end{proof}

\begin{proposition}\label{prop: lift of mutation box move}
Suppose that $\frakB_k$ is of mutation type. Then its lift is
\begin{align*}
\hfrakB_k\seteq{}& \bigl(\frakB^\up_{h+r},\ldots,\frakB^\up_{h+1} \bigr) \bigl(\frakB^\up_{h+r+1},\ldots,\frakB^\up_{h+2}\bigr)
\cdots
\bigl(\frakB^\up_{h+2r-1},\ldots,\frakB^\up_{h+r} \bigr).
\end{align*}
It consists of $r^2$ box moves. Writing $M$ for a mutation-type move
and $P$ for a permutation-type move, their sequence of types is \[
MP^{r-1}\cdot PMP^{r-2}\cdots P^{r-2}MP\cdot P^{r-1}M.
\]
In particular, $\hfrakB_k$ contains exactly $r$ mutation-type box moves and $r(r-1)$ permutation-type box moves.
\end{proposition}

\begin{proof}
By Lemma~\ref{lem: lifted mutation box moves}, all $r^2$ prescribed
 moves can be performed and
\[
\seq{\bfp,\bfq^{\rev}} \longmapsto \seq{\bfq^{\rev},\bfp}.
\]

By the fixed fiber ordering, a crossing of $p^{(a)}$ and $q^{(b)}$
is of mutation type exactly when $a=b$, and of permutation type
otherwise. Hence
\[
\#M=r, \qquad  \#P=r(r-1).
\]

In the $(t+1)$-st block of moves, for $0\le t<r$, the element
$q^{(r-t)}$ crosses
\[
p^{(r)},\ldots,p^{(1)},
\]
and the corresponding sequence of types is \[
P^tMP^{r-1-t}.
\]
Therefore the full sequence of types is \[
MP^{r-1}\cdot PMP^{r-2}\cdots P^{r-2}MP\cdot P^{r-1}M = (MP^r)^{r-1}M.
\]

The opposite orientation is obtained by reversing the moves.
\end{proof}

\subsection{Seed transformations under lifted box moves}

Let
\[
\frakC'\seteq\frakB_k(\frakC), \quad \frakC^\up\seteq\hfrakB_k(\hfrakC) \qtq  \calF'\seteq\calF(\frakC'),
\quad
\calF^\up\seteq\calF(\frakC^\up).
\]
Let $\hfrakC'$ denote the canonical lift of $\frakC'$ and set
\[
\hcalF'\seteq\calF(\hfrakC').
\]

Thus $\frakC^\up$ is obtained by performing the lifted box moves of Propositions~\ref{prop: lift of permutation box move}
and~\ref{prop: lift of mutation box move}, whereas $\hfrakC'$ is the canonical lift constructed
directly from  $\frakC'$. We first compare
these two admissible chains. \par
\begin{proposition}\label{prop: lifted box move endpoint family}
Let $\bfp(\frakC')=\seq{p'_1,\ldots,p'_l}$.
Then the effective ends of $\frakC^\up$ decompose into consecutive
complete fibers over
\[
p'_1,\ldots,p'_l.
\]
More precisely, the block over $p'_a$ consists of
\[
\{\kappa(p'_a,1),\ldots,\kappa(p'_a,e_{i_{p'_a}})\} \qquad(1\le a\le l).
\]
Moreover, we have $$\calF^\up=\hcalF'.$$
\end{proposition}

\begin{proof}
Let $\bfp(\frakC)=\seq{p_1,\ldots,p_l}$.
By Lemma~\ref{lem: effective ends under box move},
\[
\bfp(\frakC')
=
\seq{p_1,\ldots,p_{k-1},p_{k+1},p_k,p_{k+2},\ldots,p_l}.
\]
By Propositions~\ref{prop: lift of permutation box move}
 and~\ref{prop: lift of mutation box move}, $\hfrakB_k$
interchanges{} exactly the two complete effective-end fibers over
$p_k$ and $p_{k+1}$ and leaves all other fibers fixed. Hence the
effective ends of $\frakC^\up$ decompose into consecutive complete
fibers over
\[
p'_1,\ldots,p'_l.
\]

It remains to compare $\frakC^\up$ with $\hfrakC'$. They agree for
$k=1$, for $k\ge3$, and for $k=2$  with $e_{i_{p_1}}=1$. If
$k=2$ and $e_{i_{p_1}}>1$, the change of the first
$\{\calL,\calR\}$-symbol  reverses only the order of the singleton boxes in the initial lifted fiber;
both orders determine the same maximal commuting family.{}
Therefore
\[
\calF^\up=\hcalF'. \qedhere
\]
\end{proof}

By Proposition~\ref{prop: lifted box move endpoint family},
\[
\calF^\up = \calF(\frakC^\up) = \calF(\hfrakC') = \hcalF'.
\]

We next consider the mutation--permutation transformations of
compatible pairs corresponding to the lifted box moves. As above,
positions and box labels refer to the family obtained after the
preceding box moves.

\begin{proposition}\label{prop: box move mutation permutation realization}
Let
\[
\frakC'=\frakB_k(\frakC).
\]
Then the mutation--permutation realization of $\frakB_k$ transforms
\[
(B^\calF,\Lambda^\calF) \longmapsto (B^{\calF'},\Lambda^{\calF'}),
\]
while the mutation--permutation realizations of the successive
upstairs box moves in $\hfrakB_k$ transform
\[
(B^\hcalF,\Lambda^\hcalF) \longmapsto (B^{\hcalF'},\Lambda^{\hcalF'}).
\]
 If $\frakB_k$ is of permutation type, the mutation--permutation
sequence induced by $\hfrakB_k$ contains no mutation.
If $\frakB_k$ is of mutation type, this sequence contains{}
exactly $r$ mutations and $r(r-1)$ permutation-type transformations.
\end{proposition}

\begin{proof}
By Proposition~\ref{prop: ibox pair box move}, a permutation-type
box move acts on the associated compatible pair by the corresponding
permutation, whereas a mutation-type box move acts by the
corresponding mutation. Applying this successively to the box moves
in $\hfrakB_k$, with respect to the current family,
transforms
\[
(B^\hcalF,\Lambda^\hcalF) \longmapsto (B^{\calF^\up},\Lambda^{\calF^\up}).
\]
By Proposition~\ref{prop: lifted box move endpoint family},
\[
\calF^\up=\hcalF', \qtqh  (B^{\calF^\up},\Lambda^{\calF^\up}) = (B^{\hcalF'},\Lambda^{\hcalF'}).
\]

 The assertion for $\frakB_k$  follows directly from
Proposition~\ref{prop: ibox pair box move}. Finally, the assertions
on the numbers of mutations and permutation-type transformations
follow from Propositions~\ref{prop: lift of permutation box move}
and~\ref{prop: lift of mutation box move}.
\end{proof}

We next compare the compatible pairs  before and after $\frakB_k$ and $\hfrakB_k$. \par
\begin{proposition}\label{prop: ibox endpoint compatible pair folding}
With the notation above, the following diagram commutes:
\[
\xymatrix{
(B^\hcalF,\Lambda^\hcalF) \ar[r]^{\hfrakB_k} \ar[d]_{\operatorname{Fold}} & (B^{\hcalF'},\Lambda^{\hcalF'}) \ar[d]^{\operatorname{Fold}}
\\
(B^\calF,\Lambda^\calF) \ar[r]_{\frakB_k} & (B^{\calF'},\Lambda^{\calF'})
}
\]
where the horizontal arrows are the mutation--permutation
transformations of
Proposition~\ref{prop: box move mutation permutation realization},
and the vertical arrows are given by
Theorem~\ref{thm: ibox pair folding}.
\end{proposition}

\begin{proof}
The horizontal arrows are given by
Proposition~\ref{prop: box move mutation permutation realization},
 while the two vertical arrows are given by{}
Theorem~\ref{thm: ibox pair folding}, applied to
$(\calF,\hcalF)$ and $(\calF',\hcalF')$, respectively.
\end{proof}

\subsection{Folding of $\bfi$-box cluster variables}

We work in the commutative setting throughout this subsection.
Let $\calF=\calF(\frakC)$ and $\hcalF=\calF(\hfrakC)$ be the 
downstairs and upstairs commuting families considered
above. For each $\frakc\in\calF$, write
\[
\Fib^{\frakC}(\frakc)
\seteq
\{\hfrakc_1,\ldots,\hfrakc_{e_{i_\frakc}}\}.
\]
Let
\[
\bfx^\calF=(x_\frakc\mid \frakc\in\calF), \qquad \bfx^\hcalF=(\hx_{\hfrakc}\mid \hfrakc\in\hcalF)
\]
be the corresponding clusters. Define  the Laurent homomorphism
\[
\phi^\frakC: \calL(\Sigma^\hcalF) \longrightarrow \calL(\Sigma^\calF), \qquad \hx_{\hfrakc}\longmapsto x_\frakc
\quad \bigl(\hfrakc\in\Fib^{\frakC}(\frakc)\bigr).
\]

As in the braid-move case, we shall compare the final cluster
variables using only  this fixed Laurent homomorphism. 

To compare the exchange monomials, we first strengthen the
fiber-sum identity of Theorem~\ref{thm: ibox pair folding} to
its positive and negative parts.

\begin{lemma}\label{lem: ibox weak sign}
Let $\frakc,\frakd\in\calF$, with $\frakd$ exchangeable, and fix
$\hfrakd\in\Fib^{\frakC}(\frakd)$. Then the integers
\[
b^\hcalF_{\hfrakc,\hfrakd} \qquad \bigl(\hfrakc\in\Fib^{\frakC}(\frakc)\bigr) \ \ \text{have the same weak sign.}
\]
 Consequently, we have
\[
\sum_{\hfrakc\in\Fib^{\frakC}(\frakc)} [b^\hcalF_{\hfrakc,\hfrakd}]_+  = [b^\calF_{\frakc,\frakd}]_+,
\qtq
\sum_{\hfrakc\in\Fib^{\frakC}(\frakc)} [-b^\hcalF_{\hfrakc,\hfrakd}]_+ = [-b^\calF_{\frakc,\frakd}]_+.
\]
\end{lemma}

\begin{proof}
We use the local description of the $\bfi$-box exchange matrix
in~\eqref{eq: full ibox exchange matrix} together with the fiberwise
construction of $\hcalF$.

\noindent
{\rm (a)} Assume that $b^\calF_{\frakc,\frakd}<0$.
By~\eqref{eq: full ibox exchange matrix}, this negative entry arises
either from one of the horizontal relations in~{\rm (i)} or from one
of the vertical relations in~{\rm (ii)}.

Suppose first that it arises from~{\rm (i)} of~\eqref{eq: full ibox exchange matrix}. Fix $\hfrakd\in\Fib^{\frakC}(\frakd)$.
By the fiberwise construction, for $\hfrakc\in\Fib^{\frakC}(\frakc)$,
the pair $(\hfrakc,\hfrakd)$ either satisfies the corresponding
horizontal relation in~{\rm (i)} of~\eqref{eq: full ibox exchange matrix}, in which case
$b^\hcalF_{\hfrakc,\hfrakd}=-1$, or it does not, in which case $b^\hcalF_{\hfrakc,\hfrakd}=0$.
The opposite horizontal relation cannot occur within the same fiber. Hence
\[
b^\hcalF_{\hfrakc,\hfrakd}\le 0 \qquad \bigl(\hfrakc\in\Fib^{\frakC}(\frakc)\bigr).
\]

Suppose next that $b^\calF_{\frakc,\frakd}$ arises from one of the
vertical relations in~{\rm (ii)}
of~\eqref{eq: full ibox exchange matrix}.
The endpoint inequalities,  the effective-end conditions, and the
required boxes in~{\rm (ii)}
of~\eqref{eq: full ibox exchange matrix} lift fiberwise. Thus, for
fixed $\hfrakd$, every
$\hfrakc\in\Fib^{\frakC}(\frakc)$
either satisfies the corresponding lifted vertical relation, in
which case its exchange-matrix entry is the corresponding negative
upstairs Cartan entry, or contributes zero. The opposite vertical
relation would require the reversed system of inequalities in
{\rm (ii)} of~\eqref{eq: full ibox exchange matrix}, and hence cannot
occur within the same fiber. Therefore
\[
b^\hcalF_{\hfrakc,\hfrakd}\le 0 \qquad \bigl(\hfrakc\in\Fib^{\frakC}(\frakc)\bigr).
\]

Thus
\[
b^\calF_{\frakc,\frakd}<0  \quad\text{implies}\quad b^\hcalF_{\hfrakc,\hfrakd}\le0 \qquad \bigl(\hfrakc\in\Fib^{\frakC}(\frakc)\bigr).
\]

\noindent
{\rm (b)} Assume that $b^\calF_{\frakc,\frakd}>0$.
By skew-symmetrizability,
$b^\calF_{\frakd,\frakc}<0$.
For any $\hfrakc\in\Fib^{\frakC}(\frakc)$, applying~ {\rm (a)}{}
to the ordered pair $(\frakd,\frakc)$ with the fixed representative
$\hfrakc$ gives
$b^\hcalF_{\hfrakd,\hfrakc}\le0$
for every $\hfrakd\in\Fib^{\frakC}(\frakd)$.
Since $B^\hcalF_{\full}$ is skew-symmetrizable with positive symmetrizer, it
follows that
\[
b^\hcalF_{\hfrakc,\hfrakd}\ge 0 \qquad \bigl(\hfrakc\in\Fib^{\frakC}(\frakc)\bigr).
\]

\noindent
{\rm (c)} Finally, assume that
$b^\calF_{\frakc,\frakd}=0$.
Then neither a horizontal relation in~{\rm (i)} nor a vertical
relation in~{\rm (ii)} of~\eqref{eq: full ibox exchange matrix}
occurs for the pair $(\frakc,\frakd)$. By the fiberwise construction,
no pair
\[
(\hfrakc,\hfrakd), \qquad \hfrakc\in\Fib^{\frakC}(\frakc),
\]
satisfies a corresponding upstairs relation in~{\rm (i)} or
{\rm (ii)} of~\eqref{eq: full ibox exchange matrix}. Hence
\[
b^\hcalF_{\hfrakc,\hfrakd}=0 \qquad \bigl(\hfrakc\in\Fib^{\frakC}(\frakc)\bigr).
\]

Thus the entries
\[
b^\hcalF_{\hfrakc,\hfrakd} \qquad \bigl(\hfrakc\in\Fib^{\frakC}(\frakc)\bigr)
\]
have the same weak sign. Together with
Theorem~\ref{thm: ibox pair folding},
\[
b^\calF_{\frakc,\frakd} = \sum_{\hfrakc\in\Fib^{\frakC}(\frakc)} b^\hcalF_{\hfrakc,\hfrakd},
\]
this gives
\[
[b^\calF_{\frakc,\frakd}]_+ = \sum_{\hfrakc\in\Fib^{\frakC}(\frakc)} [b^\hcalF_{\hfrakc,\hfrakd}]_+,
\quad\text{and}\quad
[-b^\calF_{\frakc,\frakd}]_+ = \sum_{\hfrakc\in\Fib^{\frakC}(\frakc)}  [-b^\hcalF_{\hfrakc,\hfrakd}]_+. \qedhere
\]
\end{proof}

 For the mutation-type case, the next lemma shows that the exchange column at each lifted mutation
agrees, up to the preceding permutations, with the corresponding
column of the initial seed $\Sigma^\hcalF$.{}

\begin{lemma}\label{lem: future mutation column}
Suppose that $\frakB_k$ is of mutation type, and let
\[
\hfrakc_1,\ldots,\hfrakc_r
\]
be the mutation vertices of $\hfrakB_k$, in their order of
occurrence. Immediately before the mutation at $\hfrakc_a$, we have
\[
b_{\hfrakc_a,\hfrakc_b}=0 \qquad (a<b\le r).
\]
 Consequently, mutation at $\hfrakc_a$ leaves the
$\hfrakc_b$-column unchanged for every $b>a$.{} The
permutation-type box moves only permute these columns.
\end{lemma}

\begin{proof}
By the proof of Proposition~\ref{prop: lift of mutation box move},
the $r$ mutation-type moves occur when $p^{(a)}$ and $q^{(a)}$
cross, one for each $1\leq a\leq r$. The permutation-type box moves
only permute their positions. Hence
\[
\{i_{\hfrakc_1},\ldots,i_{\hfrakc_r}\}=I_{i_p}.
\]
By admissibility, these colors are pairwise nonadjacent.

Therefore, neither a horizontal relation in~{\rm (i)} nor a vertical
relation in~{\rm (ii)} of
\eqref{eq: full ibox exchange matrix} can occur between two distinct
mutation vertices. Thus
\[
b_{\hfrakc_a,\hfrakc_b}=0 \qquad(a\ne b)
\]
 at the corresponding stage of $\hfrakB_k$.

Now fix $a<b$. Under mutation at $\hfrakc_a$, the
$\hfrakc_b$-column changes according to
\[
b'_{\hfrakd,\hfrakc_b}
= b_{\hfrakd,\hfrakc_b} + [b_{\hfrakd,\hfrakc_a}]_+ [b_{\hfrakc_a,\hfrakc_b}]_+ - [-b_{\hfrakd,\hfrakc_a}]_+ [-b_{\hfrakc_a,\hfrakc_b}]_+.
\]
Since $b_{\hfrakc_a,\hfrakc_b}=0$,
we obtain
\[
b'_{\hfrakd,\hfrakc_b} = b_{\hfrakd,\hfrakc_b}
\]
for every row $\hfrakd$. Thus the mutation at $\hfrakc_a$ leaves the
$\hfrakc_b$-column unchanged.

Applying this successively along $\hfrakB_k$, and observing that the
permutation-type box moves only permute box labels and
seed positions, proves the assertion.
\end{proof}

We first consider the permutation-type case.

\begin{proposition}\label{prop: permutation ibox cluster folding}
Suppose that $\frakB_k$ is of permutation type. Then
\[
\phi^\frakC(\hx'_{\hfrakc'}) = x'_{\frakc'}
\qquad
\bigl( \frakc'\in\calF',\ \hfrakc'\in\Fib^{\frakC'}(\frakc') \bigr).
\]
\end{proposition}

\begin{proof}
By Proposition~\ref{prop: lift of permutation box move}, both
$\frakB_k$ and $\hfrakB_k$ act on the corresponding clusters only by
permutations. Proposition~\ref{prop: lifted box move endpoint family}
 identifies the permuted blocks with the box fibers for $\frakC'$.
 Hence the definition of $\phi^\frakC$ gives the asserted identities.{}
\end{proof}

We next consider the mutation-type case.

\begin{proposition}\label{prop: mutation ibox cluster folding}
Suppose that $\frakB_k$ is of mutation type. Then
\[
\phi^\frakC(\hx'_{\hfrakd}) =  x'_{\frakd}
\qquad
\bigl( \frakd\in\calF',\ \hfrakd\in\Fib^{\frakC'}(\frakd) \bigr).
\]
\end{proposition}

\begin{proof}
Let $\frakc\in\calF$ be the $\bfi$-box at which the downstairs
mutation occurs. By Proposition~\ref{prop: lift of mutation box move},
 $\hfrakB_k$ contains exactly $r$ mutation-type box moves,
separated by permutation-type box moves. Let
\[
\hfrakc_1,\ldots,\hfrakc_r
\]
be the corresponding mutation vertices, in their order of occurrence.

 We first consider the cluster variables not mutated by $\frakB_k$.
The variables in their corresponding fibers are not mutated along
$\hfrakB_k$; 
the permutation-type box moves only permute their
labels and positions.  Hence the definition of $\phi^\frakC$ gives{}
\[
\phi^\frakC(\hx'_{\hfrakd})=x'_{\frakd}
\]
for  every unchanged box $\frakd\in\calF$  and each
$\hfrakd\in\Fib^{\frakC'}(\frakd)$.

It remains to consider the affected fiber. Write  the exchange relation at $\frakc$  as
\[
x_\frakc x'_{\frakc'} = M^+ + M^-,
\]
where
\[
M^+ \seteq \prod_{\frakd\in\calF} x_\frakd^{[b^\calF_{\frakd,\frakc}]_+},
\qquad
M^- \seteq \prod_{\frakd\in\calF} x_\frakd^{[-b^\calF_{\frakd,\frakc}]_+},
\]
and $\frakc'\in\calF'$ is the new box replacing $\frakc$.

Fix $1\le a\le r$ and consider the mutation at $\hfrakc_a$.
By Lemma~\ref{lem: future mutation column}, after identifying the
vertices through the preceding permutation-type box moves, the
exchange column used at $\hfrakc_a$ is obtained from its corresponding
initial exchange column only by relabeling. In particular, none of
the preceding mutations changes this column.

Moreover, every previously mutated vertex $\hfrakc_b$ satisfies
\[
b_{\hfrakc_b,\hfrakc_a}=0.
\]
Hence the cluster variables created by the preceding mutations do not
occur in the exchange relation at $\hfrakc_a$. Thus, after identifying
vertices through the preceding permutations, the exchange relation
at $\hfrakc_a$ has the form
\[
\hx_{\hfrakc_a}\hx'_{\hfrakc'_a} = \hM^+_a +\hM^-_a,
\]
where
\[
\hM^+_a \seteq \prod_{\hfrakd\in\hcalF} \hx_{\hfrakd}^{[b^\hcalF_{\hfrakd,\hfrakc_a}]_+},
\qquad
\hM^-_a \seteq \prod_{\hfrakd\in\hcalF} \hx_{\hfrakd}^{[-b^\hcalF_{\hfrakd,\hfrakc_a}]_+}.
\]

We now apply  the fixed Laurent homomorphism $\phi^\frakC$.
For each $\frakd\in\calF$, we have
\[
\phi^\frakC(\hx_{\hfrakd})=x_\frakd \qquad \bigl(\hfrakd\in\Fib^{\frakC}(\frakd) \bigr).
\]
Hence
\[
\phi^\frakC(\hM^+_a) = \prod_{\frakd\in\calF} x_\frakd^{ \sum_{\hfrakd\in\Fib^{\frakC}(\frakd)} [b^\hcalF_{\hfrakd,\hfrakc_a}]_+ }.
\]
By Lemma~\ref{lem: ibox weak sign},
\[
[b^\calF_{\frakd,\frakc}]_+ = \sum_{\hfrakd\in\Fib^{\frakC}(\frakd)} [b^\hcalF_{\hfrakd,\hfrakc_a}]_+.
\]
Therefore
\[
\phi^\frakC(\hM^+_a) = \prod_{\frakd\in\calF} x_\frakd^{[b^\calF_{\frakd,\frakc}]_+} = M^+.
\]
Similarly,
\[
[-b^\calF_{\frakd,\frakc}]_+ = \sum_{\hfrakd\in\Fib^{\frakC}(\frakd)} [-b^\hcalF_{\hfrakd,\hfrakc_a}]_+,
\qtqh
\phi^\frakC(\hM^-_a)=M^-.
\]

Since $\hfrakc_a$ belongs to the fiber over $\frakc$, we have
\[
\phi^\frakC(\hx_{\hfrakc_a})=x_\frakc.
\]
Applying $\phi^\frakC$ to  the exchange relation at $\hfrakc_a$  gives
\[
x_\frakc\, \phi^\frakC(\hx'_{\hfrakc'_a}) = M^+ + M^-.
\]
Comparing this with $x_\frakc x'_{\frakc'}=M^+ + M^-$,
we obtain
\[
\phi^\frakC(\hx'_{\hfrakc'_a}) = x'_{\frakc'}.
\]

Since $a$ was arbitrary, this holds for all $r$ mutation vertices.
After the remaining permutation-type box moves, the variables
\[
\hx'_{\hfrakc'_1},\ldots,\hx'_{\hfrakc'_r}
\]
are precisely those labelled by the boxes in
\[
\Fib^{\frakC'}(\frakc').
\]
Together with the unchanged-variable case above, this proves the
assertion.
\end{proof}

\begin{corollary}\label{cor: ibox g-vector folding}
Define
\[
\Fold^{\frakC}:\Z^{\hcalF}\longrightarrow\Z^\calF,
\qquad
\bigl(\Fold^{\frakC}(\hbfg)\bigr)_{\frakc} = \sum_{\hfrakc\in\Fib^{\frakC}(\frakc)} \hbfg_{\hfrakc}.
\]
Then, for every $\frakc'\in\calF'$ and
$\hfrakc'\in\Fib^{\frakC'}(\frakc')$,
\[
\deg_{\Sigma^\calF}(x'_{\frakc'}) = \Fold^{\frakC}\left(\deg_{\Sigma^\hcalF}(\hx'_{\hfrakc'}) \right).
\]
\end{corollary}

\begin{proof}
The proof is the same as that of
Proposition~\ref{prop: g vector folding}, using
Theorem~\ref{thm: ibox pair folding} and
Propositions~\ref{prop: permutation ibox cluster folding}
and~\ref{prop: mutation ibox cluster folding}.
\end{proof}

\begin{example}\label{ex: ibox cluster folding A5 C3}
Consider the folding
\[
A_5\longrightarrow C_3,
\qquad
\widehat 1=(1,5),\quad \widehat 2=(2,4),\quad \widehat 3=(3).
\]
Take
\[
\bfi=(1,2,1,3,2,3,1),
\qquad
\calT=(\calL,\calL,\calR,\calL,\calL,\calL).
\]
The corresponding admissible chain and effective-end sequence are
\[
\frakC = \seq{[6],[5],[4,6],[7],[3,7],[2,5],[1,7]},
\qquad
\bfp(\frakC)=\seq{6,5,4,7,3,2,1}.
\]
Thus, in effective-end order,
\[
\calF = \seq{ [1,7],[2,5],[3,7],[4,6],[5],[6],[7] }.
\]

Consider the  box move $\frakB_4$. Since
\[
p_4=7,\qquad p_5=3,\qquad i_{p_4}=i_{p_5}=1,
\]
it is of mutation type and replaces
\[
\frakc=[7] \quad \text{with} \quad 
\frakc'=[3].
\]
The corresponding exchange relation is
\[
x_{[7]}x'_{[3]} = x_{[3,7]}+x_{[5]}.
\]

The fixed lift of $\bfi$ is
\[
\hii=(1,5,2,4,1,5,3,2,4,3,1,5),
\]
and the canonical lifted effective-end sequence is
\[
\seq{10,9,8,7,11,12,6,5,4,3,2,1}.
\]
The relevant part of  the lifted family $\hcalF$ before the move  is
\[
\Fib^{\frakC}([3,7]) = \seq{[5,11],[6,12]}, \ \ \Fib^{\frakC}([5])  = \seq{[8],[9]}, \ \  \Fib^{\frakC}([7]) = \seq{[11],[12]}.
\]
Hence
\begin{align*}
\phi^\frakC(\hx_{[5,11]}) &=\phi^\frakC(\hx_{[6,12]})=x_{[3,7]},\allowdisplaybreaks\\
\phi^\frakC(\hx_{[8]}) &=\phi^\frakC(\hx_{[9]})=x_{[5]},\allowdisplaybreaks\\
\phi^\frakC(\hx_{[11]}) &=\phi^\frakC(\hx_{[12]})=x_{[7]}.
\end{align*}

Here $r=2$ and $h=4$, so
Proposition~\ref{prop: lift of mutation box move} gives
\[
\hfrakB_4 = \bigl( \frakB^\up_6,\frakB^\up_5, \frakB^\up_7,\frakB^\up_6 \bigr),
\]
with the sequence of types \[
MPPM.
\]
The two mutation-type moves are
\[
[12]\longmapsto[6], \qquad [11]\longmapsto[5].
\]
Their exchange relations are
\[
\hx_{[12]}\hx'_{[6]} = \hx_{[6,12]}+\hx_{[9]}
\qtq
\hx_{[11]}\hx'_{[5]} = \hx_{[5,11]}+\hx_{[8]}.
\]
Applying $\phi^\frakC$ gives in both cases
\[
x_{[7]}\, \phi^\frakC(\hx') = x_{[3,7]}+x_{[5]} = x_{[7]}x'_{[3]}.
\]
Therefore
\[
\phi^\frakC(\hx'_{[5]}) = \phi^\frakC(\hx'_{[6]}) = x'_{[3]}.
\]

 After applying $\hfrakB_4$,
\[
\Fib^{\frakC'}([3]) = \seq{[5],[6]},
\]
so the two mutated variables are labelled by the boxes in
$\Fib^{\frakC'}([3])$.
This example also makes visible the
two ingredients in the mutation-type argument: the
permutation-type moves only permute the relevant exchange columns,
while the positive and negative exchange monomials fold separately.
\end{example}

\subsection{Folding for $\bfi$-box moves}

We conclude this section by collecting the preceding results for a
single box move.  
Let $\frakc_k$ be movable, put $\frakC'=\frakB_k(\frakC)$, and let 
$\hfrakB_k$
 be the sequence of box moves constructed in{}
Propositions~\ref{prop: lift of permutation box move}
and~\ref{prop: lift of mutation box move}.  Let \[
\calF=\calF(\frakC), \quad \calF'=\calF(\frakC') \qtq  \hcalF=\calF(\hfrakC), \quad \hcalF'=\calF(\hfrakC')
\]
be the corresponding families and their canonical lifts. \par
Let
\[
\Sigma^\calF = (\bfx^\calF,B^\calF,\Lambda^\calF),
\qquad
\Sigma^\hcalF = (\bfx^\hcalF,B^\hcalF,\Lambda^\hcalF)
\]
 be the seeds before the move, and let
\[
\Sigma^{\calF'} = (\bfx^{\calF'},B^{\calF'},\Lambda^{\calF'}),
\qquad
\Sigma^{\hcalF'} = (\bfx^{\hcalF'},B^{\hcalF'},\Lambda^{\hcalF'})
\]
 be the seeds after the move.

Recall  the fixed Laurent homomorphism
\[
\phi^\frakC: \calL(\Sigma^\hcalF) \longrightarrow \calL(\Sigma^\calF),
\qquad
\hx_{\hfrakc}\longmapsto x_\frakc \quad \bigl( \hfrakc\in\Fib^{\frakC}(\frakc) \bigr).
\]

\begin{theorem}\label{thm: ibox endpoint folding}
 The sequence $\hfrakB_k$ transforms
$\Sigma^\hcalF$ into $\Sigma^{\hcalF'}$.{}
 For
$\frakc',\frakd'\in\calF'$ and any{}
$\hfrakd'\in\Fib^{\frakC'}(\frakd')$, we have
\[
b^{\calF'}_{\frakc',\frakd'} = \sum_{\hfrakc'\in\Fib^{\frakC'}(\frakc')} b^{\hcalF'}_{\hfrakc',\hfrakd'},
\qtq
\Lambda^{\calF'}_{\frakc',\frakd'} = d_{i_{\frakc'}} \sum_{\hfrakc'\in\Fib^{\frakC'}(\frakc')} \Lambda^{\hcalF'}_{\hfrakc',\hfrakd'}.
\]
Moreover,
\[
\phi^\frakC(\hx'_{\hfrakc'}) = x'_{\frakc'}
\qquad
\bigl(\frakc'\in\calF',\ \hfrakc'\in\Fib^{\frakC'}(\frakc') \bigr).
\]
{} \end{theorem}

\begin{proof}
By Proposition~\ref{prop: lifted box move endpoint family}, the 
maximal commuting family obtained after $\hfrakB_k$ is $\hcalF'$.
 The two matrix identities follow from{} Proposition~\ref{prop: ibox endpoint compatible pair folding}.

For cluster variables, $\frakB_k$ is either of permutation or
mutation type. The desired identity follows respectively from
Proposition~\ref{prop: permutation ibox cluster folding} or
Proposition~\ref{prop: mutation ibox cluster folding}.
\end{proof}

\section{Global folding of cluster monomials} \label{sec: global folding}

We now collect the results of Sections~\ref{sec: braid} and~\ref{sec: box} and consider the
cluster algebras associated with positive braid monoid
elements.

Let $\ttb\in\ttB^+_\g$. Thus
\[
\ii=\seq{i_1,\ldots,i_l}\in\Seq(\ttb)
\]
if $\ttb=r_{i_1}\cdots r_{i_l}$.
Any two elements of $\Seq(\ttb)$ are connected by a sequence of
braid moves.

For $\ii\in\Seq(\ttb)$, choose a maximal commuting family $\calF$ of
$\bfi$-boxes associated with $\ii$, and let $\Sigma^\calF$ be the
corresponding box seed.  Let
\[
\Sigma_q^\calF=(\bfX^\calF,B^\calF,\La^\calF)
\]
be the corresponding quantum seed, where
$\bfX^\calF=(X_\frakc\mid\frakc\in\calF)$.
By the braid-move and box-move results above, the mutation class of
$\Sigma_q^\calF$ depends only on $\ttb$, and not on the choices of
$\ii$ and $\calF$.

In view of the preceding results, the following definition makes sense
for an arbitrary symmetrizable Cartan datum and an arbitrary positive
braid.

\begin{definition}\label{def: positive braid cluster algebra}
For $\ttb\in\ttB_\g^+$, we define
\[
\calA_q(\ttb) \seteq \calA_q(\Sigma_q^\calF)
\]
to be the quantum cluster algebra associated with $\ttb$.
Specializing at $q=1$, we define
\[
\calA(\ttb) \seteq \calA(\Sigma^\calF) \subseteq \calL(\Sigma^\calF), \qquad \Sigma^\calF=(\bfx^\calF,B^\calF).
\]
\end{definition}

At the commutative level, Section~\ref{sec: arbitrary positive braid folding}
identifies $\calA(\ttb)$ with the cluster algebra of Qin's all-positive
signed-word seed associated with $\ii$; see
\cite[Section~3.1]{Qin24B}.

Fix a maximal commuting family $\calF$ of $\bfi$-boxes, and choose an
admissible chain $\frakC$ realizing $\calF$.
 Let $\hii$ be the blockwise lift of $\ii$, let $\hfrakC$ be the
canonical lift of $\frakC$, and set $\hcalF=\calF(\hfrakC)$.{} Then
$\hii\in\Seq(\httb)$ for the corresponding upstairs positive braid
monoid element $\httb\in \ttB^+_\hg$, and $\hcalF$ determines the lifted
box seed
\[
\Sigma^\hcalF = (\bfx^\hcalF,B^\hcalF) \text{ of $\calA(\httb)$.}
\]

We fix the surjective Laurent polynomial homomorphism
\[
\phi^\calF: \calL(\Sigma^\hcalF) \twoheadrightarrow \calL(\Sigma^\calF),
\qquad
\hx_{\hfrakc}\longmapsto x_\frakc \quad \bigl(\hfrakc\in\Fib^\frakC(\frakc)\bigr).
\]

Our main result in this section is that this single Laurent
homomorphism simultaneously folds all cluster monomials arising from
box seeds associated with arbitrary sequences of $\ttb$.

\begin{theorem}\label{thm: global cluster monomial folding}
Let $\jj\in\Seq(\ttb)$, and let $\calG$ be any maximal commuting
family of   $\bfi$-boxes associated with $\jj$. Choose an admissible chain $\frakD$ realizing
 $\calG$. Let $\hjj$ be the blockwise lift of $\jj$, let
$\widehat{\frakD}$ be the canonical lift of $\frakD$, and set
$\hcalG=\calF(\widehat{\frakD})$.{} Then
\[
\phi^\calF\bigl(\hx^\hcalG_{\hfrakd}\bigr) = x^\calG_{\frakd}
\qquad
\bigl( \frakd\in\calG,\  \hfrakd\in\Fib^\frakD(\frakd) \bigr).
\]
In particular, $\phi^\calF$ sends every cluster monomial associated
with $\Sigma^\hcalG$ to a cluster monomial associated with
$\Sigma^\calG$.
\end{theorem}

\begin{proof}
 We first show that if $\phi^\calF$ sends each upstairs cluster variable
before a braid move or box move to the corresponding
downstairs variable, then the same assertion holds after the move.{}

 Let $\widehat{\bfu}$ be the upstairs cluster before the chosen
lifted sequence, and suppose that $\phi^\calF$ sends each of its
variables to the corresponding downstairs variable.  Set
$\bfv\seteq\phi^\calF(\widehat{\bfu})$, where
$\phi^\calF$ is applied componentwise.  Thus $\bfv$ lists the downstairs cluster variables, with repetitions
along the fibers. Let $z$ be a downstairs cluster variable after the
move, and let $\hz$ be a variable in the corresponding upstairs fiber.  By the Laurent phenomenon, there are a polynomial $P$
and a monomial $M$ with nonnegative exponents such that
\[
M(\widehat{\bfu})\hz=P(\widehat{\bfu}).
\]
Every entry of $\widehat{\bfu}$ and the variable $\hz$ belong
to $\calL(\Sigma^\hcalF)$ by the Laurent phenomenon. Applying
$\phi^\calF$ to this identity therefore gives
\[
M(\bfv)\phi^\calF(\hz)=P(\bfv).
\]
By the braid-move results of Section~\ref{sec: braid} and
Theorem~\ref{thm: ibox endpoint folding} for box moves,  the Laurent homomorphism determined by
$\widehat{\bfu}\mapsto\bfv$ sends $\hz$ to $z$.
Hence
\[
P(\bfv)=M(\bfv)z.
\]
Comparing the last two identities and cancelling the nonzero
monomial $M(\bfv)$ gives $\phi^\calF(\hz)=z$.

To apply the braid-move results, choose a sequence of box moves
from $\frakC$ to the chain with constant $\calR$-sequence, whose
family is $\calF^\ii$.
By \eqref{eq: hfrakT} and the convention for sequences of length
$1$, its canonical lift has family $\calF^\hii$: the order of the
singleton boxes in the first fiber does not affect the family.
By Theorem~\ref{thm: ibox endpoint folding} and the preceding
argument, $\phi^\calF$ identifies the cluster variables of
$\Sigma^{\calF^\hii}$ with those of $\Sigma^{\calF^\ii}$.

Choose a sequence of braid moves
\[
\ii=\ii^{(0)} \longrightarrow \ii^{(1)} \longrightarrow \cdots \longrightarrow \ii^{(N)}=\jj.
\]
By the folding results of Section~\ref{sec: braid}, together
with Proposition~\ref{prop: reduction to connected configurations},
 the identity $\phi^\calF(\hz)=z$ is preserved{} 
 under each braid move.
Induction along the chosen
braid-move sequence shows that $\phi^\calF$ identifies the cluster
variables of  $\Sigma^{\calF^\hjj}$ with those of $\Sigma^{\calF^\jj}${}.

For the fixed sequence $\jj$, connect $\calF^\jj$ to $\calG$ by box moves. 
By Theorem~\ref{thm: ibox endpoint folding},  the identity $\phi^\calF(\hz)=z$ is preserved{} under each box move.
Hence
\[
\phi^\calF\bigl(\hx^\hcalG_{\hfrakd}\bigr) =x^\calG_{\frakd}
\qquad
\bigl(\frakd\in\calG,\ \hfrakd\in\Fib^\frakD(\frakd) \bigr).
\]

Since $\phi^\calF$ is an algebra homomorphism,
\[
\phi^\calF\left(\prod_{\hfrakd\in\hcalG}(\hx_{\hfrakd}^{\hcalG})^{a_{\hfrakd}}\right)
= \prod_{\frakd\in\calG} (x_{\frakd}^{\calG})^{\sum_{\hfrakd\in\Fib^\frakD(\frakd)}a_{\hfrakd}},
\]
which is a cluster monomial associated with $\Sigma^\calG$.
\end{proof}

\begin{corollary}\label{cor: global g-vector folding}
With the notation of
Theorem~\ref{thm: global cluster monomial folding}, we have
\[
\deg_{\Sigma^\calF}(x_{\frakd}^{\calG})=\Fold^{\frakC}\left(\deg_{\Sigma^\hcalF}(\hx_{\hfrakd}^{\hcalG})\right)
\]
for every $\frakd\in\calG$ and
$\hfrakd\in\Fib^{\frakD}(\frakd)$.
The same identity holds for the corresponding cluster monomials.
\end{corollary}

\begin{proof}
The proof is the same as that of
Proposition~\ref{prop: g vector folding}, using
Theorem~\ref{thm: global cluster monomial folding}.
\end{proof}

The cluster algebras $\calA(\httb)$ and $\calA(\ttb)$ are
contained in the Laurent polynomial algebras
$\calL(\Sigma^\hcalF)$ and $\calL(\Sigma^\calF)$,
respectively.  Theorem~\ref{thm: global cluster monomial folding}
shows that the Laurent polynomial homomorphism $\phi^\calF$ is
compatible with all cluster monomials arising from box seeds.
Our main commutative folding theorem strengthens this compatibility
to the cluster algebras themselves.

\begin{theorem}\label{thm: cluster algebra folding}
The Laurent polynomial homomorphism $\phi^\calF$ restricts to a
surjective $\C$-algebra homomorphism
\[
\phi^\calF_\ttb:
\calA(\httb)\twoheadrightarrow\calA(\ttb).
\]
\end{theorem}

\begin{remark}\label{rem: comparison with earlier folding constructions}

As mentioned in the Introduction, the approach in
\eqref{eq: philosophy} has been explored in various settings.
We compare these earlier constructions with Theorem~\ref{thm: cluster algebra folding}.

\bna
\item Demonet~\cite{Demonet11} studies folding through additive
categorification for specific Weyl group elements admitting
adapted reduced expressions.
\item Huang and Li~\cite{HL18Unfolding} construct unfoldings of
acyclic sign-skew-symmetric matrices using locally finite
quivers with group actions. The existence of an acyclic seed
can be difficult to determine in general.
\item\label{item: Dupont folding}
Dupont~\cite{Dupont08} studies folding via orbit mutations
under the global foldability condition, establishing
surjectivity in finite cluster type.
\item Felikson, Shapiro, and Tumarkin~\cite{FST12Unfoldings}
study unfoldings under a compatibility condition similar
to that in~\eqref{item: Dupont folding}, requiring
compatibility under all mutation sequences.
\ee

In contrast, Theorem~\ref{thm: cluster algebra folding} applies
to arbitrary positive braids over symmetrizable Cartan data,
without assuming finite cluster type, acyclicity, an adapted
reduced expression, or global foldability.

\end{remark}

We first study two settings in which this folding homomorphism admits
quantum-group-theoretic realizations and the simply-laced upstairs
algebras admit monoidal categorifications.  In
Section~\ref{sec: reduced}, we consider positive braids represented
by reduced sequences and use quantum unipotent coordinate algebras.
In Section~\ref{sec: finite case}, we treat arbitrary positive braids
in finite type through positive-braid algebras and quantum virtual
Grothendieck rings.  We then prove
Theorem~\ref{thm: cluster algebra folding} in full generality in
Section~\ref{sec: arbitrary positive braid folding} using the
standard-monomial basis~\cite{Qin24B}.

\section{Reduced sequences and quantum unipotent coordinate algebras}
\label{sec: reduced}

Let $w\in\weyl$ and choose a reduced sequence
$\ii\in\Red(w)$.  In this setting, the quantum cluster algebra
associated with $\ii$ will be identified with the quantum unipotent
coordinate algebra $\sfA_q(\n(w))$.   The blockwise lift $\hii$ is
again reduced{} and determines the corresponding upstairs algebra
$\sfA_q(\n(\hw))$.  Both upstairs and downstairs, the dual
PBW generators are categorified by cuspidal modules over the
corresponding quiver Hecke algebras. Moreover, since the upstairs
Cartan datum is simply-laced, the upstairs cluster algebra admits a
monoidal categorification by quiver Hecke modules.

In this section, we show that, after specialization at $q=1$, the
 homomorphism $\phi^\calF$ of
Section~\ref{sec: global folding}{} restricts to a surjective algebra
homomorphism
\[
\calA(\hw)\xrightarrow{\ \phi^\calF_w\ }\calA(w).
\]
The proof uses these realizations and the polynomial
presentations by dual PBW generators.

\subsection{Quantum unipotent coordinate algebras}

We briefly recall only the notation and properties of quantum
unipotent coordinate algebras needed below.  For the definitions,
notation, and basic properties, we refer to \cite{Kimura12}
(see also~\cite{KKKO18,FHOO}).

\smallskip

Recall that $\g$ is the symmetrizable Kac--Moody algebra and let $\calU_q(\g)$ be its quantized enveloping
algebra over $\Q(q^{1/2})$.  We denote by $\calU_q^+(\g)$ the subalgebra generated by
$e_i$ $(i\in I)$, and by $\calU_\bbA^+(\g)$ its $\bbA$-form generated by
the divided powers $e_i^{(n)}$, where $\bbA = \Z[q^{\pm1/2}]$.

Set
\[
\sfA_q(\n) \seteq  \bigoplus_{\beta\in \sfQ^-}   \Hom_{\Q(q^{1/2})}\bigl(\calU_q^+(\g)_{-\beta},\Q(q^{1/2})\bigr).
\]
The multiplication dual to the coproduct on $\calU_q^+(\g)$ endows
$\sfA_q(\n)$ with an algebra structure.  We denote by
$\sfA_\bbA(\n)$ the corresponding integral form consisting
of those linear forms taking $\calU_\bbA^+(\g)$ into $\bbA$.

For $w\in \weyl$, let $\sfA_q(\n(w))$ denote the quantum
unipotent coordinate algebra associated with $w$.
Equivalently, it is the subalgebra of $\sfA_q(\n)$
corresponding to $\Phi^+\cap w\Phi^-$.

For $\ii=(i_1,\ldots,i_r)\in\Red(w)$ with $w\in\weyl$, set
\[
w_k^\ii=s_{i_1}\cdots s_{i_k}, \quad w_{<k}^\ii=s_{i_1}\cdots s_{i_{k-1}}, \quad \be_k^\ii=w_{<k}^\ii(\alpha_{i_k}).
\]
Then we have
\[
\Phi^+\cap w\Phi^- = \{\be_1^\ii,\ldots,\be_r^\ii\}.
\]

For $\lambda\in P^+$ and $\mu,\zeta\in\weyl\lambda$, we denote by
\[
D(\mu,\zeta)\in\sfA_\bbA(\n)
\]
the corresponding unipotent quantum minor.  Recall that
$D(\mu,\zeta)$ is either zero or an element of the upper global basis (or dual canonical basis) $\bfB$
of $\sfA_q(\n)$; see \cite[Section~9.1]{KKKO18}.

For $\ii=(i_1,\ldots,i_r)\in\Red(w)$ with $w\in\weyl$ and
$1\le k\le r$, let
\[
F^{\up}(\be_k^\ii)\in\sfA_q(\n(w))
\]
denote the dual PBW root vector corresponding to $\be_k^\ii$.
With a proper $q$-normalization, \cite[Section~5.1.1]{KK24}{} gives the exact identity
\[
D\bigl(w_k^\ii\varpi_{i_k},w_{<k}^\ii\varpi_{i_k}\bigr)
=F^{\up}(\be_k^\ii)
\qquad (1\le k\le r).
\]
Thus, $F^{\up}(\be_k^\ii)$ belongs to $\sfA_\bbA(\n(w))$, and
$\sfA_\bbA(\n(w))$ is generated by these elements as an
$\bbA$-algebra. The generator indexed by $k$ corresponds to the positive
root $\beta_k^\ii$.

We shall mainly work with the specialization at $q=1$.  Regarding
$\C$ as an $\bbA$-algebra via $q^{1/2}\mapsto 1$, we set
\[
\sfA(\n(w)) \seteq \C\otimes_{\bbA}\sfA_\bbA(\n(w)).
\]
For $x\in\sfA_\bbA(\n(w))$, we write
\[
x|_{q=1}\seteq 1\otimes x\in\sfA(\n(w)) \quad \text{for its specialization at $q^{1/2}=1$.}
\]
By \cite[Theorem~4.39]{Kimura12}, we have an algebra isomorphism
\[
\sfA(\n(w))\simeq \C[N(w)],
\]
where $N(w)$ is the unipotent algebraic group with Lie algebra
$$\n(w) \seteq \soplus_{\alpha\in\Phi^+\cap w\Phi^-}\g_{-\alpha}$$ and $\C[N(w)]$ denotes its coordinate algebra.

It is known that $\sfA_\bbA(\n(w))$ admits a quantum cluster
algebra structure associated with each
$\ii=(i_1,\ldots,i_r)\in\Red(w)$; see
\cite{GLS11,GY17,GY21,KKKO18,Qin20}.  The corresponding quantum seed has the compatible pair $(B^\ii,\La^\ii)$
and cluster variables, up to the standard $q$-normalization,
\[
D\bigl(w_k^\ii\varpi_{i_k},\varpi_{i_k}\bigr) \qquad (1\le k\le r).
\]
By specializing the above quantum cluster algebra structure at $q=1$,
we obtain a cluster algebra structure on $\sfA(\n(w))$ whose initial
seed is $\Sigma^\ii$. Hence
\[
\calA(w) \seteq \calA(\Sigma^\ii) \simeq \sfA(\n(w)) \simeq \C[N(w)].
\]
Henceforth, we identify $\calA(w)$ with $\sfA(\n(w))$ via the above
isomorphism.

We recall the dual PBW basis associated with
$\ii\in\Red(w)$. The  suitably normalized {}ordered monomials in the dual PBW root vectors
$F^{\up}(\be_k^\ii)$ form the dual PBW basis of
$\sfA_q(\n(w))$. In particular,
\[
\sfA_q(\n(w)) = \Q(q^{1/2}) \bigl\langle F^{\up}(\be_1^\ii),\ldots,F^{\up}(\be_r^\ii) \bigr\rangle .
\]
Moreover, each $F^{\up}(\be_k^\ii)$ belongs to (\cite[Proposition~4.23]{Kimura12})
\[
\bfB(w)\seteq \bfB\cap\sfA_\bbA(\n(w)).
\]

We call
\[
\bfB^{\rm cl}(w) \seteq \left\{ b|_{q=1}\mid b\in\bfB(w) \right\} \ \ \text{the \emph{classical dual canonical basis} of $\sfA(\n(w))$.}
\]

For $1\le k\le r$, set
\[
\sff_k^\ii \seteq F^{\up}(\be_k^\ii)|_{q=1} \in \sfA(\n(w)).
\]
Under the identification  $\sfA(\n(w))\simeq\C[N(w)]$,
the elements $\sff_1^\ii,\ldots,\sff_r^\ii$ are the coordinate functions
dual to the root vectors corresponding to
$\be_1^\ii,\ldots,\be_r^\ii$.  By the PBW theorem, their monomials
form a $\C$-basis of $\C[N(w)]$.  Hence
\[
\calA(w) \simeq \sfA(\n(w)) = \C[\sff_1^\ii,\ldots,\sff_r^\ii].
\]
In particular, $\calA(w)$ is a polynomial algebra in $r=\ell(w)$ variables.

Quantum cluster monomials of $\sfA_\bbA(\n(w))$ belong, up to $q^{\frac{1}{2}\Z}$,  to $\bfB(w)$; this follows from \cite{KKKO18}
in the symmetric case and from \cite{Qin20} in the symmetrizable case. Consequently, every cluster monomial of $\calA(w)$ belongs
to $\bfB^{\rm cl}(w)$ via the isomorphism.

\subsection{Categorification and dual PBW vectors} \label{subsec: PBW}

We briefly recall only the notation and results on quiver Hecke
algebras that are needed below.  Throughout this subsection, we follow
the notation and conventions of \cite{KKKO18,KKOP18}, where quiver
Hecke algebras, their module categories, and determinantial modules
are developed in detail.

\smallskip

Let $R(\beta)$ $(\beta\in\sfQ^+)$ be a quiver Hecke algebra
(\cite{KL1,KL2,R08}) associated with a Cartan datum.
We denote by $R(\beta)\gmod$ the category of finite-dimensional
$\mathbb Z$-graded $R(\beta)$-modules, and set
\[
R\gmod \seteq \soplus_{\beta\in\sfQ^+} R(\beta)\gmod.
\]

The convolution product makes $R\gmod$ into a monoidal category.
There is an $\bbA$-algebra isomorphism \cite{KL1,KL2,R08}
\[
\Omega: \bbA\otimes_{\Z[q^{\pm1}]} K_{\Z[q^{\pm1}]}(R\gmod) \xrightarrow{\sim} \sfA_\bbA(\n).
\]

For $\lambda\in\wl^+$ and $\mu,\zeta\in\weyl\lambda$ such that
$D(\mu,\zeta)\ne0$, we denote by $M(\mu,\zeta)$ the corresponding
\emph{determinantial} module (\cite{KKOP18}).  It is the unique \emph{self-dual} simple module satisfying
\[
\Omega([M(\mu,\zeta)])=D(\mu,\zeta).
\]

Let $w\in\weyl$ and $\ii=(i_1,\ldots,i_r)\in\Red(w)$.
For $1\le k\le r$, set
\[
S_k^\ii \seteq M\bigl(w_k^\ii\varpi_{i_k},w_{<k}^\ii\varpi_{i_k}\bigr).
\]
Then
\[
-\wt(S_k^\ii)=\be_k^\ii \qtq \Omega([S_k^\ii]) = F^{\up}(\be_k^\ii).
\]
The module $S_k^\ii$ is the $k$-th \emph{cuspidal} module associated with $\ii$.

The modules $S_k^\ii$ are the categorical counterparts of the dual
PBW root generators associated with $\ii$; see
\cite[Section~5.1]{KK24}.

Recall that, for an $\ii$-box $[x,y]$, the corresponding
determinantial module is
\[
M^\ii[x,y] = M\bigl(w_y^\ii\varpi_{i_y}, w_{<x}^\ii\varpi_{i_x}\bigr);
\]
(see \cite{KKKO18} and \cite[Section~5.1]{KK24}).  Therefore, for a singleton box
$[k]$, we have $M^\ii[k]=S_k^\ii$, and the corresponding cluster variable is precisely $\sff_k^\ii$.

We now apply the preceding discussion to the folding setting.
 Recall that $\hii$ is the blockwise lift of $\ii$.{}  Since $\ii$ is
reduced, so is $\hii$ by Corollary~\ref{cor: reduced lift}.  For $1\le k\le r$ and
$1\le s\le e_{i_k}$, set
\[
\hsff_{k,s} \seteq F^{\up}(\hbe_{k,s})|_{q=1} \in \calA(\hw),
\]
where $\hii\in\Red(\hw)$ for $\hw\in\hsfW$.

\begin{lemma}\label{lem: singleton PBW generators}
For every $1\le k\le r$, the cluster variable associated with the
singleton $\ii$-box $[k]$ is $\sff_k^\ii$.
Moreover, the fiber of $[k]$ consists of the singleton
$\bfi$-boxes associated with $\hii$
\[
[\kappa(k,1)],\ldots, [\kappa(k,e_{i_k})],
\]
whose corresponding cluster variables are
\[
\hsff_{k,1},\ldots,\hsff_{k,e_{i_k}}.
\]
Furthermore,
\[
\calA(w)=\C[\sff_1^\ii,\ldots,\sff_r^\ii] \qtq \calA(\hw) = \C[\hsff_{k,s}\mid 1\le k\le r,\ 1\le s\le e_{i_k}].
\]
\end{lemma}

\begin{proof}
By \cite{KKOP18} (see also \cite[Section~5.1]{KK24}), the module associated with an
$\ii$-box $[x,y]$ is
\[
M^\ii[x,y] = M\bigl(w_y^\ii\varpi_{i_y}, w_{<x}^\ii\varpi_{i_x}\bigr)
\qtqh
M^\ii[k] = M\bigl(w_k^\ii\varpi_{i_k}, w_{<k}^\ii\varpi_{i_k}\bigr) = S_k^\ii.
\]
By the preceding discussion, the corresponding cluster variable is
$\sff_k^\ii$.

Let $\calF$ be a maximal commuting family containing $[k]$, choose an
admissible chain $\frakC$ realizing $\calF$, and let $\hcalF$ be its
canonical lift. By
Lemma~\ref{lem: big epsilon box fiber},
\[
\Fib^\frakC([k]) = \bigl\{ [\kappa(k,1)],\ldots,[\kappa(k,e_{i_k})] \bigr\}.
\]
Since $\hii$ is reduced, applying the first assertion to $\hii$ shows
that the corresponding upstairs cluster variables are
\[
\hsff_{k,1},\ldots,\hsff_{k,e_{i_k}}.
\]

Finally, the PBW generation recalled above gives
\[
\calA(w) = \C[\sff_1^\ii,\ldots,\sff_r^\ii].
\]
Applying the same argument to $\hii\in\Red(\hw)$, we obtain
\[
\calA(\hw) = \C[\hsff_{k,s}\mid 1\le k\le r,\ 1\le s\le e_{i_k}],
\]
which completes the proof.
\end{proof}

\subsection{The folding homomorphism via PBW generators}

We now realize the folding homomorphism in terms of dual PBW
generators.  Fix a maximal commuting family $\calF$ of
$\bfi$-boxes, choose an admissible chain $\frakC$ realizing $\calF$,
and let $\hcalF$ be its canonical lift.
Recall from Section~\ref{sec: global folding} the fixed surjective Laurent polynomial
homomorphism
\[
\phi^\calF: \calL(\Sigma^{\hcalF}) \twoheadrightarrow \calL(\Sigma^\calF),
\qquad
\hx_{\hfrakc}\longmapsto x_\frakc
\quad
\bigl(\hfrakc\in\Fib^\frakC(\frakc)\bigr).
\]
Since
\[
\calA(\hw)\subset\calL(\Sigma^{\hcalF}) \qtq \calA(w)\subset\calL(\Sigma^\calF),
\]
it remains to show that $\phi^\calF$ restricts to a surjective
homomorphism from $\calA(\hw)$ onto $\calA(w)$.

\begin{proposition}\label{prop: folding PBW generators}
For every $1\le k\le r$ and $1\le s\le e_{i_k}$, we have
\[
\phi^\calF(\hsff_{k,s})=\sff_k^\ii.
\]
\end{proposition}

\begin{proof}
Fix $1\le k\le r$, and choose a maximal commuting family $\calG$
of $\bfi$-boxes containing the singleton $[k]$.  By
Lemma~\ref{lem: singleton PBW generators}, the cluster variable
associated with $[k]$ is $\sff_k^\ii$, while its lifted fiber consists
of
\[
[\kappa(k,1)],\ldots,[\kappa(k,e_{i_k})],
\]
with corresponding cluster variables
\[
\hsff_{k,1},\ldots,\hsff_{k,e_{i_k}}.
\]
Since $\ii$ remains fixed, $\calF$ and $\calG$ are connected by
$\ii$-box moves. Hence repeated application of
Theorem~\ref{thm: ibox endpoint folding} gives
\[
\phi^\calF(\hsff_{k,s})=\sff_k^\ii \qquad (1\le s\le e_{i_k}),
\]
as desired.
\end{proof}

\begin{theorem}\label{thm: reduced cluster algebra folding}
Let $w\in W$ and $\ii\in\Red(w)$.  Then
\[
\phi^\calF_w
\seteq
\left.\phi^\calF\right|_{\calA(\hw)}:
\calA(\hw)\twoheadrightarrow\calA(w)
\]
is a well-defined surjective algebra homomorphism.  Under the
polynomial presentations by dual PBW generators, it is determined by
\[
\phi^\calF_w(\hsff_{k,s})=\sff_k^\ii
\qquad
(1\le k\le r,\ 1\le s\le e_{i_k}).
\]
\end{theorem}

\begin{proof}
By Lemma~\ref{lem: singleton PBW generators}, we have
\[
\calA(\hw) = \C[\hsff_{k,s}\mid 1\le k\le r,\ 1\le s\le e_{i_k}]
\qtq
\calA(w) = \C[\sff_1^\ii,\ldots,\sff_r^\ii].
\]
By Proposition~\ref{prop: folding PBW generators},
\[
\phi^\calF(\hsff_{k,s}) = \sff_k^\ii \qquad (1\le k\le r,\ 1\le s\le e_{i_k}).
\]
It follows that
\[
\phi^\calF\bigl(\calA(\hw)\bigr) \subseteq \calA(w),
\]
so that the restriction $\phi^\calF_w$ is well-defined.

Moreover, for every $1\le k\le r$, we have
\[
\sff_k^\ii = \phi^\calF_w(\hsff_{k,1}),
\]
and hence all algebra generators of $\calA(w)$ belong to the image of
$\phi^\calF_w$. Therefore
\[
\phi^\calF_w\bigl(\calA(\hw)\bigr) = \calA(w),
\]
which proves the surjectivity.
\end{proof}

\begin{corollary}
\label{cor: classical dual canonical basis folding}
Let $\ii\in\Red(w)$, let $\calG$ be any maximal commuting family of
$\ii$-boxes, and let $\hcalG$ be its canonical lift. Then
\[
\phi^\calF_w \left(\prod_{\hfrakc\in\hcalG} (\hx_{\hfrakc}^{\hcalG})^{a_{\hfrakc}} \right)
= \prod_{\frakc\in\calG} (x_{\frakc}^{\calG})^{ \sum_{\hfrakc\in\Fib(\frakc)}a_{\hfrakc}} \in \bfB^{\rm cl}(w)
\]
for every $(a_{\hfrakc})_{\hfrakc\in\hcalG}\in\Z_{\geq0}^{\hcalG}$.
In particular, $\phi^\calF_w$ sends every cluster monomial associated
with $\Sigma^{\hcalG}$ to an element of $\bfB^{\rm cl}(w)$.
\end{corollary}

\begin{proof}
This follows from Theorem~\ref{thm: global cluster monomial folding}
and the fact that cluster monomials of $\calA(w)$ belong to
$\bfB^{\rm cl}(w)$.
\end{proof}

\begin{remark}\label{rem: canonical basis interpretation}
The preceding folding result also admits a quantum interpretation at
the level of canonical bases.  Recall that quantum cluster monomials
belong, up to the standard $q$-normalization, to $\bfB$.
In the symmetric case this follows from the monoidal categorification
by quiver Hecke algebras; see \cite{KKKO18}, while the symmetrizable
case follows from \cite{Qin20}.

Let $\hz$ be a quantum cluster monomial associated with a canonical
lifted box seed, and let $z$ be the corresponding downstairs quantum
cluster monomial.  Then

\[
\hz\in q^{\Z/2}\hbfB,\qquad z\in q^{\Z/2}\bfB.
\]
Moreover, their $g$-vectors satisfy
\[
\deg_\ii(z) = \Fold^\ii\bigl(\deg_\hii(\hz)\bigr),
\]
by the $g$-vector folding established above.
Thus Theorem~\ref{thm: global cluster monomial folding} gives a
  folding correspondence between the classical dual canonical basis
  elements $\hz|_{q=1}$ and $z|_{q=1}$.

We emphasize that this interpretation does not arise from an algebra
homomorphism between the corresponding quantum tori.  In general, the
fiberwise identification of quantum torus generators does not preserve
their $q$-commutation relations.  Rather, the quantum correspondence is
detected through the folding of $g$-vectors.  After specialization at
$q=1$, this obstruction disappears, and the correspondence is realized
by the Laurent algebra homomorphism $\phi^\calF$.
\end{remark}

\section{Positive braids of finite Cartan type and quantum virtual
Grothendieck rings}
\label{sec: finite case}

Throughout this section, we assume that $\g$ is of finite type, but
the sequence representing $\ttb\in\ttB_\g^+$ need not be reduced.
In this setting, we identify the associated quantum cluster algebra
$\calA_q(\ttb)$ with the positive-braid subalgebra
$\scrA_{\bbA}(\ttb)$ defined through the braid group action on the
bosonic extension, also known as the quantum virtual Grothendieck
ring.  For the simply-laced upstairs algebra
$\scrA_{\bbA}(\httb)$, we use its monoidal categorification by the
category $\scrC_\hg(\httb)$.

Using these realizations, we show that, after specialization at
$q=1$, the ambient Grothendieck-ring homomorphism $\uppi$ restricts
to a surjective algebra homomorphism
\[
\uppi_\ttb: \calA(\httb)\twoheadrightarrow\calA(\ttb)
\]
and agrees with $\phi^\calF_\ttb$.

\subsection{Quantum virtual Grothendieck rings}

We briefly recall only the notation and properties needed below.
Throughout this section, we follow the notation and conventions of
\cite{FHOO,JLO1,JLO2,Qin24B}, where the quantum virtual Grothendieck
rings associated with $\g$ and its simply-laced unfolding $\hg$ are
developed in detail.

\smallskip

Fix $t\in\C^\times$ which is not a root of unity.
Let $\calU_t(L\g)$ be the quantum loop algebra associated with $\g$,
and let $\scrC_\g$ be the category of finite-dimensional type~$1$
$\calU_t(L\g)$-modules.  It is a rigid monoidal category, and its
Grothendieck ring $K(\scrC_\g)$ is commutative.

The simple objects of $\scrC_\g$ are parametrized by \emph{dominant}
monomials in the variables $y_{i,a}$ $(i\in I,\ a\in\C^\times)$.
For a dominant monomial $m$, we denote by $L(m)$ the corresponding
simple module.  The $t$-character homomorphism of
Frenkel--Reshetikhin \cite{CP95,FR98,FR99}
\footnote{In the literature, this homomorphism is usually called the
$q$-character homomorphism and the corresponding Laurent variables
are usually denoted by $Y_{i,a}$.  We use $t$ and $y_{i,a}$ here to
avoid confusion with the quantum parameter and the quantum Laurent
variables used elsewhere in this paper.}
gives an injective algebra homomorphism
\[
\chi_t: K(\scrC_\g) \lhook\joinrel\longrightarrow \Z[y_{i,a}^{\pm1}\mid i\in I,\ a\in\C^\times].
\]

We first consider the simply-laced Lie algebra $\hg$.  Choose a
$\Gamma$-invariant parity function
\[
\epsilon:\hI\longrightarrow\{0,1\}, \qquad \epsilon_{\im}\equiv\epsilon_{\jm}+1\pmod 2
\quad\text{if }\im\sim\jm,
\]
and set
\[
\tDynkin_0^\hg \seteq \{(\im,p)\in\hI\times\Z \mid p\equiv\epsilon_{\im}\pmod 2\}.
\]
For $(\im,p)\in\tDynkin_0^\hg$, we write
\[
\hy_{\im,p}\seteq y_{\im,t^p}.
\]

Following Hernandez--Leclerc, let $\scrC_\hg^0$ be the Serre
subcategory of $\scrC_\hg$ whose simple objects have dominant
monomials in the variables
\[
\hy_{\im,p} \qquad ((\im,p)\in\tDynkin_0^\hg).
\]
The category $\scrC_\hg^0$ is a rigid monoidal subcategory of
$\scrC_\hg$ and is often referred to as the Hernandez--Leclerc
category, or the skeleton category.  Indeed, every simple object of
$\scrC_\hg$ can be expressed as a commuting tensor product of
spectral-parameter shifts of simple objects in $\scrC_\hg^0$
(see \cite{HL10,HL16,FHOO}).

The $t$-character homomorphism restricts to an injective algebra
homomorphism
\[
\chi_t: K(\scrC_\hg^0) \lhook\joinrel\longrightarrow
\calY_{\Z}(\hg) \seteq \Z[\hy_{\im,p}^{\pm1}\mid(\im,p)\in\tDynkin_0^\hg].
\]
Via the $t$-character homomorphism, we identify
$K(\scrC_\hg^0)$ with its image in $\calY_{\Z}(\hg)$ and write
\[
\frakK_{\Z}(\hg) \seteq \chi_t\bigl(K(\scrC_\hg^0)\bigr)
\subseteq \calY_{\Z}(\hg).
\]
We use the complexified classical rings
\[
\calY(\hg) \seteq \C\otimes_{\Z}\calY_{\Z}(\hg)
\simeq \C[\hy_{\im,p}^{\pm1}\mid(\im,p)\in\tDynkin_0^\hg],
\qquad
\frakK(\hg) \seteq \C\otimes_{\Z}\frakK_{\Z}(\hg)
\subseteq \calY(\hg).
\]
Thus $\frakK(\hg)$ is the complexified Grothendieck ring of
$\scrC_\hg^0$ realized as a ring of Laurent polynomials.

We recall the virtual counterpart associated with the possibly
non-simply-laced Lie algebra $\g$.  The $\Gamma$-invariant parity function
$\epsilon$ descends to $I$ by setting
\[
\epsilon_i\seteq\epsilon_{\im} \qquad (\im\in I_i).
\]
Accordingly, we set
\[
\tDynkin_0^\g \seteq \{(i,p)\in I\times\Z \mid p\equiv\epsilon_i\pmod 2\},
\qquad
\calY_{\Z}(\g) \seteq \Z[y_{i,p}^{\pm1}\mid(i,p)\in\tDynkin_0^\g].
\]
We use its complexification
\[
\calY(\g) \seteq \C\otimes_{\Z}\calY_{\Z}(\g)
\simeq \C[y_{i,p}^{\pm1}\mid(i,p)\in\tDynkin_0^\g].
\]
The quotient map $\hI\twoheadrightarrow I$ induces a surjective
algebra homomorphism
\[
\uppi_{\Z}: \calY_{\Z}(\hg)\twoheadrightarrow\calY_{\Z}(\g),
\qquad \hy_{\im,p}\longmapsto y_{i,p} \quad (\im\in I_i).
\]
By extension of scalars, it induces a surjective $\C$-algebra
homomorphism
\begin{align} \label{eq: uppi}
\uppi: \calY(\hg)\twoheadrightarrow\calY(\g), \qquad \hy_{\im,p}\longmapsto y_{i,p} \quad (\im\in I_i).
\end{align}

Following \cite{JLO1}, we define the virtual Grothendieck ring of
$\g$ and its complexification by
\[
\frakK_{\Z}(\g) \seteq \uppi_{\Z}\bigl(\frakK_{\Z}(\hg)\bigr)
\subseteq \calY_{\Z}(\g),
\qquad
\frakK(\g) \seteq \C\otimes_{\Z}\frakK_{\Z}(\g)
=\uppi\bigl(\frakK(\hg)\bigr) \subseteq \calY(\g).
\]
Thus $\uppi$ restricts to a surjective algebra homomorphism
\[
\uppi: \frakK(\hg)\twoheadrightarrow\frakK(\g).
\]
 See \cite{JLO1} for the construction and basic properties of $\frakK(\g)$.{}

We recall the fundamental and Kirillov--Reshetikhin
polynomials.  For $(i,p),(i,s)\in\tDynkin_0^\g$ with $p\le s$, set
\[
m^{(i)}[p,s] \seteq y_{i,p}y_{i,p+2}\cdots y_{i,s}.
\]
We set $m^{(i)}[p,s] \seteq 1$ if $p>s$. We call $m^{(i)}[p,s]$ a Kirillov--Reshetikhin \emph{$($KR$)$ monomial}.
Associated with it is a distinguished element
\[
F\bigl(m^{(i)}[p,s]\bigr)\in\frakK(\g),
\]
which we call the corresponding \emph{KR-polynomial}.  In particular,
\[
F(y_{i,p}) = F\bigl(m^{(i)}[p]\bigr) \text{ is called a \emph{fundamental polynomial}.}
\]

When $\g$ is simply laced, these elements are the $t$-characters of
the corresponding Kirillov--Reshetikhin modules.  In particular,
$F(y_{i,p})$ is the $t$-character of the fundamental module
$L(y_{i,p})$. Thus, for $\im\in\hI$, the Laurent polynomial
\[
F(\hy_{\im,p})\in\frakK(\hg) \text{ is the $t$-character of the fundamental $\calU_t(L\hg)$-module $L(\hy_{\im,p})$.}
\]

The Frenkel--Mukhin algorithm provides an algorithmic construction of
these Laurent polynomials. In particular, the KR-polynomials above are
well defined and can be obtained from their simply-laced counterparts
by folding.

More precisely, for $\im\in I_i$ and
$(\im,p),(\im,s)\in\tDynkin_0^\hg$, set
\[
\widehat{m}^{(\im)}[p,s] \seteq \hy_{\im,p}\hy_{\im,p+2}\cdots\hy_{\im,s}.
\]
Then
\[
\uppi\bigl(F(\widehat{m}^{(\im)}[p,s])\bigr) = F\bigl(m^{(i)}[p,s]\bigr).
\]
In particular,
\[
\uppi\bigl(F(\hy_{\im,p})\bigr) = F(y_{i,p})
\]
(see \cite[Corollary~3.10 and Remark~5.40]{JLO1} for instance).

We next consider the quantum deformation.  Recall that
$\sfD=\diag(d_i\mid i\in I)$ is the symmetrizer of the Cartan matrix
$\sfC=(c_{ij})_{i,j\in I}$.  For an indeterminate $t$, set
\[
c_{ij}(t) \seteq \delta_{ij}(t+t^{-1})+\delta_{i\ne j}c_{ij},
\qquad \sfC(t)\seteq(c_{ij}(t))_{i,j\in I}.
\]
We call $\sfC(t)$ the \emph{$t$-quantized Cartan matrix} of $\g$.
Set
\[
\sfB(t)\seteq \sfC(t)\sfD^{-1}.
\]
Then $\sfB(t)$ is symmetric.  Regarding $\sfB(t)$ as a matrix over
$\Q(t)$, let
\[
\tsfB(t) = \bigl(\tsfB_{ij}(t)\bigr)_{i,j\in I} \seteq \sfB(t)^{-1}.
\]
For $i,j\in I$, we write
\[
\tsfB_{ij}(t) = \sum_{u\in\Z}\tb_{ij}(u)t^u \in\Q((t)) \quad \text{ for its Laurent expansion at $t=0$.}
\]

For $(i,p),(j,s)\in\tDynkin_0^\g$, set
\[
\calN(i,p;j,s) \seteq \tb_{ij}(p-s-1)-\tb_{ij}(s-p-1) -\tb_{ij}(p-s+1)+\tb_{ij}(s-p+1).
\]
Since $\tsfB(t)$ is symmetric, we have $\calN(i,p;j,s)=-\calN(j,s;i,p)$.

Let $\calY_q(\g)$ be the $\bbA$-algebra generated by
\[
Y_{i,p}^{\pm1} \qquad ((i,p)\in\tDynkin_0^\g)
\]
subject to the relations
\[
Y_{i,p}Y_{i,p}^{-1} = Y_{i,p}^{-1}Y_{i,p} = 1 \qtq  Y_{i,p}Y_{j,s} = q^{\calN(i,p;j,s)} Y_{j,s}Y_{i,p}.
\]
We call $\calY_q(\g)$ the quantum torus associated with $\sfC(t)$. \par
Specializing $q^{1/2}$ to $1$, we recover the commutative Laurent
polynomial rings over $\Z$ introduced above.  More precisely, there are
surjective algebra homomorphisms
\[
\ev_{q=1}:\calY_q(\g)\twoheadrightarrow\calY_{\Z}(\g),
\qquad
Y_{i,p}\longmapsto y_{i,p}.
\]
Passing to the commutative cluster-algebra setting means composing
this specialization with extension of scalars from $\Z$ to $\C$.

We recall the quantum Grothendieck rings inside these quantum
tori.  When $\hg$ is of simply-laced finite type, the quantum
Grothendieck ring
\[
\frakK_q(\hg)\subseteq\calY_q(\hg)
\]
is the quantum deformation of the integral form $\frakK_{\Z}(\hg)$ constructed
geometrically by Nakajima and Varagnolo--Vasserot; see
\cite{Nak01,Nak04,VV03}.  An algebraic approach using deformed
screening operators was developed by Hernandez \cite{Her04}.

For every dominant monomial  $M$, there is a distinguished element
\[
 F_q(M)\in\frakK_q(\hg)
\]
having  $M${} as its unique dominant monomial.  The elements
\[
 \{F_q(M)\mid M\text{ is dominant}\}
\]
form an $\bbA$-basis of $\frakK_q(\hg)$.  Moreover, specialization
at $q=1$ recovers the corresponding elements of the ordinary
integral Grothendieck ring:
\[
 \ev_{q=1}\bigl(F_q(M)\bigr)=F\bigl(\ev_{q=1}(M)\bigr).
\]

For the possibly non-simply-laced Lie algebra $\g$, we consider the
\emph{quantum virtual Grothendieck ring}
\[
\frakK_q(\g)\subseteq\calY_q(\g),
\]
introduced by Kashiwara--Oh \cite{KO23} and further developed in
\cite{JLO1,JLO2}.  It admits an $\bbA$-basis
 $\{F_q(M)\mid M\text{ is dominant}\}$,
whose specialization at $q=1$ gives the $F$-basis of
$\frakK_{\Z}(\g)$ (see \cite[Corollary~5.28]{JLO1}):
\[
 \ev_{q=1}\bigl(F_q(M)\bigr) = F\bigl(\ev_{q=1}(M)\bigr).
\]

For a \emph{quantum} KR-monomial
\[
M^{(i)}[p,s] = Y_{i,p}Y_{i,p+2}\cdots Y_{i,s},
\]
where the product is understood with the standard $q$-normalization,
we call
\[
F_q\bigl(M^{(i)}[p,s]\bigr)\in\frakK_q(\g)
\]
the corresponding \emph{quantum KR-polynomial}.  In particular,
$F_q(Y_{i,p})$ is called a \emph{quantum fundamental polynomial}.
By the specialization property above,
\[
\ev_{q=1}\left(F_q\bigl(M^{(i)}[p,s]\bigr)\right) = F\bigl(m^{(i)}[p,s]\bigr).
\]
The same notation and terminology will be used for $\hg$.

The quantum Grothendieck ring $\frakK_q(\hg)$ also admits the
basis
\[
  \{L_q(M)\mid M\text{ is dominant}\},
\]
constructed by Nakajima and referred to as the \emph{canonical basis}.  Its specialization at $q=1$ recovers the
$t$-characters of the corresponding simple modules:
\[
\ev_{q=1}\bigl(L_q(M)\bigr) = \chi_t\bigl(L(m)\bigr), \qquad m=\ev_{q=1}(M).
\]
Moreover, if $M$ is a KR-monomial, then
\[
F_q(M)=L_q(M).
\]
These properties follow from Nakajima's results~\cite{Nak01}.

For a dominant monomial $M$, let $E_q(M)$ be the product of quantum
fundamental polynomials with the multiplicities prescribed by $M$,
ordered by increasing spectral label and normalized to have
coefficient~$1$ at $M$; see \cite[(5.21)]{JLO1}.
These elements form the \emph{standard basis} of $\frakK_q(\g)$.

For the quantum virtual Grothendieck ring $\frakK_q(\g)$, the
corresponding $L_q$-basis
$\{L_q(M)\mid M\text{ is dominant}\}$
was constructed in \cite{JLO1}.  It is characterized by the
bar-invariance and unitriangularity properties with respect to the
 standard basis $\{E_q(M)\mid M\text{ is dominant}\}$.  A fundamental result concerning quantum KR-polynomials is
\[
F_q\bigl(M^{(i)}[p,s]\bigr)
=
L_q\bigl(M^{(i)}[p,s]\bigr)
\]
for every KR-monomial $M^{(i)}[p,s]$.  In the simply-laced case,
this is due to Nakajima, while the general finite-type case was
proved by Qin; see \cite{Nak03} and \cite[Theorem~1.1]{Qin24B}.

The rings $\frakK(\g)$ and $\frakK_q(\g)$ are generated by their
fundamental polynomials.  More precisely,
\[
\frakK(\g) = \C\bigl[F(y_{i,p})\mid(i,p)\in\tDynkin_0^\g\bigr]
\qtq
\frakK_q(\g) = \bbA\left\langle F_q(Y_{i,p})\mid(i,p)\in\tDynkin_0^\g\right\rangle.
\]

For $i\in I$, let $\sfh$ denote the Coxeter number of $\g$, and let $i^*\in I$ be determined by
$w_0(\alpha_i)=-\alpha_{i^*}$, where $w_0$ is the longest element of
the Weyl group.
We denote by $\calD_q$ the duality automorphism, given on the ambient
variables by
\[
\calD_q(Y_{i,p})=Y_{i^*,p+\sfh}.
\]
We denote by $\calD$ its specialization at $q^{1/2}=1$, so that
\[
\calD(y_{i,p})=y_{i^*,p+\sfh}.
\]
We use the same notation for their natural extensions to the ambient
quantum torus and its skew field of fractions.

A \emph{Dynkin quiver} of $\g$ is a pair
$Q=(\Dynkin,\xi),$
where $\Dynkin$ is the Dynkin diagram of $\g$ endowed with an
orientation and
\[
\xi:I\longrightarrow\Z
\]
is a height function satisfying
\[
\xi_i=\xi_j+1 \qquad \text{whenever } i\to j.
\]
We choose $\xi$ compatibly with the parity function fixed above, so
that
\[
\xi_i\equiv\epsilon_i\pmod2 \qquad(i\in I).
\]
Set
\[
\tDynkin_{\le\xi}^\g \seteq \{(i,p)\in\tDynkin_0^\g\mid p\le\xi_i\}.
\]

For a Dynkin quiver $Q=(\Dynkin,\xi)$, set
\[
\frakK_q(\g)_{\le\xi} \seteq \bbA\left\langle F_q(Y_{i,p}) \mid (i,p)\in\tDynkin_{\le\xi}^\g \right\rangle.
\]
It is known that
\[
\frakK_q(\g)_{\le\xi} = \soplus_M  \bbA \; F_q(M),
\]
where the sum runs over all dominant monomials $M$ involving only
the variables $Y_{i,p}$ with $(i,p)\in\tDynkin_{\le\xi}^\g$; see
\cite{HL16,B21,FHOO,JLO1}.

For a Dynkin quiver $Q=(\Dynkin,\xi)$, a vertex $i\in I$ is called
a \emph{source} if
\[
\xi_i>\xi_j
\qquad\text{for every }j\in I\text{ adjacent to }i.
\]
If $i$ is a source, let $s_iQ=(\Dynkin,s_i\xi)$, where
\[
(s_i\xi)_j = \xi_j-2\delta_{ij}.
\]
We call $s_iQ$ the reflection of $Q$ at $i$.

Let $Q=(\Dynkin,\xi)$ be a Dynkin quiver.  A sequence
$\ii=\seq{i_1,i_2,\ldots}$
of finite or infinite length is said to be \emph{adapted to $Q$} if
$i_1$ is a source of $Q$ and, inductively, $i_k$ is a source of the
Dynkin quiver obtained by successively reflecting at
$i_1,\ldots,i_{k-1}$.

We shall also use the lifted Dynkin quiver associated with $Q$.
Define
\[
\hQ=(\hDynkin,\hxi), \qquad \hxi_{\im}\seteq\xi_i \quad (\im\in I_i),
\]
where each edge of $\hDynkin$ lying above an oriented edge
$i\to j$ of $\Dynkin$ is oriented from $I_i$ to $I_j$.
Then $\hxi$ is a $\Gamma$-invariant height function on $\hDynkin$.
  We call $\hQ$ the \emph{lifted Dynkin quiver} associated with $Q$.
Its quotient by $\Gamma$ is $Q$.

\begin{lemma}\label{lem: lifting adapted sequences}
Let $Q=(\Dynkin,\xi)$ be a Dynkin quiver and let
$\hQ=(\hDynkin,\hxi)$ be its   lifted Dynkin quiver.
If $\ii=\seq{i_1,i_2,\ldots}$
is a sequence of finite or infinite length adapted to $Q$,  then its
blockwise lift{}
\[
\hii = \seq{\im_{1,1},\ldots,\im_{1,e_{i_1}},  \im_{2,1},\ldots,\im_{2,e_{i_2}},\ldots},
\qquad
I_{i_k} = \seq{\im_{k,1},\ldots,\im_{k,e_{i_k}}},
\]
is adapted to $\hQ$.
\end{lemma}

\begin{proof}
If $i$ is a source of $Q$, then every vertex of $I_i$ is a source
of $\hQ$.  Moreover, the vertices in $I_i$ are pairwise
nonadjacent, so the corresponding source reflections commute.
After reflecting successively at all vertices of $I_i$, the height function is
\[
\hxi'_{\im} = \hxi_{\im}-2\delta_{\im\in I_i},
\]
which is precisely the lift of the height function obtained from
$\xi$ by reflecting $Q$ at $i$.
The assertion now follows inductively.
\end{proof}

Every infinite adapted sequence considered below is assumed to
contain each color infinitely often.

For each position $k$ of $\ii$, set
\[
p_k \seteq \xi_{i_k} - 2\#\{1\le a<k\mid i_a=i_k\}.
\]
If $\ii$ is infinite, then
\[
\{(i_k,p_k)\mid k\ge1\} = \tDynkin_{\le\xi}^\g.
\]

\begin{theorem}[{\cite{HL16,B21,FHOO,JLO1,KKOP24A}}]
\label{thm: cluster structure of K xi}
Let $Q=(\Dynkin,\xi)$ be a Dynkin quiver and let
$\ii=(i_1,i_2,\ldots)$ be an infinite sequence adapted to $Q$.
Then $\frakK_q(\g)_{\le\xi}$ has a quantum cluster algebra
structure whose initial quantum seed is
\[
\left( \left\{F_q\bigl(M^{(i_k)}[p_k,\xi_{i_k}]\bigr) \mid k\ge1 \right\}, B^\ii,\La^\ii \right).
\]
\end{theorem}

\begin{theorem}[{\cite{KKOP24A,Qin20,Qin24B}}]
\label{thm: cluster monomials and global basis}
Every quantum cluster monomial of
$\frakK_q(\g)_{\le\xi}$ is an element of the $L_q$-basis of
$\frakK_q(\g)$.
\end{theorem}

\begin{remark}\label{rem: Lq global basis}
Under the identification of $\frakK_q(\g)$ with the corresponding
bosonic extension established in \cite{HL15,JLO2}, the $L_q$-basis
is identified with the normalized global basis $\tbfG$ of \cite{KKOP24} (see~\eqref{eq: nor global basis} below).
Indeed, the two bases are uniquely characterized by the same
bar-invariance and unitriangularity conditions with respect to the
 standard basis $\{E_q(M)\}$ and its image in the bosonic extension.  Consequently, every quantum
cluster monomial corresponds to a normalized global basis element.
\end{remark}

By the Kazhdan--Lusztig characterization of the $L_q$-basis, its
transition matrix from the  standard basis $\{E_q(M)\}$ is unitriangular over
$\Z[q^{\pm1}]$.  Specializing at $q=1$ therefore gives a
unitriangular transition matrix over $\Z$.  Hence
\[
\sfL\seteq
\{\ev_{q=1}(L_q(M))\mid M\text{ is dominant}\}
\]
forms a $\Z$-basis of $\frakK_{\Z}(\g)$ and, after complexification,
a $\C$-basis of $\frakK(\g)$, which we call the \emph{classical global basis}.\footnote{This basis is also referred
to as the \emph{canonical basis}.}
By Remark~\ref{rem: Lq global basis}, it is the specialization at
$q=1$ of the normalized global basis $\tbfG$.
Set
\[
\frakK(\g)_{\le\xi} \seteq \C\bigl[ F(y_{i,p}) \mid (i,p)\in\tDynkin_{\le\xi}^\g \bigr].
\]
Specializing Theorem~\ref{thm: cluster structure of K xi} at
$q=1$, we obtain the following.

\begin{corollary}\label{cor: cluster structure of K xi classical}
Let $Q=(\Dynkin,\xi)$ be a Dynkin quiver and let
$\ii=(i_1,i_2,\ldots)$ be an infinite sequence adapted to $Q$.
Then $\frakK(\g)_{\le\xi}$ has a cluster algebra structure whose
initial seed is
\[
\left(\left\{ F\bigl(m^{(i_k)}[p_k,\xi_{i_k}]\bigr) \mid k\ge1 \right\},  B^\ii \right).
\]
Moreover, every cluster monomial of $\frakK(\g)_{\le\xi}$ belongs to the $L$-basis $\sfL$ of
$\frakK(\g)$.
\end{corollary}

\begin{proof}
This follows by specializing Theorem~\ref{thm: cluster structure of K xi} and Theorem~\ref{thm: cluster monomials and global basis}
at $q=1$.
\end{proof}

\subsection{Cluster structures on Grothendieck subrings}

For $i\in I$, set $q_i\seteq q^{d_i}$ and
\[
[n]_{q_i} \seteq \frac{q_i^n-q_i^{-n}}{q_i-q_i^{-1}},
\qquad
[n]_{q_i}! \seteq \prod_{k=1}^n[k]_{q_i},
\qquad
\begin{bmatrix} n\\ k\end{bmatrix}_{q_i}
\seteq
\frac{[n]_{q_i}!}{[k]_{q_i}![n-k]_{q_i}!}.
\]

\begin{definition} \cite{HL15,FHOO,JLO2}
The \emph{bosonic extension} $\scrA_q(\g)$ associated with $\g$ is
the $\Q(q^{1/2})$-algebra generated by $\{f_{i,p} \ | \  (i\in I,\ p\in\Z) \}$, subject to the relations
\begin{align*}
&\sum_{k=0}^{1-c_{ij}} (-1)^k
\begin{bmatrix}
1-c_{ij}\\ k
\end{bmatrix}_{q_i}
f_{i,p}^{\,k}f_{j,p}f_{i,p}^{\,1-c_{ij}-k} =0 \qquad (i\ne j)\\
&
f_{i,m}f_{j,p} = q^{(-1)^{p-m+1}(\alpha_i,\alpha_j)} f_{j,p}f_{i,m} + \delta_{(j,p),(i,m+1)}(1-q_i^2)
\qquad (m<p).
\end{align*}

For $-\infty\le a\le b\le\infty$, let
$\scrA_q(\g)[a,b]$ be the subalgebra generated by
$f_{i,p}$ with $a\le p\le b$.  We set
\[
\scrA_q(\g)_{\ge m} \seteq \scrA_q(\g)[m,\infty], \qquad \scrA_q(\g)_{\le m} \seteq \scrA_q(\g)[-\infty,m],
\]
and
$\scrA_q(\g)_{>m}$, $\scrA_q(\g)_{<m}$ similarly. We also abbreviate $\scrA_q(\g)[m]\seteq\scrA_q(\g)[m,m]$.

\end{definition}

When there is no danger of confusion, we omit $\g$ from the notation
and simply write
\[
\scrA_q,\qquad \scrA_q[a,b],\qquad \scrA_q[m],\qquad \scrA_{q,\ge m},\qquad \scrA_{q,\le m},
\]
and similarly for the other subalgebras.

Recall that there are anti-involutions  
$\omega$ and $\overline{(\,\cdot\,)}$ of $\scrA_q(\g)$ defined by 
\[
\omega(q^{1/2})=q^{1/2}, \quad \omega(f_{i,p})=f_{i,-p},
\qtq 
\overline{q^{1/2}}=q^{-1/2}, \quad \overline{f_{i,p}}=f_{i,p},
\]
respectively (\cite[(2.3)]{KKOP24B}).  They commute, and hence
$\Upomega \seteq \omega\circ\overline{(\,\cdot\,)}$ is a multiplicative
$\mathbb Q$-algebra automorphism satisfying
\[
\Upomega (q^{1/2})=q^{-1/2}, \qquad
\Upomega (f_{i,p})=f_{i,-p}. 
\]

Set
\[
\kappa_i\seteq q_i^{-1/2}(1-q_i^2), \qquad F_{i,p}\seteq\kappa_i^{-1}f_{i,p}, \qquad F_{i,p}^{(r)} \seteq \frac{F_{i,p}^r}{[r]_{q_i}!}.
\]

\begin{definition} \cite{KKOP21B,JLO2,KKOP24B}
For each $i\in I$, define the braid symmetry
\[
\TT_i:\scrA_q\longrightarrow\scrA_q
\]
on the generators by
\[
\TT_i(f_{j,p}) =
\begin{cases}
f_{i,p+1}, & j=i,\\[2mm]
\displaystyle \sum_{\substack{r,s\ge0\\ r+s=-c_{ij}}} (-q_i)^s F_{i,p}^{(r)}\,f_{j,p}\,F_{i,p}^{(s)}, & j\ne i.
\end{cases}
\]
\end{definition}

The maps $\TT_i$ are  $\Q(q^{1/2})$-algebra automorphisms and satisfy the braid
relations.  Hence they define an action of the braid group
$\ttB_\g$ on $\scrA_q$.  We denote the inverse of $\TT_i$ by
$\TT_i^\star$, so that
\[
\TT_i^\star \TT_i=\TT_i\TT_i^\star={\rm id}.
\]

\begin{remark}
For each $m\in\Z$, the local subalgebra
$\scrA_q[m]$ is isomorphic to $\calU_q^-(\g)$.
The braid symmetries $\TT_i$ on $\scrA_q$ extend Lusztig's
braid symmetries on these local pieces; see
\cite{KKOP24,KKOP24B}.
\end{remark}

For $\ttb\in\ttB_\g^+$ and
$\ii=(i_1,\ldots,i_l)\in\Seq(\ttb)$, set
\[
\TT_{\ttb} \seteq \TT_{i_1}\cdots \TT_{i_l}.
\]
By the braid relations, $\TT_{\ttb}$ is independent of the choice of
$\ii\in\Seq(\ttb)$.

\begin{definition} \cite{OP24,KKOP24B}
For $\ttb\in\ttB_\g^+$, define
\[
\scrA_q(\ttb) \seteq \scrA_{q,\ge0} \cap \TT_{\ttb}\bigl(\scrA_{q,<0}\bigr).
\]
We call $\scrA_q(\ttb)$ the \emph{positive-braid subalgebra}
associated with $\ttb$.
\end{definition}

Let $\ii=(i_1,\ldots,i_l)\in\Seq(\ttb)$.  For $1\le k\le l$, set
\[
P_k^\ii \seteq \TT_{i_1}\cdots \TT_{i_{k-1}} \bigl(q_{i_k}^{1/2}f_{i_k,0}\bigr),
\]
and call it the \emph{cuspidal element} of $\ii$ at $k$.
For $m\in\Z_{\ge0}$, set
\[
P_k^{\ii,\{m\}} \seteq q_{i_k}^{m(m-1)/2}(P_k^\ii)^m.
\]
For $\bfu=(u_1,\ldots,u_l)\in\Z_{\ge0}^l$, define
\[
P^\ii(\bfu) \seteq \oprod_{k\in[1,l]} P_k^{\ii,\{u_k\}}.
\]
The product is taken in decreasing order of $k$.
We shall also use the integral form of the bosonic extension.
Let
\[
\scrA_{\bbA} \subseteq \scrA_q
\]
be its $\bbA$-integral form. The integral form $\scrA_{\bbA}$ is preserved by
$\TT_i^{\pm1}$.  For $\ttb\in\ttB_\g^+$, set
\[
\scrA_{\bbA}(\ttb) \seteq \scrA_q(\ttb)\cap\scrA_{\bbA}.
\]

For finite type $\g$ and each Dynkin quiver $Q$ of $\g$, the $\bbA$-integral form
$\scrA_{\bbA}(\g)$ is isomorphic to the quantum virtual
Grothendieck ring $\frakK_q(\g)$ introduced above
\cite[Theorem~6.2]{KKOP25M}; see also \cite{HL15} and
\cite[Theorems~5.2 and~7.2]{JLO2}.  In the simply-laced case, this recovers the
quantum Grothendieck ring of the Hernandez--Leclerc category.

More precisely, for every Dynkin quiver $Q=(\Dynkin,\xi)$, there is a
$\mathbb Z$-algebra isomorphism
\[
\Upphi_Q:\scrA_{\bbA}\xrightarrow{\sim}\frakK_q(\g),
\qquad
\Upphi_Q(q^{1/2})=q^{-1/2}
\]
(see \cite[Theorem 6.2]{KKOP25M}). 
It restricts to a $\mathbb Z$-isomorphism
\[
\Upphi_Q:
(\scrA_{\bbA})_{\le0}
\xrightarrow{\sim}
\frakK_q(\g)_{\le\xi}.
\]
We set
\[
\Psi_Q\seteq\Upphi_Q\circ\Upomega.
\]
Since both $\Upphi_Q$ and $\Upomega$ invert the quantum parameter $q$,
$\Psi_Q$ is an $\bbA$-algebra isomorphism
\[
\Psi_Q:\scrA_{\bbA}\xrightarrow{\sim}\frakK_q(\g),
\]
whose restriction gives an $\bbA$-algebra isomorphism
\[
\Psi_Q:
(\scrA_{\bbA})_{\ge0}
\xrightarrow{\sim}
\frakK_q(\g)_{\le\xi}.
\]

For later use, we also recall the substitution isomorphisms of
\cite{FHOO2,LO24}.  Let $\bbF_q(\g)$ denote the skew field of fractions
of the ambient quantum torus, and set
\[
\bbF(\g)\seteq\C(y_{i,p}\mid(i,p)\in\tDynkin_0^\g),
\]
the corresponding field of rational functions after specializing
$q^{1/2}=1$.

For two Dynkin quivers $Q$ and $Q'$ of $\g$, set
\[
\Psi_{Q,Q'}
\seteq
\Psi_Q\circ\Psi_{Q'}^{-1}
=
\Upphi_Q\circ\Upphi_{Q'}^{-1},
\]
where the second equality follows from the definition of $\Psi_Q$.
Thus $\Psi_{Q,Q'}$ is an $\bbA$-algebra automorphism of
$\frakK_q(\g)$.
By \cite{FHOO2,LO24}, this automorphism extends to an
$\bbA$-algebra automorphism
\[
\tPsi_{Q,Q'}:
\bbF_q(\g)\xrightarrow{\sim}\bbF_q(\g),
\]
called the \emph{substitution isomorphism}, which commutes with the
duality automorphism $\calD_q$. We denote its specialization at $q^{1/2}=1$ by the
$\C$-algebra automorphism
\[
\tPsi^{\rm cl}_{Q,Q'}:
\bbF(\g)\xrightarrow{\sim}\bbF(\g).
\]
 
Let $\calB(\infty)$ denote the crystal of $\calU_q^-(\g)$,
and let $1$ denote its highest weight element.
Recall that Kashiwara--Park \cite{KP22} introduced the
\emph{extended crystal}
\[
\hcalB(\infty) \seteq \left\{ (b_m)_{m\in\Z}\in\prod_{m\in\Z}\calB(\infty) \;\middle|\; b_m=1\text{ for all but finitely many }m \right\}.
\]
The bosonic extension $\scrA_q$ admits the normalized global basis
\begin{align} \label{eq: nor global basis}
\tbfG   = \{\tG(b)\mid b\in\hcalB(\infty)\}
\end{align}
which is compatible with the integral form and is preserved by the
braid symmetries $\TT_i$; see
\cite{KKOP24,KKOP24B}.

\begin{theorem}[{\cite{OP24,KKOP24B}}]
\label{thm: PBW basis bosonic braid}
Let $\ttb\in\ttB_\g^+$ with $\ell(\ttb)=l$ and let
$\ii\in\Seq(\ttb)$.  Then
\[
\calP^\ii \seteq \left\{ P^\ii(\bfu) \mid \bfu\in\Z_{\ge0}^l \right\}
\]
forms an $\bbA$-basis of $\scrA_{\bbA}(\ttb)$.
In particular, $\scrA_q(\ttb)$ is generated by the cuspidal elements
\[
P_1^\ii,\ldots,P_l^\ii.
\]
Moreover, the global basis is compatible with
$\scrA_q(\ttb)$; namely,
\[
\tbfG(\ttb) \seteq \tbfG\cap\scrA_q(\ttb) \ \ \ \text{forms an $\bbA$-basis of $\scrA_{\bbA}(\ttb)$.}
\]
\end{theorem}
In particular,
\[
\scrA_\bbA(\ttb) = \bbA\langle P_1^\ii,\ldots,P_l^\ii\rangle,
\]
with PBW basis $\calP^\ii$.

We recall a consequence of the Garside normal form that will be
used to realize arbitrary finite sequences inside adapted ones.
Let
\[
\pi:\ttB_\g\longrightarrow \sfW
\]
be the canonical projection, and let $\Updelta\in\ttB_\g^+$ be the
positive braid lift of the longest element $w_0\in \sfW$, characterized by
\[
\pi(\Updelta)=w_0, \qquad \ell(\Updelta)=\ell(w_0).
\]
A sequence is called \emph{locally reduced} if every consecutive
subsequence of length at most $\ell(w_0)$ is reduced
\cite[Definition~2.3(iii)]{KKOP25M}. 
By \cite[Remark~5.4]{KKOP25M}, for every finite sequence
\[
\ii=(i_1,\ldots,i_l),
\]
there exist $m\gg0$ and a sequence
\[
\tii\in\Seq(\Updelta^m) \quad  \text{such that} \quad \tii_{[1,l]}=\ii.
\]
Moreover, $\tii$ can be obtained from a locally reduced
sequence in $\Seq(\Updelta^m)$ by finitely many commutation and braid
moves, where the locally reduced sequence may be chosen to be adapted
to a suitable Dynkin quiver.

Now let $\ttb\in\ttB_\g^+$ and
$\ii=(i_1,\ldots,i_l)\in\Seq(\ttb)$, and choose
$\tii\in\Seq(\Updelta^m)$ as above.  Thus $\ii$
occurs as the initial segment of $\tii$:
\[
\tii_{[1,l]}=\ii.
\]

Choose a locally reduced sequence
$\jj\in\Seq(\Updelta^m)$
adapted to a Dynkin quiver
$Q=(\Dynkin,\xi)$
such that $\tii$ is obtained from $\jj$ by finitely many commutation
and braid moves.  Let
$\hQ=(\hDynkin,\hxi)$
be the   lifted Dynkin quiver associated with $Q$.

 By Lemma~\ref{lem: lifting adapted sequences}, the blockwise lift{}
$\hjj$ of $\jj$ is adapted to $\hQ$.  Moreover, by the lifting
results of the preceding sections, the sequence of commutation and
braid moves from $\jj$ to $\tii$ lifts to a sequence of moves from
 $\hjj$ to the blockwise lift $\htii$ of $\tii$.{}

Since $\ii$ is the initial segment of $\tii$, the corresponding
 initial block of $\htii$ is precisely the blockwise lift $\hii$ of{}
$\ii$.  {}  Writing
\[
\tii=\ii\star\tau, \qquad \htii=\hii\star\htau,
\]
 where $\htau$ is the blockwise lift of $\tau$,{} we obtain
\[
\begin{array}{ccc}
\hjj\text{ adapted to }\hQ &\longrightarrow& \htii=\hii\star\htau \\[1mm]
\downarrow && \downarrow \\[1mm]
\jj\text{ adapted to }Q &\longrightarrow& \tii=\ii\star\tau,
\end{array}
\]
where the horizontal arrows are given by the corresponding sequences
of commutation and braid moves.

We now apply this construction to the simply-laced upstairs
Grothendieck ring $\frakK_q(\hg)$.  Let
$\httb\in\ttB_\hg^+$ be the positive braid monoid element represented
 by the blockwise lift $\hii$ of $\ii$.{}  By the monoidal
categorification of \cite{KKOP25M}, the positive-braid subalgebra
$\scrA_{\bbA}(\httb)$ admits a categorical realization by the
monoidal subcategory $\scrC_\hg(\httb)$ of the Hernandez--Leclerc
category $\scrC_\hg^0$.  In particular, admissible chains of $\bfi$-boxes give
monoidal seeds whose quantum cluster algebra realizes
$\scrA_{\bbA}(\httb)$.

For a simple module $L$ in $\scrC_\hg^0$, let $[L]_q$ denote
its $(q,t)$-character, with deformation parameter $q$, and set
\[
L_q\seteq\Upphi_{\hQ}^{-1}([L]_q)\in\scrA_\bbA(\hg).
\]
The modules and the map $\Upphi_{\hQ}$ are associated with the
same Dynkin quiver $\hQ$.
For an admissible chain $\hfrakC$, the exchange and commutation
matrices in \cite[Definition~9.5 and Proposition~5.20]{KKOP25M}
are $-B^{\hfrakC}$ and $-\Lambda^{\hfrakC}$, respectively.
Thus \cite[Theorem~9.7]{KKOP25M}, with $t^{1/2}=q^{-1/2}$
and the exchange matrix negated, gives the quantum seed below.

\begin{theorem}[{\cite[Theorems~9.6 and~9.12]{KKOP25M}}]
\label{thm: monoidal seed for positive braid}
Let $\httb\in\ttB_\hg^+$ and
$\hii\in\Seq(\httb)$.  For each $\bfi$-box $\hfrakc$, let
$\hM[\hfrakc]$ be the corresponding affine determinantial module.
For a singleton box $[k]$, $\hM[k]$ is the corresponding
affine cuspidal module.

For every admissible chain
$\hfrakC=(\hfrakc_1,\ldots,\hfrakc_{\hl})$
of $\bfi$-boxes, the family
$\{\hM[\hfrakc_k]\mid 1\le k\le\hl\}$
forms a monoidal seed of $\scrC_\hg(\httb)$, and
\begin{align} \label{eq: quantum seed}
\left( \{\hM[\hfrakc_k]_q\}_{1\le k\le\hl}, B^{\hfrakC}, \Lambda^{\hfrakC} \right) \ \
\text{is a quantum seed of $\scrA_{\bbA}(\httb)$.}
\end{align}
\end{theorem}

\begin{corollary}
\label{cor: upstairs cluster algebra bosonic}
For every $\httb\in\ttB_\hg^+$, there is an
$\bbA$-algebra isomorphism
\[
\calA_q(\httb) \xrightarrow{\sim} \scrA_{\bbA}(\httb).
\]
Under this isomorphism, the cluster variable associated with an
$\bfi$-box $\hfrakc$ is identified with
$\hM[\hfrakc]_q$.
\end{corollary}

\begin{proof}
Fix $\hii\in\Seq(\httb)$.  By
Theorem~\ref{thm: monoidal seed for positive braid}, for every
admissible chain
$\hfrakC=(\hfrakc_1,\ldots,\hfrakc_{\hl})$
of $\bfi$-boxes,
\[
\left( \{\hM[\hfrakc_k]_q\}_{1\le k\le\hl}, B^{\hfrakC}, \Lambda^{\hfrakC} \right)
\]
is a quantum seed of $\scrA_{\bbA}(\httb)$.
On the other hand, the $\bfi$-box seed defining
$\calA_q(\httb)$ has the same compatible pair
$(B^{\hfrakC},\Lambda^{\hfrakC})$.

Moreover, under box moves, the two families of seeds are related by
the same mutation--permutation transformations.  Hence the
identification of the initial cluster variables
\[
X_{\hfrakc}\longmapsto \hM[\hfrakc]_q \qquad (\hfrakc\in\hfrakC)
\]
extends throughout the corresponding mutation classes and induces
the desired isomorphism
\[
\calA_q(\httb)
\xrightarrow{\sim}
\scrA_{\bbA}(\httb). \qedhere
\]
\end{proof}

\begin{theorem}[{\cite[Theorem~9.12 and its proof]{KKOP25M}}]
\label{thm: affine cuspidal and bosonic PBW}
Let $\httb\in\ttB_\hg^+$ and write
$\hii=\seq{\im_1,\ldots,\im_{\hl}}\in\Seq(\httb)$.
Then, for every $1\le k\le\hl$,
\[
\hM[k]_q=q_{\im_k}^{-1/2}P_k^\hii.
\]
\end{theorem}

We compare the quantum seeds associated with $\hii$ with those
obtained from the ambient sequence $\htii=\hii\star\htau$.

Extend $\hjj$ to an infinite $\hQ$-adapted sequence and take
the KR seed of Theorem~\ref{thm: cluster structure of K xi}.
Freeze the last occurrence of each color in $\hjj$ and delete
all vertices outside $\hjj$.
By \eqref{eq: GLS exchange matrix}, the deleted rows have zero
entries in every retained exchangeable column, and
\eqref{eq: GLS Lambda matrix} restricts to $\Lambda^{\hjj}$.
This gives a finite quantum seed in $\frakK_q(\hg)_{\le\hxi}$
with compatible pair $(B^{\hjj},\Lambda^{\hjj})$.
The same construction applies to $\jj$ and $Q$ downstairs.
The sequence of braid moves from $\hjj$ to $\htii$ is realized by
the corresponding mutation--permutation sequence.  Hence  the quantum seed associated with $\htii$ belongs to the same mutation
class.

\begin{proposition}[{\cite[Proposition~8.10]{KKOP25M}}]
\label{prop: cutting positive braid seed}
Let $\htii=\hii\star\htau$ be as above, and let
$\hfrakC$ be an admissible chain of $\bfi$-boxes associated with $\hii$.
Extend $\hfrakC$ to an admissible chain $\htfrakC$ of
$\htii$-boxes associated with $\htii$.  Assume that $\hfrakC$ has extent $[1,\hl]$ and that
$\htfrakC$ has length equal to that of $\htii$.  Write
\[
\hJ=\hJ_{\rm ex}\sqcup\hJ_{\rm fr}
\]
for the vertices corresponding to the $\hii$-part, and let
$\hJ_{\rm tail}$ denote the remaining vertices.  Then
\[
b^{\htfrakC}_{ab}=0 \quad (a\in\hJ_{\rm tail},\ b\in\hJ_{\rm ex}),
\qtq
B^{\htfrakC}\big|_{\hJ\times\hJ_{\rm ex}} = B^{\hfrakC},
\quad
\Lambda^{\htfrakC}\big|_{\hJ\times\hJ} = \Lambda^{\hfrakC}.
\]
Consequently, declaring the vertices in $\hJ_{\rm fr}$ frozen
and deleting the vertices in $\hJ_{\rm tail}$ from  the quantum
seed in \eqref{eq: quantum seed} for $\htfrakC$ gives
\[
\Sigma^{\hfrakC} = \left(\{\hM[\hfrakc]_q\}_{\hfrakc\in\hfrakC},B^{\hfrakC},\Lambda^{\hfrakC}\right),
\]
the quantum seed of $\scrA_{\bbA}(\httb)$ in~\eqref{eq: quantum seed}.
\end{proposition}

We return to the downstairs sequence
$\ii\in\Seq(\ttb)$.  Let
$\frakC=(\frakc_1,\ldots,\frakc_l)$
be an admissible chain of $\bfi$-boxes, and let $\hfrakC$ be its
canonical lift.  By Theorem~\ref{thm: ibox pair folding},
the compatible pair of the upstairs quantum seed
$\Sigma^{\hfrakC}$ folds to
$(B^{\frakC},\Lambda^{\frakC}).$
More precisely, for $\frakc,\frakd\in\frakC$ and
$\hfrakd\in\Fib^{\frakC}(\frakd)$,
\[
b^{\frakC}_{\frakc,\frakd} = \sum_{\hfrakc\in\Fib^{\frakC}(\frakc)} b^{\hfrakC}_{\hfrakc,\hfrakd},
\qquad
\Lambda^{\frakC}_{\frakc,\frakd} = d_{i_{\frakc}} \sum_{\hfrakc\in\Fib^{\frakC}(\frakc)} \Lambda^{\hfrakC}_{\hfrakc,\hfrakd}.
\]

Thus the exchange and quantum commutation data of the upstairs seed
fold precisely to the pair
$(B^{\frakC},\Lambda^{\frakC})$
associated with the downstairs admissible chain $\frakC$.

From now on, we distinguish the positive-braid subalgebra
$\scrA_{\bbA}(\ttb)$ of the bosonic extension from the quantum
cluster algebra associated with $\ttb$, which will be denoted by
$\calA_q(\ttb)$.

We consider the cluster variables on the downstairs side.
We first identify the cluster variables associated with the singleton boxes.
For this reason, we use the following result of Qin.
For the adapted sequences in \cite[Section~5.1]{Qin24B},
the quantization matrix is $\Lambda^\ii$ by
\cite[(7.1) and~(7.5)]{KO23} and \eqref{eq: GLS Lambda matrix}.
The equality is preserved by braid moves, by
Proposition~\ref{prop: mutation permutation braid moves}.
Restriction to an initial segment preserves the root sequence
and box vectors, and hence also this equality.
The quantization in \cite[Section~6.3]{Qin24B} is obtained by
these operations, so it agrees with the present quantization.

\begin{theorem}[{\cite[Theorem~6.6]{Qin24B}}]
\label{thm: Qin fundamental variables}
Let $\ii=(i_1,\ldots,i_l)$ be an arbitrary finite sequence.  Under the bosonic realization,
the cluster variable associated with the singleton $\ii$-box $[k]$ is represented by
\[
\TT_{i_1}\cdots \TT_{i_{k-1}}(f_{i_k,0}) \qquad (1\le k\le l).
\]
\end{theorem}

Here the braid symmetries in
Theorem~\ref{thm: Qin fundamental variables} agree with
the automorphisms $\TT_i$ defined above. 
 Comparing the formula in Theorem~\ref{thm: Qin fundamental variables} with{}
\[
P_k^\ii = \TT_{i_1}\cdots\TT_{i_{k-1}} \bigl(q_{i_k}^{1/2}f_{i_k,0}\bigr),
\]
we obtain
the correspondence
\[
X_{[k]}\longmapsto q_{i_k}^{-1/2}P_k^\ii
\qquad(1\le k\le l).
\]

\begin{corollary}
\label{cor: positive braid cluster algebra}
Let $\ttb\in\ttB_\g^+$ and
$\ii=\seq{i_1,\ldots,i_l}\in\Seq(\ttb)$.
The assignment
\[
X_{[k]}\longmapsto q_{i_k}^{-1/2}P_k^\ii
\qquad(1\le k\le l)
\]
extends to an $\bbA$-algebra isomorphism
\[
\Upsilon_\ii:\calA_q(\ttb)\xrightarrow{\sim}\scrA_\bbA(\ttb).
\]
For the empty sequence, $\Upsilon_\emptyset=\operatorname{id}_{\bbA}$.
\end{corollary}

\begin{proof}
The assertion is immediate for $\ttb=1$.
For $\ttb\ne1$, Theorem~\ref{thm: Qin fundamental variables}
gives the stated isomorphism after extension of scalars to
$\Q(q^{1/2})$.

By \cite[Theorem~8.10(1)]{Qin24A}, $\calA_q(\ttb)$ is
generated over $\bbA$ by the $X_{[k]}$.
By Theorem~\ref{thm: PBW basis bosonic braid},
$\scrA_\bbA(\ttb)$ is generated over $\bbA$ by the $P_k^\ii$.
Since $q_{i_k}^{-1/2}\in\bbA^\times$, the isomorphism restricts
to the asserted $\bbA$-algebra isomorphism.
\end{proof}
We identify $\calA_q(\ttb)$ with $\scrA_\bbA(\ttb)$ via
$\Upsilon_\ii$.  Under this identification, the KR realization
constructed above is the restriction of $\Psi_Q$.
Indeed, with the generator correspondence in
\cite[Section~6.1 and footnote~3]{Qin24B}, the bosonic
 map{} used there sends
\[
f_{a,m}\longmapsto\Upphi_Q(f_{a,-m})=\Psi_Q(f_{a,m})
\qquad(a\in I,\ m\in\Z);
\]
see \cite[(7.2) and Theorem~7.2]{JLO2}.
Thus this  map{} agrees with $\Psi_Q$ on every generator.

For the chosen $Q$-adapted word $\jj=\seq{j_1,\ldots,j_N}$,
put $p_r=\xi_{j_r}-2\#\{a<r\mid j_a=j_r\}$.  The initial
variables satisfy
\[
(\Psi_Q\circ\Upsilon_\jj)
\bigl(X^{\calF^\jj}_{[r_{\min},r]}\bigr)
=F_q\bigl(M^{(j_r)}[p_r,\xi_{j_r}]\bigr)
\qquad(1\le r\le N).
\]
By \cite[Lemma~5.10]{Qin24B}, this identity holds for repeated
Coxeter words when $\xi$ takes at most two consecutive values
on each connected component.  Commutations and
restriction to prefixes give it for the chosen adapted word.
For a source $i$ of $Q$, \cite[(8.16)]{JLO1} gives the mutation
from the KR seed for $Q$ to that for $s_iQ$, while
\cite[(6.3)]{Qin24B} gives
\[
\Psi_{s_iQ}=\Psi_Q\circ\TT_i.
\]
These identities preserve the initial-variable correspondence
under source reflections.  Together with even shifts of the
height function on each connected component, they give the
same correspondence for arbitrary $Q$.
After the chosen braid moves from $\jj$ to $\tii$ and
restriction to the prefix $\ii$,
Theorem~\ref{thm: Qin fundamental variables} gives the image
$\Psi_Q(q_{i_k}^{-1/2}P_k^\ii)$ for each singleton variable
$X_{[k]}$ in the KR realization.
Since these variables generate $\calA_q(\ttb)$, this
realization is $\Psi_Q\circ\Upsilon_\ii$.
The same argument applies upstairs.

\begin{remark}
The algebra $\calA_q(\ttb)$ is defined combinatorially for an
arbitrary symmetrizable Cartan datum in
Definition~\ref{def: positive braid cluster algebra}, while
$\scrA_{\bbA}(\ttb)$ is the positive-braid algebra introduced in
\cite{KKOP24B}. In finite type, these algebras are isomorphic, as
recalled in Corollary~\ref{cor: positive braid cluster algebra}.
We expect that
\[
\calA_q(\ttb)\simeq\scrA_{\bbA}(\ttb)
\]
for arbitrary symmetrizable Cartan data.
\end{remark}

\subsection{Surjectivity for arbitrary sequences}
By Corollary~\ref{cor: positive braid cluster algebra}, after
specialization at $q=1$ and extension of scalars to $\C$, we identify
\begin{align} \label{eq: identi}
\calA(\httb)\subset\frakK(\hg)
\qquad\text{and}\qquad
\calA(\ttb)\subset\frakK(\g)
\end{align}
with the corresponding specialized positive-braid algebras through
$\Psi_{\hQ}$ and $\Psi_Q$, respectively.
Under these
identifications, we shall show that $\uppi$ restricts to a surjective
homomorphism
\[
\uppi_\ttb:\calA(\httb)\twoheadrightarrow\calA(\ttb).
\]

{}

Let $\ttb\in\ttB_\g^+$ and let
$\ii=(i_1,\ldots,i_l)\in\Seq(\ttb)$.  Fix a maximal commuting
family $\calF$ of $\bfi$-boxes, choose an admissible chain $\frakC$
realizing $\calF$, and let $\hcalF$ be its canonical lift.

\begin{theorem}
\label{thm: surjectivity arbitrary sequences}
Assume that $\g$ is of finite type.  The ambient homomorphism $\uppi$
restricts to a surjective algebra homomorphism
\[
\uppi_\ttb:
\calA(\httb)\twoheadrightarrow\calA(\ttb).
\]
Moreover, for any $\ii\in\Seq(\ttb)$ and any maximal commuting family
$\calF$ of $\ii$-boxes, under the Laurent realizations associated with
its canonical lift $\hcalF$ and $\calF$, we have
\[
\uppi_\ttb=\phi^\calF_\ttb.
\]
In particular, the folding homomorphisms associated with all choices
of sequences and maximal commuting families are induced by the single
ambient homomorphism $\uppi$.
\end{theorem}

\begin{proof}
 We first compare $\uppi$ with the Laurent folding map on the
adapted KR seeds for $\hjj$ and $\jj$ constructed above.
By the folding formula for KR-polynomials, for each position $r$ of
$\jj$, $\uppi$ sends every variable at a corresponding position of
$\hjj$ to the variable at position $r$ of $\jj$. Hence the two maps agree on the
Laurent rings of these seeds.

Theorem~\ref{thm: global cluster monomial folding} then gives, for
every $1\le k\le l$ and every
$\hfrakc\in\Fib^\frakC(\frakc_k)$,
\[
\uppi(\hx_{\hfrakc}) = \phi^\calF(\hx_{\hfrakc}) = x_{\frakc_k}.
\]
By Corollary~\ref{cor: positive braid cluster algebra} and
specialization at $q=1$, the cluster algebras
$\calA(\httb)$ and $\calA(\ttb)$ are generated by the cluster
variables corresponding to singleton boxes.  It follows that
\[
\uppi\bigl(\calA(\httb)\bigr)=\calA(\ttb).
\]
Since $\uppi$ and $\phi^\calF$ agree on these generators, their
restrictions to $\calA(\httb)$ coincide.  This proves the assertion.
\end{proof}

\begin{corollary}
\label{cor: classical global basis folding}
Let $\ii\in\Seq(\ttb)$, let $\calG$ be any maximal commuting family
of $\ii$-boxes, and let $\hcalG$ be its canonical lift. Then
\[
\uppi_{\ttb}\left(\prod_{\hfrakc\in\hcalG}(\hx_{\hfrakc}^{\hcalG})^{a_{\hfrakc}} \right)
= \prod_{\frakc\in\calG} (x_{\frakc}^{\calG})^{ \sum_{\hfrakc\in\Fib(\frakc)}a_{\hfrakc}} \in \sfL
\]
for every
$(a_{\hfrakc})_{\hfrakc\in\hcalG} \in\mathbb Z_{\geq0}^{\hcalG}$.
In particular, $\uppi_{\ttb}$ sends every cluster monomial associated
with $\Sigma^{\hcalG}$ to an element of the classical global basis
$\sfL$.
\end{corollary}

\begin{proof}
This follows from Theorem~\ref{thm: global cluster monomial folding} and
the fact that cluster monomials of $\frakK(\g)_{\leq\xi}$ belong to the classical global basis $\sfL$.
\end{proof}

\begin{remark}
Theorem~\ref{thm: surjectivity arbitrary sequences} is a commutative statement obtained after
specializing $q^{1/2}=1$.  At the quantum level, we do not claim the
existence of a folding homomorphism between the corresponding
quantum cluster algebras.

  Quantum cluster monomials correspond, up to the
standard $q$-normalization, to global basis elements, and the
folding of their $g$-vectors gives a correspondence between the
associated upstairs and downstairs global basis elements.  More
precisely, if $\hz$ and $z$ are corresponding quantum cluster
monomials, then
\[
\deg_{\ii}(z)
=
\Fold^\ii\bigl(\deg_{\hii}(\hz)\bigr).
\]
Thus the quantum correspondence is detected by $g$-vector folding,
rather than by an algebra homomorphism between quantum tori.
\end{remark}

\subsection{Folding of substitution isomorphisms}
\label{subsec: folding substitution isomorphisms}

Let
\[
\calS_\uppi\seteq \{f\in\calY(\hg)\mid \uppi(f)\ne0\}.
\]
Then $\uppi$ extends uniquely to
\[
\uppi:\calS_\uppi^{-1}\calY(\hg)\longrightarrow\bbF(\g),
\qquad
\frac{f}{g}\longmapsto\frac{\uppi(f)}{\uppi(g)}.
\]
We call the elements of $\calS_\uppi^{-1}\calY(\hg)$
\emph{$\uppi$-specializable}.

\begin{corollary}
\label{cor: substitution folding}
Let $Q,Q'$ be Dynkin quivers of $\g$, and let
$\hQ,\hQ'$ be their   lifted Dynkin quivers.
For every $(i,p)\in\tDynkin^\g_0$ and every $\im\in I_i$, the
rational function
\[
\tPsi^{\rm cl}_{\hQ,\hQ'}(\hy_{\im,p}) \ \ \text{is $\uppi$-specializable}
\qtq
\uppi\left( \tPsi^{\rm cl}_{\hQ,\hQ'}(\hy_{\im,p}) \right) = \tPsi^{\rm cl}_{Q,Q'}(y_{i,p}).
\]
\end{corollary}

\begin{proof}
We first assume that $(i,p)\in\tDynkin_{\leq\xi'}^\g$.
Choose a reduced expression of $w_0$ adapted to $Q'$ and extend
it to an infinite sequence $\ii=(i_1,i_2,\ldots)$ by
$i_{r+\ell(w_0)}=i_r^*$ for $r\ge1$.
This sequence is $Q'$-adapted; see \cite[(5.11)]{LO24}.
Then $(i,p)=(i_k,p_k)$ for some $k$, where
\[
p_k = \xi'_{i_k} -2\#\{1\leq a<k\mid i_a=i_k\}.
\]
 By Lemma~\ref{lem: lifting adapted sequences}, its blockwise lift{}
$\hii$ is $\hQ'$-adapted.

By Theorem~\ref{thm: cluster structure of K xi} and specialization
at $q=1$, the corresponding initial cluster variables of the
$Q'$- and $\hQ'$-adapted seeds are the KR polynomials
\[
F\bigl(m^{(i)}[p,\xi'_i]\bigr) \qtq  F\bigl(\hhm^{(\im)}[p,\hxi'_\im]\bigr),
\]
respectively.
On the other hand, under the ambient Laurent realizations of these
same adapted seeds, the corresponding initial seed variables are
represented by
\[
m^{(i)}[p,\xi'_i] \qtq \hhm^{(\im)}[p,\hxi'_\im],
\]
respectively; see~\cite{FHOO2} and \cite[Proposition~5.7 and
Theorem~5.12]{LO24}.

Under the monomial realizations, $\uppi$ identifies the initial variables
along their fibers, since
\[
\uppi\bigl(\hhm^{(\im)}[p,\hxi'_\im]\bigr)
=m^{(i)}[p,\xi'_i]
\qquad(\im\in I_i).
\]

Choose the analogous infinite $Q$-adapted sequence.
By \cite[Theorem~6.3 and diagram~(7.1)]{FHOO2} and
\cite[Theorem~5.21 and diagram~(5.19)]{LO24}, the substitution
isomorphism between the two monomial realizations is obtained by
repeating a finite braid sequence between their reduced words
for $w_0$.
 For the variable at position $k$, use the finite transformation
between sufficiently long initial segments of length a multiple of $\ell(w_0)$. The finite words obtained in this way represent the same power of  $\Updelta$.
Every mutated position
has a later occurrence of the same color   within the subword to which the braid move applies, so all mutations are at exchangeable vertices.
Theorem~\ref{thm: global cluster monomial folding}, applied to
this finite braid sequence and its  lift, together with
the initial monomial folding above, gives
\[
\uppi\left(\tPsi^{\rm cl}_{\hQ,\hQ'}\bigl(\hhm^{(\im)}[p,\hxi'_\im]\bigr)\right)
= \tPsi^{\rm cl}_{Q,Q'} \bigl(m^{(i)}[p,\xi'_i]\bigr).
\]
The same argument with $p$ replaced by $p+2$ gives
\[
\uppi\left( \tPsi^{\rm cl}_{\hQ,\hQ'} \bigl(\hhm^{(\im)}[p+2,\hxi'_\im]\bigr) \right)
= \tPsi^{\rm cl}_{Q,Q'} \bigl(m^{(i)}[p+2,\xi'_i]\bigr).
\]

Using the usual empty-product convention when $p=\xi'_i$, we have
\[
y_{i,p} = m^{(i)}[p,\xi'_i]\, m^{(i)}[p+2,\xi'_i]^{-1},
\]
and, since $\hxi'_\im=\xi'_i$,
\[
\hy_{\im,p} = \hhm^{(\im)}[p,\hxi'_\im]\, \hhm^{(\im)}[p+2,\hxi'_\im]^{-1}.
\]
The second folding identity shows that the denominator of
$\tPsi^{\rm cl}_{\hQ,\hQ'}(\hy_{\im,p})$ has nonzero image under
$\uppi$.  Hence
$\tPsi^{\rm cl}_{\hQ,\hQ'}(\hy_{\im,p})$ is
$\uppi$-specializable, and
\[
\uppi\left(\tPsi^{\rm cl}_{\hQ,\hQ'}(\hy_{\im,p})\right)
=\frac{\tPsi^{\rm cl}_{Q,Q'}\bigl(m^{(i)}[p,\xi'_i]\bigr)}{\tPsi^{\rm cl}_{Q,Q'}\bigl(m^{(i)}[p+2,\xi'_i]\bigr)}
=\tPsi^{\rm cl}_{Q,Q'}(y_{i,p}).
\]

Finally, let $\widehat{\sfh}$ denote the Coxeter number of $\hg$.
For the finite-type admissible foldings considered here, we have
\[
\widehat{\sfh}=\sfh,
\qquad I_i^*=I_{i^*}.
\] 
Hence
\[
\uppi\bigl(\widehat{\calD}(\hy_{\im,p})\bigr)
= \calD\bigl(\uppi(\hy_{\im,p})\bigr) \qquad (\im\in I_i).
\]
Since the classical substitution isomorphisms commute with the
duality automorphisms,  \cite[Theorem~6.3]{FHOO} and \cite[Proposition~5.16]{LO24} imply that
the identity proved above for
$(i,p)\in\tDynkin_{\leq\xi'}^\g$
extends to every $(i,p)\in\tDynkin_0^\g$.
\end{proof}

In particular, since the substitution isomorphisms preserve the
classical Grothendieck rings, Corollary~\ref{cor: substitution folding}
gives the commutative diagram
\[
\xymatrix@C=55pt@R=32pt{
\frakK(\hg)
  \ar[r]^-{\tPsi^{\rm cl}_{\hQ,\hQ'}}
  \ar@{->>}[d]_{\uppi}
&
\frakK(\hg)
  \ar@{->>}[d]^{\uppi}
\\
\frakK(\g)
  \ar[r]_-{\tPsi^{\rm cl}_{Q,Q'}}
&
\frakK(\g).
}
\]

The settings considered in
Sections~\ref{sec: reduced} and~\ref{sec: finite case} admit
quantum-group-theoretic realizations, and their simply-laced upstairs
cluster algebras admit monoidal categorifications.  We now turn to
arbitrary positive braids for general symmetrizable Cartan data,
where no such representation-theoretic realization is used.

\section{Folding for positive braids over arbitrary symmetrizable Cartan data}
\label{sec: arbitrary positive braid folding}

We prove Theorem~\ref{thm: cluster algebra folding} for an arbitrary
symmetrizable Cartan datum and an arbitrary positive braid.  We first
show that the commutative cluster algebra is freely generated by the
cluster variables attached to singleton $\bfi$-boxes.

Fix a nonempty sequence
\[
\ii=\seq{i_1,\ldots,i_l}\in\Seq(\ttb).
\]
We take $\Sigma^{\calF^\ii}$ as the initial seed of $\calA(\ttb)$,
where $\calF^\ii=\{[k_{\min},k]\mid1\le k\le l\}$.
By Subsection~\ref{subsec: ibox compatible pair}, the exchangeable set
is $J^\ii_{\ex}$ and the exchange matrix is $B^\ii$, given by
\eqref{eq: GLS exchange matrix}.

We recall the definition of interval variables from \cite{Qin24A}
for the commutative cluster algebra $\calA(\ttb)$.
For $p\in[1,l]$, set $p[0]=p$ and, for each integer $d\ge0$ with
$p[d]$ finite, define
\[
p[d+1]\seteq(p[d])_+.
\]
For $1\le k\le l$, set
\[
o_+^\ii(k)\seteq\#\{p\in[1,l]\mid k<p,\ i_p=i_k\}
\]
and define the mutation block
\begin{equation} \label{eq: mutation block}
\bmu^{(k)}\seteq
\begin{cases}
\mu_{k_{\min}[o_+^\ii(k)-1]}\cdots\mu_{k_{\min}[0]},
&o_+^\ii(k)>0,\\
\operatorname{id},&o_+^\ii(k)=0,
\end{cases}    
\end{equation}
where the rightmost mutation is performed first.  If
$o_+^\ii(k)>0$, then $k_{\min}[s]\in J^\ii_{\ex}$ for every integer
$s$ with $0\le s<o_+^\ii(k)$; hence every mutation in $\bmu^{(k)}$ is
performed at an exchangeable vertex.

For a seed $\Sigma$ indexed by $J^\ii$ and $k\in J^\ii$, let
$x_k(\Sigma)$ denote its cluster variable at vertex $k$.
An empty composition of mutation blocks is understood to be the
identity.

\begin{definition}[{\cite[Definition~8.3]{Qin24A}}]
\label{def: Qin interval and fundamental variables}
Let $[a,b]$ be an $\bfi$-box, and let $d\ge0$ be the unique integer
such that $b=a[d]$.  The \emph{interval variable} attached to $[a,b]$
is defined by
\[
W_{[a,b]}\seteq
x_{a_{\min}[d]}\!\left(
\bmu^{(a-1)}\cdots\bmu^{(1)}(\Sigma^{\calF^\ii})
\right).
\]

For $1\le k\le l$, set
\[
W_k\seteq W_{[k,k]}=W_{[k]}.
\]
The variables $W_1,\ldots,W_l$ are called the
\emph{fundamental variables}.
We also set $W_\emptyset\seteq1$.
\end{definition}

Each block $\bmu^{(r)}$ mutates only vertices whose labels $p$
satisfy $i_p=i_r$.
For $1\le k\le l$, one has $i_r\ne i_k$ whenever $r<k_{\min}$, so
none of these blocks mutates the vertex $k$.  Hence,
\begin{equation}
W_{[k_{\min},k]}=x^{\calF^\ii}_{[k_{\min},k]}
\qquad
(1\le k\le l).
\label{eq: Qin initial interval variables}
\end{equation}

The following proposition identifies the cluster variables attached
to $\bfi$-boxes with the corresponding interval variables.

\begin{proposition}
\label{lem: Qin box interval identification}
Let $\calG$ be a maximal commuting family of $\bfi$-boxes associated
with $\ii$.  Then
\[
x^{\calG}_{[a,b]}=W_{[a,b]}
\qquad ([a,b]\in\calG).
\]
Consequently, the cluster variable attached to an $\bfi$-box is
independent of the maximal commuting family containing it.
\end{proposition}

\begin{proof}
The exchangeable set and exchange matrix in
\cite[Section~3.1, equation~(3.1)]{Qin24B}, with sign $+$ attached
to every letter of $\ii$, agree with those of $\Sigma^{\calF^\ii}$.
The mutation blocks $\bmu^{(1)},\ldots,\bmu^{(l)}$ also agree with
those in \cite[Section~3.3]{Qin24B}, and
\eqref{eq: Qin initial interval variables} identifies the variables
of the initial seed. Let $x<y$ be positions of the same color $i$.
By \cite[Proposition~8.6]{Qin24A} at $q^{1/2}=1$, the interval
variables satisfy
\begin{equation}
W_{[x_+,y]}W_{[x,y_-]}
=W_{[x,y]}W_{[x_+,y_-]}+
\prod_{\substack{j\in I\setminus\{i\},\ c_{j,i}<0\\ x(j)^+\le y(j)^-}}
W_{[x(j)^+,y(j)^-]}^{-c_{j,i}}.
\label{eq: Qin interval classical T system}
\end{equation}
For each color $j\ne i$, the interval in Qin's product is
$[x(j)^+,y(j)^-]$, when nonempty, and its exponent is $-c_{j,i}$.
Thus \eqref{eq: Qin interval classical T system} agrees with
\eqref{eq: box Qin classical T system}; the conventions
$W_\emptyset=x^\calG_\emptyset=1$ cover the empty intervals.

Every maximal commuting family is connected to $\calF^\ii$ by box
moves.  We prove by induction on the length of a box-move sequence
starting at $\calF^\ii$ that the seed $\Sigma^\calG$ satisfies
\[
x^\calG_{[a,b]}=W_{[a,b]}
\qquad ([a,b]\in\calG).
\]
The initial case is \eqref{eq: Qin initial interval variables}.
Suppose that the assertion holds for $\calG$, and let
$\calG'=\frakB_k\calG$ be obtained by one further box move.
A permutation-type move preserves the variable attached to each box,
so the assertion holds for $\calG'$.
For a mutation-type move, Proposition~\ref{lem: box classical T system}
and \eqref{eq: Qin interval classical T system} give the same exchange
relation. The induction hypothesis identifies all common variables
and the variable being mutated. Since the mutated variable is nonzero
in the ambient field, the exchange relation uniquely determines the
new variable. Solving the exchange relation for the other factor
gives the reverse mutation-type move.
All other cluster variables are unchanged by the mutation, so the
assertion holds for every box of $\calG'$.  This completes the
induction.
Since $W_{[a,b]}$ is defined from $\Sigma^{\calF^\ii}$ and $[a,b]$
alone, the asserted independence follows.
\end{proof}

We henceforth write $x_{[a,b]}$ for this family-independent variable.
In particular,
\begin{equation}
x_{[k]}=W_k
\qquad(1\le k\le l).
\label{eq: Qin singleton identification}
\end{equation}

The singleton box variables give the following polynomial presentation
of $\calA(\ttb)$.

\begin{proposition}
\label{prop: singleton polynomial presentation}
For every finite sequence $\ii=\seq{i_1,\ldots,i_l}\in\Seq(\ttb)$,
\[
\calA(\ttb)=\C[x_{[1]},\ldots,x_{[l]}].
\]
More precisely, the monomials
\[
x_{[1]}^{a_1}\cdots x_{[l]}^{a_l}
\qquad ((a_1,\ldots,a_l)\in\Z_{\ge0}^l)
\]
form a $\C$-basis of $\calA(\ttb)$.  In particular,
$x_{[1]},\ldots,x_{[l]}$ are algebraically independent.
\end{proposition}

\begin{proof}
For the empty sequence, the algebra is $\C$ and the asserted basis
is $\{1\}$.  Suppose that $l>0$.
In the commutative case, \cite[Theorem~3.6(1)]{Qin24B} states that
the monomials $W_1^{a_1}\cdots W_l^{a_l}$, with
$(a_1,\ldots,a_l)\in\Z_{\ge0}^l$, form a $\Z$-basis of the
$\Z$-algebra generated by the cluster variables obtained from
$\Sigma^{\calF^\ii}$.
After extension of scalars to $\C$, these monomials form a $\C$-basis
of $\calA(\ttb)$.  Equation~\eqref{eq: Qin singleton identification}
gives the asserted basis in the singleton box variables.  This proves
both generation and algebraic independence.
\end{proof}

We now prove Theorem~\ref{thm: cluster algebra folding} using
Proposition~\ref{prop: singleton polynomial presentation}.

Let $\ii=\seq{i_1,\ldots,i_l}$, $\calF$, and $\frakC$ be the sequence,
maximal commuting family, and admissible chain used to define
$\phi^\calF$ in Section~\ref{sec: global folding}.
Recall from \eqref{eq: kappa k,r} that $\kappa(k,r)$ is the position
in $\hii$ of the $r$th entry of the sequence replacing $i_k$.

\begin{proof}[Proof of Theorem~\ref{thm: cluster algebra folding}]
The assertion is immediate for the identity braid, for which both
cluster algebras are $\C$.  Suppose that $\ttb\ne1$.

By
Proposition~\ref{prop: singleton polynomial presentation},
\begin{equation} \label{eq: singleton polynomial rings}
\begin{aligned}
\calA(\ttb) &=\C[x_{[k]}\mid1\le k\le l], \\
\calA(\httb) &=\C[\hx_{[\kappa(k,r)]}\mid 1\le k\le l,\ 1\le r\le e_{i_k}].
\end{aligned}
\end{equation}

We determine the images of the generators $\hx_{[\kappa(k,r)]}$
under $\phi^\calF$.
Fix $k\in[1,l]$ and choose a maximal commuting family $\calG$ containing
$[k]$.
Let $\frakD$ be an admissible chain realizing $\calG$.
By Lemma~\ref{lem: big epsilon box fiber}, the $r$th box in
$\Fib^\frakD([k])$, in the fixed fiber order, has box vector
$\hve_{k,r}$ for $1\le r\le e_{i_k}$.
By \eqref{eq: big epsilon box}, each endpoint of a box indexes
a nonzero coordinate of its box vector.
Since the $r$th box has box vector $\hve_{k,r}$, both endpoints
are $\kappa(k,r)$, and hence this box is $[\kappa(k,r)]$.
Theorem~\ref{thm: global cluster monomial folding} therefore gives
\begin{equation}
\phi^\calF\bigl(\hx_{[\kappa(k,r)]}\bigr)=x_{[k]}
\qquad(1\le k\le l,\ 1\le r\le e_{i_k}).
\label{eq: singleton generator folding}
\end{equation}

The restriction $\phi^\calF_\ttb$ is well defined, since
\eqref{eq: singleton polynomial rings} and
\eqref{eq: singleton generator folding} show that every generator of
$\calA(\httb)$ is mapped into $\calA(\ttb)$.
It is surjective because every generator $x_{[k]}$ of $\calA(\ttb)$
is the image of $\hx_{[\kappa(k,1)]}$.
\end{proof}

\section{Twist automorphisms and folding for positive braids}
\label{sec: arbitrary positive word DT twists}

Kimura--Oya constructed quantum twist automorphisms of localized
quantum unipotent coordinate algebras and proved that they permute
the localized dual canonical bases
\cite[Theorem~6.1]{KimuraOya21}.
We express their twist formula in terms of $\bfi$-box variables
and box vectors, and use this expression to extend the construction
to arbitrary positive braids.
The automorphism is a twist automorphism of
Donaldson--Thomas type in the sense of Kimura--Qin--Wei
\cite{KimuraQinWei}.
It is independent of the sequence representing the braid and
permutes the common triangular basis.
We prove that the commutative specializations of these twists
commute with the folding homomorphisms.

\subsection{Quantum twists and box vectors}
\label{subsec: reduced twists and box vectors}

Fix a positive braid $\ttb\in\ttB_\g^+$ and a sequence
\[
\ii=\seq{i_1,\ldots,i_l}\in\Seq(\ttb).
\]
We use $\Sigma_q^{\calF^\ii}$ as the initial quantum seed of
$\calA_q(\ttb)$.
For a maximal commuting family $\calG$ of $\bfi$-boxes, write
$X^\calG_{\frakd}$ for the cluster variable labelled by
$\frakd\in\calG$ in $\Sigma_q^\calG$.

By Proposition~\ref{lem: Qin box interval identification}, the
specialization of $X^\calG_{[a,b]}$ at $q^{1/2}=1$ is independent
of $\calG$.  Since specialization is injective on quantum cluster
monomials \cite[Lemma~A.4]{FHOO2}, the same holds for
$X^\calG_{[a,b]}$.
We write $X_{[a,b]}$ for this variable and set $X_\emptyset=1$.

For $1\le k\le l$, let $k_{\max}$ be the last occurrence of
$i_k$ in $\ii$.  Set
\begin{equation}
\frakc_k\seteq[k_{\min},k],
\qquad
\frakf_k\seteq[k_{\min},k_{\max}],
\qquad
\frakn_k\seteq
\begin{cases}
[k_+,k_{\max}],&k<k_{\max},\\
\emptyset,&k=k_{\max}.
\end{cases}
\label{eq: twist initial frozen suffix boxes}
\end{equation}
Thus $\frakc_k$ labels the initial variable at $k$,
$\frakf_k$ is the frozen box of color $i_k$, and
$\frakn_k$ contains the occurrences of this color after $k$.
Their box vectors satisfy
\begin{equation}
\calE_{\frakf_k}=\calE_{\frakc_k}+\calE_{\frakn_k}
\qquad(1\le k\le l).
\label{eq: twist box vector decomposition}
\end{equation}

Let $\tcalA_q(\ttb)$ be the $\bbA$-subalgebra of the skew
field of fractions of $\calT(\Lambda^\ii)$ generated by
\[
\calA_q(\ttb)
\quad\text{and}\quad
X_{\frakf_k}^{-1}\qquad(k=k_{\max}).
\]
For the identity braid, set $\tcalA_q(1)=\bbA$.

Suppose now that $\ii\in\Red(w)$ for $w\in\weyl$.
Let $\tsfA_\bbA(\n(w))$ be the localization of
$\sfA_\bbA(\n(w))$ obtained by adjoining the inverses of
$D(w\varpi_i,\varpi_i)$ for the colors occurring in $\ii$.
Let
\[
\eta_{w,q}:\tsfA_\bbA(\n(w))
\xrightarrow{\sim}\tsfA_\bbA(\n(w))
\]
denote the quantum twist automorphism of Kimura--Oya
\cite[Theorems~4.13 and~6.1]{KimuraOya21}.
For $1\le k\le l$, put $\omega=\varpi_{i_k}$.  Then
\begin{equation}
\eta_{w,q}\bigl(D(w_k^\ii\omega,\omega)\bigr)
=q^{(\omega,\omega-w_k^\ii\omega)}
D(w\omega,\omega)^{-1}D(w\omega,w_k^\ii\omega).
\label{eq: reduced quantum twist minor formula}
\end{equation}

For an $\bfi$-box $\frakd=[a,b]$, set
\[
\beta_{\frakd}\seteq
w_{<a}^\ii\varpi_{i_a}-w_b^\ii\varpi_{i_a},
\qquad
D_{\frakd}\seteq
D(w_b^\ii\varpi_{i_a},w_{<a}^\ii\varpi_{i_a}).
\]

\begin{lemma}[{\cite[Theorem~7.3]{GY21} (see also~\cite[\S 8]{KKKO18})}]
\label{lem: reduced quantum box realization}
Let $\ii\in\Red(w)$.
The assignment
\[
X_{\frakd}\longmapsto
q^{-(\beta_{\frakd},\beta_{\frakd})/4}D_{\frakd}
\qquad(\frakd\text{ an }\bfi\text{-box})
\]
extends to an $\bbA$-algebra isomorphism
\[
\calA_q(w)\xrightarrow{\sim}
\sfA_\bbA(\n(w)).
\]
This isomorphism extends to an $\bbA$-algebra isomorphism
\[
\tcalA_q(w)\xrightarrow{\sim}
\tsfA_\bbA(\n(w)).
\]
\end{lemma}

\begin{proof} The first assertion follows from \cite[Theorem~7.3]{GY21}, after
applying the normalization isomorphism
\[
\sfA'_\bbA(\n(w))\xrightarrow{\sim}\sfA_\bbA(\n(w))
\]
of \cite[Remark~8.2.3 and (8.3)]{KKKO18}.
Under this isomorphism, the quantum cluster variable associated with
an $\bfi$-box $\frakd$ is sent to
\[
q^{-(\beta_{\frakd},\beta_{\frakd})/4}D_{\frakd}.
\]
The isomorphism extends to the localizations because each
frozen variable is sent to a unit in $\bbA$ times the corresponding
frozen minor.
\end{proof}

We identify $\tcalA_q(\ttb)$ with
$\tsfA_\bbA(\n(w))$ by this isomorphism and write
\begin{equation}
X_{\frakd}=
q^{-(\beta_{\frakd},\beta_{\frakd})/4}D_{\frakd}.
\label{eq: reduced quantum box normalization}
\end{equation}
Set $\beta_\emptyset=0$ and $D_\emptyset=1$.

The following proposition expresses the formula
\eqref{eq: reduced quantum twist minor formula} in terms of the
quantum cluster variables $X_{\frakd}$ and the
pairing $\lambda^\ii$ on box vectors.

\begin{proposition}
\label{prop: reduced quantum twist box vectors}
For $1\le k\le l$, one has
\begin{equation}
\eta_{w,q}(X_{\frakc_k})
=q^{\frac12\lambda^\ii(\calE_{\frakf_k},\calE_{\frakn_k})}
X_{\frakf_k}^{-1}X_{\frakn_k}.
\label{eq: reduced quantum twist box vectors}
\end{equation}
In particular,
$\eta_{w,q}(X_{\frakf_k})=X_{\frakf_k}^{-1}$.
\end{proposition}

\begin{proof}
Fix $k$ and put $\omega=\varpi_{i_k}$.
For $j\ne i_k$, the simple reflection $s_j$ fixes $\omega$.
Hence
\begin{align*}
 D_{\frakc_k}&=D(w_k^\ii\omega,\omega),&
 \beta_{\frakc_k}&=\omega-w_k^\ii\omega,\allowdisplaybreaks\\
 D_{\frakf_k}&=D(w\omega,\omega),&
 \beta_{\frakf_k}&=\omega-w\omega,\allowdisplaybreaks\\
 D_{\frakn_k}&=D(w\omega,w_k^\ii\omega),&
 \beta_{\frakn_k}&=w_k^\ii\omega-w\omega.
\end{align*}
The last line also holds for $k=k_{\max}$, using
$D_\emptyset=1$ and $\beta_\emptyset=0$.
In particular,
\[
 \beta_{\frakf_k}
 =
 \beta_{\frakc_k}+\beta_{\frakn_k}.
\]
Since $(w_k^\ii\omega,w_k^\ii\omega)=(\omega,\omega)$, we have
\[
 (\beta_{\frakc_k},\beta_{\frakc_k})
 =2(\omega,\omega)-2(\omega,w_k^\ii\omega)
 =2(\omega,\beta_{\frakc_k}).
\]

Substituting \eqref{eq: reduced quantum box normalization}
into \eqref{eq: reduced quantum twist minor formula},
we obtain
\begin{align*}
 \eta_{w,q}(X_{\frakc_k})
 =
 q^{(\omega,\beta_{\frakc_k})
 -\frac14(\beta_{\frakc_k},\beta_{\frakc_k})
 -\frac14(\beta_{\frakf_k},\beta_{\frakf_k})
 +\frac14(\beta_{\frakn_k},\beta_{\frakn_k})}
 X_{\frakf_k}^{-1}X_{\frakn_k}
 =
 q^{-\frac12(\beta_{\frakc_k},\beta_{\frakn_k})}
 X_{\frakf_k}^{-1}X_{\frakn_k}.
\end{align*}

It remains to express this exponent in terms of box vectors.
For an $\bfi$-box $\frakd=[a,b]$, one has
\[
 \beta_{\frakd}  =  \sum_{p=a}^{b}
 \bigl(w_{<p}^\ii\varpi_{i_a}-w_p^\ii\varpi_{i_a}\bigr)
 =  \sum_{\substack{a\le p\le b\\i_p=i_a}}\be_p^\ii.
\]
Every position contributing to $\calE_{\frakc_k}$ precedes
every position contributing to $\calE_{\frakn_k}$.
Thus \eqref{eq: lambda pairing},
\eqref{eq: twist box vector decomposition}, and skew-symmetry give
\begin{align*}
 \lambda^\ii(\calE_{\frakf_k},\calE_{\frakn_k})
 =  \lambda^\ii(\calE_{\frakc_k},\calE_{\frakn_k}) 
 =  -(\beta_{\frakc_k},\beta_{\frakn_k}).
\end{align*}
Substituting this into the formula for
$\eta_{w,q}(X_{\frakc_k})$ proves
\eqref{eq: reduced quantum twist box vectors}.

When $k=k_{\max}$, we have
$\frakc_k=\frakf_k$ and $\frakn_k=\emptyset$.
Since $X_\emptyset=1$ and $\calE_\emptyset=0$, the formula
reduces to
$\eta_{w,q}(X_{\frakf_k})=X_{\frakf_k}^{-1}$.
\end{proof}

\subsection{Quantum twists for arbitrary positive braids}
\label{subsec: terminal box seed arbitrary word}
\label{subsec: arbitrary positive word DT construction}
The $\bfi$-box variables and the pairing on box vectors in
\eqref{eq: reduced quantum twist box vectors} are defined for
arbitrary sequences.  Motivated by this formula, we define a quantum
twist automorphism for an arbitrary positive braid.

Fix a nonidentity positive braid $\ttb\in\ttB_\g^+$ and a sequence
\[
\ii=\seq{i_1,\ldots,i_l}\in\Seq(\ttb).
\]

For a quantum torus $R$, we write $\bbF(R)$
for its skew field of fractions.

\begin{lemma}
\label{lem: positive word twist coordinates}
The assignment
\begin{equation}
X_{\frakc_k}\longmapsto
q^{\frac12\lambda^\ii(\calE_{\frakf_k},\calE_{\frakn_k})}
X_{\frakf_k}^{-1}X_{\frakn_k}
\qquad(1\le k\le l)
\label{eq: quantum twist box coordinates}
\end{equation}
extends uniquely to a $\Q(q^{1/2})$-algebra automorphism of
$\bbF(\calT(\Lambda^\ii))$.
\end{lemma}

\begin{proof}
Set
\[
\calG^\ii\seteq\{[k,k_{\max}]\mid1\le k\le l\}.
\]
This family is associated with the constant $\calL$-sequence.
We compare the quantum seeds
\begin{align*}
\Sigma_q^{\calF^\ii}
=(\bfX^{\calF^\ii},B^\ii,\Lambda^\ii) \qtq 
\Sigma_q^{\calG^\ii}
=(\bfX^{\calG^\ii},B^{\calG^\ii},\Lambda^{\calG^\ii}),
\end{align*}
whose based quantum tori are
$\calT(\Lambda^\ii)$ and $\calT(\Lambda^{\calG^\ii})$,
respectively.
The cluster of the latter seed is
\[
\{X_{\frakn_k}\mid k<k_{\max}\}
\ \cup\
\{X_{\frakf_k}\mid k=k_{\max}\}.
\]

Write $T_k$ for the image prescribed by
\eqref{eq: quantum twist box coordinates}.
Each $T_k$ is a based monomial in
$\calT(\Lambda^{\calG^\ii})$.
By \eqref{eq: twist box vector decomposition},
\[
\calE_{\frakn_k}-\calE_{\frakf_k}
=-\calE_{\frakc_k}.
\]
Hence \eqref{eq: ibox Lambda matrix} gives
\begin{align*}
\lambda^\ii \bigl(\calE_{\frakn_k}-\calE_{\frakf_k}, \calE_{\frakn_p}-\calE_{\frakf_p}\bigr) 
= \lambda^\ii(-\calE_{\frakc_k},-\calE_{\frakc_p}) =\Lambda^\ii_{k,p}.
\end{align*}
Thus $T_kT_p=q^{\Lambda^\ii_{k,p}}T_pT_k$, and the assignment
defines an $\bbA$-algebra homomorphism
\[
\Phi:\calT(\Lambda^\ii)\longrightarrow
\calT(\Lambda^{\calG^\ii}),
\qquad
X_{\frakc_k}\longmapsto T_k.
\]

For each $k$, we have
$\frakf_{k_{\max}}=\frakf_k$ and
$T_{k_{\max}}=X_{\frakf_k}^{-1}$.
Consequently,
\begin{align*}
X_{\frakf_k}
&=T_{k_{\max}}^{-1}
\quad (1\le k\le l) \qtq 
X_{\frakn_k}
=
q^{-\frac12\lambda^\ii(\calE_{\frakf_k},\calE_{\frakn_k})}
T_{k_{\max}}^{-1}T_k
\quad (k<k_{\max}).
\end{align*}
These identities show that the exponent vectors of
$T_1,\ldots,T_l$ form a $\Z$-basis of the exponent lattice of
$\calT(\Lambda^{\calG^\ii})$.
Thus $\Phi$ maps the based monomial basis of $\calT(\Lambda^\ii)$
bijectively onto that of $\calT(\Lambda^{\calG^\ii})$ and is
therefore an $\bbA$-algebra isomorphism.

Since $\Sigma_q^{\calF^\ii}$ and $\Sigma_q^{\calG^\ii}$ are
mutation-equivalent, we have
\[
\bbF(\calT(\Lambda^\ii))
=
\bbF(\calT(\Lambda^{\calG^\ii})).
\]
Hence $\Phi$ extends uniquely to a $\Q(q^{1/2})$-algebra
automorphism of this common skew field.
\end{proof}

\begin{definition}
\label{def: arbitrary positive word twist}
The \emph{quantum twist automorphism} associated with $\ii$ is
the $\Q(q^{1/2})$-algebra automorphism $\tau_{\ii,q}$ of
$\bbF(\calT(\Lambda^\ii))$ determined by
\eqref{eq: quantum twist box coordinates}, namely,
\[
\tau_{\ii,q}(X_{\frakc_k})
=
q^{\frac12\lambda^\ii(\calE_{\frakf_k},\calE_{\frakn_k})}
X_{\frakf_k}^{-1}X_{\frakn_k}
\qquad(1\le k\le l).
\]
\end{definition}

We recall the notion of a twist automorphism of
Donaldson--Thomas type from \cite{KimuraQinWei}, restricting to
integral powers of the frozen variables.

Let $\Sigma_q=(\bfX,B,\Lambda)$ be a quantum seed indexed by $J$,
with exchangeable set $J_\ex$, and put $J_{\fr}=J\setminus J_\ex$.
If $\Xi_q$ is a quantum seed indexed by $J$, we denote by
$X_j(\Xi_q)$ its quantum cluster variable at $j\in J$.
Let $\bmu$ be a sequence of mutations and write
$\Sigma_q'=\bmu(\Sigma_q)=(\bfY,B',\Lambda')$.
Write $(\ve_j)_{j\in J}$ and $(\ve_k)_{k\in J_\ex}$ for the
standard bases of $\Z^J$ and $\Z^{J_\ex}$, respectively.
Let $\sigma$ be a permutation of $J_\ex$ preserving the
skew-symmetrizer.  Assume, as in \cite[Section~2.12]{KimuraQinWei}, 
that there are integers $a_{j,k}$
$(j\in J_{\fr},\ k\in J_\ex)$ such that
\begin{equation}
\deg_{\Sigma_q}(Y_{\sigma(k)})
=-\ve_k+\sum_{j\in J_{\fr}}a_{j,k}\ve_j
\qquad(k\in J_\ex).
\label{eq: KQW twist degree hypothesis}
\end{equation}
For the fixed seeds $\Sigma_q$ and $\Sigma_q'$, condition
\eqref{eq: KQW twist degree hypothesis} determines $\sigma$
uniquely.  Indeed, the restrictions of $\deg_{\Sigma_q}(Y_j)$,
$j\in J_\ex$, to the exchangeable coordinates are the distinct
vectors $-\ve_k$, $k\in J_\ex$, each occurring once, so
$\sigma(k)$ is the unique index for which this restriction equals
$-\ve_k$.

\begin{definition}[{\cite[Definitions~3.3 and~3.15 and Section~5.1]{KimuraQinWei}}]
\label{def: DT type twist integral frozen variables}
A $\Q(q^{1/2})$-algebra automorphism $\psi$ of the ambient skew
field is a \emph{twist automorphism of Donaldson--Thomas type}
passing through $\Sigma_q'$ if
\[
\text{there are vectors
$\bfv_j\in\Z^J$ $(j\in J)$ \emph{supported on $J_{\fr}$} }
\] such that
\[
\psi(X_k)=Y^{\ve_{\sigma(k)}+\bfv_k}\quad(k\in J_\ex),
\qquad
\psi(X_j)=Y^{\bfv_j}\quad(j\in J_{\fr}),
\]
and
\[
\psi\bigl(X^{B\ve_k}\bigr)=Y^{B'\ve_{\sigma(k)}}
\qquad(k\in J_\ex).
\]
\end{definition}

A based monomial in a quantum seed is called a
\emph{localized quantum cluster monomial} if its exponents at all
exchangeable vertices are nonnegative.

\begin{theorem}
\label{thm: arbitrary positive word DT twist}
The automorphism $\tau_{\ii,q}$ is a twist automorphism of
Donaldson--Thomas type and
restricts to an $\bbA$-algebra automorphism
\[
\tau_{\ii,q}:\tcalA_q(\ttb)
\xrightarrow{\sim}\tcalA_q(\ttb).
\]
Both $\tau_{\ii,q}$ and $\tau_{\ii,q}^{-1}$ send every quantum
cluster variable to a localized quantum cluster monomial.
\end{theorem}

\begin{proof} Recall the mutation blocks $\bmu^{(k)}$ defined in
\eqref{eq: mutation block}. We use the same notation for the
corresponding mutation sequences acting on quantum seeds. 
Set
\[
\Sigma_q=\Sigma_q^{\calF^\ii}=(\bfX,B^\ii,\Lambda^\ii),
\qquad
\bmu_\ii=\bmu^{(l)}\cdots\bmu^{(1)},
\qquad
\Sigma_q'=\bmu_\ii(\Sigma_q)=(\bfY,B',\Lambda'),
\]
and write
$J^\ii_{\fr}=J^\ii\setminus J^\ii_\ex
=\{j\in J^\ii\mid j=j_{\max}\}$.
If $J^\ii_\ex=\emptyset$, then $\tau_{\ii,q}$ inverts every
cluster variable, and the assertion is immediate.

We first identify the variables of $\Sigma_q'$.
For each color occurring at $p_0<\cdots<p_{m-1}$ with $m\ge2$,
define
\begin{equation}
\sigma(p_r)=p_{m-2-r}\qquad(0\le r\le m-2).
\label{eq: arbitrary word terminal permutation}
\end{equation}
This permutation preserves colors and hence the skew-symmetrizer.

Fix $k=p_r\in J^\ii_\ex$.
By Definition~\ref{def: Qin interval and fundamental variables},
$W_{[k_+,k_{\max}]}$ occurs at $\sigma(k)=p_{m-2-r}$
 after applying $\bmu^{(k_+-1)}\cdots\bmu^{(1)}$.  
Each subsequent mutation block of the same color is $\bmu^{(p_s)}$
with $s\ge r+1$ and mutates only vertices $p_t$ with
\[
0\le t<m-1-s\le m-2-r.
\]
Blocks of other colors do not mutate $\sigma(k)$ either.
Thus $Y_{\sigma(k)}$ specializes to $W_{[k_+,k_{\max}]}$,
which is also the specialization of $X_{\frakn_k}$ by
Proposition~\ref{lem: Qin box interval identification}.
Since quantum cluster variables are uniquely determined by
their classical specializations, we have
$Y_{\sigma(k)}=X_{\frakn_k}$.
Together with \cite[Proposition~8.7]{Qin24A}, this gives
\begin{equation}
Y_{\sigma(k)}=X_{\frakn_k},
\qquad
\deg_{\Sigma_q}(Y_{\sigma(k)})
=-\ve_k+\ve_{k_{\max}}
\qquad(k\in J^\ii_\ex).
\label{eq: terminal suffix box arbitrary word}
\end{equation}
Frozen variables are unchanged by mutation, so
\begin{equation}
Y_j=X_{\frakf_j},
\qquad
\deg_{\Sigma_q}(Y_j)=\ve_j
\qquad(j\in J^\ii_{\fr}).
\label{eq: terminal frozen box arbitrary word}
\end{equation}
In particular, \eqref{eq: KQW twist degree hypothesis} holds
with $a_{j,k}=\delta_{j,k_{\max}}$.

We now express $\tau_{\ii,q}$ in the variables of $\Sigma_q'$.
 By \eqref{eq: terminal suffix box arbitrary word},
\eqref{eq: terminal frozen box arbitrary word}, and
\eqref{eq: ibox Lambda matrix}, we have{}
\[
\Lambda'_{k_{\max},\sigma(k)}
=
\lambda^\ii(\calE_{\frakf_k},\calE_{\frakn_k})
\qquad(k\in J^\ii_\ex).
\]
Hence \eqref{eq: quantum twist box coordinates}, with the
normalization of based monomials, becomes
\begin{equation}
\tau_{\ii,q}(X_{\frakc_k})=
\begin{cases}
Y^{\ve_{\sigma(k)}-\ve_{k_{\max}}},&k\in J^\ii_\ex,\\
Y^{-\ve_k},&k\in J^\ii_{\fr}.
\end{cases}
\label{eq: terminal monomial coordinates arbitrary word}
\end{equation}

It remains to verify 
\begin{align} \label{eq: exchange-monomial identity}
\tau_{\ii,q}\bigl(X^{B\ve_k}\bigr)
=
Y^{B'\ve_{\sigma(k)}}
\qquad(k\in J^\ii_{\ex}).    
\end{align}

Let $G$ be the matrix whose $j$-th column is
$\deg_{\Sigma_q}(Y_j)$.\footnote{This is the $G$-matrix of
$\Sigma_q'$ with respect to $\Sigma_q$, including the frozen coordinates.}
The degree formulas above show that $G$ is unimodular and that
\[
G(\ve_{\sigma(k)}-\ve_{k_{\max}})=-\ve_k
\quad(k\in J^\ii_\ex),
\qquad
G(-\ve_j)=-\ve_j
\quad(j\in J^\ii_{\fr}).
\]
Consequently,
\[
\tau_{\ii,q}(X^{\bfa})=Y^{-G^{-1}\bfa}
\qquad(\bfa\in\Z^{J^\ii}).
\]

Since the exchangeable parts of the columns of $G$ indexed by 
$\sigma(k)$ are $-\ve_k$ and $\sigma$ preserves the
skew-symmetrizer,
\cite[equation~(2.4) and Lemma~2.15(1)]{KimuraQinWei} give
\[
GB'\ve_{\sigma(k)}=-B^\ii\ve_k
\qquad(k\in J^\ii_{\ex}).
\]
Therefore
\[
\tau_{\ii,q}(X^{B^\ii\ve_k})
=
Y^{-G^{-1}B^\ii\ve_k}
=
Y^{B'\ve_{\sigma(k)}},
\]
which verifies the remaining condition in
Definition~\ref{def: DT type twist integral frozen variables}.

Finally, \eqref{eq: terminal monomial coordinates arbitrary word}
gives
\begin{align} \label{eq: monomial correspondences}
\tau_{\ii,q}^{-1}(Y_{\sigma(k)})
=X^{\ve_k-\ve_{k_{\max}}}
\quad(k\in J^\ii_\ex),
\qquad
\tau_{\ii,q}^{-1}(Y_j)=X^{-\ve_j}
\quad(j\in J^\ii_{\fr}).    
\end{align}
Applying $\tau_{\ii,q}^{-1}$ to~\eqref{eq: exchange-monomial identity}
shows that $\tau_{\ii,q}^{-1}$ satisfies the same condition
with the two seeds interchanged. 
By \cite[Propositions~3.12 and~3.14]{KimuraQinWei}, the identities in~\eqref{eq: monomial correspondences},  with integral frozen exponents,
are preserved under mutations at corresponding vertices.
Iterating along a mutation sequence therefore shows that both
$\tau_{\ii,q}$ and $\tau_{\ii,q}^{-1}$ send every exchangeable
quantum cluster variable to a cluster variable multiplied by
a frozen Laurent monomial, with the based normalization.
These are localized quantum cluster monomials in seeds
mutation-equivalent to $\Sigma_q$.

Thus both maps send every quantum cluster variable into
$\tcalA_q(\ttb)$.
They also invert the frozen variables and fix $\bbA$.
Since $\tcalA_q(\ttb)$ is generated over $\bbA$ by the quantum
cluster variables and the inverses of the frozen variables,
both maps preserve $\tcalA_q(\ttb)$.
Their restrictions are mutually inverse, proving the assertion.
\end{proof}

\begin{lemma}
\label{lem: normalized DT twist every seed}
Let $\ii\in\Seq(\ttb)$.
For every quantum seed $\Sigma_q$ mutation-equivalent to
$\Sigma_q^{\calF^\ii}$, we have
\begin{equation}
\deg_{\Sigma_q}
\tau_{\ii,q}\bigl(X_j(\Sigma_q)\bigr)=-\ve_j
\qquad(j\in J^\ii).
\label{eq: normalized DT twist every seed}
\end{equation}
\end{lemma}

\begin{proof}
By Theorem~\ref{thm: arbitrary positive word DT twist}, the
images under $\tau_{\ii,q}$ of quantum cluster variables are
localized quantum cluster monomials.
Quantum cluster monomials are pointed in every seed, as recalled
in Section~\ref{sec: compatible pair folding}.
Multiplication by a frozen Laurent monomial shifts the degree in
\eqref{eq: pointed element} by its exponent vector.
Hence localized quantum cluster monomials are also pointed in
every seed, and their degrees satisfy \eqref{eq: g-vector mutation}.

We prove \eqref{eq: normalized DT twist every seed} by induction
on the length of a mutation sequence from $\Sigma_q^{\calF^\ii}$.
For the initial seed, \eqref{eq: terminal suffix box arbitrary word}
and \eqref{eq: terminal monomial coordinates arbitrary word} give
\[
\deg_{\Sigma_q^{\calF^\ii}}
\tau_{\ii,q}(X_{\frakc_k})
=
(-\ve_k+\ve_{k_{\max}})-\ve_{k_{\max}}
=
-\ve_k
\qquad(k\in J^\ii_\ex).
\]
The assertion also holds at frozen vertices, since
$\tau_{\ii,q}$ inverts the frozen variables.

Suppose that \eqref{eq: normalized DT twist every seed} holds
for $\Sigma_q=(\bfX,B,\Lambda)$.
Let $k\in J^\ii_\ex$ and
$\Sigma_q'=\mu_k(\Sigma_q)=(\bfX',B',\Lambda')$.
Set
\[
\bfb^\pm
=\sum_{j\in J^\ii}[\pm b_{j,k}]_+\ve_j.
\]
The quantum exchange relation has the form
\[
X_k(\Sigma_q)X_k(\Sigma_q')
=
q^{r_+}X^{\bfb^+}
+
q^{r_-}X^{\bfb^-},
\qquad
r_\pm\in\tfrac12\Z,
\]
where the based monomials on the right are taken in $\Sigma_q$.
Applying $\tau_{\ii,q}$ gives
\[
\begin{aligned}
\tau_{\ii,q}\bigl(X_k(\Sigma_q)\bigr)  \tau_{\ii,q}\bigl(X_k(\Sigma_q')\bigr) 
= q^{r_+}\tau_{\ii,q}(X^{\bfb^+}) + q^{r_-}\tau_{\ii,q}(X^{\bfb^-}).
\end{aligned}
\]

By the induction hypothesis and additivity of degrees under
multiplication, the two summands on the right are pointed at
$-\bfb^+$ and $-\bfb^-$, respectively.
Since
\[
-\bfb^-=-\bfb^++B\ve_k,
\]
their expansions in \eqref{eq: pointed element} take the form
\begin{align*}
q^{r_+}\tau_{\ii,q}(X^{\bfb^+})
&=
q^{s_+}X^{-\bfb^+}
+
\sum_{\bfn\in\N^{J^\ii_\ex}\setminus\{0\}}
c_{\bfn}X^{-\bfb^++B\bfn},\\
q^{r_-}\tau_{\ii,q}(X^{\bfb^-})
&=
\sum_{\bfn\in\N^{J^\ii_\ex}}
d_{\bfn}X^{-\bfb^++B(\ve_k+\bfn)}.
\end{align*}
Compatibility implies that $B$ has full column rank.
Thus $B(\ve_k+\bfn)\ne0$ for every
$\bfn\in\N^{J^\ii_\ex}$, so the second summand contains
no term with exponent $-\bfb^+$.
The leading term $q^{s_+}X^{-\bfb^+}$ therefore cannot
be cancelled, and the sum is pointed at $-\bfb^+$.

Both factors on the left are pointed in $\Sigma_q$, and
the induction hypothesis gives
\[
\deg_{\Sigma_q}
\tau_{\ii,q}\bigl(X_k(\Sigma_q)\bigr)=-\ve_k.
\]
Comparing degrees on both sides, we obtain
\[
-\ve_k+
\deg_{\Sigma_q}
\tau_{\ii,q}\bigl(X_k(\Sigma_q')\bigr)
=
-\bfb^+.
\]
Write
\[
\begin{aligned}
\bfg=(g_j)_{j\in J^\ii}
&\seteq
\deg_{\Sigma_q}
\tau_{\ii,q}\bigl(X_k(\Sigma_q')\bigr)
=
\ve_k-\bfb^+,\\
\bfg'=(g'_j)_{j\in J^\ii}
&\seteq
\deg_{\Sigma_q'}
\tau_{\ii,q}\bigl(X_k(\Sigma_q')\bigr).
\end{aligned}
\]
Since $b_{k,k}=0$, the coordinates of $\bfg$ are
\[
g_k=1,
\qquad
g_j=-[b_{j,k}]_+
\quad(j\ne k).
\]
The $k$-th equation in \eqref{eq: g-vector mutation} gives
\[
g'_k=-g_k=-1.
\]
For $j\ne k$, the case $g'_k\le0$ of the same formula gives
\[
g_j
=
g'_j+[-b'_{j,k}]_+g'_k
=
g'_j-[b_{j,k}]_+,
\]
where we used $b'_{j,k}=-b_{j,k}$.
Substituting $g_j=-[b_{j,k}]_+$, we obtain
\[
g'_j = -[b_{j,k}]_++[b_{j,k}]_+ = 0 \qquad(j\ne k).
\]
Thus $\bfg'=-\ve_k$, so
\[
\deg_{\Sigma_q'} \tau_{\ii,q}\bigl(X_k(\Sigma_q')\bigr)
= -\ve_k.
\]

For $j\ne k$, we have $X_j(\Sigma_q')=X_j(\Sigma_q)$.
Hence the induction hypothesis gives
\[
\begin{aligned}
\deg_{\Sigma_q} \tau_{\ii,q}\bigl(X_j(\Sigma_q')\bigr)
= \deg_{\Sigma_q} \tau_{\ii,q}\bigl(X_j(\Sigma_q)\bigr)
=-\ve_j.
\end{aligned}
\]
Since the $k$-th coordinate of $-\ve_j$ is zero,
\eqref{eq: g-vector mutation} gives
\[
\deg_{\Sigma_q'}
\tau_{\ii,q}\bigl(X_j(\Sigma_q')\bigr)
=
-\ve_j
\qquad(j\ne k).
\]
Together with the calculation for $j=k$, this proves
\eqref{eq: normalized DT twist every seed} for $\Sigma_q'$,
completing the induction.
\end{proof}

\begin{proposition}
\label{prop: normalized DT twist independence}
The automorphism $\tau_{\ii,q}$ is independent of the choice of
$\ii\in\Seq(\ttb)$.
\end{proposition}

\begin{proof}
Let $\ii,\ii'\in\Seq(\ttb)$ and put
$\Sigma_q=\Sigma_q^{\calF^{\ii'}}$.
By Proposition~\ref{prop: mutation permutation braid moves},
$\Sigma_q$ is obtained from $\Sigma_q^{\calF^\ii}$ by mutations
and permutations.
Relabelling vertices also relabels the coordinates of degree
vectors, so Lemma~\ref{lem: normalized DT twist every seed}
applies and gives
\[
\begin{aligned}
\deg_{\Sigma_q}
\tau_{\ii,q}\bigl(X_j(\Sigma_q)\bigr)
&=-\ve_j,\\
\deg_{\Sigma_q}
\tau_{\ii',q}\bigl(X_j(\Sigma_q)\bigr)
&=-\ve_j
\end{aligned}
\qquad(j\in J^{\ii'}).
\]
Both images are localized quantum cluster monomials by
Theorem~\ref{thm: arbitrary positive word DT twist}.
Multiplying both images by a common frozen monomial and
normalizing by the same power of $q^{1/2}$ gives quantum cluster
monomials of the same degree.
By \cite[Theorem~A.8]{FHOO2}, these monomials are equal.
Cancelling the common factors, we obtain
\[
\tau_{\ii,q}\bigl(X_j(\Sigma_q)\bigr)
=
\tau_{\ii',q}\bigl(X_j(\Sigma_q)\bigr)
\qquad(j\in J^{\ii'}).
\]
The two automorphisms agree on all variables of $\Sigma_q$
and on $\Q(q^{1/2})$, hence on their common skew field of
fractions.
\end{proof}
We write $\tau_{\ttb,q}\seteq\tau_{\ii,q}$ and set
$\tau_{1,q}=\operatorname{id}_{\bbA}$.

If $\ttb$ is the positive braid lift of $w\in\weyl$, we identify
$\tcalA_q(\ttb)$ with $\tsfA_\bbA(\n(w))$ by
Lemma~\ref{lem: reduced quantum box realization}.

\begin{corollary}
\label{cor: reduced quantum twist recovered}
Let $w\in\weyl$, and let $\ttb$ be its positive braid lift.  Then
\[
\tau_{\ttb,q}=\eta_{w,q}.
\]
\end{corollary}

\begin{proof}
Choose $\ii\in\Red(w)$.
Proposition~\ref{prop: reduced quantum twist box vectors} and
\eqref{eq: quantum twist box coordinates} show that
$\tau_{\ttb,q}$ and $\eta_{w,q}$ agree on all initial quantum
cluster variables.  These generate the ambient skew field over
$\Q(q^{1/2})$, so the two automorphisms agree.
\end{proof}

By \cite[Theorems~6.13(1) and~6.14]{Qin24A}, for $\ttb\ne1$
the algebra $\tcalA_q(\ttb)$ admits a common triangular
$\bbA$-basis in the sense of
\cite[Definitions~6.2.1--6.2.2]{Qin20}.
This basis is unique by \cite[Lemma~6.2.4]{Qin20}.
We denote it by $\bfL_\ttb$ and set $\bfL_1=\{1\}$.

For $\ii\in\Red(w)$, the isomorphism in
Lemma~\ref{lem: reduced quantum box realization} identifies
$\bfL_\ttb$ with the localized dual canonical basis:
the products of an element of $\bfB(w)$ with inverses of the
frozen minors, normalized by powers of $q^{1/2}$ so that their
leading coefficients in $\Sigma_q^{\calF^\ii}$ are $1$; see
\cite[Theorems~9.3.2 and~9.5.1]{Qin20}.
The following corollary extends the basis-preservation property
of the Kimura--Oya quantum twist $\eta_{w,q}$ \cite[Theorem~6.1]{KimuraOya21}
to arbitrary positive braids.

\begin{corollary}
\label{prop: positive braid twist common triangular basis}
The quantum twist automorphism $\tau_{\ttb,q}$ permutes
$\bfL_\ttb$:
\[
\tau_{\ttb,q}(\bfL_\ttb)=\bfL_\ttb.
\]
\end{corollary}

\begin{proof}
The assertion is immediate for $\ttb=1$.
For $\ttb\ne1$, it follows from
Theorem~\ref{thm: arbitrary positive word DT twist},
\cite[Theorem~6.14]{Qin24A}, and
\cite[Section~6.1 and Theorem~6.2(2)]{KimuraQinWei}.
\end{proof}

In finite Cartan type, the preceding corollary also gives preservation
of the normalized global basis of Kashiwara--Kim--Oh--Park
\cite{KKOP24} after localization.
We first identify this basis with the common triangular basis.

Assume that $\sfC$ is of finite type, and fix
$\ii=\seq{i_1,\ldots,i_l}\in\Seq(\ttb)$.
Recall from Section~\ref{sec: finite case} that the cuspidal elements are
\[
P_k^\ii=\TT_{i_1}\cdots\TT_{i_{k-1}}
\bigl(q_{i_k}^{1/2}f_{i_k,0}\bigr)
\qquad(1\le k\le l),
\]
where $q_i=q^{d_i}$, the $f_{i,p}$ are the generators of the
bosonic extension $\scrA_q$, and the $\TT_i$ are its braid
symmetries.
Recall also the normalized global basis $\tbfG(\ttb)$ of
$\scrA_\bbA(\ttb)$ from
Theorem~\ref{thm: PBW basis bosonic braid}.

Corollary~\ref{cor: positive braid cluster algebra} gives the
$\bbA$-algebra isomorphism
\[
\Upsilon_\ii:\calA_q(\ttb)\xrightarrow{\sim}\scrA_\bbA(\ttb).
\]
Set
\[
\widetilde P_k^\ii\seteq\Upsilon_\ii(X_{[k]})
=q_{i_k}^{-1/2}P_k^\ii
\qquad(1\le k\le l).
\]

\begin{proposition}
\label{prop: KKOP common triangular basis identification}
The isomorphism $\Upsilon_\ii$ identifies the two bases:
\[
\Upsilon_\ii\bigl(\bfL_\ttb\cap\calA_q(\ttb)\bigr)
=\tbfG(\ttb).
\]
\end{proposition}

\begin{proof}
The assertion is immediate for $\ttb=1$. Suppose that $\ttb\ne1$.
The variables $X_{[k]}$ are bar-invariant, and their images
$\widetilde P_k^\ii$ are bar-invariant by
\cite[(2.4) and Lemma~4.17]{KKOP24B}. Since these elements generate the two algebras by
Corollary~\ref{cor: positive braid cluster algebra} and
Theorem~\ref{thm: PBW basis bosonic braid}, the isomorphism
$\Upsilon_\ii$ commutes with the bar anti-involutions.

For $\bfu\in\Z_{\ge0}^l$, let $M(\bfu)$ and
$\widetilde P^\ii(\bfu)$ denote the standard monomial of
\cite[Section~8.3]{Qin24A} and the normalized PBW element of
\cite[(2.4) and (4.7)]{KKOP24B}, respectively. The definitions,
\eqref{eq: GLS Lambda matrix}, \eqref{eq: lambda pairing}, and
\cite[Lemma~8.4]{Qin24A} give
\[
\Upsilon_\ii(M(\bfu))=\overline{\widetilde P^\ii(\bfu)}.
\]
The bar-invariance and PBW triangularity in
\cite[Definition~4.14 and Lemmas~4.16--4.17]{KKOP24B}, after applying
the bar anti-involution, give the standard monomial characterization in
\cite[Definition~2.20 and Theorem~8.10(4)]{Qin24A}. The uniqueness in
that characterization applies because the bi-lexicographic order
refines the reverse lexicographic order and normalization preserves
the coefficients of terms of the same weight. Hence
$\Upsilon_\ii^{-1}(\tbfG(\ttb))=\bfL_\ttb\cap\calA_q(\ttb)$.
\end{proof}

We extend the basis identification to the frozen localizations.
For $j\in J^\ii_{\fr}$, put
$H_j=\Upsilon_\ii(X_{\frakf_j})$, where
$\frakf_j=[j_{\min},j_{\max}]$ is the frozen $\bfi$-box of
color $i_j$.
Let $\widetilde{\scrA}_\bbA(\ttb)$ be the localization of
$\scrA_\bbA(\ttb)$ at the elements $H_j$.
Then $\Upsilon_\ii$ extends to an isomorphism
\[
\Upsilon_\ii:\tcalA_q(\ttb)
\xrightarrow{\sim}\widetilde{\scrA}_\bbA(\ttb).
\]
Let $\tbfG^{\loc}(\ttb)$ be the set obtained from the products
\[
\left(\prod_{j\in J^\ii_{\fr}}H_j^{-a_j}\right)\bfg,
\qquad a_j\in\Z_{\ge0},\quad \bfg\in\tbfG(\ttb),
\]
by normalizing their leading coefficients in the image of
$\Sigma_q^{\calF^\ii}$ to $1$.
The product is taken in increasing order of $j$.
Normalized multiplication by frozen variables and their inverses
preserves the common triangular basis
\cite[Definition~2.20]{Qin24A}.
Clearing the frozen denominators and applying the proposition gives
\[
\Upsilon_\ii(\bfL_\ttb)=\tbfG^{\loc}(\ttb).
\]
In particular, $\tbfG^{\loc}(\ttb)$ is an $\bbA$-basis of
$\widetilde{\scrA}_\bbA(\ttb)$.
Write
\[
\tau^{\mathrm{bos}}_{\ii,q}
\seteq\Upsilon_\ii\circ\tau_{\ttb,q}\circ\Upsilon_\ii^{-1}.
\]

\begin{corollary}
\label{cor: positive braid twist KKOP global basis}
If $\sfC$ is of finite type, then
\[
\tau^{\mathrm{bos}}_{\ii,q}
\bigl(\tbfG^{\loc}(\ttb)\bigr)=\tbfG^{\loc}(\ttb).
\]
\end{corollary}

\begin{proof}
This follows from
Corollary~\ref{prop: positive braid twist common triangular basis} and Proposition~\ref{prop: KKOP common triangular basis identification}.
\end{proof}

\subsection{Compatibility with folding}
\label{subsec: arbitrary word DT folding}

We prove that the commutative specializations of the quantum
twists commute with the folding homomorphisms.

Let $S_\ttb$ be the multiplicative subset of $\calA(\ttb)$
generated by its frozen variables, and set
\[
\tcalA(\ttb)\seteq S_\ttb^{-1}\calA(\ttb).
\]
Evaluation at $q^{1/2}=1$ induces a surjective homomorphism
\[
\C\otimes_{\bbA}\tcalA_q(\ttb)
\twoheadrightarrow\tcalA(\ttb),
\]
since quantum cluster variables and inverses of frozen variables
specialize to the corresponding generators of
$\tcalA(\ttb)$.
The elements of $\bfL_\ttb$ are pointed at distinct degrees
with leading coefficient $1$ in $\Sigma_q^{\calF^\ii}$.
Their specializations at $q^{1/2}=1$ are therefore linearly
independent.  Since $\bfL_\ttb$ is an $\bbA$-basis of
$\tcalA_q(\ttb)$, this homomorphism is also injective and hence
an isomorphism.
We denote the specialization of $\tau_{\ttb,q}$ by
\[
\tau_\ttb:\tcalA(\ttb)\xrightarrow{\sim}\tcalA(\ttb).
\]
For the identity braid, $S_1=\{1\}$, $\tcalA(1)=\C$, and
$\tau_1=\operatorname{id}_\C$.
By \eqref{eq: quantum twist box coordinates}, for every
$\ii\in\Seq(\ttb)$ one has
\begin{equation}
\tau_\ttb(x_{\frakc_k})
=x_{\frakf_k}^{-1}x_{\frakn_k}
\qquad(1\le k\le l),
\label{eq: arbitrary positive word DT box formula}
\end{equation}
where $x_\emptyset=1$.
In particular, $\tau_\ttb(x_{\frakf_k})=x_{\frakf_k}^{-1}$.
By Corollary~\ref{prop: positive braid twist common triangular basis},
$\tau_\ttb$ permutes the basis $\bfL_\ttb^{\mathrm{cl}}$
obtained by specializing $\bfL_\ttb$ at $q^{1/2}=1$.

Fix a maximal commuting family $\calF$ of $\bfi$-boxes and an
admissible chain $\frakC$ realizing $\calF$, and let
$\phi^\calF$ be the Laurent homomorphism defined in
Section~\ref{sec: global folding}.
Let $\hii\in\Seq(\httb)$  be the blockwise lift of $\ii$.{}
Recall from \eqref{eq: kappa k,r} that $\kappa(k,a)$ is the
position of the $a$th lifted letter over $i_k$ in $\hii$.
Theorem~\ref{thm: cluster algebra folding} gives the surjective homomorphism
\[
\phi^\calF_\ttb:\calA(\httb)\twoheadrightarrow\calA(\ttb).
\]

For $1\le k\le l$ and $1\le a\le e_{i_k}$, set
\[
\hfrakc_{k,a}\seteq
[\kappa(k_{\min},a),\kappa(k,a)],
\qquad
\hfrakf_{k,a}\seteq
[\kappa(k_{\min},a),\kappa(k_{\max},a)],
\]
and set
\[
\hfrakn_{k,a}\seteq
\begin{cases}
[\kappa(k_+,a),\kappa(k_{\max},a)],&k<k_{\max},\\
\emptyset,&k=k_{\max}.
\end{cases}
\]
Set $\hx_{\emptyset}=1$.
Lemma~\ref{lem: big epsilon box fiber} identifies the $a$th box
in a fiber over $[u,v]$ as $[\kappa(u,a),\kappa(v,a)]$.
Applying Theorem~\ref{thm: global cluster monomial folding} and
Proposition~\ref{lem: Qin box interval identification},
we obtain
\begin{equation}
\phi^\calF_\ttb(\hx_{\hfrakc_{k,a}})=x_{\frakc_k},
\qquad
\phi^\calF_\ttb(\hx_{\hfrakf_{k,a}})=x_{\frakf_k},
\qquad
\phi^\calF_\ttb(\hx_{\hfrakn_{k,a}})=x_{\frakn_k}.
\label{eq: aligned box variable folding arbitrary word DT}
\end{equation}

By \eqref{eq: aligned box variable folding arbitrary word DT},
$\phi^\calF_\ttb(S_{\httb})=S_\ttb$.
Thus $\phi^\calF_\ttb$ extends to a surjective homomorphism
\begin{equation}
\phi^{\calF,\loc}_\ttb:
\tcalA(\httb)
\twoheadrightarrow
\tcalA(\ttb).
\label{eq: localized arbitrary word folding for DT}
\end{equation}

\begin{theorem}
\label{thm: arbitrary positive word DT folding compatibility}
The twist automorphisms satisfy
\begin{equation}
\phi^{\calF,\loc}_\ttb\circ\tau_{\httb}
=\tau_\ttb\circ\phi^{\calF,\loc}_\ttb.
\label{eq: arbitrary positive word DT folding compatibility}
\end{equation}
\end{theorem}

\begin{proof}
For $\ttb=1$, all maps are the identity of $\C$, so the assertion
holds.  Assume that $\ttb\ne1$.

The two homomorphisms agree on
$\hx_{\hfrakc_{k,a}}$ for $1\le k\le l$ and $1\le a\le e_{i_k}$.
The occurrences of the
lifted color at $\kappa(k,a)$ are precisely the positions
$\kappa(p,a)$ with $i_p=i_k$.  Hence the three boxes in
\eqref{eq: arbitrary positive word DT box formula} for this lifted
position are $\hfrakc_{k,a}$, $\hfrakf_{k,a}$, and
$\hfrakn_{k,a}$.  By \eqref{eq: arbitrary positive word DT box formula}
and \eqref{eq: aligned box variable folding arbitrary word DT},
\[
\begin{aligned}
\phi^{\calF,\loc}_\ttb\tau_{\httb}
 (\hx_{\hfrakc_{k,a}})
=\phi^{\calF,\loc}_\ttb
 (\hx_{\hfrakf_{k,a}}^{-1}\hx_{\hfrakn_{k,a}})
=x_{\frakf_k}^{-1}x_{\frakn_k}
=\tau_\ttb(x_{\frakc_k})
=\tau_\ttb\phi^{\calF,\loc}_\ttb
 (\hx_{\hfrakc_{k,a}}).
\end{aligned}
\]

The variables $\hx_{\hfrakc_{k,a}}$, over all $k$ and $a$, form the
cluster of $\Sigma^{\calF^\hii}$, and the two homomorphisms have
the same nonzero values on them.
By the Laurent phenomenon, every element of
$\tcalA(\httb)$ is a Laurent polynomial in these
variables.  Multiplying by a monomial in the initial variables
clears the denominators.  The two homomorphisms agree on the polynomial, and the common image of the monomial is
nonzero.  Cancelling this image shows that the two homomorphisms
agree on
$\tcalA(\httb)$.
\end{proof}

For $w\in\weyl$, the specialization of $\eta_{w,q}$ at
$q^{1/2}=1$ is an automorphism of
$\C\otimes_\bbA\tsfA_\bbA(\n(w))$ by
\cite[Corollary~6.2]{KimuraOya21}.
If $\ttb$ is the positive braid lift of $w$, we identify
$\C\otimes_\bbA\tsfA_\bbA(\n(w))$ with
$\tcalA(\ttb)$ by specializing the isomorphism in
Lemma~\ref{lem: reduced quantum box realization}.

\begin{corollary}
\label{cor: reduced twist recovered by arbitrary word DT}
Let $w\in\weyl$, set $\hw=\iota(w)$, and let $\ttb$ be the
positive braid lift of $w$.
Then
\[
\phi^{\calF,\loc}_\ttb\circ
\left.\eta_{\hw,q}\right|_{q=1}
=\left.\eta_{w,q}\right|_{q=1}
\circ\phi^{\calF,\loc}_\ttb.
\]
\end{corollary}

\begin{proof}
Since $\ii\in\Seq(\ttb)$  is reduced, its blockwise lift $\hii${}
belongs to $\Red(\hw)$ by Corollary~\ref{cor: reduced lift}.
Specializing Corollary~\ref{cor: reduced quantum twist recovered}
for $w$ and $\hw$ gives
\[
\tau_\ttb=\left.\eta_{w,q}\right|_{q=1},
\qquad
\tau_{\httb}=\left.\eta_{\hw,q}\right|_{q=1}.
\]
The assertion now follows from
Theorem~\ref{thm: arbitrary positive word DT folding compatibility}.
\end{proof}

\appendix

\section{\texorpdfstring{  Additional calculations for affine braid moves}{Additional calculations for affine braid moves}} \label{app: affine braid folding}

 The braid-move calculations in Section~\ref{sec: braid} treat
the subdiagrams on the two color orbits arising in finite-type
foldings. Here we give the affine cases with orbit sizes
$(3,3)$ and $(4,2)$.{}   The exchange columns are computed as in Section~\ref{sec: braid}.
We describe the lifted paths and the corresponding mutation--permutation
sequences, and prove \eqref{eq: folding of cluster var}.
Throughout, the homomorphism
$\phi^\ii$ is the one defined on  $\calL(\Sigma^\hii)$.

\subsection{\texorpdfstring{  The $3$-move with orbit sizes $(3,3)$}{The 3-move with orbit sizes (3,3)}}
\label{subsec: affine braid 33}

Use the Sage numbering of $E_6^{(1)}$, whose three arms are
\[
 0-2-4,\qquad 1-3-4,\qquad 6-5-4.
\]
The automorphism
\[
 \sigma=(0\ 1\ 6)(2\ 3\ 5),\qquad \sigma(4)=4,
\]
is admissible, with ordered orbits
\[
 I_i=(a_1,a_2,a_3)=(0,1,6),\qquad
 I_j=(b_1,b_2,b_3)=(2,3,5).
\]
The quotient matrix is
\[
 \begin{pmatrix}
 2&-1&0\\ -1&2&-3\\ 0&-1&2
 \end{pmatrix},
\]
the Cartan matrix of $D_4^{(3)}$.  In particular,
$c_{i,j}=c_{j,i}=-1$ and $e_i=e_j=3$.  The only cross edges between
$I_i$ and $I_j$ are the three mutually disjoint edges $a_r-b_r$
$(1\le r\le3)$.

The blockwise lift of $iji\longrightarrow jij$ is
\[
 \seq{a_1,a_2,a_3,b_1,b_2,b_3,a_1,a_2,a_3}
 \longrightarrow
 \seq{b_1,b_2,b_3,a_1,a_2,a_3,b_1,b_2,b_3}.
\]
 The following sequence of $2$- and $3$-moves is applied from
left to right at the current positions:{}
\begin{equation}
\label{eq: affine lifted 3-move 33}
\begin{aligned}
 &(2;\kappa(k,1)+2),(2;\kappa(k,1)+1),(2;\kappa(k,1)+3),(2;\kappa(k,1)+5),(2;\kappa(k,1)+4),\\
 &(2;\kappa(k,1)+3),(2;\kappa(k,1)+2),(2;\kappa(k,1)+6),(2;\kappa(k,1)+5),\\
 &(3;\kappa(k,1)),(3;\kappa(k,1)+3),\\
 &(2;\kappa(k,1)+2),(2;\kappa(k,1)+1),(2;\kappa(k,1)+3),(3;\kappa(k,1)+6),\\
 &(2;\kappa(k,1)+5),(2;\kappa(k,1)+4),(2;\kappa(k,1)+3),(2;\kappa(k,1)+2), \\
 &(2;\kappa(k,1)+6),(2;\kappa(k,1)+5).
\end{aligned}
\end{equation}
After the first nine commuting moves,   the nine-letter subword has the form  \[
 \underline{a_1b_1a_1}\,
 \underline{a_2b_2a_2}\,
 \underline{a_3b_3a_3}.
\]
The three underlined factors are changed to $b_ra_rb_r$; the remaining
moves  give the word\linebreak
$b_1b_2b_3a_1a_2a_3b_1b_2b_3$.  Thus
\eqref{eq: affine lifted 3-move 33} consists of eighteen $2$-moves and
exactly three simply-laced $3$-moves.  Its mutation--permutation
realization is obtained by replacing $(2;p)$ with $P_p$ and $(3;p)$
with $\seq{\mu_p,P_{p+1}}$, with every index read in the current
seed.

Let $D^\pm$ denote the contributions of the rows indexed by
$p\notin\{k,k+1,k+2\}$ to the positive and negative exchange
monomials at the unique downstairs mutation.
For $r=1,2,3$, let $U_r^\pm$ denote the contributions of the rows
with indices outside $[\kappa(k,1),\kappa(k+2,3)]$ to the positive
and negative exchange monomials, computed immediately before
the $r$th upstairs mutation in \eqref{eq: affine lifted 3-move 33}.

For $r=1,2,3$, the $r$th upstairs mutation in
\eqref{eq: affine lifted 3-move 33} occurs at
$\hk_r\seteq\kappa(k,1)+3(r-1)$.
Let $\hB^{(r)}=(\hb^{(r)}_{u,v})$ be the exchange matrix
immediately before the mutation at $\hk_r$.

\begin{lemma}
\label{lem: coefficient folding three move 33}
For $r=1,2,3$,
\[
 \phi^\ii(U_r^\pm)=D^\pm.
\]
More precisely, if   $p\notin\{k,k+1,k+2\}$  and
$\delta_p=b^\ii_{p,k}$,
then 
\[
 \sum_s[\mathord\pm \hb^{(r)}_{\kappa(p,s),\hk_r}]_+
   =[\mathord\pm\delta_p]_+.
\]
\end{lemma}

\begin{proof}
The verification is parallel to the size-$(2,2)$ calculation in
Lemma~\ref{lem: coefficient folding three move}.  Here the  subdiagram on $I_i\sqcup I_j${}
has the three mutually disjoint edges $a_r-b_r$.
Fix an outside position $p$ with $\delta_p\ne0$.
In the $r$th mutation column, only the row $\kappa(p,s)$ whose color
lies in the component containing $a_r$ has a nonzero entry.
This entry has the same sign as $\delta_p$.
The exchange matrix entries between the three mutation vertices
vanish. Hence mutating at one of these vertices does not change
the columns indexed by the other two. The $2$-moves only relabel
the vertices.
\end{proof}

\begin{proposition}
\label{prop: cluster variable folding three move 33}
Suppose that $\ii'$ is obtained from $\ii$ by a $3$-move
\[
 iji\longrightarrow jij,\qquad e_i=e_j=3,
\]
in the affine folding $E_6^{(1)}\longrightarrow D_4^{(3)}$ above.
Then \eqref{eq: folding of cluster var} holds.
\end{proposition}

\begin{proof}
  Write $x_1,x_2,x_3$ for the initial cluster variables at positions $k,k+1,k+2$, respectively.    The final downstairs cluster variables are
\[
 x'_1=\frac{D^+x_2+D^-x_3}{x_1},\qquad
 x'_2=x_3,\qquad x'_3=x_2.
\]
  For the cluster variables at positions $\kappa(k,1),\ldots,\kappa(k+2,3)$ in $\Sigma^\hii$, write 
\[
 \hx_1,\hx_2,\hx_3;\quad
 \hx_4,\hx_5,\hx_6;\quad
 \hx_7,\hx_8,\hx_9
\]
 with three variables for each successive position of $iji$.  The three mutations create
variables $z_r$ satisfying
\[
 \hx_rz_r=U_r^+\hx_{3+r}+U_r^-\hx_{6+r}
 \qquad(1\le r\le3).
\]
Lemma~\ref{lem: coefficient folding three move 33} gives
\[
 \phi^\ii(z_1)=\phi^\ii(z_2)=\phi^\ii(z_3)=x'_1.
\]
 After the final permutations, the nine upstairs cluster variables are{}
\[
 \seq{z_1,z_2,z_3,\hx_7,\hx_8,\hx_9,
       \hx_4,\hx_5,\hx_6}.
\]
 Consequently, $\phi^\ii$ sends the variables corresponding to the successive positions of $jij$ respectively to{}
$x'_1,x'_2,x'_3$, proving \eqref{eq: folding of cluster var}.  
\end{proof}

\subsection{\texorpdfstring{  The $4$-move with orbit sizes $(4,2)$}{The 4-move with orbit sizes (4,2)}}
\label{subsec: affine braid 42}

For $D_{2n+2}^{(1)}$ with $n\ge2$, put $N=2n+2$ and consider the order-four
automorphism whose orbits are
\[
 I_0=(0,1,N-1,N),\qquad
 I_q=(q+1,N-1-q)\ (1\le q<n),\qquad
 I_n=(n+1).
\]
The action is admissible.  Direct orbit summation gives the Cartan
matrix of $A_{2n}^{(2)}$; in particular,
\[
 (c_{0,1},c_{1,0})=(-2,-1),\qquad(e_0,e_1)=(4,2).
\]
We take $i=0$, $j=1$ and write
\[
 I_i=(a_1,a_2,a_3,a_4),\qquad I_j=(b_1,b_2).
\]
The  subdiagram on $I_i\sqcup I_j${} has two connected components, each consisting of one
vertex of $I_j$ adjacent to two vertices of $I_i$.
 The size-$(2,4)$ case is obtained by reversing the path below.

 The blockwise lifts of $ijij$ and $jiji$ are, respectively
\[
\begin{aligned}
 &\seq{a_1,a_2,a_3,a_4,b_1,b_2,
       a_1,a_2,a_3,a_4,b_1,b_2},\\
 &\seq{b_1,b_2,a_1,a_2,a_3,a_4,
       b_1,b_2,a_1,a_2,a_3,a_4}.
\end{aligned}
\]
 The following sequence of $2$- and $3$-moves connects the two displayed words:{}
\begin{equation}
\label{eq: affine lifted 4-move 42}
\begin{aligned}
&(2;\kappa(k,1)),(2;\kappa(k,1)+2),(2;\kappa(k,1)+3),(2;\kappa(k,1)+2),
 (2;\kappa(k,1)+5) \\
& (2;\kappa(k,1)+4),(2;\kappa(k,1)+3),(2;\kappa(k,1)+6),(2;\kappa(k,1)+5),
  (2;\kappa(k,1)+4), \\
& (2;\kappa(k,1)+9),(2;\kappa(k,1)+8), (2;\kappa(k,1)+7),(2;\kappa(k,1)+6),
  (2;\kappa(k,1)+5), \\
& (3;\kappa(k,1)+7),(3;\kappa(k,1)+9),(2;\kappa(k,1)+8), (3;\kappa(k,1)+6),
  (3;\kappa(k,1)+8), \\
& (3;\kappa(k,1)+1),(3;\kappa(k,1)+3), (2;\kappa(k,1)+2),(3;\kappa(k,1)),
  (3;\kappa(k,1)+2), \\
& (2;\kappa(k,1)+5),(2;\kappa(k,1)+6),(2;\kappa(k,1)+7), (2;\kappa(k,1)+8),
  (2;\kappa(k,1)+4), \\
& (2;\kappa(k,1)+5), (2;\kappa(k,1)+6),(2;\kappa(k,1)+7),(2;\kappa(k,1)+3),
  (2;\kappa(k,1)+4), \\
& (2;\kappa(k,1)+5), (2;\kappa(k,1)+4), (2;\kappa(k,1)+1),(2;\kappa(k,1)+2),
  (2;\kappa(k,1)+1).
\end{aligned}
\end{equation}
It contains thirty-two $2$-moves and eight simply-laced $3$-moves.
Replacing the former by adjacent permutations and each $(3;p)$ by
$\seq{\mu_p,P_{p+1}}$ gives eight mutations and forty permutations,
again with the labels at the current step.

For $a=1,2,3$, let $D_a^\pm$ denote the contributions of the rows
indexed by $p\notin\{k,k+1,k+2,k+3\}$ to the positive and negative
exchange monomials at the $a$th downstairs mutation in the $4$-move
realization.
For $q=1,\ldots,8$, let $U_q^\pm$ denote the contributions of the rows
with indices outside $[\kappa(k,1),\kappa(k+3,2)]$ to the positive
and negative exchange monomials, computed immediately before
the $q$th upstairs mutation in \eqref{eq: affine lifted 4-move 42}.

\begin{lemma}
\label{lem: coefficient folding four move 42}
 For the $4$-move $ijij\longrightarrow jiji$ at position $k$ in
$\ii$, the coefficient factors satisfy
\[
(D_2^+,D_2^-)=(1,D_1^-),\qquad
(D_3^+,D_3^-)=(1,D_1^+),
\]
 and
\[
\bigl(\phi^\ii(U_q^+),\phi^\ii(U_q^-)\bigr)=
\begin{cases}
(D_1^+,D_1^-),&q=1,3,5,7,\\
(1,1),&q=2,6,\\
(1,D_1^+),&q=4,8.
\end{cases}
\]
 Consequently, $\phi^\ii$ sends the variables corresponding to each position of $jiji$ to the corresponding downstairs cluster variable, as in \eqref{eq: folding of cluster var}.
\end{lemma}

\begin{proof}
The verification uses the same calculation of the factors $D_a^\pm$ and $U_q^\pm$ from exchange matrix entries as for the $(e_i,e_j)=(2,1)$ $4$-move 
in Lemma~\ref{lem: tropical four move 21},
applied independently to the two  connected{} components described above.
Moves between the components are
$2$-moves and hence only permute the current labels.   Substituting these coefficient identities into the calculation in
Proposition~\ref{prop: cluster folding four move 21} gives the
cluster-variable identities for each component.
\end{proof}

\begin{proposition}
\label{prop: cluster variable folding four move 42}
Suppose that $\ii'$ is obtained from $\ii$ by the affine $4$-move
\[
 ijij\longrightarrow jiji,\qquad(e_i,e_j)=(4,2),
\]
above, or by  the reverse move.  Then
\eqref{eq: folding of cluster var} holds.
\end{proposition}

\begin{proof}
We first compute  the final cluster variables  with all factors $D_a^\pm$ and
$U_q^\pm$ set to $1$. Write $x_1,x_2,x_3,x_4$ for  the four initial downstairs variables  and put
\[
 A=x_1^2x_4+(x_2+x_3)^2,\qquad
 C=x_1^2x_4+x_2x_3+x_3^2.
\]
 The three downstairs mutations and final permutations give
\[
 \seq{x'_1,x'_2,x'_3,x'_4}
 =\seq{\frac{A}{x_1^2x_2},
        \frac{C}{x_1x_2},x_4,x_3}.
\]
  Let $\hx_1,\ldots,\hx_{12}$ be the cluster variables at positions $\kappa(k,1),\ldots,\kappa(k+3,2)$ in $\Sigma^\hii$, in this order.   Thus $\phi^\ii$ sends
\[
 \hx_1,\ldots,\hx_4\mapsto x_1,\quad
 \hx_5,\hx_6\mapsto x_2,\quad
 \hx_7,\ldots,\hx_{10}\mapsto x_3,\quad
 \hx_{11},\hx_{12}\mapsto x_4.
\]
The eight upstairs mutations in
\eqref{eq: affine lifted 4-move 42} produce
\[
\begin{aligned}
 z_1&=\frac{\hx_1\hx_2\hx_{11}+\hx_5^2+
 \hx_5\hx_7+\hx_5\hx_8+\hx_7\hx_8}
 {\hx_1\hx_2\hx_5},\\
 z_2&=\frac{\hx_3\hx_4\hx_{12}+\hx_6^2+
 \hx_6\hx_9+\hx_6\hx_{10}+\hx_9\hx_{10}}
 {\hx_3\hx_4\hx_6},\\
 w_1&=\frac{\hx_1\hx_2\hx_{11}+\hx_5\hx_7+\hx_7\hx_8}
 {\hx_2\hx_5},
&w_2&=\frac{\hx_1\hx_2\hx_{11}+\hx_5\hx_8+\hx_7\hx_8}
 {\hx_1\hx_5},\\
 w_3&=\frac{\hx_3\hx_4\hx_{12}+\hx_6\hx_9+\hx_9\hx_{10}}
 {\hx_4\hx_6},
&w_4&=\frac{\hx_3\hx_4\hx_{12}+\hx_6\hx_{10}+\hx_9\hx_{10}}
 {\hx_3\hx_6}.
\end{aligned}
\]
 After the final permutations, the twelve upstairs cluster variables are{}
\[
 \seq{z_1,z_2,w_1,w_2,w_3,w_4,
       \hx_{11},\hx_{12},
       \hx_7,\hx_8,\hx_9,\hx_{10}}.
\]
Direct substitution gives
\[
\begin{aligned}
 \phi^\ii(z_1)=\phi^\ii(z_2)&=\frac{A}{x_1^2x_2}=x'_1,\\
 \phi^\ii(w_1)=\cdots=\phi^\ii(w_4)&=\frac{C}{x_1x_2}=x'_2,\\
 \phi^\ii(\hx_{11})=\phi^\ii(\hx_{12})&=x_4=x'_3,\\
 \phi^\ii(\hx_7)=\cdots=\phi^\ii(\hx_{10})&=x_3=x'_4.
\end{aligned}
\]
 The four groups correspond to the successive positions of $jiji$
and have cardinalities $2,4,2,4$. For a general ambient sequence, the same calculation is performed 
with the factors $D_a^\pm$ and $U_q^\pm$ retained. By Lemma~\ref{lem: coefficient folding four move 42}, after applying
$\phi^\ii$ the exchange relations agree componentwise with those of the $(e_i,e_j)=(2,1)$ calculation in
Proposition~\ref{prop: cluster folding four move 21}. Hence the same elimination proves
  \eqref{eq: folding of cluster var}.
 For the reverse move, we use  the sequence \eqref{eq: affine lifted 4-move 42} in reverse order.
\end{proof}

Together with the  braid-move calculations in Section~\ref{sec: braid}{}, these two
propositions complete \eqref{eq: folding of cluster var} for every
 finite braid move in the admissible finite and affine
foldings considered here.

 For a general admissible folding, the subdiagram on the two
color orbits involved in the braid move may have more connected
components.{}  In the next subsection, we show that the verification of
\eqref{eq: folding of cluster var} reduces to the case in which the
subdiagram on
$I_i\sqcup I_j$
is connected.  Thus no new calculation is required   when the subdiagram on $I_i\sqcup I_j$ has several connected components.

\subsection{\texorpdfstring{  Reduction to connected subdiagrams}{Reduction to connected subdiagrams}}

{}  Let $i,j\in I$
be the vertices involved in a downstairs braid move.  The subdiagram
of $\hDynkin$ on $I_i\sqcup I_j$ decomposes into isomorphic connected
components
\[
\hDynkin_1\sqcup\cdots\sqcup\hDynkin_r.
\]
For $1\le s\le r$, let $\hJ_s$ be the positions in the lifted braid
block whose colors belong to $\hDynkin_s$.
When the vertex labels are permuted, the sets $\hJ_s$ are
relabelled in the same way; mutation leaves these sets unchanged.

\begin{lemma}\label{lem: componentwise braid reduction}
 At every stage of the mutation--permutation sequence realizing the $2$- and $3$-moves in Proposition~\ref{prop: canonical finite braid lift},
\[
b_{ab}=0
\qquad
(a\in\hJ_s,\ b\in\hJ_t,\ s\ne t)
\]
whenever $b$ is exchangeable.
Consequently, mutations supported on distinct $\hJ_s$ commute, and
mutations supported on $\hJ_s$ do not change the mutation columns
 indexed by{} $\hJ_t$ for $s\ne t$.
\end{lemma}

\begin{proof}
 For the initial exchange matrix,
\[
b_{ab}=0
\qquad
(a\in\hJ_s,\ b\in\hJ_t,\ s\ne t).
\]
If $k\in\hJ_s$ is mutated and $b\in\hJ_t$, $s\ne t$, then
$b_{kb}=0$.  Hence
\[
b'_{ab} = b_{ab} +[b_{ak}]_+[b_{kb}]_+ -[-b_{ak}]_+[-b_{kb}]_+ = b_{ab} \quad (a\ne k)
\qtq 
b'_{kb}=-b_{kb}=0.
\]
Thus  the entries between distinct components remain zero.  The permutations in
 this sequence only permute the labels, and the assertion follows
inductively.
\end{proof}

\begin{proposition}\label{prop: reduction to connected configurations}
 If \eqref{eq: folding of cluster var} holds for each connected
component $\hDynkin_s$ of the subdiagram on $I_i\sqcup I_j$,
then it holds for the original braid move.
\end{proposition}

\begin{proof}
In the construction of Proposition~\ref{prop: canonical finite braid lift},
$2$-moves first group the letters by connected component.
For each $s$, let $\hsfT_s$ be the mutation--permutation sequence
for $\hDynkin_s$; it acts only at positions in $\hJ_s$.
Further $2$-moves then give the letter order of $\hii'$.

By Lemma~\ref{lem: componentwise braid reduction}, mutations supported
on distinct components commute, and the exchange columns occurring
along $\hsfT_s$, including their entries in   rows indexed outside $\bigcup_{u=1}^r\hJ_u$, are unchanged by mutations supported on $\hJ_t$, $t\ne s$.   
Thus the calculations for the connected components
may be performed independently.

 Assume that \eqref{eq: folding of cluster var} holds for every connected component.{}  Then, for each $s$,
\[
\phi^\ii\bigl(\hx'_{\kappa'(p,u)}\bigr)=x'_p
\]
whenever the position $\kappa'(p,u)$ belongs to $\hJ_s$.

The $2$-moves before and after the sequences $\hsfT_s$
only permute the position labels.
Hence
\[
\phi^\ii\bigl(\hx'_{\kappa'(p,u)}\bigr)=x'_p
\qquad
(1\le p\le l,\ 1\le u\le e_{i'_p}),
\]
which proves the assertion.
\end{proof}

Thus the verification of \eqref{eq: folding of cluster var} reduces to the calculations for the connected subdiagrams treated in Section~\ref{sec: braid}.

\bibliographystyle{amsplain}
\bibliography{ref}{}
\end{document}